\documentclass[10pt]{amsart}
\pdfoutput=1
\usepackage[a4paper,bottom=35mm]{geometry}
\usepackage[OT2,OT1]{fontenc}
\usepackage{amsmath,amssymb}
\usepackage[UKenglish]{babel}
\usepackage{enumerate,mathrsfs,mathtools}
\usepackage{graphicx,tikz,tikz-cd}
    \graphicspath{ {./gfx/} }
\usepackage{caption}
\usepackage[colorlinks=true,urlcolor=blue,linkcolor=blue,citecolor=blue]{hyperref}
    \newtheorem{theorem}{Theorem}
    \newtheorem{proposition}[theorem]{Proposition}
    \newtheorem{corollary}[theorem]{Corollary}
    \newtheorem{lemma}[theorem]{Lemma}
\theoremstyle{definition}
    \newtheorem{definition}[theorem]{Definition}
    \newtheorem{example}[theorem]{Example}
    \newtheorem{question}[theorem]{Question}
    \newtheorem{conjecture}[theorem]{Conjecture}
    \newtheorem{remark}[theorem]{Remark}
    
\numberwithin{equation}{section}
\numberwithin{theorem}{section}
\numberwithin{figure}{section}

\newcommand\switchToCyrilic{%
    \renewcommand\rmdefault{wncyr}%
    \renewcommand\encodingdefault{OT2}%
    \normalfont\selectfont}
\DeclareTextFontCommand{\textcyr}{\switchToCyrilic}
\newcommand{\Lob}{\mathord{{\textcyr{L}}}}

\DeclareMathOperator{\Vol}{Vol}
\DeclareMathOperator{\Arg}{Arg}
\DeclareMathOperator{\Real}{Re}
\DeclareMathOperator{\Ad}{Ad}

\DeclareMathOperator{\Imag}{Im}
\DeclareMathOperator{\Id}{Id}

\DeclareMathOperator{\hol}{hol}
\DeclareMathOperator{\hyperF}{F}
\DeclareMathOperator{\Ker}{Ker}
\DeclareMathOperator{\Obs}{Obs}
\DeclareMathOperator{\linspan}{span}

\DeclareMathOperator{\diag}{diag}

\DeclareMathOperator{\pned}{pned}
\DeclareMathOperator{\sgn}{sgn}
\DeclareMathOperator*{\Res}{Res}

\newcommand{\abs}[1]{\left|#1\right|}
\newcommand{\ahol}[1]{\measuredangle_{#1}}
\newcommand{\aholM}{A_\mu}
\newcommand{\aholL}{A_\lambda}

\newcommand{\Areal}{\mathcal{A}^{\RR}}
\newcommand{\Astrict}{\mathcal{A}^{(0,\pi)}}
\newcommand{\Betaweight}{\mathcal{B}_{\omega,\alpha}^{\tet}}
\newcommand{\Betasigma}[1]{\mathcal{B}^{#1}}
\newcommand{\Beta}{\mathrm{B}}
\newcommand{\bd}{\partial}
\newcommand{\CC}{\mathbb{C}}
\newcommand{\circLob}{\widetilde{\Lob}}
\newcommand{\coldotsvec}[2]{\begin{pmatrix}#1\\\vdots\\#2\end{pmatrix}}
\newcommand{\contribreal}{\mathcal{C}_\mathbb{R}}
\newcommand{\contribnonreal}{\mathcal{C}_*}
\newcommand{\customBetaweight}[1]{\mathcal{B}_{\omega,#1}^{\tet}}
\newcommand{\del}{\partial}
\newcommand{\dHaar}{\mathrm{dHaar}}
\newcommand{\edges}{\mathcal{E}}
\newcommand{\eps}{\varepsilon}
\newcommand{\eqdef}{\overset{\text{def}}{=\joinrel=}}
\newcommand{\FF}{\mathbb{F}}
\newcommand{\geo}{{\mathrm{geo}}}
\newcommand{\gluvar}{\mathcal{V}_{\triang}}

\newcommand{\hyperinv}{\beta}
\newcommand{\hypersymbol}[2]{\sideset{_{#1}\mkern-0.5mu}{\mkern-4mu_{#2}}{\hyperF}}
\newcommand{\I}[3]{\mathcal{I}_{#1,(#2,#3)}}
\def\ID{I_{\Delta}}
\newcommand{\ind}{3D\nobreakdash-index}
\newcommand{\intersection}{\iota}
\newcommand{\Ipre}[3]{\sideset{_{\mathrm{pre}}}{_{#1,(#2,#3)}}{\mathop{\mathcal{I}}}}

\newcommand{\justI}{\mathcal{I}}
\newcommand{\justIpre}{\sideset{_{\mathrm{pre}}}{}{\mathop{\mathcal{I}}}}
\newcommand{\Li}{\mathrm{Li}_{2}}
\newcommand{\Lis}[1]{\mathrm{Li}_{#1}}

\newcommand{\maillink}[1]{\href{mailto:#1}{#1}}
\newcommand{\Mat}{\mathcal{M}}
\newcommand{\MB}{I_{\mathcal{MB}}}

\newcommand{\NN}{\mathbb{N}}
\newcommand{\nonrealrep}{non-real}
\newcommand{\nword}[2]{\ensuremath{#1}\nobreakdash-#2}
\newcommand{\oA}{\overline{A}}
\newcommand{\oB}{\overline{B}}
\newcommand{\oC}{\overline{C}}
\newcommand{\oneloop}{\tau_{1\text{-loop}}}
\newcommand{\PachA}{\mathcal{P}_{2\mathord{\rightarrow}3}}
\newcommand{\PachD}{\mathcal{P}_{3\mathord{\rightarrow}2}}

\newcommand{\pFq}[5]{\hypersymbol{#1}{#2}\mkern-4mu\left(%
    \genfrac{}{}{0pt}{0}{#3}{#4}%
    \middle|\,#5\right)\mkern-4mu}
\newcommand{\iiiFii}[5]{\pFq{3}{2}{{#1}\mkern20mu{#2}\mkern20mu{#3}}{{#4}\mkern20mu{#5}}{1}}

\newcommand{\preAholo}[1]{\mathcal{P}_{#1}^{(0,\pi)}}
\newcommand{\preAreal}{\mathcal{P}^{\RR}}
\newcommand{\preArealsub}[1]{\mathcal{P}^{\RR}_{#1}}
\newcommand{\preAsub}[1]{\mathcal{P}^{(0,\pi)}_{#1}}
\newcommand{\preMB}{\sideset{_{\mathrm{pre}}}{}{\mathop{I_{\mathcal{MB}}}}}
\newcommand{\pvint}{\vp\mkern-6mu\int}
\newcommand{\Proofpart}[2]{\par\noindent\textbf{#1~#2}\ }
\newcommand{\PSLC}{PSL(2,\CC)}
\newcommand{\PSLR}{PSL(2,\RR)}
\newcommand{\qPoch}[2]{\left(#1;#2\right)_\infty}

\newcommand{\realrep}{real}

\newcommand{\RR}{\mathbb{R}}
\newcommand{\SAS}{\mathcal{S\mkern-2.5mu A}}
\newcommand{\SASdistinguished}{\mathcal{S\mkern-2.5mu A}^\geo}
\newcommand{\SASnonflat}{\mathcal{S^+\mkern-5mu A}}
\newcommand{\Shear}{\mathscr{S}}
\newcommand{\SLC}{SL(2,\CC)}
\newcommand{\Step}[1]{\Proofpart{Step}{#1}}
\newcommand{\Sval}{\nword{S^1}{valued}}
\newcommand{\TAS}{T\mkern-2mu\mathcal{AS}}
\newcommand{\tet}{\Delta}
\newcommand{\T}{T}
\newcommand{\Tet}{\mathbb{T}}
\newcommand{\threetwo}{\ensuremath{3\mathord{\rightarrow}2}}
\newcommand{\torsion}{\mathbb{T}_{\Ad}}
\newcommand{\torsionNormalised}{\mathbb{T}_{0}}
\newcommand{\transp}{\top} % Transpose symbol; it always goes in the superscript
\newcommand{\triang}{\mathcal{T}}
\newcommand{\twothree}{\ensuremath{2\mathord{\rightarrow}3}}
\newcommand{\Volhyp}{\Vol_{\mathbb{H}}}
\newcommand{\vp}{\mathrm{v.p.}}
\newcommand{\Ztaut}{\nword{\ZZ_2}{taut}}
\newcommand{\Ztautset}{\Omega_{\text{\rm taut}}}

\newcommand{\ZZ}{\mathbb{Z}}
\newcommand{\Z}{\mathbb{Z}}

\newcommand{\Case}[1]{\Proofpart{Case}{#1}}

\newcommand{\pheq}{\phantom{{}={}}}

\begin{document}

% Define the title:
\title[On the asymptotics of the meromorphic 3D-index]{On the asymptotics of the meromorphic 3D-index}

% Define authors
\author[C. Hodgson]{Craig D. Hodgson}
\author[A. Kricker]{Andrew J. Kricker}
\author[R. Siejakowski]{Rafa\l\ M. Siejakowski}

% MSC classification and keywords:

% Addresses of authors and their emails:
\address{% Craig's address
School of Mathematics and Statistics\endgraf
The University of Melbourne\endgraf
Parkville, Victoria 3010\endgraf
Australia
}
\email{\maillink{craigdh@unimelb.edu.au}}
\address{% Andrew's address
School of Physical and Mathematical Sciences\endgraf
Nanyang Technological University\endgraf
21~Nanyang Link\endgraf
Singapore~637371
}
\email{\maillink{ajkricker@ntu.edu.sg}}
\address{% Rafael's address
Instituto de Matem\'atica e Estat\'\i{}stica\endgraf
Universidade de S\~ao Paulo\endgraf
Rua do Mat\~ao 1010, 05508-090 S\~ao Paulo, SP\endgraf
Brazil
}
%\email{\maillink{rafal@ime.usp.br}}
\email{\maillink{rs@rs-math.net}}

% The abstract of the paper
\begin{abstract}
In their recent work, Garoufalidis and Kashaev extended the 3D-index of an
ideally triangulated 3-manifold with toroidal boundary to a well-defined
topological invariant which takes the form of a meromorphic function of $2$
complex variables per boundary component and which depends in addition on a
quantisation parameter $q$.
In this paper, we study
asymptotics of this invariant as $q$ approaches 1 and
develop a conjectural asymptotic approximation in the form of a sum of 
contributions associated to conjugacy classes of certain boundary parabolic $PSL(2,\mathbb{C})$ representations of the fundamental group.
Furthermore, we study the coefficients appearing in these contributions, which 
include the hyperbolic volume, the `1-loop invariant' of Dimofte and Garoufalidis, 
as well as a new topological invariant of 3-manifolds with torus boundary, 
which we call the `beta invariant'.
The technical heart of our analysis is the expression of the state-integral 
of the Garoufalidis--Kashaev invariant as an integral over 
one 
connected
component of the space of circle-valued angle structures introduced by Luo. 
Our stationary phase analysis of the asymptotics of this integral reveals many connections to the theory of angle structures and volume optimization.

This investigation was motivated by extensive numerical experiments. In addition, we prove a variety of theorems about the quantities appearing in the analysis which support the overall conjectural picture.

\end{abstract}

% Print the title and TOC:
\maketitle
\tableofcontents

%==============================================================================
% Section 1: Introduction
%==============================================================================
\section{Introduction}

The \emph{\ind}\ is a fascinating invariant of ideal triangulations of oriented 3-manifolds with toroidal
boundary which associates a $q$-series with integer coefficients to each first homology class of the boundary.
It was originally obtained by Dimofte, Gaiotto and Gukov in a mathematical physics context from a certain
superconformal field theory \cite{dimofte-gaiotto-gukov 3D index, dimofte-gaiotto-gukov gauge theories}.
Garoufalidis provided an axiomatic formulation of the \ind\ and studied convergence of this invariant,
which is only defined for some ideal triangulations \cite{stavros}.

Although the physics predicts that the \ind\ should be a topological invariant,
there is no \emph{a priori} mathematical reason that this should be so.
While the \ind\ is invariant under Pachner 2-3 and 0-2 
moves on ideal triangulations
(whenever it is defined on both sides of the move),
it is not known whether all topological ideal triangulations on which the \ind\
is defined can be connected by such moves, so topological invariance does not follow.
Nonetheless, Garoufalidis, Rubinstein, Segerman and the first named author showed in
\cite{garoufalidis-hodgson-rubinstein-segerman} that the \ind\ 
does indeed give a topological invariant
in the special case of cusped hyperbolic 3-manifolds.

This paper concerns 
an extension of the \ind\ that was introduced by Garoufalidis and Kashaev
\cite{garoufalidis-kashaev} and that we will call the \emph{meromorphic \ind}.
The extension associates to \emph{every} ideal triangulation of an oriented
3-manifold with toroidal boundary a meromorphic function which has 2 variables
for every boundary component.
The meromorphic \ind\ is invariant under \emph{all} 2-3 Pachner moves, and is
hence a topological invariant.

The meromorphic \ind\ is an extension of the $q$-series \ind\ in the following sense. In the case that the ideal
triangulation has a strict angle structure, there is a Laurent expansion of the meromorphic \ind\ whose
coefficients are the $q$-series appearing in the original \ind. Therefore, we can think of the meromorphic
\ind\ as a kind of generating function for these $q$-series, whenever they are defined. (See \cite{garoufalidis-kashaev} and Appendix \ref{recover-qseries-appendix}.)
In particular, the $q$-series \ind\ is a true topological invariant when restricted to
ideal triangulations admitting strict angle structures.

The meromorphic \ind\ introduced in \cite{garoufalidis-kashaev} depends on a complex parameter $q$
satisfying $0<|q|<1$.
The main goal of this paper is to study the asymptotics of this invariant when $q\to 1$.
More precisely, we express the meromorphic \ind\ in terms
of a parameter $\hbar\in (-\infty,0)$ which determines a $q\in(0,1)$ by $q=e^{\hbar}$.
The meromorphic \ind\ evaluated on the ideal triangulation~$\triang$ of a 3-manifold $M$ 
with a single torus boundary component 
is then denoted $\I{\triang}{\aholM}{\aholL}(\hbar)$ and is a meromorphic function of
two complex variables $\aholM,\aholL$.
These variables can be interpreted as the \emph{angle-holonomies,} along 
peripheral curves $\mu$ and $\lambda$ giving a basis for $H_1(\del M;\ZZ)$, 
of a complexified angle structure on $\triang$, in a sense that will be made more precise later.

Of particular interest is the behaviour of this invariant at the point~$\aholM=\aholL=0$, corresponding
to trivial peripheral holonomies, since the value at this point does not depend on the choice of the
curves $\mu$, $\lambda$.
Moreover, in the case the manifold is hyperbolic, this point corresponds to an angle structure
derived from the complete hyperbolic structure.

In this paper, we mostly focus on the case of a one-cusped hyperbolic manifold $M$ and study
the asymptotic behaviour of the meromorphic \ind\ for $\hbar=-1/\kappa$,
which corresponds to $q=\exp(-1/\kappa)$:
\begin{equation}\label{limit-we-study}
    \I{M}{0}{0}(-1/\kappa),\quad
    \kappa>0,\quad
    \kappa\rightarrow\infty.
\end{equation}
In Section~\ref{the-conjecture}, we provide a conjectural asymptotic approximation for the above
quantity, in the form of a sum over conjugacy classes of irreducible, boundary-parabolic
representations $\rho:\pi_1(M)\rightarrow\PSLC$ with a certain fixed obstruction class
to lifting the representation to $\SLC$. 
Moreover, we will use heuristic saddle-point methods to develop explicit formulae for the leading terms of these contributions.
Our asymptotic approximations show strong agreement with numerical experiments,
which we describe in Section~\ref{sec:numerical}.

Although the overall picture is conjectural and based on a heuristic analysis,
we will present many rigorous theorems about the various expressions that
arise as coefficients in our asymptotic approximations.
In particular, we expect some of these coefficients to be given by classical quantities,
such as the volume of a representation and the Reidemeister torsion associated to its adjoint.
Moreover, we use the total predicted contribution of boundary-parabolic \emph{real}
representations $\rho:\pi_1(M)\rightarrow\PSLR$ to construct a new topological invariant,
which we call the \emph{beta invariant}.
The beta invariant is defined in terms of state-integrals on ideal triangulations built from the Euler beta function
and shown, in Section~\ref{sec:Mellin-Barnes}, to be fully invariant under Pachner 2-3 moves.

\subsection{A sketch of the meromorphic 3D-index}\label{sketch-of-the-index}

Here we'll give a brief sketch of the construction of meromorphic \ind\
in the special case of a 3-manifold with a \emph{single} torus boundary component.
Section~\ref{section-index} below gives a detailed review.

To be precise, we need to consider two cases: the case when the triangulation admits a strict angle structure,
and the case when it does not. (A strict angle structure is an assignment of strictly positive real numbers to
normal quadrilateral types which sums to $\pi$ in each tetrahedron and to $2\pi$ around each edge.)

In the case the triangulation does admit a strict angle structure, we will write down a \emph{state-integral}
expression which depends on the triangulation and a strict angle structure on it.
The strict angle structure ensures convergence of the integral, but the value of the integral
depends on the strict angle structure.
The dependence factors through the peripheral angle-holonomies~$\aholM$ and~$\aholL$.

Thus, the state-integral determines a germ of a meromorphic function defined initially
at suitable points $(\aholM,\aholL)\in\RR^2$.
The full function $\I{\triang}{\aholM}{\aholL}(\hbar)$ is then obtained as the analytic continuation 
of this germ into the complex domain.
The other case, when the triangulation does not admit a strict angle structure, is detailed in
Section~\ref{section-index} and involves a more intricate analytic
continuation argument, presented in \cite{garoufalidis-kashaev}.

In order to describe the state-integral depending on a given oriented ideal triangulation~$\triang$
with $N\geq 2$ tetrahedra, we shall denote by $Q(\triang)$ the set of normal quadrilateral types in $\triang$
and by $\edges(\triang)$ the set of edges of $\triang$.
Suppose $\alpha: Q(\triang)\to(0,\pi)$ is a strict angle structure on $\triang$.
For every edge~$E\in\edges(\triang)$, there will be one integration variable~$y_E$
taking values in the unit circle~$S^1\subset\CC$.
The integrand will be a product of \emph{tetrahedral weights}, one for each tetrahedron.
The weight of each tetrahedron is a function of the integration variables appearing on the edges of that tetrahedron.
Moreover, each tetrahedral weight consists of 
four factors -- one for each normal quad type in the tetrahedron, and one additional normalisation factor.
So the general shape of the state integral is
\begin{equation}\label{eq:general-state-integral}
c_q^N \int_{(S^1)^{\edges(\triang)}} \prod_{\square\in Q(\triang)} F(\square)\, \bigwedge_{E\in \edges(\triang)} \mathrm{d}y_E.
\end{equation}%
\begin{figure}%
    \centering%
    \begin{tikzpicture}[baseline=(BASE)]
        \node[anchor=south west,inner sep=0] (image) at (0,0)
            {\includegraphics[width=0.2\columnwidth]{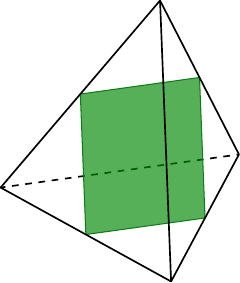}};
        \begin{scope}[x={(image.south east)},y={(image.north west)},every node/.style={inner sep=1pt}]
            \node[anchor=south east] at (0.34, 0.67) {$y_a$};
            \node[anchor=north west] at (0.85, 0.23) {$y_b$};
            \node[anchor=south west] at (0.84, 0.75) {$y_c$};
            \node[anchor=north east] at (0.3, 0.195) {$y_d$};
            \node (BASE) at (1,1) {};
         \end{scope}
    \end{tikzpicture}
    \hfill
    \parbox[t][][t]{\dimexpr0.8\columnwidth-5mm\relax}{%
        \caption{\label{fig:GK-factor}%
        The edge variables which appear in the factor in the integrand of the state-integral
        corresponding to a quad type. See Equation~\eqref{eq:stateIntegralFactor}.
        }
    }
\end{figure}
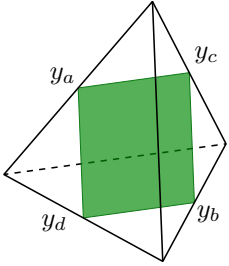%
]The factor $F(\square)$ corresponding to the quad type $\square$ shown in Figure~\ref{fig:GK-factor} is
given by
\begin{equation}\label{eq:stateIntegralFactor}
    F(\square)=G_q\Bigl(e^{\alpha(\square)\hbar/\pi}e^{i\alpha(\square)}\frac{y_ay_b}{y_cy_d}\Bigr),
\end{equation}
where $q=e^\hbar$ and $\hbar<0$. 
In the above formulae, $c_q=c(q)$ and $G_q=G_q(z)$ denote certain meromorphic functions defined in the punctured
plane and related to Kirillov's
$q$-dilogarithm~\cite{kirillov}; see Section~\ref{section-index} below for precise definitions.

The ``states'' we are integrating over in the
state-integral~\eqref{eq:general-state-integral} are assignments~$y_E$ of a complex
number of modulus 1 to every edge $E\in\edges(\triang)$.
In the factor~\eqref{eq:stateIntegralFactor}, corresponding to some quad type~$\square$,
the argument of the function~$G_q$ is given by the following product
\[
    e^{\alpha(\square)\hbar/\pi} \times \left(e^{i\alpha(\square)}\frac{y_ay_b}{y_cy_d}\right),
\]
in which the first factor is a real number in $(0,1)$, whereas the
second term is a complex number of modulus 1.

In this paper, we will rewrite the integral in terms of the variables
\begin{equation}\label{change-to-angle-structure-variables}
    \omega(\square)=e^{i\alpha(\square)}\frac{y_a y_b}{y_c y_d},
\end{equation}
which replaces the integration over states~$y:\edges(\triang)\rightarrow S^1$ with an
integral over a space of \emph{$S^1$-valued angle structures} $\omega: Q(\triang)\rightarrow S^1$.

\subsection{\texorpdfstring{$S^1$-valued angle structures}{Circle-valued angle structures}}

Our study of the meromorphic \ind\ in this paper rests on a space of structures on an ideal triangulation called \emph{\Sval\ angle structures} or \emph{circle-valued angle structures}, first introduced by Luo \cite{luo-volume}. 

In the companion work \cite{hodgson-kricker-siejakowski} we describe the spaces of \Sval\ angle structures on a given triangulation, and present theorems which count the number of connected components of these spaces in terms of the cohomology of the underlying manifold~$M$. 
To do this we construct a canonical correspondence between the topological components of the space of 
\Sval\ angle structures and the group of obstruction classes associated to the problem of lifting $\PSLC$ representations of $\pi_1(M)$ to $\SLC$. This correspondence is in a precise sense natural with respect to Thurston's gluing and completeness equations. Here we will briefly recall from \cite{hodgson-kricker-siejakowski} the details we will need in the current paper.

 There are two versions of \Sval\ angle structures on an ideal triangulation: structures with and without peripheral constraints.
 First we introduce the structures without any peripheral constraints.
 An \Sval\ angle structure on $\triang$ is an assignment $\omega:Q(\triang)\rightarrow S^1$ of a
 complex number of modulus 1 to every quad type, satisfying the two sets of equations which say that:
 \begin{enumerate}[(i)]
     \item
     For each tetrahedron, $\omega(\square)\omega(\square')\omega(\square'')=-1$,
     where $\square$, $\square'$ and $\square''$ are the three quad types associated to the
     tetrahedron -- see Figure~\ref{fig:quads}.
     \item
     For each edge $E\in\edges(\triang)$, $\prod_{\square\in Q(\triang)}\omega(\square)^{G(E,\square)}=1$,
     where $G(E,\square)$ counts the number of times the quad type $\square$ faces the edge~$E$.
 \end{enumerate}
 The manifold of such structures will be denoted $\SAS(\triang)$;
 see also Definition~\ref{angle structures} below, which follows Luo~\cite{luo-volume}.

The second version of an \Sval\ angle structure has additional peripheral constraints.
As before, we fix peripheral curves $\mu$ and $\lambda$ giving a basis for the first homology of the boundary.
An \Sval\ angle structure on $\triang$ associates to each of $\mu$ and $\lambda$ an \Sval\
\emph{angle-holonomy} -- see equation~\eqref{peripheral_S1_holonomy_defn}.
 Of special importance to our study is the space~$\SAS_0(\triang)$
 of \emph{peripherally trivial} \Sval\ angle structures, where we prescribe the peripheral
 holonomies to both be $1$.

The importance of \Sval\ angle structures in our analysis of the meromorphic \ind\ is the fact that by making the 
substitution~\eqref{change-to-angle-structure-variables}
one can express the defining state integral as an integral over \emph{precisely one} of the topological components of the space 
of \Sval\ angle structures with corresponding prescribed peripheral angle holonomy.
To be precise:
 \begin{equation}\label{rewritten-integral}
     \I{\triang}{\aholM}{\aholL}(\hbar)
     =
     \frac{c_q^N \bigl|H_1(\hat{M};\FF_2)\bigr|}{(2\pi)^{N-1}}
     \int\limits_{\Omega} \prod_{\square\in Q(\triang)}
     G_q\bigl(e^{\alpha(\square)\hbar/\pi}\omega(\square)\bigr)\;d\omega,
 \end{equation}
  where $\hat{M}$ is the end compactification of $M$.
 In the above formula,
 $\Omega$ denotes the connected component of the circle-valued angle structure~$\exp(i\alpha)$
 in the space $\SAS_{(\exp(i\aholM), \exp(i\aholL))}(\triang)$ of \Sval\ angle structures
 with prescribed peripheral angle-holonomies equal to $\bigl(\exp(i\aholM), \exp(i\aholL)\bigr)\in S^1\times S^1$.
 We obtain this integral from the definition in Section \ref{state-integral-to-circle-valued}, and the precise meaning of this expression is discussed there.

In \cite{hodgson-kricker-siejakowski} we present the following theorem classifying the topological components of the space of \Sval\ angle structures in terms of obstruction theory. The statement uses the circle-valued angle structure associated to an algebraic solution $z : Q(\triang)\rightarrow\CC$ of Thurston's gluing equations, which is defined by the formula
 \begin{equation}\label{associated-S1-structure}
     \omega(\square)=\frac{z(\square)}{|z(\square)|},
     \quad\square\in Q(\triang).
 \end{equation}

\begin{theorem}\label{obstruction-theorem}
There exist maps
\[
    \Phi:\SAS(\triang)\to H^2(M;\FF_2)
    \quad\text{and}\quad
    \Phi_0:\SAS_0(\triang)\to H^2(M,\partial M;\FF_2)
\]
satisfying the following conditions:
\begin{enumerate}[(i)]
\item\label{item:continuity}
$\Phi$ and $\Phi_0$ are continuous, i.e., constant on the
connected components of their domains.
\item\label{obstructions-bijectivity}
The induced map on the connected components,
\[
    \pi_0(\Phi_0):\pi_0\bigl(\SAS_0(\triang)\bigr)\to H^2(M,\partial M;\FF_2),
\]
is a bijection. Moreover, the induced map
\[
    \pi_0(\Phi):\pi_0\bigl(\SAS(\triang)\bigr)\to H^2(M;\FF_2)
\]
is an injection whose image coincides with the kernel of the map
\begin{equation}\label{restriction-to-boundary-H2F2}
    \iota^* : H^2(M;\FF_2) \to H^2(\partial M;\FF_2)
\end{equation}
induced by the inclusion of the boundary into $M$.
\item\label{item:correct-obstructions}
If $z$ is an algebraic solution of Thurston's gluing equations on $\triang$
defining a representation $\rho:\pi_1(M)\to\PSLC$
and $\omega(\square)=z(\square)/|z(\square)|$
is the associated circle-valued angle structure, then $\Phi(\omega)\in H^2(M;\FF_2)$ is the
cohomological obstruction $\Obs(\rho)$ 
to lifting $\rho$ to an $\SLC$ representation.
Moreover, if $\rho$ is boundary-parabolic,
then $\Phi_0(\omega)\in H^2(M,\partial M;\FF_2)$ is the cohomological obstruction $\Obs_0(\rho)$
to lifting $\rho$ to a boundary-unipotent $\SLC$ representation.
\end{enumerate}
\end{theorem}

\subsection{Asymptotic conjecture for the meromorphic 3D-index}\label{the-conjecture}

Kashaev was the first to conjecture that an asymptotic limit of a quantum invariant contains geometric
information \cite{kashaev-volume1,kashaev-volume2}.
He conjectured that the hyperbolic volume of a knot complement
is the exponential growth rate of a certain sequence
of knot invariants obtained from quantum dilogarithms.
Murakami and Murakami identified Kashaev's invariant in terms of coloured Jones polynomials, and this
influential conjecture has since become known as the Volume Conjecture \cite{murakami-murakami}.
By now there exist numerous asymptotic conjectures for quantum invariants, such as \cite{chen-yang}.
It is worth noting that Andersen and Kashaev \cite{kashaev-andersen-TQFT} have an asymptotic volume conjecture for a state-integral invariant that is in a class related to the meromorphic \ind, although 
we do not currently know if there is a relationship with our work.

In this paper, we conjecture an asymptotic formula for the limit \eqref{limit-we-study}
in terms of a 
variety of classical geometric quantities associated to the manifold,
some of which appear for the first time here.
In order to discuss the general shape of our conjecture, we assume that $M$ 
is an orientable hyperbolic 3-manifold with a single cusp 
and introduce the following notations:
\begin{itemize} 
    \item
    Let $\rho^\geo:\pi_1(M)\rightarrow\PSLC$ denote the boundary-parabolic representation (defined up
    to conjugation only) associated to the 
    complete hyperbolic structure of finite volume.
    \item
    A representation $\rho:\pi_1(M)\rightarrow \PSLC$ will be called 
     \emph{\realrep}\ if it can
    be conjugated to lie in the subgroup $\PSLR\subset \PSLC$. Otherwise it will be called 
       \emph{\nonrealrep}.
    \item
    $\mathcal{X}^\geo$ will denote the set of conjugacy classes of irreducible, boundary-parabolic
    representations $\rho:\pi_1(M)\to\PSLC$ with the same obstruction class as $\rho^\geo$, i.e.,
    satisfying $\Obs_0(\rho)=\Obs_0(\rho^\geo)$.
       \item
    $\mathcal{X}^\geo_\RR\subset\mathcal{X}^\geo$ denotes the subset consisting of \realrep\ representations.
\end{itemize}
In some cases, $\mathcal{X}^\geo$ can be infinite; see section \ref{example-infinite-family} below. 
Assuming $\mathcal{X}^\geo$ is {\em finite},  
we conjecture an asymptotic expansion for the meromorphic \ind\ of the following general form:
\begin{equation}\label{general-approximation}
    \I{M}{0}{0}(-1/\kappa)
    = \sum_{[\rho]\in\mathcal{X}^\geo_{\RR} }\contribreal([\rho],\kappa)
    + \sum_{[\rho]\in\mathcal{X}^\geo\setminus\mathcal{X}^\geo_{\RR} } \contribnonreal([\rho],\kappa)
    + o(1),\ \ \text{as}\ \kappa\to\infty.
\end{equation}
Here $\contribreal([\rho],\kappa)$ denotes a contribution associated to a conjugacy class of \realrep\
representation~$\rho$, which we postulate to be a linear function of the parameter $\kappa$,
\begin{equation}\label{real-contribution}
    \contribreal([\rho],\kappa) = |H_1(\hat{M};\FF_2)|\cdot\MB(\rho)\,\kappa,
\end{equation}
where $\MB(\rho)\in\RR$ is a 
 new conjectural invariant of (boundary-parabolic)
$\PSLR$ representations on $M$, formulated more explicitly below in terms of a certain multivariate
contour integral of Mellin--Barnes type.

Meanwhile, $\contribnonreal([\rho],\kappa)$ denotes a contribution associated to a conjugacy class of
a \nonrealrep\ representation $\rho$.
If $\rho$ is such a representation and $\overline{\rho}$ is its complex conjugate representation,
then we have
\begin{equation}\label{combined-contribution-generic}
    \contribnonreal([\rho],\kappa)+\contribnonreal([\overline{\rho}],\kappa)
    =
    C_\rho\sqrt{\kappa}\cos\left(\kappa\Vol(\rho) + \frac{\pi n_\rho}{4}\right),
\end{equation}
where $C_\rho > 0$, $n_\rho\in\ZZ$, and 
$\Vol(\rho)\in\RR$ denotes the volume of
the representation $\rho$; see \cite{francaviglia} for a discussion of this volume.

We remark that our asymptotic conjecture is substantially different from the Volume Conjecture of 
Kashaev and Murakami--Murakami \cite{murakami-murakami}, since the latter
predicts that the leading-order asymptotics of the coloured Jones polynomials at roots
of unity is determined by the volume of the manifold with its complete hyperbolic structure.
In contrast, our conjecture includes the algebraic volumes of
boundary-parabolic $\PSLC$ representations other than the one arising from hyperbolic geometry.
Moreover, the leading-order term in our approximation is a sum of contributions which appear to
originate from $\PSLR$ representations (which have zero volume).

We shall now detail what these contributions are and, in particular, provide 
explicit expressions for the quantities $\contribnonreal([\rho],\kappa)$ and $\MB(\rho)$, 
when the representation
$\rho$ arises from a solution to Thurston's gluing equations for an ideal triangulation $\triang$.

\subsection{Contributions to the asymptotics from \nonrealrep\ representations}

We start with the generic case of a boundary-parabolic
representation~$\rho : \pi_1(M)\rightarrow\PSLC$ whose conjugacy class $[\rho]$ lies in
$\mathcal{X}^\geo\setminus\mathcal{X}^\geo_\RR$.
Assume that $\rho$ arises from an algebraic solution of Thurston's edge and completeness
equations $z : Q(\triang)\rightarrow\CC$, satisfying the condition that $\Imag(z(\square))\neq0$
for all $\square\in Q(\triang)$.
Then the conjectured contribution to the asymptotics, obtained in Section~\ref{nonreal-contrib}, is: 
\begin{equation}\label{intro:contribution-of-smooth-crit}
    \contribnonreal([\rho],\kappa)=
    |H_1(\hat{M};\FF_2)|\,\tau(z)\,\exp\left[i\tfrac{\pi}{4}\bigl(N_+ - N_- +\Sigma(\omega)\bigr)\right]
    \sqrt{2\pi\kappa}\,e^{i\kappa\Vol(\rho)}.
\end{equation}
As before, $\hat{M}$ denotes the end compactification of $M$ and moreover:
\begin{itemize}
    \item
    $\tau(z)$ is a positive real number defined as a function of the solution $z$,
    in terms of the determinant of the Hessian of the hyperbolic volume on a space of $S^1$-valued angle structures, with a certain trigonometric correction.
    The next section (\ref{tau-invt-sect}) discusses the rich properties of $\tau(z)$, including the fact it is invariant
    under shaped 2-3 and 3-2 Pachner moves involving no real shapes.
    \item 
    $N_+-N_-+\Sigma(\omega)$ is an integer obtained from the signature of
    the aforementioned Hessian matrix associated to the solution $z$ and $\omega=z/|z|$. 
    It equals $+1$ for a positively oriented geometric solution and $-1$ for a negatively oriented solution.
\end{itemize}
Furthermore, the contribution from the complex conjugate representation $\overline{\rho}$ is the complex
conjugate of the contribution from $\rho$, so the sum of the two contributions is of the
form \eqref{combined-contribution-generic}.

We remark that there is generally no guarantee that a given \nonrealrep\ representation can 
be obtained from a shape parameter solution without real shapes.
Even assuming that a \nonrealrep\ representation $\rho$ is recovered from
a shaped ideal triangulation $(\triang,z)$ without real shapes,
it is not known in general whether the set of all such shaped triangulations is
connected under shaped 2-3 Pachner moves.
Consequently, our proof of invariance of the quantity~$\tau(z)$ under shaped Pachner moves does not imply
that it is a topological invariant associated to the representation, although it provides evidence for that belief.

Instead, we will study the components of our expression \eqref{intro:contribution-of-smooth-crit} as functions
of algebraic solutions to Thurston's equations and prove their various properties from that starting point.
We conjecture that there are well-defined topological invariants which specialise to these expressions in the case
the representation does in fact come from a suitable solution to Thurston's equations.

\subsection{\texorpdfstring{The $\tau$-invariant}{The tau-invariant}}
\label{tau-invt-sect}
Here we sketch the $\tau$-invariant~$\tau(z)$, an intriguing function of 
algebraic solutions of Thurston's equations which arises from our asymptotic analysis.
Theorem \ref{th:tau=oneloop}, which is one of the main theorems of this paper,
connects the $\tau$-invariant to the \emph{1-loop invariant} of Dimofte and Garoufalidis \cite{tudor-stavros},
conjectured to be equal to the adjoint Reidemeister torsion associated to the
holonomy representation.

The formula that defines the $\tau$-invariant $\tau(z)$ depends on the solution~$z$ through the
corresponding \Sval\ angle structure given by $\omega(\square)=z(\square)/|z(\square)|$.
It also involves a generalised angle structure $\alpha : Q(\triang)\rightarrow\RR$ on $\triang$ with
vanishing peripheral angle-holonomy, as defined in Section~\ref{angle_structures_section}. 
Because $\rho$ lies in $\mathcal{X}^\geo$, it turns out there exists such an angle structure $\alpha$
lifting $\omega$, in the sense that $\omega(\square)=\exp(i\alpha(\square))$.
We will prove that $\tau$ does not depend on the choice of $\alpha$.
Explicitly, $\tau(z)$ can be computed from the formula
\[
    \tau(z) = \tau(z,\alpha)
        = \frac{
            \prod\limits_{\square\in Q(\triang)}
            |\Imag \omega(\square)|^{\frac{\alpha(\square)}{\pi}-\frac{1}{2}}
          }%
          {
            \sqrt{\left|\det\bigl(\mathbf{L_*}
            \diag\bigl\{\frac{\Real\omega(\square)}{\Imag\omega(\square)}\bigr\}
            \mathbf{L_*^\transp} \bigr) \right|}
          }.
\]
Here $\mathbf{L_*}$ denotes an $(N-1)\times 3N$ integer matrix which results
from deleting a row from the matrix of coefficients of Tollefson's Q-matching equations~\cite{tollefson}.
Equivalently, the rows of $\mathbf{L_*}$ can be viewed as the
\emph{leading-trailing deformations} \cite{futer-gueritaud} associated to all edges of $\triang$ except one;
see equation \eqref{definition-L_star} for a precise definition.
Meanwhile, the matrix
$\diag\bigl\{\frac{\Real\omega(\square)}{\Imag\omega(\square)}\bigr\}$ is a $3N\times 3N$ diagonal
matrix whose diagonal entries $\frac{\Real\omega(\square)}{\Imag\omega(\square)}$
can be interpreted as cotangents of the dihedral angles in $\triang$.

The expression above appears in our asymptotic analysis through a stationary phase approximation of an integral
where the phase function is the hyperbolic volume functional on the space of \Sval\ angle structures $\SAS_0(\triang)$,
as discussed by Luo~\cite{luo-volume}.
In particular, the matrix
$\mathbf{L_*}\diag\bigl\{\frac{\Real\omega(\square)}{\Imag\omega(\square)}\bigr\}\mathbf{L_*^\transp}$
arises as (minus) the Hessian matrix of the volume functional with respect to leading-trailing coordinates;
see equation~\eqref{computation-of-Hessian} for more details.

In Section \ref{tau-exponent-invariance} we prove that $\tau(z,\alpha)$ does not depend on the choice of the
angle structure $\alpha$ and then, in Section \ref{pachner-invariance}, we establish the following theorem.
\begin{theorem}
Let $(\triang_2, z^{(2)})$ and $(\triang_3, z^{(3)})$ be two shaped
ideal triangulations, where $z^{(2)}$ and $z^{(3)}$ are algebraic solutions of Thurston's edge and completeness
equations in $\CC\setminus\RR$.
Suppose that $(\triang_2, z^{(2)})$ and $(\triang_3, z^{(3)})$ are related by a shaped Pachner 2-3 move,
so that both solutions correspond to the same conjugacy class of a boundary-parabolic representation $\rho$.
Assume moreover that $\Obs_0(\rho)=\Obs_0(\rho^\geo)$.
Then $\tau\bigl(z^{(2)}\bigr)=\tau\bigl(z^{(3)}\bigr)$.
\end{theorem}

Our main asymptotic conjecture implies that $\tau(z)$ should be a topological invariant
associated to the pair consisting of the manifold~$M$ and the boundary-parabolic
representation~$\rho$ determined by $z$ (up to conjugation).
While the above theorem provides evidence in favour of this claim, it does not explain the potential
topological meaning of the quantity~$\tau(z)$.
Nonetheless, after further study of $\tau(z)$, we establish Theorem~\ref{th:tau=oneloop} in
which we connect $\tau(z)$ to the \nword{1}{loop} invariant~$\oneloop(z, \gamma)$
of Dimofte and Garoufalidis \cite{tudor-stavros} in the case when the shape
parameter solution~$z$ represents the complete hyperbolic structure of finite volume on $M$.

The \nword{1}{loop} invariant has its origins in a quantisation of an $\SLC$ Chern-Simons gauge theory
due to Dimofte \cite{dimofte}, based on quantum dilogarithms and the Neumann-Zagier symplectic structure
of the gluing equations, and yielding a non-rigorous state-integral expression over shape parameters.
A formal asymptotic expansion of this expression, which the gauge theory predicts should be an asymptotic
expansion of the Kashaev invariant, was computed and studied in \cite{tudor-stavros}.
The \nword{1}{loop} invariant arises as the sub-leading term of the aforementioned formal expansion.
It depends on both an algebraic solution~$z$ of Thurston's equations yielding the holonomy representation
of the hyperbolic structure, and also a choice of an essential simple closed peripheral curve $\gamma$.

Assuming that the solution $z$ does not involve real shapes and recovers the geometric
representation $\rho^\geo$, we prove, in Theorem \ref{th:tau=oneloop}, that
\begin{equation}\label{tau=1loop-into}
    \tau(z) = \frac{L(\gamma)}{\sqrt{2A}}\left|\oneloop(z,\gamma)\right|^{-1},
\end{equation}
where $A$ is the Euclidean area of a horospherical cross-section of the cusp of $M$
and $L(\gamma)$ is the length of a geodesic representative of the free homotopy class of
$\gamma$ computed in the same cross-section.

One of the main conjectures of \cite{tudor-stavros}, known as the ``\nword{1}{loop} conjecture'', predicts that
\begin{equation}\label{oneloop-conjecture-eq-intro}
    \oneloop(z, \gamma) = \pm\torsion(M,\rho_z, \gamma),
\end{equation}
where $\torsion(M,\rho_z,\gamma)\in\CC^*/\{\pm1\}$ is the adjoint Reidemeister
torsion introduced by Porti~\cite{porti1997}, corresponding to the holonomy representation~$\rho_z$
of an algebraic solution~$z$, which is either $\rho^\geo$ or a small deformation of it.
In particular, combining \eqref{tau=1loop-into} and \eqref{oneloop-conjecture-eq-intro}
we obtain a conjectural topological interpretation of our $\tau$-invariant
in terms of the adjoint Reidemeister torsion associated to the geometric representation $\rho^\geo$.

\subsection{Contributions to the asymptotics from \realrep\ representations and beta invariant} 
In our conjecture for the asymptotics of the meromorphic \ind, there is a linear contribution,
 obtained in Section~\ref{sec:derivation-of-MB}, 
of the form~\eqref{real-contribution} for each $[\rho]\in\mathcal{X}^\geo_\RR$,
a conjugacy class of boundary-parabolic real representations having the geometric obstruction class and detected by a solution of Thurston's gluing equations over the reals.
We conjecture that there exists a topological invariant~$\MB(\rho)$ associated to the class
$[\rho]$ which determines the slope of the corresponding linear term~\eqref{real-contribution}.
To study this slope, in this paper we will make the assumption that the representation $\rho$ is determined
by a solution of Thurston's edge and completeness equations {\em lying in the real numbers},
$z: Q(\triang)\rightarrow(\RR\setminus\{0,1\})\subset\CC$.
We give an explicit integral formula for the coefficient $\MB(\rho)$ for such a representation,
and conjecture that the defining integrals always converge and that there is a topological invariant which specialises to the given formula.

As with the $\tau$-invariant discussed in the last section,
this contribution will factor through the \Sval\ angle structure corresponding to the solution:
$\omega(\square) = z(\square)/|z(\square)| = \sgn z(\square)$.
The contribution arises from a neighbourhood of $\omega$ in the integration domain $\Omega$ of the state-integral
for the meromorphic \ind .

As $z$ takes real values, the \Sval\ angle structure $\omega$ becomes a function $\omega:Q(\triang)\rightarrow\{+1,-1\}$.
Such angle structures were first studied by Luo \cite{luo-volume},
who called them \emph{\Ztaut\ angle structures}.
A \Ztaut\ angle structure is an assignment of one of $\{+1,-1\}$ to every quad type such that two conditions are met.
The first condition is that in every tetrahedron two of the quad types are assigned $+1$ while the third is assigned $-1$.
The second condition is that the product of the signs around every edge is $+1$.
See also parts (iv) and (v) of Definition~\ref{angle structures}.

In Section~\ref{sec:Mellin-Barnes}, we introduce a real-valued 
multivariate contour integral of Mellin--Barnes type which
we will denote $\MB(\triang,\omega,\alpha)$; see Definition \ref{definition-MB-integral-at-omega}. 
This integral can be viewed as a state integral of Turaev--Viro type on ideal triangulations, with charged tetrahedral weights given by Euler's beta function, similar to the integrals discussed in \cite{kashaev-luo-vartanov}.
The integral~$\MB(\triang,\omega,\alpha)$ can be written down explicitly using the gluing matrices of the triangulation $\triang$,
a \Ztaut\ angle structure $\omega$, and a strict angle structure $\alpha$ on $\triang$ with
vanishing peripheral angle-holonomy.
We will call $\MB(\triang,\omega,\alpha)$ the \emph{Mellin--Barnes integral} associated to this data.

Subsequently, we prove that the Mellin--Barnes integral does not in fact depend on the choice or
existence of a peripherally-trivial strict angle structure~$\alpha$. 
Therefore, $\MB(\triang,\omega)=\MB(\triang,\omega,\alpha)$ is fully determined
by the \Ztaut\ angle structure $\omega$ on $\triang$.
If the Mellin--Barnes integrals $\MB(\triang,\omega)$ converge in the sense of Cauchy's principal value for all \Ztaut\ angle
structures $\omega$ in the connected component of $\SAS_0(\triang)$ corresponding to the geometric obstruction
class, then we call their sum~$\hyperinv(\triang)\in\RR$  
the \emph{beta invariant} of the triangulation $\triang$.
\begin{theorem}
When $\triang$, $\triang'$ are two ideal triangulations of the same connected, orientable, non-compact
3-manifold $M$ whose ideal boundary is a torus,
and the beta invariant is defined on both triangulations,
then $\hyperinv(\triang)=\hyperinv(\triang')$.
In other words, the beta invariant is a topological invariant $\beta(M)$ of $M$ whenever it is defined on at least one ideal triangulation of $M$.
\end{theorem}
The above theorem is a consequence of the invariance of $\hyperinv(\triang)$
under 2-3 Pachner moves, which we establish in Theorem~\ref{beta-invariance-theorem}.
This shows that $\hyperinv(M)$ constitutes a previously unknown topological invariant of 3-manifolds
with torus boundary.
Although the existence of such a state-integral invariant was previously suggested by Kashaev \cite{kashaev-beta},
its topological meaning and detailed properties remain a mystery. 

The aforementioned results strengthen our belief 
that the coefficient $\MB(\rho)$ in the
conjectural asymptotic contribution~\eqref{real-contribution}
associated to a real solution $z$ of Thurston's equations is given by the Mellin--Barnes
integral $\MB(\triang,\omega)$ where $\omega(\square)=\sgn z(\square)$.
Curiously, our asymptotic analysis predicts such a contribution from \emph{every}
\Ztaut\ angle structure $\omega$ in the domain $\Omega$ of the state integral \eqref{rewritten-integral}.
Some of those points will correspond to real solutions of Thurston's equations but some will not.
In Conjecture \ref{MB-vanishing-conjecture}, we therefore predict that at points $\omega$ not arising
from solutions to Thurston's gluing equations, the corresponding Mellin--Barnes integral is zero,
i.e., that such \Ztaut\ angle structures do not contribute to the beta invariant.
In Section~\ref{sec:example-4_1}, we establish this property for the standard two-tetrahedron
triangulation of the figure-eight knot complement.

\subsection{Example: the SnapPea census manifold `m011'}
We will now illustrate our conjectured asymptotic approximation \eqref{general-approximation}
in the case of the SnapPea census manifold $M=\mathtt{m011}$ \cite{snappy}.
This orientable one-cusped manifold has a 3-tetrahedron
ideal triangulation $\triang$ with Regina \cite{regina} isomorphism signature \texttt{dLQacccjsjb}.

Using Regina, we compute $H_1(M;\ZZ)\cong\ZZ$,
which implies $H^2(M,\partial M;\FF_2)\cong\FF_2$ via the Universal Coefficient Theorem and
Poincar\'e duality. So there are two possible obstruction classes in this case.
The Ptolemy database \cite{ptolemy} records the conjugacy classes of boundary-parabolic
$\PSLC$-representations of $\pi_1(M)$ and groups them by corresponding obstruction class.
The entry for $\mathtt{m011}$ confirms that both elements of $H^2(M,\partial M;\FF_2)$
arise as obstructions to lifting certain boundary-parabolic $\PSLC$-representations of $\pi_1(M)$ to
boundary-unipotent representations in $\SLC$.
Furthermore, by \cite[Proposition~9.19]{garoufalidis-thurston-zickert}, the unique non-zero class
must correspond to the holonomy representation $\rho^\geo$.

By Theorem~\ref{obstruction-theorem}, the
two distinct elements of $H^2(M,\partial M;\FF_2) = \{0, \Obs_0(\rho^\geo)\}$
correspond to two connected components of the space $\SAS_0(\triang)$ of peripherally
trivial \Sval\ angle structures.
According to the Ptolemy database, the representations with trivial obstruction class
include a single \realrep\ representation and a complex conjugate pair of \nonrealrep\ representations of
approximate volume $\pm0.9427$.
However, our theory predicts that these three representations \emph{do not} contribute to the asymptotics
of the meromorphic \ind.

So we shall now focus exclusively on the component $\Omega=\Phi_0^{-1}\bigl(\Obs_0(\rho^\geo)\bigr)$.
A list of boundary-parabolic representations with the geometric obstruction class is given in the
Ptolemy database \cite{ptolemy} under the ``obstruction index $1$''.
In particular, the component~$\Omega$ contains \Sval\ angle structures corresponding to shape
parameter solutions for the following representations:
\begin{enumerate}[(i)]
    \item The complex conjugate pair $\{\rho^\geo, \overline{\rho^\geo}\}$ given by the holonomy
    representations of the complete hyperbolic structure on $M$ with the two possible orientations,
    \item A $\PSLR$-representation $[\rho_1]\in\mathcal{X}^\geo_\RR$ with Chern-Simons invariant $\approx0.313$,
    \item A $\PSLR$-representation $[\rho_2]\in\mathcal{X}^\geo_\RR$ with Chern-Simons invariant $\approx0.681$.
\end{enumerate}

The geometric solution $z^\geo$ of Thurston's edge consistency and completeness equations
on $\triang$ induces a strict angle structure $\alpha^\geo$ given by
\begin{align*}
    \alpha^\geo \approx
    (&0.956349, 0.899471, 1.285772;\\
     &0.899471, 1.285772, 0.956349;\\
     &0.658845, 1.855821, 0.626927)
\end{align*}
in SnapPy's default ordering of normal quads for the census triangulation \texttt{m011}.
Using the above values and the formula \eqref{tau-recalled}, we find the corresponding
$\tau$-invariant to be
\[
    \tau(z^\geo) \approx 0.188233922367388.
\]
Moreover, $H_1(\hat{M};\FF_2)=0$, so from (\ref{intro:contribution-of-smooth-crit}) the total contribution of the hyperbolic
holonomy representation and its mirror image is expected to take the form
\begin{align*}
      \contribnonreal([\rho^\geo],\kappa)
    + \contribnonreal([\overline{\rho^\geo}],\kappa)
    &= \tau(z^\geo)\cdot 2\cos\left(\kappa\Volhyp(M) + \tfrac{\pi}{4}\right) \sqrt{2\pi\kappa}\\
    &\approx 0.9436649441\sqrt{\kappa}\cos\left(\kappa\Volhyp(M) + \tfrac{\pi}{4}\right),
\end{align*}
where $\Volhyp(M) \approx 2.781833912396$ is the 
volume of the complete hyperbolic structure on $M$.

The two \realrep\ representations $\rho_1$, $\rho_2$ are detected by real solutions to
Thurston's equations on $\triang$ which project, under \eqref{associated-S1-structure}, to
the following \Ztaut\ angle structures:
\[
    \omega_1  = ( 1,  1, -1;\;
                  1, -1,  1;\;
                  1,  1, -1),\quad
    \omega_2  = (-1,  1,  1;\;
                  1,  1, -1;\;
                  1, -1,  1).
\]
As explained in detail in Section~\ref{sec:Mellin-Barnes} below, the slopes of the corresponding
linear terms \eqref{real-contribution} can be expressed by the following two-variable
contour integrals:
\begin{align*}
    \MB(\triang,\omega_1) &=
    \frac{1}{(2\pi i)^2} \int_{(i\RR)^2}
        \Beta^2\bigl(\tfrac{3}{16}-s_1, \tfrac{9}{16} + 2s_1-s_2\bigr)
        \Beta\bigl(\tfrac18 + 2s_2, \tfrac34 + s_1 - s_2\bigr)
    \,ds_1 ds_2,
    \\
    \MB(\triang,\omega_2) &=
    \frac{1}{(2\pi i)^2} \int_{(i\RR)^2}
        \Beta^2\bigl(\tfrac{9}{16} + 2s_1 - s_2, \tfrac14 - s_1 + s_2\bigr)
        \Beta\bigl(\tfrac18 + 2s_2, \tfrac18 - s_1 - s_2\bigr)
    \,ds_1 ds_2,
\end{align*}
where $\Beta(x,y)=\Gamma(x)\Gamma(y)/\Gamma(x+y)$ is the Euler beta function.
Note that the contour of integration is symmetric with respect to simultaneous
conjugation of both variables, so the above integrals are real.
By numerical quadrature, we obtain the approximate values
$\MB(\triang,\omega_1)\approx 0.6616$ and $\MB(\triang,\omega_2)\approx 0.5050$,
although it is difficult to establish error estimates owing to the expected conditional
convergence of such integrals.

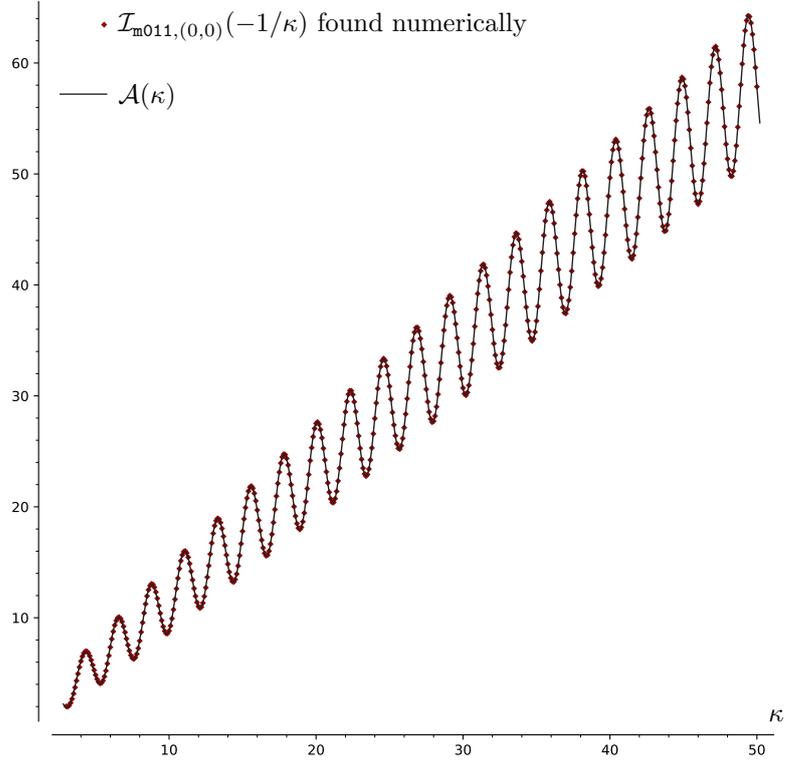
\begin{figure}%
    \centering
    \begin{tikzpicture}
        \node[anchor=south west,inner sep=0] (image) at (0,0)
            {\includegraphics[width=0.7\columnwidth]{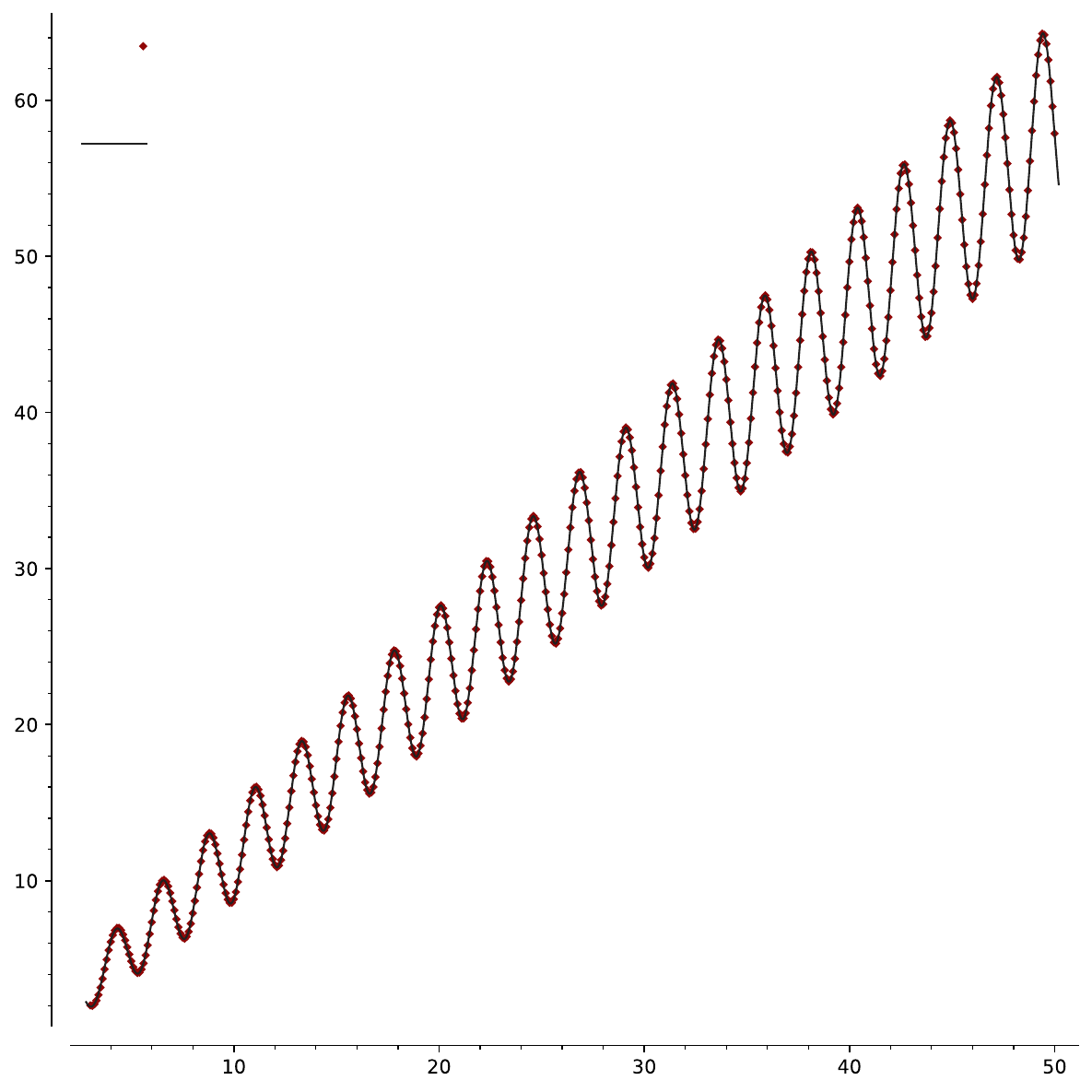}};
        \begin{scope}[x={(image.south east)},y={(image.north west)},
                every node/.style={inner sep=0, outer sep=0}]
            \node[anchor=south east] at (1, 0.06)  {$\kappa$};
            \node[anchor=west] at (0.15, 0.955) {$\I{\mathtt{m011}}{0}{0}(-1/\kappa)$ found numerically};
            \node[anchor=west] at (0.15, 0.865) {$\mathcal{A}(\kappa)$};
        \end{scope}
    \end{tikzpicture}
    \caption{\label{asym-m011}%
    Comparison between the values $\I{M}{0}{0}(-1/\kappa)$ of the meromorphic \ind\ for the
    SnapPy census manifold $M=\mathtt{m011}$ with the approximation $\mathcal{A}(\kappa)$
    of \eqref{predicted-approx-m011}.
    }
\end{figure}

Gathering the contributions of all conjugacy classes in $\mathcal{X}^\geo$,
our conjectured asymptotic approximation from (\ref{real-contribution}) and (\ref{intro:contribution-of-smooth-crit}) becomes:
\begin{align}\label{predicted-approx-m011}
    \mathcal{A}(\kappa) &=  \contribreal([\rho_1],\kappa)
              + \contribreal([\rho_2],\kappa) + \contribnonreal([\rho^\geo],\kappa)
              + \contribnonreal([\overline{\rho^\geo}],\kappa)
              \\
    &\approx 1.166578\,\kappa
             + 0.9436649441\sqrt{\kappa}\cos\left(\kappa\Volhyp(M) + \tfrac{\pi}{4}\right).\notag
\end{align}
In order to compare the above approximation with the values of the meromorphic \ind\ at
$\hbar=-1/\kappa$, we used the state-integral \eqref{rewritten-integral} to
express $\I{\mathtt{m011}}{0}{0}(-1/\kappa)$ and evaluated the integral numerically
for $\kappa\in\{3, 3.1, 3.2,\dotsc,50\}$.
We observe very good agreement between the two quantities, as illustrated in Figure~\ref{asym-m011}.

We performed numerical quadrature of the state-integral of the meromorphic \ind\
for many other census triangulations, obtaining in each case
a very good agreement with the asymptotic approximations conjectured here.
Section~\ref{sec:numerical} contains more examples of such computations
and a description of the computer program that we created for this purpose.
\par\medskip\noindent
\textbf{Acknowledgements.}
The authors would like to thank Rinat Kashaev and Blake Dadd for helpful
conversations.  Blake Dadd was involved at the beginning of this collaboration
in the second half of 2017.  This research has been partially supported by
grants DP160104502 and DP190102363 from the Australian Research Council.  The
project also received support through the AcRF Tier~1 grant RG~32/17 from the
Singapore Ministry of Education.  Rafa\l\ Siejakowski was supported by grant
\char`#2018/12483-0 of the S\~ao Paulo Research Foundation (FAPESP).  Our work
has benefited from the support of Nanyang Technological University and the
University of Melbourne.  We would like to thank these institutions for their
hospitality.

We would also like to thank Stavros Garoufalidis and Campbell Wheeler for discussing their recent work \cite{garoufalidis-wheeler}, which studies asymptotics of a related, but different invariant. (See Example~\ref{appendix-fig-eight-example} from Appendix \ref{recover-qseries-appendix}.)

%==============================================================================
% Section 2: Preliminaries
%==============================================================================
\section{Preliminaries on ideal triangulations and angle structures}
\label{sec:preliminaries}
Throughout the paper, we assume that $M$ is an orientable non-compact
\nword{3}{manifold} homeomorphic to the interior of a compact connected
\nword{3}{manifold}~$\overline{M}$ whose boundary~$\del\overline{M}$
consists of $k\geq1$ tori.
Whenever there is no risk of confusion, we will omit the bars and write simply
$\del M$ for this toroidal boundary.

%------------------------------------------------------------------------------
\subsection{Notations and conventions for ideal triangulations}
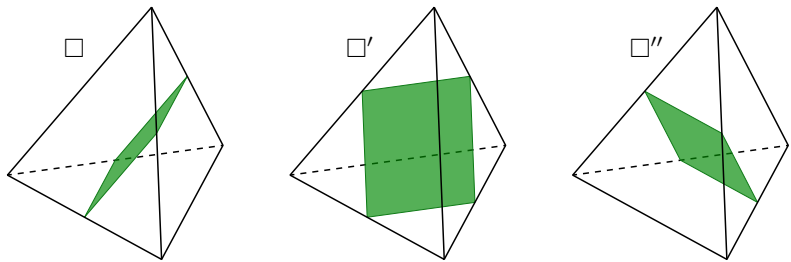
\begin{figure}[tbh]
    \centering
    \begin{tikzpicture}
        \node[anchor=south west,inner sep=0] (image) at (0,0)
            {\includegraphics[width=0.7\columnwidth]{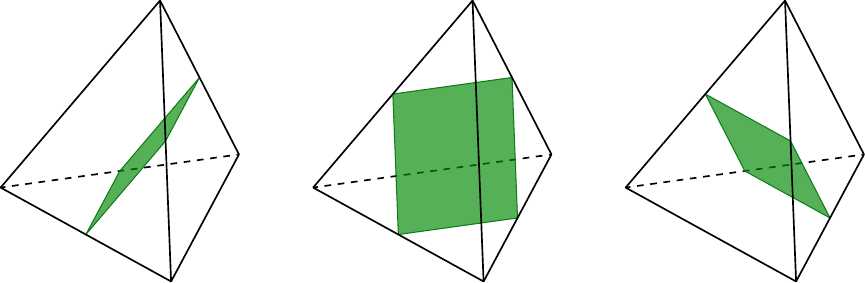}};
        \begin{scope}[x={(image.south east)},y={(image.north west)},
                every node/.style={inner sep=0, outer sep=0}]
            \node[anchor=south east] at (0.1, 0.8)  {$\square$};
            \node[anchor=south east] at (0.47, 0.8)  {$\square'$};
            \node[anchor=south east] at (0.84, 0.8) {$\square''$};
        \end{scope}
    \end{tikzpicture}
    \caption{%
    The cyclic ordering of normal quadrilaterals $\square\to\square'\to\square''\to\square$.
    }\label{fig:quads}
\end{figure}
In what follows, we understand an \emph{ideal triangulation} of $M$ to be a
collection $\{\tet_1,\dotsc,\tet_N\}$ of $N\geq2$ distinct, oriented topological ideal tetrahedra
equipped with a system of
face pairings, with the property that the resulting quotient space is homeomorphic to $M$.
The face pairings are orientation-reversing
homeomorphisms identifying pairs of distinct faces of the tetrahedra~$\tet_1,\dotsc,\tet_N$
in such a way that no face remains unglued.
In particular, we allow identifications between different faces of the same tetrahedron, as well
as multiple face gluings between the same tetrahedra.

When $\triang$ is an ideal triangulation,
we denote by $Q(\triang)$ the set of \emph{normal quadrilateral types} in $\triang$,
and by $\edges(\triang)$ the set of edges of $\triang$.
As there are three normal quadrilateral types in each tetrahedron,
we have $|Q(\triang)|=3N$.
Since all connected components of $\partial M$ are tori,
a simple Euler characteristic calculation \cite{youngchoi} shows that $|\edges(\triang)|=N$.

We will use the term \emph{edges} to refer to the elements of $\edges(\triang)$
and the term \emph{tetrahedral edges} when speaking about the edges of the unglued tetrahedra.
The gluing pattern of $\triang$ induces identifications among the tetrahedral edges,
so that each edge~$E\in\edges(\triang)$
can be thought of as an equivalence class of tetrahedral edges.

The three normal quads
contained in an oriented tetrahedron~$\tet$ admit a cyclic ordering
$\square\to\square'\to\square''\to\square$ depicted in Figure~\ref{fig:quads};
note that reversing the orientation of $\tet$ reverses this ordering.
We shall write $\square\subset\tet$ when
$\square$ is a normal quadrilateral type in a tetrahedron~$\tet$.

%--------------------------------------------------------------------------------------------------
\subsection{Thurston's gluing equations}\label{sec:gluing-equations}
Every normal quad~$\square\subset\tet$ determines a pair of opposite tetrahedral edges of $\tet$,
namely the two edges \emph{facing} the quadrilateral~$\square$ (and disjoint from it).
For any $E\in\edges(\triang)$ and any $\square\in Q(\triang)$,
let $G(E,\square)\in\{0,1,2\}$ be
the number of tetrahedral edges facing the quadrilateral~$\square$ in the edge class of $E$.

Following Neumann and Zagier \cite{neumann-zagier}, it is common to organise the incidence
numbers $G(E,\square)$ into matrices.
To this end, we fix a numbering of both the tetrahedra and the edges with integers $\{1,\dotsc,N\}$,
so that $\edges(\triang)=\{E_1,\dotsc,E_N\}$.
For every $j\in\{1,\dotsc,N\}$, we choose a distinguished normal
quadrilateral~$\square_j\subset\tet_j$ in the $j$th tetrahedron.
With these choices, the \emph{gluing matrix} of $\triang$ is defined to be the matrix
\begin{align}\label{def-gluing-matrix}
    \mathbf{G} &= [\,G\;|\;G'\;|\;G''\,] \in\Mat_{N\times 3N}(\ZZ),\ \text{where for any}\
    m,n\in\{1,\dotsc,N\}\ \text{we have}\\
    G_{m,n} &= G(E_m,\square_n),\quad
    G'_{m,n} = G(E_m,\square'_n),\quad
    G''_{m,n} = G(E_m,\square''_n).\notag
\end{align}
In the above formula, and throughout the paper, the symbol $\Mat_{m\times n}(\mathcal{R})$
denotes the space of $m\times n$ matrices with entries in a ring $\mathcal{R}$.

In what follows, the term \emph{peripheral curve} stands for
an oriented, homotopically non-trivial simple closed curve contained in $\del M$.
By placing a peripheral curve~$\gamma$ in normal position with respect to
the triangulation~$\triang$, we may
associate to $\gamma$ the gluing equations coefficients $G(\gamma, \square)$, $\square\in Q(\triang)$
in the usual way \cite{thurston-notes,neumann-zagier}.
By numbering the connected components of $\del M=T_1\sqcup\dotsb\sqcup T_k$,
and fixing, for every $1\leq j\leq k$, a pair of peripheral curves~$\mu_j, \lambda_j\subset T_j$
generating $H_1(T_j,\ZZ)$, we obtain the \emph{peripheral gluing matrix}~$\mathbf{G_\del}$,
\begin{align}\label{peripheral-gluing-matrix}
    \mathbf{G_\del} &= [\,G_\del\;|\;G'_\del\;|\;G''_\del\,] \in\Mat_{2k\times 3N}(\ZZ),\ \text{where}\\
        (G_\del)_{2j-1,n} &= G(\mu_j,\square_n),\quad\notag
        (G_\del)_{2j,n} = G(\lambda_j,\square_n),\quad 1\leq j \leq k;
\end{align}
$G'_\del$, $G''_\del$ are defined analogously.
With the above notations, Thurston's hyperbolic gluing equations,
as presented by Neumann and Zagier~\cite{neumann-zagier}, have the form
\begin{equation}\label{edge-consistency-equations}
    \prod_{\square\in Q(\triang)} {z(\square)}^{G(E_m,\square)}
        =\prod_{n=1}^N {z_n}^{G_{m,n}}{z'_n}^{G'_{m,n}}{z''_n}^{G''_{m,n}}
                = 1\ \text{for all}\ m,
\end{equation}
where the shape parameter variables $z$ satisfy $z_n = z(\square_n)$ and
\begin{equation}\label{z-relations}
    z(\square') = \frac{1}{1-z(\square)}\ \text{for all}\
                  \square\in Q(\triang),
\end{equation}
i.e., $z'_n = 1/(1-z_n)$ and $z''_n = (z_n-1)/z_n$ for all $n$.

The \emph{gluing variety}~$\gluvar$ of the triangulation $\triang$
is defined to be the affine algebraic set cut out from $\CC^{Q(\triang)}$ by the edge
consistency equations~\eqref{edge-consistency-equations} and \eqref{z-relations}.

After taking logarithms, we may rewrite equation~\eqref{edge-consistency-equations}
in matrix form as
\begin{equation}\label{edge-consistency-general}
    GZ + G'Z' + G''Z'' \in 2\pi i \ZZ^{\edges(\triang)},
\end{equation}
where $Z=(\log z_1, \dotsc,\log z_N)$, $Z'=(\log z'_1,\dotsc,\log z'_N)$,
and $Z''=(\log z''_1,\dotsc,\log z''_N)$.
Of particular importance is the geometric case, where we impose
the equations
\begin{equation}\label{edge-consistency-logarithmic}
    GZ + G'Z' + G''Z'' = 2\pi i\, (1,1,\dotsc,1)^\transp,
\end{equation}
which can be used in the search for hyperbolic structures, cf.~\cite[eq.~4.2.2]{thurston-notes}.

An \emph{algebraic solution} $z=\exp(Z)$ is any solution of
\eqref{edge-consistency-equations} subject to \eqref{z-relations}
with $z(\square)\in\CC\setminus\{0,1\}$ for all $\square\in Q(\triang)$.
A \emph{geometric solution} of the gluing equations is a solution
$z$ of \eqref{edge-consistency-equations}--\eqref{z-relations} in which
all shapes additionally satisfy
$\Imag z(\square)>0$ and \eqref{edge-consistency-logarithmic} holds with the principal
branch of the logarithms.
In this case, the principal arguments $\Arg z(\square)\in(0,\pi)$ can be interpreted as dihedral
angles of positively oriented hyperbolic ideal tetrahedra.
By gluing such geometric tetrahedra via hyperbolic isometries, according
to the face identification pattern of $\triang$, one
obtains a finite-volume hyperbolic metric on $M$~\cite{thurston-notes}.

For a fixed peripheral curve $\gamma\subset\del M$,
let $\log\hol(\gamma)$ denote its \emph{log-holonomy}, defined by the formula
\begin{equation}\label{loghol}
    \log\hol(\gamma) = \sum_{\square\in Q(\triang)} G(\gamma, \square) \log z(\square).
\end{equation}
The hyperbolic structure defined by a geometric solution~$z$ is complete iff $\log\hol(\gamma)$
vanishes for every $\gamma$.
It is sufficient to impose this condition
on the chosen generator curves $(\mu_j, \lambda_j)$, $1\leq j\leq k$,
which leads to the finite system of \emph{completeness equations}
\begin{equation}\label{completeness-equations}
    G_\del Z + G'_\del Z' + G''_\del Z'' = 0.
\end{equation}
In the case of a single cusp, we customarily write $(\mu,\lambda)=(\mu_1,\lambda_1)$
and use the short notation from Neumann and Zagier~\cite{neumann-zagier},
\begin{equation}\label{def:uv}
    u = \log\hol(\mu),\quad v = \log\hol(\lambda).
\end{equation}

More generally, Yoshida~\cite{yoshida-ideal} showed
that any algebraic solution~$z=\exp(Z)\in(\CC\setminus\{0,1\})^{Q(\triang)}$
of the gluing equations determines a conjugacy
class of a representation~$\rho\colon\pi_1(M)\to\PSLC$.
When the log-parameters~$Z$ additionally satisfy
\begin{equation}\label{completeness-equations-mod-2pi}
    G_\del Z + G'_\del Z' + G''_\del Z'' \in 2\pi i \ZZ^{2k},
\end{equation}
then $\rho$ is boundary-parabolic.

Suppose that $z\in\gluvar$ is a (positively oriented) geometric solution inducing a hyperbolic
structure of finite volume on $M$, not necessarily complete.
Choi~\cite{youngchoi} 
showed that 
a neighbourhood of $z$ is parametrised bi-holomorphically
by the local coordinate $u = \log\hol(\mu)$ (or, alternatively, $v=\log\hol(\lambda)$).

%--------------------------------------------------------------------------------------------------
\subsection{Angle structures on ideal triangulations}\label{angle_structures_section}
We recall the definitions of several different types of angle structures on ideal triangulations.

\begin{definition}
\label{angle structures}
    \begin{enumerate}[(i)]
    \item
    A \emph{generalised angle structure} on the ideal triangulation $\triang$
    is a function $\alpha : Q(\triang)\to\RR$ satisfying
    \begin{align}
        \alpha(\square)+\alpha(\square')+\alpha(\square'') \label{AS:tetrah_sum=pi}
            &= \pi\ \text{for every}\ \square\in Q(\triang),\ \text{and}\\
        \sum_{\square\in Q(\triang)} G(E,\square)\,\alpha(\square) \label{AS:edge_sum=2pi}
            &= 2\pi\ \text{for every}\ E\in\edges(\triang).
    \end{align}
    The set of all generalised angle structures on $\triang$ is denoted by $\Areal(\triang)$.
    \item
    A \emph{strict angle structure} is a generalised angle structure satisfying
    additionally $0 < \alpha(\square) < \pi$ for all $\square\in Q(\triang)$.
    The set of all strict angle structures on $\triang$ is denoted by $\Astrict(\triang)$.
    \item
    An \emph{\Sval\ angle structure} on the triangulation $\triang$ is a function
    $\omega : Q(\triang)\to S^1$ satisfying
    \begin{align}
        \omega(\square)\omega(\square')\omega(\square'') \label{SAS:tetrah_prod=-1}
            &= -1\ \text{for every}\ \square\in Q(\triang),\ \text{and}\\
        \prod_{\square\in Q(\triang)} \omega(\square)^{G(E,\square)} \label{SAS:edge_prod=1}
            &= 1\ \text{for every}\ E\in\edges(\triang).
    \end{align}
    The set of all \Sval\ angle structures on $\triang$ is denoted by
    $\SAS(\triang)$.
    \item
    A \emph{\nword{\ZZ_2}{valued} angle structure} is an \Sval\ angle structure $\omega$
    satisfying $\omega(\square)\in\{+1,-1\}$ for all $\square\in Q(\triang)$.
    \item
    A \emph{\Ztaut\ angle structure} is a \nword{\ZZ_2}{valued} angle structure which takes
    the value $+1$ on exactly two normal quads in each tetrahedron and $-1$ on the remaining quad.
    \end{enumerate}
\end{definition}
Since the columns of the matrices $\mathbf{G}$, $\mathbf{G_\del}$ correspond to normal quadrilaterals,
we may simply write \eqref{AS:edge_sum=2pi} as $\mathbf{G}\alpha = (2\pi,\dotsc,2\pi)^\transp$ with
the understanding that the normal quads are ordered in accordance with the ordering of tetrahedra, i.e.,
\begin{equation}\label{quad-ordering}
    Q(\triang) = \{ \square_1, \square_2, \dotsc, \square_N;\;
            \square'_1, \square'_2, \dotsc, \square'_N;\;
            \square''_1, \square''_2, \dotsc, \square''_N\}.
\end{equation}

The space $\SAS(\triang)$ of \Sval\ angle structures was first introduced and studied
by Feng Luo in \cite{luo-volume}.
In particular, Luo showed that $\SAS(\triang)$ is a smooth closed manifold of dimension $N+k$
and remarked that it is disconnected in general.
Our Theorem~\ref{obstruction-theorem} classifies the connected components of $\SAS(\triang)$
and relates them to cohomological obstructions to lifting $\PSLC$ representations of the fundamental
group to $\SLC$. The structure of these spaces is described in more detail in \cite{hodgson-kricker-siejakowski}.

Observe that the formula~$\omega(\square)=e^{i\alpha(\square)}$ defines
an ``exponential map''~$\exp : \Areal(\triang)\to\SAS(\triang)$.
As $\Areal(\triang)\subset\RR^{Q(\triang)}$ is a real affine space (and is therefore connected),
the image of the exponential map lies entirely in a single connected component of $\SAS(\triang)$,
which we shall call the \emph{geometric component}
and denote $\SASdistinguished(\triang)$.
More generally, all of $\SAS(\triang)$ can be viewed as the set of exponentiated
real-valued \emph{pseudo-angles}, i.e., solutions of the angle
equations~\eqref{AS:tetrah_sum=pi}, \eqref{AS:edge_sum=2pi} modulo $2\pi\ZZ$.

When $\alpha \in \Areal(\triang)$ is a generalised angle structure
and $\gamma$ is a peripheral curve, the \emph{angle-holonomy} of
$\alpha$ along $\gamma$ is given by
\begin{equation}\label{def-angle-holonomy}
    \ahol{\alpha}(\gamma)= \sum_{\square\in Q(\triang)} G(\gamma, \square) \alpha(\square)\in\RR.
\end{equation}
We say that $\alpha$ has \emph{vanishing peripheral angle-holonomy},
or that $\alpha$ is \emph{peripherally trivial},
if $\ahol{\alpha}(\mu_j)=\ahol{\alpha}(\lambda_j)=0$ for all $j\in\{1,\dotsc,k\}$,
i.e., when
\(
    \mathbf{G_\del} \alpha = 0
\).
Similarly, we may define the multiplicative angle-holonomy of an \Sval\ angle structure
$\omega$ along $\gamma$ by the formula
\begin{equation}\label{peripheral_S1_holonomy_defn}
    \ahol{\omega}(\gamma) =
        \prod_{\square\in Q(\triang)} {\omega(\square)}^{G(\gamma,\square)}
    \in S^1.
\end{equation}

We shall often need to consider spaces of angle structures with
prescribed peripheral angle-holonomies.
For simplicity, suppose that $M$ has only one torus boundary component,
which we orient in accordance with the orientation of $M$, using the convention of
``outward-pointing normal vector in the last position''.
We may further assume that the peripheral curves~$\mu$, $\lambda$ have
intersection number~$\intersection(\mu,\lambda)=1$ with respect to this orientation. 
%cdh: this is opposite to usual convention for meridian longitude of knot complement
For any $(x,y)\in\RR^2$, we define
\begin{equation}\label{def-Axy}
    \Areal_{(x,y)}(\triang) = \{\alpha\in\Areal(\triang)\ |\
        \bigl(\ahol{\alpha}(\mu), \ahol{\alpha}(\lambda)\bigr)=(x,y)\}.
\end{equation}
Likewise, for two elements $\xi,\eta\in S^1$ we shall consider
\Sval\ angle structures with prescribed multiplicative angle-holonomy,
\begin{equation}\label{def-SASxy}
    \SAS_{(\xi,\eta)}(\triang) = \{\omega\in\SAS(\triang)\ |\
    \bigl(\ahol{\omega}(\mu), \ahol{\omega}(\lambda)\bigr)=(\xi,\eta)\}.
\end{equation}
An important special case is presented by the spaces of \emph{peripherally trivial}
angle structures.
In order to lighten the notation, we set
\begin{equation*}
    \Areal_0(\triang) := \Areal_{(0,0)}(\triang),
    \qquad
    \SAS_0(\triang) := \SAS_{(1,1)}(\triang).
\end{equation*}
Note that the space $\SAS_0(\triang)$ is independent of the choice of the curves~$\mu$ and~$\lambda$.

\begin{definition}\label{def:Zvals}
\begin{enumerate}[(i)]
\item
An \emph{integer-valued angle structure} on the triangulation $\triang$ is a
vector $f\in\ZZ^{Q(\triang)}$ such that
$\pi f \in \Areal(\triang)$.
\item
An integer-valued angle structure~$f$ is \emph{peripherally trivial} if $\pi f\in\Areal_0(\triang)$.
\end{enumerate}
\end{definition}
By a result of Neumann~\cite[Theorem~2]{neumann1990}, an ideal triangulation
of a manifold with toroidal boundary always admits a peripherally trivial integer-valued angle structure.
As a consequence, using the exponential map, the set $\SAS_0(\triang)$ is always non-empty.
We refer to \cite{luo-volume,garoufalidis-hodgson-hoffman-rubinstein, hodgson-kricker-siejakowski} for more details on angle
structures of various types.

%--------------------------------------------------------------------------------------------------
\subsection{Tangential angle structures}\label{subsec:TAS}
\begin{definition}
The space of \emph{tangential angle structures}~$\TAS(\triang)\subset\RR^{Q(\triang)}$ is the vector
space consisting of vectors~$w$ satisfying the equations
\begin{align}
    w(\square) + w(\square') + w(\square'')&=0,\quad \square\in Q(\triang),\label{TAS:tet_sum=0}\\
    \label{TAS:edge_sum=0}
    \mathbf{G}\,w&=0.
\end{align}
Moreover, define $\TAS_0(\triang)\subset\TAS(\triang)$ to be the subspace cut out by the equations
\(
    \mathbf{G_\del}\,w = 0
\).
\end{definition}
Since the above equations are homogenous versions of the angle structure
equations \eqref{AS:tetrah_sum=pi}--\eqref{AS:edge_sum=2pi},
$\TAS(\triang)$ is canonically identified with the tangent space to
the affine subspace $\Areal(\triang)\subset\RR^{Q(\triang)}$ at every one of its points,
and by using the exponential map, it is also identified
with the tangent space to $\SAS(\triang)$.
By the same reasoning, $\TAS_0(\triang)$ is the tangent space to the subspaces
$\Areal_{(x,\,y)}(\triang)$ and $\SAS_{(\xi,\,\eta)}(\triang)$ with fixed
peripheral angle-holonomy; see~\cite[p.~307]{luo-volume} for more details.

Based on the treatment of Futer and Gu\'eritaud~\cite[\S4]{futer-gueritaud},
we introduce explicit bases for $\TAS(\triang)$ and $\TAS_0(\triang)$.
For every $j\in\{1,\dotsc,N\}$, we define the \emph{leading--trailing deformation about
the edge}~$E_j\in\edges(\triang)$
to be the vector $l_j\in\ZZ^{Q(\triang)}$ given by the formula
\begin{equation}\label{definition-leading-trailing}
    l_j(\square) = G(E_j,\square'') - G(E_j,\square').
\end{equation}
It is easy to see that $l_j\in\TAS_0(\triang)$;
in fact, the vectors $\{l_1,\dotsc,l_N\}$ span $\TAS_0(\triang)$,
cf.~\cite[Proposition~4.6]{futer-gueritaud}.
Using the matrix notation of \eqref{def-gluing-matrix}, we may define
\begin{equation}\label{definition-L}
    \begin{aligned}
        \mathbf{L} &= [\,L\;|\;L'\;|\;L''\,] \in\Mat_{N\times 3N}(\ZZ),\ \text{where}\\
        L &= G''-G',\quad L' = G-G'', \quad L'' = G'-G.
    \end{aligned}
\end{equation}
With this notation, the vector $l_j$ can be understood as the \nword{j}{th} row of $\mathbf{L}$.
Note that we always have $-l_N = l_1+\dotsb+l_{N-1}$.
In the case of a single boundary component,
this is the only nontrivial linear relation
among the vectors $l_j$, so that
after removing the last vector ($j=N$), the set~$\{l_1,\dotsc,l_{N-1}\}$ becomes
a basis of $\TAS_0(\triang)$, cf. \cite[\S4]{futer-gueritaud}.
In this case, it is convenient to define
\begin{equation}\label{definition-L_star}
    \begin{aligned}
    \mathbf{L_*} &= [\,L_*\;|\;L'_*\;|\;L''_*\,] \in\Mat_{(N-1)\times 3N}(\ZZ),\\
    [\mathbf{L_*}]_{j,\square} &= [\mathbf{L}]_{j,\square},\quad 1\leq j\leq N-1,
    \end{aligned}
\end{equation}
so that the rows of $\mathbf{L_*}$ form an explicit ordered basis of $\TAS_0(\triang)$.

More generally, Futer and Gu\'eritaud~\cite{futer-gueritaud} consider
leading-trailing deformations along oriented peripheral curves in general position with
respect to $\triang$.
When $\gamma$ is such peripheral curve, we may define
\begin{equation}\label{peripheral-LTD}
    l_\gamma(\square) = G(\gamma,\square'') - G(\gamma,\square').
\end{equation}
It suffices for our purposes to consider fixed meridian-longitude
pairs~$(\mu_j,\lambda_j)$, $1\leq j \leq k$ and the associated
leading-trailing deformation vectors $\{l_{\mu_j}, l_{\lambda_j}\}_{j=1}^k$.
These vectors arise as the rows of the matrix
\begin{equation}\label{def-L-peripheral}
    \begin{aligned}
    \mathbf{L}_\del &= [\,L_\del\;|\;L'_\del\;|\;L''_\del\,]\in\Mat_{2k\times 3N}(\ZZ),\\
    L_\del &= G''_\del - G'_\del,\quad
    L'_\del = G_\del - G''_\del,\quad
    L''_\del = G'_\del - G_\del.
    \end{aligned}
\end{equation}
In this notation, we have $\TAS_0(\triang)=\mathbf{L}\mkern-2mu^\transp(\RR^N)=\mathbf{L}_*^\transp(\RR^{N-1})$ and
$\TAS(\triang) = \TAS_0(\triang) + \mathbf{L}^\transp_\del(\RR^{2k})$, where $^\transp$ denotes transpose; cf.~\cite[Proposition~4.6]{futer-gueritaud}.
Moreover, the dimension of $\TAS_0(\triang)$ is always $N-k$,
whereas the dimension of $\TAS(\triang)$ is $N+k$.
We refer to \cite[Proposition~3.2]{futer-gueritaud} and \cite[Proposition~2.6]{luo-volume} for more
details and proofs.

\begin{remark}\label{TAS-and-connected-components}
As a consequence of the properties of leading-trailing deformations discussed above,
two \Sval\ angle structures~$\omega_1$ and $\omega_2$ lie in the same connected component of
$\SAS_0(\triang)$ if and only if there exists a
vector~$r\in\Imag\mathbf{L}^\transp=\mathbf{L}^\transp(\RR^N)$ such that
\begin{equation}\label{TAS-connected-SAS}
    \omega_2(\square) = \omega_1(\square) e^{ir(\square)}\ \text{for all}\ \square\in Q(\triang).
\end{equation}
Likewise, $\omega_1$ and $\omega_2$ belong to the same connected component of $\SAS(\triang)$
if and only if a vector~$r$ satisfying \eqref{TAS-connected-SAS}
can be found in the larger space~$\Imag\mathbf{L}^\transp+\Imag\mathbf{L}_\del^\transp$.
\end{remark}

%--------------------------------------------------------------------------------------------------
\subsection{Differentiating the peripheral angle-holonomies along leading-trailing deformations.}
\label{shear-section}

Here we briefly recall the effect of the leading-trailing deformations along peripheral curves
on the peripheral angle-holonomies of angle structures
as described by Lemma~4.4 of Futer--Gu\'eritaud~\cite{futer-gueritaud}.
For a peripheral curve $\gamma$, we shall treat $l_\gamma$ as a tangent vector
to $\Areal(\triang)$ at any of its points and denote by $\nabla_{l_\gamma}$ the
directional derivative in the direction of $l_\gamma$.
Then, for any two peripheral curves $\gamma$, $\theta$ and any $\alpha \in \Areal(\triang)$, we have
\begin{align}\label{futer-gueritaud-peripheral}
    \nabla_{l_\gamma}\ahol{\alpha}(\theta) &= 2\intersection(\gamma,\theta),\\
    \nabla_{l_m}\ahol{\alpha}(\theta) &= 0,\quad 1\leq m \leq N,\label{futer-gueritaud-nonperipheral}
\end{align}
where $\intersection(\,\cdot\,,\,\cdot\,)$ is the intersection pairing on $\partial\overline{M}$
and where $l_m$ is the leading trailing-deformation about the $m$-th edge, $1\leq m \leq N$.

%--------------------------------------------------------------------------------------------------
\subsection{Variational approach to the gluing equations}
The variational approach to solving Thurston's equations~\eqref{edge-consistency-equations}
was introduced by Casson and Rivin, with an excellent exposition provided
by Futer and Gu\'eritaud in~\cite{futer-gueritaud}.
Here, we give a short summary based on the treatment due to Luo~\cite{luo-volume}, in terms of
\Sval\ angle structures.

Denote by $\Lob$ the \emph{Lobachevsky function},
\begin{equation*}
    \Lob : \RR\to\RR,\quad
    \Lob(\theta) = -\int_0^\theta \log|2\sin x|\,dx.
\end{equation*}
The function~$\Lob$ is a variant of Clausen's function~$\mathrm{Cl}_2$ first introduced
in \cite{clausen}
and is related to Euler's dilogarithm function~$\Li(z)$ by the identity
$2\Lob(\theta)=\Imag\Li(e^{2i\theta})$, cf.~\cite{neumann-zagier}.
In particular, Lobachevsky's function is continuous for all $\theta\in\RR$
and differentiable whenever $\theta\not\in\pi\ZZ$.
Since $\Lob$ is periodic with period~$\pi$, the function
\begin{equation}\label{circle-Lob}
\circLob : S^1\to\RR,\quad
\circLob(e^{i\theta}) = \Lob(\theta)
\end{equation}
is well-defined and continuous on the unit circle.

As shown by Milnor in~\cite[Chapter~7]{thurston-notes},
the volume of an ideal hyperbolic tetrahedron with dihedral angles
$\alpha$, $\alpha'$, $\alpha''$ is given by $\Lob(\alpha)+\Lob(\alpha')+\Lob(\alpha'')$.
Following Luo~\cite{luo-volume}, we define the \emph{volume function}
on the space of \Sval\ angle structures of a triangulation~$\triang$ by
\begin{equation}\label{volume-on-SAS}
    \Vol : \SAS(\triang)\to\RR,\quad
    \Vol(\omega) = \sum_{\square\in Q(\triang)}\circLob\bigl(\omega(\square)\bigr).
\end{equation}
The above volume function may fail to be differentiable when $\omega(\square)=\pm1$
for some $\square\in Q(\triang)$.

We are now going to review the relationship between the smooth critical points of $\Vol$ 
and complex solutions to Thurston's gluing equations. 
For simplicity, we limit ourselves to the case of a single torus boundary component
and we choose a meridian-longitude pair~$(\mu,\lambda)$.
For a fixed peripheral angle-holonomy~$(x,y)\in S^1\times S^1$,
consider the restriction $\Vol\bigr|_{(x,y)}$ on $\SAS_{(x,y)}(\triang)$.
The smooth 
critical points of volume satisfy the
equations~$\nabla_{l_j}\Vol\bigr|_{(x,y)}=0$ for $1\leq j \leq N$, where $\nabla_{l_j}$ denotes the directional derivative in the direction of the leading trailing deformation $l_j$ about edge $E_j$.
Since $\Lob'(\theta) = -\log|\sin\theta|-\log2$ for $\theta\not\in\pi\ZZ$,
an elementary calculation shows that these equations are equivalent
to the equations
\begin{equation}\label{shears-with-circle-valued-angles}
    \prod_{\square\in Q(\triang)}\abs{\frac{\Imag\omega(\square')}{\Imag\omega(\square'')}}^{G(E_j,\square)} = 1.
\end{equation}
Given an \Sval\ angle structure~$\omega$ satisfying
$\omega(\square)\neq\pm1$ for all $\square\in Q(\triang)$,
we may therefore define the complex numbers
\begin{equation}\label{from-omega-to-z}
    z(\square) = \left|\frac{\Imag\omega(\square')}{\Imag\omega(\square'')}\right|\omega(\square),
       \quad \square\in Q(\triang).
\end{equation}
As observed by Luo~\cite[Lemma~3.3]{luo-volume}, and as is easily seen from
\eqref{shears-with-circle-valued-angles}, when $\omega$ is a smooth critical
point of $\Vol\bigr|_{(x,y)}$, the above complex numbers satisfy the edge consistency
equations~\eqref{edge-consistency-equations}. However, additional conditions on $\omega$
are needed so that $z$ also satisfies the relations~\eqref{z-relations} required for the shape parameters of an ideal tetrahedron.

\begin{lemma}
Suppose that $\tet$ is a tetrahedron with \Sval\ angles $\omega=\omega(\square)$, $\omega'=\omega(\square')$, $\omega''=\omega(\square)''$,
none of which is real.
Use \eqref{from-omega-to-z} to define the complex numbers
\begin{equation}\label{tentatively-defining-z-from-omega}
    z = \abs{\frac{\Imag\omega'}{\Imag\omega''}}\omega,
    \quad
    z' = \abs{\frac{\Imag\omega''}{\Imag\omega}}\omega',
    \quad
    z'' = \abs{\frac{\Imag\omega}{\Imag\omega'}}\omega''.   
\end{equation}
Then these define shape parameters of an ideal hyperbolic tetrahedron
if and only if 
the imaginary parts of $\omega,\omega',\omega''$  all have the same sign.
\end{lemma}

\begin{proof} 
By definition, the numbers $z,z',z''$ are all in $\CC\setminus\RR$,
but in order for them to define cross-ratios
of a non-degenerate hyperbolic ideal tetrahedron,
the relation 
$$z' = \frac{1}{1-z}$$ must be satisfied.
In particular, the imaginary parts of both sides of this 
equality must have the same sign, which implies that $\sgn(\Imag z')=\sgn(\Imag z)$, hence
$\sgn(\Imag\omega')=\sgn(\Imag\omega)$.
Repeating this reasoning with $z'$ and $z''$ shows that in fact the imaginary parts of all three
numbers $\omega$, $\omega'$, $\omega''$ must have identical signs. 

Conversely, when the imaginary parts of $\omega$, $\omega'$, $\omega''$ are all
positive, we can write $\omega=e^{i\alpha}$, $\omega'=e^{i\alpha'}$ and $\omega''=e^{i\alpha''}$
for unique $\alpha,\alpha',\alpha''\in(0,\pi)$ satisfying $\alpha+\alpha'+\alpha''=\pi$.
Then the shapes $z$, $z'$ and $z''$, defined as the cross-ratios of the positively oriented
hyperbolic ideal tetrahedron with dihedral angles~$\alpha$, $\alpha'$ and $\alpha''$,
automatically satisfy the relations~\eqref{tentatively-defining-z-from-omega}.
The case of imaginary parts all negative follows by complex conjugation
and corresponds to a negatively oriented tetrahedron.
\end{proof}

Combining this with the previous discussion gives the following.

\begin{proposition}\label{vol_crit_point_to_sol_gluing_eqns}
Let $\triang$ be an ideal triangulation of a orientable 3-manifold $M$ whose boundary consists of a single torus.
Let $\omega$ be a critical point of the volume function $\Vol\bigr|_{(x,y)}$ restricted to $\SAS_{(x,y)}(\triang)$ such that $\omega(\square)\ne \pm 1$ for all $\square\in Q(\triang)$. Assume that in each tetrahedron,
the imaginary parts of the three angles~$\omega(\square)$,
$\omega(\square')$ and $\omega(\square'')$ all have the same sign.
Then the complex numbers $z$ defined by equation (\ref{from-omega-to-z})
give an algebraic solution of Thurston's gluing
equations whose log-holonomy $(u,v)$ satisfies $x = \exp(i\Imag u)$ and $y = \exp(i\Imag v)$.
\end{proposition}

% ==============================================================================
% Section 3: Review of the meromorphic 3D-index
% ==============================================================================
\section{Review of the meromorphic 3D-index}\label{section-index}

In this section we assume, for simplicity, that $\overline{M}$ has a single
torus boundary component with a fixed meridian-longitude pair~$(\mu,\lambda)$
of oriented curves whose homology classes form a basis of
$H_1(\partial\overline{M};\ZZ)$.
We shall now provide a brief summary of the original definition of the
meromorphic \ind\ given by Garoufalidis and Kashaev in
\cite{garoufalidis-kashaev}.
This invariant is obtained as an analytic continuation
of a state integral of Turaev--Viro type associated to an
angled ideal triangulation of $M$.
In analogy with the new formulation of the ``Teichm\"uller TQFT''
of Andersen and Kashaev~\cite{kashaev-andersen-new},
the integration domain can be identified with a certain space of \Sval\
angle structures on the triangulation.

Recall that the \emph{\nword{q}{Pochhammer} symbol} is defined by
\begin{equation}\label{q-Pochhammer}
    \qPoch{\;\cdot\;}{q}\colon\CC\to\CC,\qquad
    \qPoch{z}{q} = \prod_{n=0}^\infty (1-q^nz),
\end{equation}
where $0<|q|<1$ is a complex parameter.
Denote by $\Li(z;q)$ the \emph{\nword{q}{dilogarithm}} function
given for $|z|<1$ by the power series
\[
    \Li(z;q) = \sum_{n=1}^\infty \frac{z^n}{n(1-q^n)},
    \quad 0<|q|<1.
\]
The \nword{q}{Pochhammer} symbol is related to the \nword{q}{dilogarithm}
by the identity
\[
    \qPoch{z}{q}=\exp\bigl(-\Li(z;q)\bigr)
\]
(cf. e.g. Kirillov~\cite[\S2.5]{kirillov}, Zagier~\cite{zagier-dilog}).
The building block of the tetrahedral weights for the meromorphic \ind\ is
the complex function~$G_q$ given by
\begin{equation}\label{def:G_q}
    G_q(z) = \frac{\qPoch{-qz^{-1}}{q}}{\qPoch{z}{q}}
    = \exp\bigl(\Li(z;q)-\Li(-qz^{-1};q)\bigr),\quad 0<|q|<1.
\end{equation}
As stated in \cite[Lemma~2.2]{garoufalidis-kashaev},
the function $G_q(z)$ is analytic for all $z$
except for the simple poles at $z\in\{q^0, q^{-1}, q^{-2}, \dotsc\}$ and
an essential singularity at $z=0$.
In particular, $G_q(z)$ is always holomorphic on the punctured unit
disc $0<|z|<1$.

%-------------------------------------------------------------------------------
\subsection{State integral definition for triangulations with strict angle structures}
\label{sec:3D-index-with-strict-angles}
We choose parameters $(\aholM, \aholL)\in\RR^2$ and
assume that $\triang$ is an ideal triangulation of $M$ admitting a strict
angle structure~$\alpha$ with peripheral
angle-holonomy prescribed by the chosen parameters, i.e.,
$\bigl(\ahol{\alpha}(\mu), \ahol{\alpha}(\lambda)\bigr)=(\aholM,\aholL)$.
In the special case of $(\aholM,\aholL)=(0,0)$,
this amounts to the requirement that $\triang$ admit a strict angle structure
with vanishing peripheral angle-holonomy.
Note that such triangulations always exist when $M$ is the complement of a
hyperbolic knot in $S^3$, as shown by Rubinstein, Segerman and the first
author, cf.~\cite[Corollary~1.2]{hodgson-rubinstein-segerman}.

We start by fixing a constant $\hbar\in\CC$ with a negative real part and
setting $q=e^\hbar$, so that $0<|q|<1$.
\begin{figure}[tb]
    \centering%
    \begin{tikzpicture}[baseline=(BASE)]
        \node[anchor=south west,inner sep=0] (image) at (0,0)
            {\includegraphics[width=0.2\columnwidth]{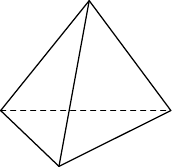}};
        \begin{scope}[x={(image.south east)},y={(image.north west)},
                every node/.style={outer sep=0}]
            \node[anchor=south] at (0.52, 1)  {$0$};
            \node (BASE) at (1,1) {};
            \node[anchor=east] at (0, 0.348)   {$1$};
            \node[anchor=north] at (0.34, 0)  {$2$};
            \node[anchor=west] at (1, 0.34)  {$3$};
            \node[anchor=south east] at (0.3, 0.65)  {$\alpha_{01}$};
            \node[anchor=west] at (0.405, 0.5)  {$\alpha_{02}$};
            \node[anchor=west] at (0.75, 0.67)  {$\alpha_{03}$};
        \end{scope}
    \end{tikzpicture}%
    \hspace{5mm}
    \parbox[t][][t]{\dimexpr0.8\columnwidth - 2cm\relax}{\caption{%
        A labelled tetrahedron~$\Delta$ equipped with dihedral angles.
    }\label{fig:tetrahedron}}%
\end{figure}
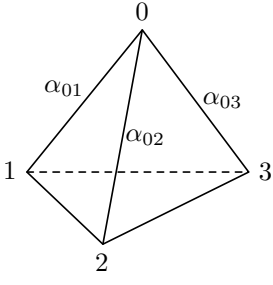
Let $\tet$ be an ideal tetrahedron of the triangulation $\triang$.
We label the ideal vertices of $\tet$ with the integers~$\{0,1,2,3\}$
and denote the tetrahedral edge of $\tet$ connecting the vertices $m$ and $n$
by $\Delta([mn])$ for all $0\leq m<n\leq 3$.
Let $\alpha_{mn}$ be the value of $\alpha$ on the normal quad facing the
tetrahedral edge $\Delta([mn])$, as depicted in Figure~\ref{fig:tetrahedron}.
With these notations, we define the Boltzmann weight of the angled tetrahedron~$\tet$ to be
\begin{align}
    B^\tet_\alpha\!%
        \left(\begin{smallmatrix} x_{01} & x_{02} & x_{03}\\
                                  x_{23} & x_{13} & x_{12}\end{smallmatrix}\right)
        = c_q
        &G_q\Bigl(e^{\alpha_{01}(i+\hbar/\pi)}\frac{x_{03}x_{12}}{x_{02}x_{13}}\Bigr)\notag\\
        &G_q\Bigl(e^{\alpha_{02}(i+\hbar/\pi)}\frac{x_{01}x_{23}}{x_{03}x_{12}}\Bigr)\label{Boltzmann-GK}\\
        &G_q\Bigl(e^{\alpha_{03}(i+\hbar/\pi)}\frac{x_{02}x_{13}}{x_{01}x_{23}}\Bigr)\notag,
        \hspace{2ex}x_{mn}\in S^1,
\end{align}
where the normalisation constant~$c_q\in\CC$ is given by the formula
\begin{equation}\label{definition-of-c(q)}
    c_q = \frac{\qPoch{q}{q}^2}{\qPoch{q^2}{q^2}}.
\end{equation}
The invariance of $B^\tet_\alpha$
under orientation-preserving tetrahedral symmetries is evident from~\eqref{Boltzmann-GK},
and an angled pentagon identity was established by
Garoufalidis and Kashaev~\cite[eq.~(50)]{garoufalidis-kashaev}.

The Boltzmann weight of the triangulation $\triang$ is defined as the product
\[
    B^\triang_\alpha(y_1,\dotsc,y_N) =
    \prod_{\tet\ \text{in}\ \triang} B^\tet_\alpha
    \Bigl|_{\textstyle
        x_{mn}=y_j\ \text{iff}\ \Delta([mn])\in E_j\in\edges(\triang)
        }
\]
in which each variable~$x_{mn}$
is replaced by the new variable~$y_j$ corresponding to the edge class~$E_j\in\edges(\triang)$
of the tetrahedral edge $\Delta([mn])$.
The state integral of the meromorphic \ind\ is obtained by
integrating over the edge variables $y_j\in S^1$
with respect to the probabilistic Haar measure on $(S^1)^{\edges(\triang)}$,
\begin{equation}\label{state-integral-formula}
    \I{\triang}{\aholM}{\aholL}(\hbar) =
    \int_{{(S^1)}^{\edges(\triang)}}
    B^\triang_\alpha(y_1,\dotsc,y_N)\;\dHaar(y_1,\dotsc,y_N).
\end{equation}
The above integral can be made more explicit if we
use the spanning set ~$\{l_1,\dotsc,l_N\}\subset \TAS_0(\triang)$ of
leading-trailing deformations
given in \eqref{definition-leading-trailing}.
In this notation, we have
\begin{equation}\label{3Dindex-integral-definition-strict-angles}
    \I{\triang}{\aholM}{\aholL}(\hbar)
    = c_q ^N
    \int_{(S^1)^N} \prod_{\square\in Q(\triang)} G_q \Bigl(
    e^{(\hbar/\pi+i)\alpha(\square)} \prod_{j=1}^N y_j^{l_j(\square)}
    \Bigr)\,\dHaar(y_1,\dotsc,y_N).
\end{equation}
Following Garoufalidis and Kashaev, we may eliminate the last variable
by making the substitution $y_j = w_j  y_N$, $1\leq j \leq N-1$,
and observing that
\[
    \prod_{j=1}^N y_j^{l_j(\square)}
    = y_N^{l_N(\square)}
    \prod_{j=1}^{N-1} {w_j}^{l_j(\square)}
    \prod_{j=1}^{N-1} y_N^{l_j(\square)}
    = \prod_{j=1}^{N-1} {w_j}^{l_j(\square)},
\]
where the last equality is a consequence of the relation $l_1+\dotsb+l_{N-1} = -l_N$.

In this way, we obtain the formula
\begin{equation}\label{3Dindex_N-1}
    \I{\triang}{\aholM}{\aholL}(\hbar)
    = c_q ^N
    \int_{(S^1)^{N-1}} \prod_{\square\in Q(\triang)} G_q \Bigl(
    e^{(\hbar/\pi+i)\alpha(\square)} \prod_{j=1}^{N-1} w_j^{l_j(\square)}
    \Bigr)\,\dHaar(w_1,\dotsc,w_{N-1}).
\end{equation}

The above integral can be interpreted in terms of \Sval\ angle structures
on $\triang$.
To this end, we take $t=(t_1,\dotsc,t_{N-1})\in[0,2\pi]^{N-1}$ and put
\begin{equation}\label{ltd-parametrization-of-SAS}
    \omega_t(\square) = \exp\Bigl(i \bigl(\alpha(\square) +
        \sum_{j=1}^{N-1} t_j l_j(\square)\bigr) \Bigr).
\end{equation}
Setting $a(\square) = \alpha(\square)/\pi\in(0,1)$ for all $\square\in Q(\triang)$,
we may now write
\begin{align}
    \I{\triang}{\aholM}{\aholL}(\hbar)
    &= \frac{c_q^N}{(2\pi)^{N-1}}
    \int\limits_{[0,2\pi]^{N-1}} \prod_{\square\in Q(\triang)} G_q \Bigl(
    e^{a(\square)\hbar}
    \exp\Bigl[i \bigl(\alpha(\square) + \sum_{j=1}^{N-1}t_j l_j(\square)\bigr) \Bigr]
    \Bigr)\,dt\label{3Dindex-int-2}\\
    &= \frac{c_q^N}{(2\pi)^{N-1}}
    \int\limits_{[0,2\pi]^{N-1}} \prod_{\square\in Q(\triang)}
    G_q\bigl(e^{a(\square)\hbar}\omega_t(\square)\bigr)\,dt.\label{integral-wrt-ltd-coords}
\end{align}
Observe that the scaling constants $\exp\bigl(a(\square)\hbar\bigr)$
ensure that each factor $G_q(z)$ is only evaluated in the punctured unit disc $0<|z|<1$,
where it is holomorphic.
Interpreting the integral~\eqref{3Dindex_N-1} as a
contour integral in several complex variables, Garoufalidis and Kashaev~\cite{garoufalidis-kashaev}
show that changing the strict angle structure~$\alpha$
amounts to merely scaling the radii of the contour.
As a result, the state integral~\eqref{3Dindex_N-1} is independent of
the choice of the initial strict angle structure~$\alpha$.

Moreover, since the integrand depends analytically on the
peripheral parameters~$(\aholM, \aholL)$, we obtain a germ of a two-variable
meromorphic function
\begin{equation*}
    \justI_\triang(\hbar) : (\aholM, \aholL)\mapsto
    \I{\triang}{\aholM}{\aholL}(\hbar)
\end{equation*}
called the \emph{meromorphic \ind\ of $\triang$}.

\begin{remark}
In the original definition of \cite{garoufalidis-kashaev},
the meromorphic \ind\ is treated as a function of
the variables~$e_\mu = \exp\bigl(\aholM(\hbar/\pi + i)\bigr)$ and
$e_\lambda = \exp\bigl(\aholL(\hbar/\pi + i)/2\bigr)$.
We prefer to work directly with the peripheral
angle-holonomies~$(\aholM,\aholL)=\bigl(\ahol{\alpha}(\mu), \ahol{\alpha}(\lambda)\bigr)$ due
to their more immediate geometric meaning.
\end{remark}

%------------------------------------------------------------------------------
\subsection{General definition of the meromorphic 3D-index}
\label{3D-index-general}
In general, an ideal triangulation $\triang$ may not admit a strict angle
structure, in which case the construction of the preceding section fails.
One can work around this issue by using analytic continuation, as described
by Andersen and Kashaev~\cite{kashaev-andersen-new},
Kashaev~\cite{kashaev-beta}, and Garoufalidis-Kashaev~\cite{garoufalidis-kashaev}.
The resulting invariant is still a meromorphic function of $(\aholM,\aholL)$,
but in the weaker sense of an analytic function $\CC^2\to\CC{P}^1$,
which includes the constant function equal to $\infty$.

Consider the set of \emph{generalised pre-angle structures},
\[
    \preAreal(\triang) = \{\alpha: Q(\triang)\to\RR\;|\;
                           \alpha(\square)+\alpha(\square')+\alpha(\square'')=\pi
                           \ \text{for all}\ \square\in Q(\triang)\}
\]
which is a real affine subspace of $\RR^{Q(\triang)}$ of dimension $2N$.
The peripheral angle-holonomy extends to a map
\begin{equation}\label{peripheral-holonomy-map-on-preangle-structures}
    \ahol{}(\mu,\lambda): \preAreal(\triang)\to\RR^2,\quad
    \ahol{}(\mu,\lambda)(\alpha)=\bigl(\ahol{\alpha}(\mu), \ahol{\alpha}(\lambda)\bigr),
\end{equation}
which is full rank.
Therefore, the preimage of an arbitrarily fixed point $(\aholM, \aholL)\in\RR^2$,
denoted by $\preArealsub{(\aholM,\aholL)}(\triang)$, is an affine subspace of dimension $2N-2$.
We now introduce $N-1$ independent functions
\begin{equation*}
    \epsilon_j: \preArealsub{(\aholM,\aholL)}(\triang)\to\RR,\quad
    \epsilon_j(\alpha) = \Bigl(\sum_{\square\in Q(\triang)} G(E_j,\square)\alpha(\square)\Bigr)-2\pi,\quad
    (1\leq j\leq N-1),
\end{equation*}
which are called \emph{angle excesses} along the edges $E_1,\dotsc,E_{N-1}$ of $\triang$.
In particular, writing $\epsilon=(\epsilon_1,\dotsc,\epsilon_{N-1})$ reduces the angle structure
equations~\eqref{AS:edge_sum=2pi} to the linear system $\epsilon=0$.

Observe that the space of tangential angle structures~$\TAS_0(\triang)$ acts freely
on $\preArealsub{(\aholM,\aholL)}(\triang)$ by vector addition.
Since the angle excesses are invariant under this action,
$\epsilon:\preArealsub{(\aholM,\aholL)}(\triang)\to\RR^{N-1}$ is
a principal \nword{\TAS_0(\triang)}{bundle}.

We fix the coordinates~$r\in\RR^{N-1}$ on $\TAS_0(\triang)$ by using the basis~$\{l_1,\dotsc,l_{N-1}\}$
of leading-trailing deformation vectors.
Then, an affine trivialisation of $\preArealsub{(\aholM,\aholL)}(\triang)$
equips it with coordinates $(r,\epsilon)$ such that
$\alpha(\square)$ is an affine function in the variables $(r_1,\dotsc,r_{N-1};\epsilon_1,\dotsc,\epsilon_{N-1})$
for all $\square\in Q(\triang)$.
Suppose that such a trivialisation has been chosen arbitrarily
and restrict it to the subspace
$\preAholo{(\aholM,\aholL)}(\triang)$ consisting of \emph{strict} pre-angle structures.
In other words, we impose the condition $0<\alpha(\square)<\pi$ for all $\square\in Q(\triang)$.
At this stage, we must assume that the peripheral parameters~$(\aholM, \aholL)$ are chosen in such a way that $\preAholo{(\aholM,\aholL)}(\triang)$ still has the full dimension $2N-2$,
so that the variables $(r,\epsilon)$ are restricted to a certain open polyhedron in $\RR^{N-1}\times\RR^{N-1}$.

Given a strict pre-angle structure $\alpha=\alpha(r,\epsilon)\in\preAholo{(\aholM,\aholL)}(\triang)$,
we define the associated state integral
\begin{equation}\label{pre-index}
    \Ipre{\triang}{\aholM}{\aholL}(\hbar)
    =
    \frac{c_q^N}{(2\pi)^{N-1}}
    \int\limits_{[0,2\pi]^{N-1}} \prod_{\square\in Q(\triang)} G_q \Bigl(
    e^{a(\square)\hbar}
    \exp\Bigl[i \bigl(\alpha(\square) + \sum_{j=1}^{N-1}t_j l_j(\square)\bigr) \Bigr]
    \Bigr)\,dt,
\end{equation}
where, as before, $a(\square)=\alpha(\square)/\pi$.
The above formula coincides with \eqref{3Dindex-int-2},
but this time $\justIpre_{\triang}$ depends non-trivially on $\alpha=\alpha(r,\epsilon)$.
By a contour shift argument analogous to that of the preceding section, we see that
$\Ipre{\triang}{\aholM}{\aholL}(\hbar)$ is gauge invariant, in the sense that it does not depend on
the variables $(r_1,\dotsc,r_{N-1})$.
Moreover, since the dependence on $\epsilon$ is analytic,
$\Ipre{\triang}{\aholM}{\aholL}(\hbar)$ yields a germ of an analytic function of the
remaining variables~$\epsilon=(\epsilon_1,\dotsc,\epsilon_{N-1})$.
We now wish to
define the meromorphic \ind\ of $\triang$ at $(\aholM,\aholL)$ as the value of the
analytic continuation of this germ to $\epsilon=0$.
A difficulty arises if this analytic continuation happens to have a singularity at $\epsilon=0$,
but in this case we simply set $\I{\triang}{\aholM}{\aholL}(\hbar)=\infty$.

\begin{definition}
When $\triang$ is an ideal triangulation of $M$ and a meridian-longitude
pair~$(\mu,\lambda)$ is fixed,
the \emph{meromorphic \ind\ of M,} denoted~$\justI_{M}(\hbar)$,
is the germ of the two-variable meromorphic function
assigning to $(\aholM,\aholL)$
the value of the analytic continuation of $\Ipre{\triang}{\aholM}{\aholL}(\hbar)$
to the point $\eps=0$.
\end{definition}
By using an analytic continuation of the pentagon identity, Garoufalidis and Kashaev
show in \cite{garoufalidis-kashaev} that $\justI_M(\hbar)$,
viewed as a function $\CC^2\to\CC{P}^1$ in the variables $(\aholM,\aholL)$,
is a topological invariant of the manifold~$M$ equipped with the basis $(\mu,\lambda)$
of $H_1(\partial\overline{M};\ZZ)$.

%--------------------------------------------------------------------------------------------------
\subsection{The state integral in terms of circle-valued angle structures}
\label{state-integral-to-circle-valued}
As mentioned previously, the state integral~\eqref{integral-wrt-ltd-coords} can be interpreted
in terms of the space of circle-valued angle structures on $\triang$.
By equation~\eqref{ltd-parametrization-of-SAS}, when the variables $t=(t_1,\dotsc,t_{N-1})$
vary in the cube $[0,2\pi]^{N-1}$, the \Sval\ angle structures~$\omega_t$ cover the
entire connected component of $\SAS_{(\exp(i\aholM), \exp(i\aholL))}(\triang)$
containing the initial angle structure $\exp(i\alpha)$.
In terms of the circle coordinates of~\eqref{3Dindex_N-1}, we may define the
associated locally diffeomorphic parametrisation,
\begin{align}\label{parametrization-of-SAS-with-circles}
    F_\triang: \bigl(S^1\bigr)^{N-1} &\to \SAS_{(\exp(i\aholM), \exp(i\aholL))}(\triang),\\
    F_\triang(w)=F_\triang(w_1,\dotsc,w_{N-1}) &=
    \Bigl(e^{i\alpha(\square)}\prod_{m=1}^{N-1}w_m^{l_m(\square)}\Bigr)_{\square\in Q(\triang).}
    \notag
\end{align}
Denote by $d_\triang$ the covering degree of $F_\triang$.
Exponentiating the Lebesgue measure endows the Lie group~$\bigl(S^1\bigr)^{N-1}$ with the canonical
invariant measure of total mass $(2\pi)^{N-1}$, which is a scalar multiple of the probabilistic (unit mass)
Haar measure.
Pushing this measure forward under the map~$F_\triang$
and rescaling by $d_\triang$, we obtain a well-defined
finite measure~``$d\omega$'' on $\SAS_{(\exp(i\aholM), \exp(i\aholL))}(\triang)$
with a total mass of $\dfrac{(2\pi)^{N-1}}{d_\triang}$.
In this way, the integral~\eqref{integral-wrt-ltd-coords} can be rewritten as
\begin{equation}\label{3Dindex_integral_over_SAS}
    \I{\triang}{\aholM}{\aholL}(\hbar)
    =
    \frac{d_\triang c_q^N}{(2\pi)^{N-1}}
    \int\limits_{\Omega} \prod_{\square\in Q(\triang)}
    G_q\bigl(e^{a(\square)\hbar}\omega(\square)\bigr)\;d\omega,
\end{equation}
where $\Omega$ is the connected component of $\SAS_{(\exp(i\aholM), \exp(i\aholL))}(\triang)$
containing the circle-valued angle structure~$\exp(i\alpha)$.

\begin{remark}\label{hyperbolic-Omega}
In the special case of $\aholM=\aholL=0$, Theorem~\ref{obstruction-theorem}
establishes a bijective correspondence between the
connected components of $\SAS_0(\triang)$
and cohomology classes in $H^2(M,\partial M;\FF_2)$.
If the triangulation $\triang$ admits a (positively oriented) geometric
solution~$z$ of edge consistency and completeness equations and we take $\alpha=\Arg(z)$
to be the hyperbolic dihedral angles, then we have
$\Omega=\Phi_0^{-1}\bigl(\Obs_0(\rho^\geo)\bigr)$, where $\Obs_0(\rho^\geo)\in H^2(M,\partial M;\FF_2)$
is the obstruction class to lifting the holonomy representation $\rho^\geo$ of the complete hyperbolic
structure on $M$ to a boundary-unipotent $\SLC$ representation.
This provides a topological interpretation of the component~$\Omega$
appearing as the domain of integration in \eqref{3Dindex_integral_over_SAS}.
\end{remark}

We shall now compute the covering degree $d_\triang$ and show that it too
can be defined purely topologically.
Let $\hat{M}$ be the end compactification of $M$, i.e., the
compact topological space obtained from
$M$ by adding a point at infinity to the toroidal end of $M$.
In terms of the ideal triangulation~$\triang$,
$\hat{M}$ can be constructed as the simplicial complex~$\hat{\triang}$
resulting from filling in the ideal vertex in $\triang$.
\begin{proposition}\label{d=homology_order}
The degree $d_\triang$ of the map \eqref{parametrization-of-SAS-with-circles}
is equal to the order of the homology group $H_1(\hat{M};\FF_2)$.
In particular, $d_\triang$ does not depend on the choice of the triangulation $\triang$.
\end{proposition}

\begin{proof}
The map $F_\triang$ of \eqref{parametrization-of-SAS-with-circles} can be understood in terms of
the linear transformation given by the transpose of
the integer matrix~$\mathbf{L_*}$ of eq.~\eqref{definition-L_star}, with respect to the standard bases.
To this end, we choose an arbitrary angle structure $\omega\in\SAS_0(\triang)$
and identify $S^1$ with the quotient $\RR/\ZZ$,
obtaining the following commutative diagram:
\begin{equation}\label{commie-diagram}
    \begin{tikzcd}[column sep=5em]
        &\ZZ^{N-1}\arrow{r}{\mathbf{L}_*^\transp\bigr|_{\ZZ^{N-1}}}\arrow[hook]{d}
        &\TAS_0(\triang)\cap\ZZ^{Q(\triang)}\arrow[hook]{d}\arrow[hook]{r}
        &\ZZ^{Q(\triang)}\arrow[hook]{d}
        \\
        &\RR^{N-1}\arrow[hook,two heads]{r}{\mathbf{L}_*^\transp}\arrow[two heads]{d}{\exp(2\pi i\;\cdot\;)}
        &\TAS_0(\triang)\arrow[hook]{r}\arrow{d}{f}
        &\RR^{Q(\triang)}\arrow[two heads]{d}{\exp(2\pi i\;\cdot\;)}
        \\
        &\bigl(S^1\bigr)^{N-1}\arrow{r}{F_\triang}
        &\SAS_0(\triang)\arrow[hook]{r}
        &\bigl(S^1\bigr)^{Q(\triang)},
    \end{tikzcd}
\end{equation}
where the middle vertical map~$f=f_\omega$ is defined by
$f(v)(\square) = \omega(\square)\exp\bigl(2\pi i v(\square)\bigr)$
and the unlabelled maps are the obvious ones.

Note that the leftmost and the rightmost columns of the diagram~\eqref{commie-diagram} are
short exact sequences of abelian groups.
The central column is also ``exact'' in the following sense: the fibres of the map $f$
are in a one-to-one correspondence with integer lattice points in $\TAS_0(\triang)$.
To make this observation precise, it suffices to focus on the chosen basepoint~$\omega$ and to write
\[
    f^{-1}(\{\omega\})
    = \bigl\{v\in\TAS_0(\triang): \exp\bigl(2\pi iv(\square)\bigr)=1\ \text{for all}\ \square\in Q(\triang)\bigr\}
    = \TAS_0(\triang)\cap\ZZ^{Q(\triang)}.
\]

For an $m\times n$ integer matrix $A$, denote by $\pned(A)\in\NN$ the product of non-zero
elementary divisors of $A$, i.e., the product of the non-zero diagonal entries
in the Smith normal form of $A$.
It is well known that  $\pned(A)$ is equal to the index of the
sublattice $A(\ZZ^n)$ in the lattice $A(\RR^n)\cap\ZZ^m$.
Since $\mathbf{L}_*^\transp$ is an isomorphism whose image coincides with the image
of $\mathbf{L}^\transp$ both over $\ZZ$ and over $\RR$, we have
\begin{align*}
    d_\triang
    &= \left|F_\triang^{-1}(\{\omega\})\right|\\
    &= \bigl[f^{-1}(\{\omega\}) : \mathbf{L}_*^\transp(\ZZ^{N-1})\bigr]\\
    &= \bigl[\TAS_0(\triang)\cap\ZZ^{Q(\triang)} : \mathbf{L}_*^\transp(\ZZ^{N-1})\bigr]\\
    &= \bigl[\mathbf{L}_*^\transp(\RR^{N-1})\cap\ZZ^{Q(\triang)} : \mathbf{L}_*^\transp(\ZZ^{N-1})\bigr]\\
    &= \bigl[\mathbf{L}^\transp(\RR^N)\cap\ZZ^{Q(\triang)} : \mathbf{L}^\transp(\ZZ^N)\bigr]\\
    &= \pned(\mathbf{L}^\transp).
\end{align*}
Neumann proves as part of Theorem 4.2 in \cite{neumann1990} that 
that $H_1(\hat{\triang};\FF_2) \cong \Ker_\ZZ p / \mathbf{L}(\ZZ^{Q(\triang)})$.
Since this homology group is finite, we have
$\mathbf{L}(\RR^{Q(\triang)})=\Ker_\RR p$, so that
\begin{align*}
    d_\triang
    &= \pned(\mathbf{L})\\
    &= \bigl[ \mathbf{L}(\RR^{Q(\triang)})\cap\ZZ^N : \mathbf{L}(\ZZ^{Q(\triang)})\bigr]\\
    &= \bigl[ \Ker_\ZZ p : \mathbf{L}(\ZZ^{Q(\triang)}) \bigr]\\
    &= |H_1(\hat{\triang};\FF_2)|. \qedhere
\end{align*}
\end{proof}

%==============================================================================
% Section 4: Asymptotic toolbox
%==============================================================================
\section{Asymptotic toolbox}
In this section, we summarise the fundamental asymptotic expansions
that will become the essential tools in the study of the $\hbar\to0$ asymptotic
behaviour of the meromorphic \nword{3D}{index}~$\justI_M(\hbar)$.
For simplicity, we shall restrict our attention to $\hbar$ real and negative,
but some of the results of this section hold for complex $\hbar$ with $\Real\hbar<0$.

%--------------------------------------------------------------------------------------------------
\subsection{\texorpdfstring{Asymptotic expansions of the $q$-dilogarithm}{Asymptotic expansions of the q-dilogarithm}}
We shall start by describing the asymptotic behaviour of
the \nword{q}{dilogarithm} $\Li(z;q)$ when
$z$ approaches a point on the unit circle from within the unit disc and $q$ approaches $1$.
More precisely, we set $q=e^{\hbar}$ and fix complex numbers $a,\omega$ with
$|\omega|=1$.
As suggested by the expression~\eqref{3Dindex_integral_over_SAS},
we need to study the asymptotics  of $G_q(\omega e^{a\hbar})$ when
$\hbar\to0^-$.

For any non-negative integer $k$, denote by $B_k(x)$ the \nword{k}{th} Bernoulli
polynomial~\cite[\S1.13]{erdelyi-bateman}, defined by the equation
\[
    \frac{t e^{xt}}{e^t-1} = \sum_{k=0}^\infty B_k(x)\frac{t^n}{n!},\quad |t|<2\pi.
\]
The value $B_k = B_k(0)$ is called the \nword{k}{th} Bernoulli number.
For any integer $k$, let $\Lis{k}(z)$ be the \nword{k}{th} polylogarithm function;
when $k>0$, it is given by the power series
\[
    \Lis{k}(z) = \sum_{n=1}^\infty \frac{z^n}{n^k},\quad |z|<1
\]
and thus satisfies the functional equation~$\Lis{k}(z)=z\frac{d}{dz}\Lis{k+1}(z)$.
Using this functional equation, we may define $\Lis{k}(z)$ for $k\leq0$
recursively; in particular
\begin{equation}\label{definition-polylogarithm}
    \Lis{0}(z) = \frac{z}{1-z}, \quad \Lis{k}(z) = z\frac{d}{dz} \Lis{k+1}(z),\quad k<0.
\end{equation}
As a consequence, for $k\leq 0$ the polylogarithms $\Lis{k}(z)$ are rational functions with
a pole at $z=1$ and analytic for $z\neq1$.

\begin{lemma}\label{fundamental-asymptotic-lemma}
    In terms of the variable $\hbar < 0$, the following asymptotic expansions hold.
    \begin{enumerate}[(i)]
        \item\label{item:Zhangs-expansion}
        When $a\in\CC$ is arbitrarily fixed, we have
        \begin{equation}\label{Zhang}
            \Li(e^{a\hbar}; e^\hbar) =
            -\frac{\pi^2}{6} \hbar^{-1}
            + (a-\tfrac{1}{2})\log(-\hbar) + \log\frac{\Gamma(a)}{\sqrt{2\pi}}
            + \frac{B_2(a)}{4}\hbar\\
            + \sum_{k=1}^\infty\frac{B_{2k+1}(a)B_{2k}}{2k(2k+1)!}\hbar^{2k}
        \end{equation}
        as $\hbar\to 0^-$, uniformly for $a$ in any compact subset of $\CC$.
        \item\label{item:Andrews-expansion}
        When $\omega\neq1$, $|\omega|=1$, and the parameter $a\in\CC$
        satisfies $\Real a>0$, we have
        \begin{equation}\label{Andrews-expansion}
            \Li(\omega e^{a\hbar}; e^\hbar) =
            -\Li(\omega)\hbar^{-1}
            + (a-{\textstyle\frac12})\log(1-\omega)
            - \sum_{k=1}^{\infty}\frac{B_{k+1}(a)}{(k+1)!}\Lis{1-k}(\omega)\,\hbar^k
        \end{equation}
        as $\hbar\to0^-$, where ``$\log$'' refers to the standard branch of
        the logarithm.
    \end{enumerate}
\end{lemma}
The first part of the above lemma can be found in McIntosh~\cite{mcintosh-qfactorials}
and Zhang~\cite[Theorem~2]{zhang-qexp}.
An expansion analogous to \eqref{Andrews-expansion} was derived by Bouzeffour and
Halouani~\cite{bh-qdilog} for $\Real a>1$.
The extension to $\Real a>0$ is necessary in order to allow for $a\in(0,1)$,
as in \eqref{3Dindex-int-2},
and therefore we provide an independent proof of the asymptotic expansion~\eqref{Andrews-expansion}
in Appendix~\ref{asymptotic-appendix}.

\begin{remark}
Since the Bernoulli numbers~$B_k$ vanish for odd $k\geq3$, we may replace
the infinite asymptotic series of \eqref{Zhang} with
\begin{equation}\label{Zhang'}
    \sum_{k=2}^\infty\frac{B_{k+1}(a)B_{k}}{k(k+1)!}\,\hbar^k ,
\end{equation}
which is in agreement with  McIntosh~\cite[Theorem~2]{mcintosh-qfactorials} after setting
$\hbar=-t$, $a=c$.
\end{remark}

\begin{lemma}\label{simplifying-lemma}
When $|\omega|=1$ and $\omega\neq\pm1$, we have the identities
\begin{align}
    \Li(\omega)-\Li(-\bar{\omega})
    &= \frac{\pi^2}{4} - \frac{\pi}{2}|\Arg\omega| + i\circLob(\omega),
    \label{combining-dilogarithms-into-Lobachevsky}\\
    \Lis{1-k}(\omega) + (-1)^k\Lis{1-k}(-\bar{\omega})
    &= 2^k\Lis{1-k}(\omega^2), \quad k\geq2
    \label{simplifying-polylogarithms}
\end{align}
where $\Arg\omega\in(-\pi,0)\cup(0,\pi)$ is the principal argument of $\omega$ and $\circLob$
is the Lobachevsky function on the circle, defined in \eqref{circle-Lob}.
\end{lemma}
\begin{proof}
In order to demonstrate \eqref{combining-dilogarithms-into-Lobachevsky},
we choose two arguments $\theta\in[0,2\pi]$ and $\theta'\in[-\pi,\pi]$ so that
$\omega=e^{i\theta}=e^{i\theta'}$.
We first write
\[
      \Li(\omega)-\Li(-\bar{\omega})
    = \bigl(\Li(e^{i\theta})-\Li(1)\bigr)-\bigl(\Li(e^{i(\pi-\theta')})-\Li(1)\bigr),
\]
and then apply the identity
\begin{equation}\label{Milnor-functional-equation}
    \Li(e^{i\alpha})-\Li(1)
    = -\frac{\alpha}{2}\left(\pi-\frac{\alpha}{2}\right)+2i\Lob\left(\frac{\alpha}{2}\right),
    \quad\text{where}\ \alpha\in[0,2\pi]
\end{equation}
which is stated as eq.~7.3 in Chapter~7 (by J.~Milnor) of Thurston's notes \cite{thurston-notes}.
Looking at the real parts, we see
\begin{equation*}
    \Real\bigl[\Li(\omega)-\Li(-\bar{\omega})\bigr]
    = \tfrac{1}{2}\left[
        \theta\bigl(\tfrac{\theta}{2}-\pi\bigr)
        + (\pi-\theta')\bigl(\tfrac{\theta'}{2}+\tfrac{\pi}{2}\bigr)\right].
\end{equation*}
When $0\leq\theta\leq\pi$, then $\theta'=\theta$ and the right-hand side
simplifies to $\pi^2/4 - \pi\theta'/2$.
When $\pi\leq\theta\leq2\pi$, the expression becomes $\pi^2/4+\pi\theta'/2$ instead.
Since $\theta' = \Arg\omega$, we always have
\begin{equation*}
    \Real\bigl[\Li(\omega)-\Li(-\bar{\omega})\bigr]
    = \frac{\pi^2}{4}-\frac{\pi}{2}|\Arg\omega|,
\end{equation*}
which establishes the real part of \eqref{combining-dilogarithms-into-Lobachevsky}.
At the same time, the imaginary parts transform into
\begin{equation*}
    \Imag\bigl[\Li(\omega)-\Li(-\bar{\omega})\bigr]
    = 2\left(\Lob\bigl(\tfrac{\theta}{2}\bigr)-\Lob\bigl(\tfrac{\pi}{2}-\tfrac{\theta}{2}\bigr)\right)
    = \Lob(\theta),
\end{equation*}
the last equality being the $n=2$ case of Lemma~7.1.4 of~\cite{thurston-notes}.

In order to establish \eqref{simplifying-polylogarithms},
we utilise the following functional equations for the polylogarithm
$\Lis{n}(z)$ with $n\in\ZZ$, $n<0$:
\[
    \Lis{n}\bigl(\tfrac{1}{z}) = (-1)^{n+1}\Lis{n}(z),\qquad
    \Lis{n}(z) + \Lis{n}(-z) = 2^{1-n}\Lis{n}(z^2).
\]
The formula on the left is known as the `inversion formula' and is stated by
Erd\'elyi et al.~\cite{erdelyi-bateman} as eq.~1.11(17).
The `duplication formula' on the right is well known and can be easily
proved from \eqref{definition-polylogarithm} by induction on $n$.
Note that since the polylogarithm $\Lis{n}(z)$ is analytic everywhere except for
the pole at $z=1$, the duplication formula always makes sense at
$z=\bar{\omega}\neq\pm1$, yielding
\[
    \Lis{1-k}(-\bar{\omega})
    = 2^k\Lis{1-k}\bigl(\tfrac{1}{\omega^2}\bigr) - \Lis{1-k}\bigl(\tfrac{1}{\omega}\bigr)
    = (-1)^k \bigl(2^k\Lis{1-k}(\omega^2)-\Lis{1-k}(\omega)\bigr)
\]
and \eqref{simplifying-polylogarithms} follows immediately.
\end{proof}

%--------------------------------------------------------------------------------------------------
\subsection{Asymptotics of the Boltzmann weights of the meromorphic 3D-index}
The goal of this section is to establish asymptotic approximations
of the quantities $c_q$ and $G_q(\omega e^{a\hbar})$ occurring in \eqref{3Dindex_integral_over_SAS}.
Here, we fix an $\omega$ with $|\omega|=1$ and set $q=e^\hbar$,
where $\hbar$ is a negative real parameter.
Since the expansions in Lemma~\ref{fundamental-asymptotic-lemma}
are fundamentally different for $\omega=1$ and for $\omega\neq 1$,
the form of the subsequent expansions for $\omega=\pm1$ will also differ from the
expansions for $\omega\neq\pm1$.

\begin{theorem}\label{th:asym-c_q-and-G_q}
Let $m$ be an arbitrary non-negative integer and let $\hbar$ be a negative real parameter.
For any fixed $a\in\CC$, we have the following asymptotic approximations as $\hbar\to0^-$.
\begin{align}
    c_q = c_{\exp(\hbar)}
    &= \exp\Bigl[
    \frac{\pi^2}{4}\hbar^{-1} - \frac{1}{2}\log(-\hbar) + \log\sqrt{4\pi} + O(\hbar^\infty)
    \Bigr], \label{asym-c}\\
    G_q(q^a) = G_{e^\hbar}(e^{a\hbar})
    &= \exp\Bigl[ -\frac{\pi^2}{4}\hbar^{-1} + (a-\tfrac{1}{2})\log(-\hbar)
        + \log\frac{\Gamma(a)}{\sqrt{2\pi}} + (a-\tfrac{1}{2})\log2 \label{asym-G-omega=+1}\\
    &\phantom{=\exp x} + \sum_{k=1}^m
        \frac{B_{2k}B_{2k+1}(a)}{k(2k+1)!}2^{2k}\,
        \hbar^{2k}
    + O(\hbar^{2(m+1)})\Bigr],\notag\\
    G_q(-q^a) = G_{e^\hbar}(-e^{a\hbar})
    &= \exp\Bigl[ \frac{\pi^2}{4}\hbar^{-1} + (a-\tfrac{1}{2})\log(-\hbar)
        - \log\frac{\Gamma(1-a)}{\sqrt{2\pi}} + (a-\tfrac{1}{2})\log 2 \label{asym-G-omega=-1}\\
    &\phantom{=\exp x} + \sum_{k=1}^m\frac{B_{2k} B_{2k+1}(a)}{2k(2k+1)!} 2^{2k}\,\hbar^{2k}
        + O(\hbar^{2(m+1)}) \Bigr],\notag
\end{align}
where the term~$O(\hbar^\infty)$ is $O(\hbar^m)$ for all $m\geq0$.
Moreover, \eqref{asym-G-omega=+1}--\eqref{asym-G-omega=-1} hold uniformly for $a$ on compact subsets of $\CC$.
Further assuming that $\omega\in\CC\setminus\RR$ satisfies $|\omega|=1$ and that $a\in(0,1)$, we have
\begin{align}
    G_{e^\hbar}(\omega e^{a\hbar})
    &= \exp\Bigl[
        \bigl(\tfrac{\pi}{2}|\Arg\omega|-\tfrac{\pi^2}{4}-i\circLob(\omega)\bigr)\hbar^{-1}\notag\\
    &\phantom{= \exp\Bigl[}
       + \bigl(a-\tfrac{1}{2}\bigr)
        \bigl(\log|2\Imag\omega|-i\tfrac{\pi}{2}\sgn\Imag\omega\bigr)\label{asym-G-general}
       + \frac{B_2(a)\Real\omega}{2i\Imag\omega}\,\hbar\\
       &\phantom{= \exp\Bigl[}
       - \sum_{k=2}^m\frac{2^k B_{k+1}(a)}{(k+1)!}
        \Lis{1-k}(\omega^2)\,\hbar^k
        +O(\hbar^{m+1})
        \Bigr]\ \text{as}\ \hbar\to0^-.\notag
\end{align}
\end{theorem}

\begin{proof}
In order to derive the asymptotic expansion~\eqref{asym-c} of the constant $c_q=c_{\exp(\hbar)}$ of
\eqref{definition-of-c(q)},
we may write
\[
    c_{\exp(\hbar)}
    = \exp\left[\Li(e^{2\hbar};e^{2\hbar}) - 2\Li(e^\hbar; e^\hbar)\right]
\]
and apply \eqref{Zhang} with $a=1$ to both terms on the right-hand side.
Recall the reflection formula for Bernoulli polynomials~\cite[eq.~1.13(12)]{erdelyi-bateman},
\begin{equation}\label{reflection-Bernoulli}
    B_k(1-a) = (-1)^k B_k(a).
\end{equation}
Since the higher odd Bernoulli number vanish, i.e., $B_{2k+1}=0$ for all $k\geq1$,
we have
\(
    B_{2k+1}(1)B_{2k} = -B_{2k}B_{2k+1} = 0
\)
for all $k\geq 1$.
In this way, we see
that the higher \nword{\hbar}{order} error terms vanish and we are left with
\begin{align}
    c_{\exp(\hbar)}
    &= \exp\left[\Li(e^{2\hbar};e^{2\hbar}) - 2\Li(e^\hbar; e^\hbar)\right]\notag\\
    &= \exp\Bigl[
    \frac{\pi^2}{4}\hbar^{-1} - \frac{1}{2}\log(-\hbar) + \log\sqrt{4\pi} + O(\hbar^\infty)
    \Bigr].
\end{align}

In order to establish equation~\eqref{asym-G-omega=+1}, we use
the duplication formula for the \nword{q}{dilogarithm} to obtain
\[
    \Li(-e^{a\hbar}; e^\hbar) = \Li(e^{2a\hbar}; e^{2\hbar}) - \Li(e^{a\hbar}; e^\hbar)
\]
and then apply \eqref{Zhang} to both terms on the right-hand side.
This results in the asymptotic expansion
\begin{align}
    \Li(-e^{a\hbar}; e^\hbar) &=
            -\frac{\pi^2}{6}\left(\frac{1}{2\hbar}-\frac{1}{\hbar}\right)
            + (a-\tfrac{1}{2})\left(\log(-2\hbar)-\log(-\hbar)\right)
            + \frac{B_2(a)}{4}(2\hbar-\hbar)\notag\\
            & \phantom{{}={}} + \sum_{k=1}^{m}\frac{B_{2k+1}(a)B_{2k}}{2k(2k+1)!}
                \bigl((2\hbar)^{2k}-\hbar^{2k}\bigr)
            + O(\hbar^{2(m+1)})\notag\\
    &=\label{nicer-minus-one-expansion}
    \frac{\pi^2}{12}\hbar^{-1}
    + (a-\tfrac{1}{2})\log 2
    + \frac{B_2(a)}{4}\hbar\\
    & \phantom{{}={}} + \sum_{k=1}^{m}\frac{B_{2k+1}(a)B_{2k}}{2k(2k+1)!}(2^{2k}-1)\hbar^{2k}
    + O(\hbar^{2(m+1)}),\notag
\end{align}
valid uniformly for $a$ on compact subsets as $\hbar\to0^-$.
Using \eqref{nicer-minus-one-expansion} and once again \eqref{Zhang},
we compute
\begin{align}
    G_{e^\hbar}(e^{a\hbar})
    &= \exp\left[ \Li(e^{a\hbar}; e^\hbar) - \Li(-e^{(1-a)\hbar}; e^\hbar)\right]\notag\\
    &= \exp\Bigl[ -\frac{\pi^2}{4}\hbar^{-1} + (a-\tfrac{1}{2})\log(-\hbar)
        + \log\frac{\Gamma(a)}{\sqrt{2\pi}} + (a-\tfrac{1}{2})\log2 +
        \frac{B_2(a)-B_2(1-a)}{4}\hbar\notag\\
    &\phantom{=\exp x} + \sum_{k=1}^m
        \frac{B_{2k}\bigl(B_{2k+1}(a)-B_{2k+1}(1-a)\bigr)}{2k(2k+1)!}2^{2k}\,
        \hbar^{2k}
    + O(\hbar^{2(m+1)})
    \Bigr]\notag\\
    &= \exp\Bigl[ -\frac{\pi^2}{4}\hbar^{-1} + (a-\tfrac{1}{2})\log(-\hbar)
        + \log\frac{\Gamma(a)}{\sqrt{2\pi}} + (a-\tfrac{1}{2})\log2 \label{omega-plus}\\
    &\phantom{=\exp x} + \sum_{k=1}^m
        \frac{B_{2k}B_{2k+1}(a)}{k(2k+1)!}2^{2k}\,
        \hbar^{2k}
    + O(\hbar^{2(m+1)})\Bigr],\notag
\end{align}
the last equality being a consequence of the reflection
formula~\eqref{reflection-Bernoulli}. Thus \eqref{asym-G-omega=+1} is established.

In order to prove \eqref{asym-G-omega=-1}, we proceed analogously and obtain
\begin{align}
    G_{e^\hbar}(-e^{a\hbar}) \notag
    &= \exp\left[ \Li(-e^{a\hbar}; e^\hbar) - \Li(e^{(1-a)\hbar}; e^\hbar)\right]\notag\\
    &= \exp\Bigl[ \frac{\pi^2}{4}\hbar^{-1} + (a-\tfrac{1}{2})\log(-\hbar)
        - \log\frac{\Gamma(1-a)}{\sqrt{2\pi}} + (a-\tfrac{1}{2})\log 2 \\
    &\phantom{=\exp x} + \sum_{k=1}^m\frac{B_{2k} B_{2k+1}(a)}{2k(2k+1)!} 2^{2k}\,\hbar^{2k}
        + O(\hbar^{2(m+1)}) \Bigr],\notag
\end{align}
as required.
It remains to study the case when $|\omega|=1$ but $\omega\neq\pm1$.
In this case, we use the asymptotic expansion~\eqref{Andrews-expansion} to obtain
\begin{align}
    G_{e^\hbar}(\omega e^{a\hbar})
    &= \exp\left[\Li(\omega e^{a\hbar};e^\hbar) - \Li(-\bar{\omega}e^{(1-a)\hbar};e^\hbar)\right]\notag\\
    &= \exp\Bigl[\bigl(\Li(-\bar{\omega})-\Li(\omega)\bigr)\hbar^{-1}
       + \bigl(a-\tfrac{1}{2}\bigr)\bigl(\log(1-\omega)+\log(1+\bar{\omega})\bigr)\label{smooth-preasym}\\
       &\phantom{\exp\Bigl[}
       + \bigl(B_2(1-a)\Lis0(-\bar\omega) - B_2(a)\Lis0(\omega)\bigr)\frac{\hbar}{2!}\notag\\
       &\phantom{\exp\Bigl[}
       - \sum_{k=2}^m\bigl( B_{k+1}(a)\Lis{1-k}(\omega)-B_{k+1}(1-a)\Lis{1-k}(-\bar{\omega})\bigr)
        \frac{\hbar^k}{(k+1)!}
     + O(\hbar^{m+1})
        \Bigr].\notag
\end{align}
We may simplify the coefficients in the above expansion as follows.
Firstly, the coefficient of $\hbar^{-1}$ has an expression given by
\eqref{combining-dilogarithms-into-Lobachevsky}.
Secondly, setting $\omega=e^{i\theta}$, where $\theta=\Arg\omega$, we may simplify the constant
($\hbar^0$) term as follows:
\begin{align*}
    \log(1-\omega)+\log(1+\bar{\omega})
    &= \log\bigl((1-e^{i\theta})(1+e^{-i\theta})\bigr)
    = \log(-2i\sin\theta)
    = \log|2\sin\theta|-i\tfrac{\pi}{2}\sgn\theta\\
    &= \log|2\Imag\omega|-i\tfrac{\pi}{2}\sgn\Imag\omega.
\end{align*}
Thirdly, the coefficient of $\hbar^1$ can be transformed with the help of the
definition \eqref{definition-polylogarithm} and the reflection
formula~\eqref{reflection-Bernoulli}, yielding the expression
\[
    \tfrac{1}{2}\bigl(B_2(1-a)\Lis0(-\bar\omega) - B_2(a)\Lis0(\omega)\bigr)
    = \frac{B_2(a)}{2}\, \frac{\omega+\bar\omega}{\omega-\bar\omega}
    = \frac{B_2(a)}{2i}\,\frac{\Real\omega}{\Imag\omega}.
\]
Lastly, we apply \eqref{simplifying-polylogarithms} to simplify
the coefficient of $\hbar^k$ for $k\geq2$.
The equation~\eqref{asym-G-general} follows.
\end{proof}

%--------------------------------------------------------------------------------------------------
\subsection{Regions of pointwise exponential suppression of Boltzmann weights}

Suppose that an \Sval\ angle structure $\omega\in\SAS(\triang)$ is fixed.
For arbitrarily fixed constants $a(\square)\in(0,1)$, we may consider the
quantity
\begin{equation}\label{integrand-as-function-of-hbar}
    b_\omega(\hbar) =
    c_q^N \prod_{\square\in Q(\triang)}
    G_q\bigl(\omega(\square)e^{a(\square)\hbar}\bigr)
\end{equation}
which corresponds to the value at $\omega$ of the integrand of \eqref{3Dindex_integral_over_SAS}.
In this section, we characterise \Sval\ angle structures $\omega$ for which $b_\omega(\hbar)$
decays exponentially as $\hbar\to0^-$.

Let $\tet$ be a tetrahedron of $\triang$ and denote by $\square$, $\square'$, $\square''$
the three normal quadrilaterals in $\tet$.
The contribution of $\tet$ to $b_\omega(\hbar)$ is then
\[
    b^\tet_\omega(\hbar) = c_q\,G_q\bigl(\omega(\square  )e^{a(\square  )\hbar}\bigr)\,
                                G_q\bigl(\omega(\square' )e^{a(\square' )\hbar}\bigr)\,
                                G_q\bigl(\omega(\square'')e^{a(\square'')\hbar}\bigr),
\]
so that $b_\omega(\hbar) = \prod_{\tet} b^\tet_\omega(\hbar)$.
\begin{definition}\label{definition-suppressing}
    We say that the angled
    tetrahedron~$\bigl(\tet,\;\omega(\square), \omega(\square'),\omega(\square'')\bigr)$
    is \emph{non-suppressing} if at least one of the following conditions is true:
    \begin{enumerate}[(i)]
        \item\label{cond:angle=1}
            $1 \in \bigl\{\omega(\square), \omega(\square'), \omega(\square'')\bigr\}$, or
        \item\label{cond:same_signs}
            the imaginary parts of $\omega(\square)$, $\omega(\square')$, $\omega(\square'')$
            are either all positive or all negative.
    \end{enumerate}
    If neither of the conditions \eqref{cond:angle=1}--\eqref{cond:same_signs} holds,
    we say that $\tet$ is \emph{suppressing}.
\end{definition}

\begin{remark} Condition (ii) arises geometrically in Proposition \ref{vol_crit_point_to_sol_gluing_eqns}, guaranteeing that a critical point of hyperbolic volume gives an algebraic solution to Thurston's gluing equations.
\end{remark}

\begin{proposition}\label{proposition-suppression}
    Let $\hbar = -1/\kappa$ where $\kappa>0$ is a real parameter.
    For a fixed $\omega\in\SAS(\triang)$, the value of $b_\omega(\hbar)$ decays
    exponentially fast as $\kappa\to+\infty$ if and only if $\triang$ contains a suppressing tetrahedron.
\end{proposition}

In order to prove this proposition, we must analyse the coefficients of $\hbar^{-1}$ in
the asymptotic approximations given in Theorem~\ref{th:asym-c_q-and-G_q}.
Since $\hbar < 0$, exponential suppression will occur when the real parts of these coefficients are positive.
Observe that since $\Arg\omega(\square)\in[-\pi,\pi]$,
the expansions~\eqref{asym-G-omega=+1}, \eqref{asym-G-omega=-1} and \eqref{asym-G-general}
imply
\[
    \abs{G_{e^\hbar}\bigl(\omega(\square)e^{a(\square)\hbar}\bigr)}
    = \exp\Bigl[
        \bigl(\tfrac{\pi}{2}|\Arg\omega(\square)|-\tfrac{\pi^2}{4}\bigr)\hbar^{-1}
        + O\bigl(\log(-\hbar)\bigr)\Bigr],
    \ \hbar\to0^-
\]
for any $\omega(\square)\in S^1$.
Using additionally the expansion \eqref{asym-c}, we get
\begin{equation}\label{suppression-analysis-pertet}
    \abs{b^\tet_\omega(\hbar)}
    = \exp\Bigl[\tfrac{\pi}{2}
        \bigl(
        |\Arg\omega(\square)|+
        |\Arg\omega(\square')|+
        |\Arg\omega(\square')|
        -\pi\bigr)\hbar^{-1} + O\bigl(\log(-\hbar)\bigr)\Bigr].
\end{equation}
We now characterise the sign of the overall coefficient of $\hbar^{-1}$ in the above approximation.
\begin{lemma}\label{suppression-inequality-lemma}
    For a tetrahedron $\tet$ with \Sval\
    angles~$\omega(\square)$, $\omega(\square')$, $\omega(\square'')$, we have
    \begin{equation}\label{tetrahedral-suppression-inequality}
        |\Arg\omega(\square)|+
        |\Arg\omega(\square')|+
        |\Arg\omega(\square')|
        \geq \pi.
    \end{equation}
    Moreover, the inequality is strict if and only if $\tet$ is suppressing.
\end{lemma}
It is easy to see that the above lemma implies Proposition~\ref{proposition-suppression}.
Indeed, when $\tet$ is a suppressing tetrahedron, the coefficient of $\hbar^{-1}$ in
\eqref{suppression-analysis-pertet} is positive, so that $b_\omega(\hbar)$ decays exponentially
as $-1/\hbar=\kappa\to+\infty$.
Conversely, if all tetrahedra are non-suppressing, then all inequalities
\eqref{tetrahedral-suppression-inequality} are saturated, implying that the
coefficients of $\hbar^{-1}$ in \eqref{suppression-analysis-pertet} vanish for all $\tet$.
In this case, the behaviour of the integrand $b_\omega(\hbar)$ is oscillatory at the leading order.

\begin{proof}[Proof of Lemma~\ref{suppression-inequality-lemma}]
Fix a tetrahedron $\tet$ with quads $\square$, $\square'$ and $\square''$ and write
$(\omega, \omega', \omega'') = \bigl(\omega(\square), \omega(\square'), \omega(\square'')\bigr)$ for short.
We first show that if $\tet$ is non-suppressing, then the
inequality~\eqref{tetrahedral-suppression-inequality} is in fact an equality.
Since $\omega\omega'\omega''=-1$, when
$(\omega, \omega', \omega'')$ all have positive imaginary parts, the absolute value signs
can be omitted and the equality holds trivially.
The negative case follows by complex conjugation, since $|\Arg\bar{\omega}| = |\Arg\omega|$.

In case $\tet$ satisfies condition~\eqref{cond:angle=1} in Definition~\ref{definition-suppressing},
we may assume, without loss of generality, that $\omega=1$.
If $\omega',\omega''\in\RR$, then $\{\omega', \omega''\}=\{1,-1\}$,
so that exactly one term in \eqref{tetrahedral-suppression-inequality} is $\pi$
and the two remaining ones are zeros, and equality holds.
Next, suppose that $\omega=1$ and $\omega'\not\in\RR$.
The equality $(-\omega')\omega''=1$ now implies $\Imag\omega'=\Imag\omega''$.
When both of these imaginary parts are positive, we must have $\Arg\omega'=\pi-\Arg\omega''$, so that
\[
    |\Arg\omega|+
    |\Arg\omega'|+
    |\Arg\omega''|
    =\pi + |\Arg\omega| = \pi.
\]
The case of $\omega=1$ and imaginary parts of $\omega'$ and $\omega''$ both
negative follows by complex conjugation.

We now show that whenever $\tet$ is suppressing, we have a strict inequality
in \eqref{tetrahedral-suppression-inequality}.
Consider first the case when one of the angles (without loss of generality,
$\omega$) is real, so that we have $\omega=-1$.
Since condition~\eqref{cond:angle=1} does not hold, neither $\omega'$ nor
$\omega''$ can equal $1$ either,
which implies that $|\Arg\omega'|, |\Arg\omega''| > 0$ and
\[
    |\Arg\omega|+
    |\Arg\omega'|+
    |\Arg\omega''|
    > |\Arg\omega| = \pi.
\]
It remains to consider the case when none of the angles is real but condition
\eqref{cond:same_signs} in Definition~\ref{definition-suppressing} fails.
It suffices to consider $\Imag\omega>0$, $\Imag\omega'<0$, $\Imag\omega''<0$;
the case of one negative and two positive imaginary parts will then follow by
complex conjugation.
Since $(-\omega)\omega'\omega''=1$, we have $\Imag(\omega'\omega'')=\Imag\omega>0$.
Hence, $\Arg\omega', \Arg\omega''\in(-\pi, -\frac{\pi}{2}]$ and
\[
    |\Arg\omega|+
    |\Arg\omega'|+
    |\Arg\omega''|
    \geq |\Arg\omega| + \tfrac{\pi}{2} + \tfrac{\pi}{2} > \pi.\qedhere
\]
\end{proof}

\begin{corollary}\label{cor:nonsupp-closed}
The subset of $\SAS(\triang)$ consisting of non-suppressing angle structures is closed.
\end{corollary}

%==============================================================================
% Section 5: Asymptotic analysis of the meromorphic 3D-index
%==============================================================================
\section {Asymptotic analysis of the meromorphic 3D-index}

In this section, we study the behaviour of the state integral \eqref{3Dindex_N-1}
defining the meromorphic \ind~$\I{\triang}{\aholM}{\aholL}(\hbar)$ of
a triangulation~$\triang$ admiting a strict angle structure~$\alpha$ with the
prescribed peripheral angle-holonomy
$(\aholM,\aholL)=\bigl(\ahol{\alpha}(\mu), \ahol{\alpha}(\lambda)\bigr)$.
We set $\hbar = -1/\kappa$, where $\kappa$ is a positive real parameter and we let
$\kappa\to+\infty$.

The equation~\eqref{3Dindex_integral_over_SAS}
allows us to treat $\I{\triang}{\aholM}{\aholL}(\hbar)$ as
an integral over the
domain~$\Omega$ arising as the connected component of
$\SAS_{(\exp(i\aholM),\exp(i\aholL))}(\triang)$
containing the angle structure~$\exp(i\alpha)$.
By Corollary~\ref{cor:nonsupp-closed}, \Sval\ angle structures without
suppressing tetrahedra form a closed subset $\Omega_+\subset\Omega$.
It is easy to see that its complement~$\Omega_- := \Omega\setminus\Omega_+$
has an exponentially small contribution
to the $\kappa\to+\infty$ asymptotics of the meromorphic \ind.
Indeed, consider an angle structure $\omega_-\in\Omega_-$.
Since $\Omega_-$ is open, there exists an $\eps>0$ and an open neighbourhood
$U\ni\omega_-$ with $\overline{U}\subset\Omega_-$ and such that
\[
        \sum_{\square\in Q(\triang)} |\Arg\omega(\square)| > N\pi + \eps
\]
for all $\omega\in\overline{U}$.
Using \eqref{suppression-analysis-pertet},
we estimate
\begin{align*}
    \abs{
    \int_{\overline{U}} c_q^N\mkern-8mu\prod_{\square\in Q(\triang)}\mkern-5mu
    G_q\bigl(\omega(\square)e^{a(\square)\hbar}\bigr)\,d\omega}
    &\leq
    \int_{\overline{U}}
    \exp\Bigl[
        \Bigl(
        -N\pi+\mkern-8mu\sum_{\square\in Q(\triang)}\mkern-5mu|\Arg\omega(\square)|
    \Bigr)\frac{\pi}{2\hbar} + O\bigl(\log(-\hbar)\bigr)\Bigr]d\omega\\
    &\leq
    \int_{\overline{U}} \exp\Bigl[\frac{-\pi\eps}{2}\kappa + O(\log\kappa)\Bigr]\,d\omega
    \longrightarrow 0\ \text{as}\ \kappa\to+\infty.
\end{align*}
As a consequence, any nonzero contribution
to the asymptotics of the integral~$\I{\triang}{\aholM}{\aholL}(\hbar)$
must come from the non-suppressed closed region $\Omega_+$,
including small neighbourhoods of points on its boundary.

%--------------------------------------------------------------------------------------------------
\subsection{The contribution of circle-valued angle structures without real angles}\label{nonreal-contrib}

Suppose that $D\subset\overline{D}\subset\Omega_+$ is a non-empty open
subset consisting entirely of \Sval\ angle structures
satisfying condition~\eqref{cond:same_signs} of Definition~\ref{definition-suppressing}.
In other words, for every $\omega\in D$ and every tetrahedron~$\tet$ of $\triang$
with normal quads $\square$, $\square'$, $\square''$,
we require the three imaginary parts~$\Imag\omega(\square)$, $\Imag\omega(\square')$ and $\Imag\omega(\square'')$
to be either all positive or all negative,
but we do not require different tetrahedra to be oriented consistently.
The goal of this section is to study the $\kappa\to+\infty$ asymptotic behaviour of the
localised integral
\begin{equation}\label{integral-over-D-smooth}
    I_D(\kappa) =
    \frac{d_\triang c_q^N}{(2\pi)^{N-1}}
    \int_D \prod_{\square\in Q(\triang)}
    G_q\bigl(\omega(\square)e^{a(\square)\hbar}\bigr)\;d\omega
    =
    \frac{d_\triang}{(2\pi)^{N-1}} \int_D b_\omega(-1/\kappa)\;d\omega,
\end{equation}
where we have used the notation of \eqref{integrand-as-function-of-hbar} and where $\hbar = -1/\kappa$.
Applying Theorem~\ref{th:asym-c_q-and-G_q} and Lemma~\ref{suppression-inequality-lemma},
we obtain
\begin{multline}\label{asym-approx-smooth-integrand}
    b_\omega(-1/\kappa)
    =
    \exp\Bigl[
    i\kappa\Vol(\omega)
    +\tfrac{N}{2}\log\kappa
    +N\log\sqrt{4\pi}\\
        +\sum_{\square\in Q(\triang)}\bigl(a(\square)-\tfrac{1}{2}\bigr)
        \bigl(\log|2\Imag\omega(\square)|-i\tfrac{\pi}{2}\sgn\Imag\omega(\square)\bigr)
    + O\bigl(\tfrac{1}{\kappa}\bigr)
    \Bigr],
\end{multline}
where $\Vol(\omega)$ refers to the volume function of~\eqref{volume-on-SAS}.
Ignoring the error term $O\bigl(\frac{1}{\kappa}\bigr)$ and exponentiating, we
obtain an approximate expression~$\widetilde{b_\omega}(\kappa)$ given by the formula
\begin{multline}\label{initial-expression-for-smooth-approx}
    \widetilde{b_\omega}(\kappa)
    =
    \left(\sqrt{4\pi\kappa}\right)^N
    \prod_{\square\in Q(\triang)}
    \abs{2\Imag\omega(\square)}^{a(\square)-\frac{1}{2}}\\
    \times \exp\Bigl[
    -i\tfrac{\pi}{2}\sum_{\square\in Q(\triang)}\sgn\Imag\omega(\square)\bigl(a(\square)-\tfrac{1}{2}\bigr)
    \Bigr]
    e^{i\kappa\Vol(\omega)}.
\end{multline}
Let $N_+$ be the number of ``positively oriented'' tetrahedra,
i.e., tetrahedra~$\tet$ with $\Imag\omega(\square)>0$ for all $\square\subset\tet$,
and correspondingly let $N_-=N-N_+$ be the number of ``negatively oriented'' tetrahedra.
Since $a(\square)+a(\square')+a(\square'')=1$, we have
\[
    \sum_{\square\in Q(\triang)}\sgn\Imag\omega(\square)\bigl(a(\square)-\tfrac{1}{2}\bigr)
    = -\tfrac{1}{2}(N_+-N_-),
\]
which allows us to rewrite \eqref{initial-expression-for-smooth-approx} as
\begin{equation}\label{nice-approx-smooth-integrand}
    \widetilde{b_\omega}(\kappa)
    = \left(\sqrt{2\pi\kappa}\right)^N
    \exp\bigl[i\tfrac{\pi}{4}(N_+-N_-)\bigr]
    \prod_{\square\in Q(\triang)}
    \abs{\Imag\omega(\square)}^{a(\square)-\frac{1}{2}}
    e^{i\kappa\Vol(\omega)}.
\end{equation}
We are now going to insert the above approximation into \eqref{integral-over-D-smooth},
and it is at this point that our reasoning becomes merely heuristic, since the integral
of an approximation may not yield an approximation of the integral.
Explicitly, our heuristic approximation has the form
\begin{equation}\label{heuristic-smooth-integral}
    \widetilde{I_D}(\kappa) =
    \frac{d_\triang}{(2\pi)^{N-1}} \int_D \widetilde{b_\omega}(\kappa)\;d\omega.
\end{equation}
Since the real-valued volume function $\Vol(\omega)$ occurring as the phase
factor in \eqref{nice-approx-smooth-integrand} is smooth on $D$, a $\kappa\to+\infty$
asymptotic approximation of the above integral can be found
with the help of the principle of stationary phase.

Suppose that $\omega\in D$ is a non-degenerate critical point of $\Vol\bigr|_{\Omega_+}$.
By the work of Luo~\cite{luo-volume}, the shape parameters defined via \eqref{from-omega-to-z}
satisfy complex gluing equations, and therefore produce a conjugacy class of a
representation $\pi_1(M)\to\PSLC$ with peripheral angle-holonomy $(\aholM,\aholL)$.
Hence, leading order asymptotic terms shall originate from conjugacy classes of representations.
We are now going to write these terms down explicitly.

Using equation~\eqref{3Dindex-int-2}, we can analyse the asymptotic contribution
of a small neighbourhood of the critical point $\omega$ to $\widetilde{I_D}(\kappa)$ with the help of
the localised integral
\begin{equation}\label{epsilon-integral}
    I_{\omega,\eps}(\kappa)=
    \frac{d_\triang}{(2\pi)^{N-1}}
    \int\limits_{|t|\leq\eps}
    \widetilde{b_{\omega_t}}(\kappa)\,dt,\quad
     \omega_t(\square) = \omega(\square) e^{i(t\mathbf{L}_*)(\square)},
\end{equation}
where $t = (t_1,\dotsc,t_{N-1})$, $\mathbf{L}_*$ is given in equation~\eqref{definition-L_star},
 and $\eps>0$ is small enough that $\omega_t$ is a smooth point of $\Vol\bigr|_{\Omega_+}$
 and that $\omega=\omega_0$ is the only critical point in the domain of integration.

We are now going to approximate the integral $\eqref{epsilon-integral}$
with the stationary phase contribution of the critical point.
Recall from \eqref{definition-L_star} that the $(N-1)\times3N$ matrix $\mathbf{L}_*$
contains the quad entries of the
leading-trailing deformation vectors~$\{l_1,\dotsc,l_{N-1}\}$, where we are canonically
excluding the redundant last vector~$l_N$.
Since the Lobachevsky function satisfies
\[
    \Lob'(\theta) = -\log|2\sin\theta|, \quad \Lob''(\theta) = -\cot\theta\quad (\theta\not\in\pi\ZZ),
\]
the Hessian matrix of the volume function $\Vol(\omega_t)$ at $t=0$ is given explicitly as
\begin{align}
    \frac{\del^2\Vol(\omega_t)}{\del t_m \del t_n}(0)
        &= -\sum_{\square\in Q(\triang)} l_m(\square)l_n(\square) \cot(\arg \omega(\square))\notag\\
        &= -\bigl[\mathbf{L_*}\diag\{\cot\arg\omega(\square)\}_\square\mathbf{L_*^\transp}\bigr]_{m,n},
        \label{computation-of-Hessian}
\end{align}
where $\diag\{\cot\arg\omega(\square)\}$ is a $3N\times3N$ real matrix
whose diagonal entries are the cotangents of the angles.
Denote by $\Sigma(\omega)$ the signature of the Hessian matrix \eqref{computation-of-Hessian}.
As $\kappa\to+\infty$, we have the stationary phase approximation
\begin{align}
    I_{\omega,\eps}(\kappa)&\sim
    \frac{% numerator:
        d_\triang\left(\sqrt{2\pi\kappa}\right)^N\mkern-3mu
        \exp\bigl[i\tfrac{\pi}{4}\bigl(N_+ - N_- + \Sigma(\omega)\bigr)\bigr]
        \left(\sqrt{\frac{2\pi}{\kappa}}\right)^{N-1}\mkern-5mu
        \prod\limits_{\square\in Q(\triang)}
            \abs{\Imag\omega(\square)}^{a(\square)-\frac{1}{2}}
    }
    {% denominator:
        (2\pi)^{N-1}
        \sqrt{\abs{\det\bigl[\mathbf{L_*}\diag\{\cot\arg\omega(\square\}_\square\mathbf{L_*^\transp}\bigr]}}
    }
    e^{i\kappa\Vol(\omega)}\notag\\
    &= \label{messy-tau}
        \frac{ d_\triang
        \prod\limits_{\square\in Q(\triang)}
            \abs{\Imag\omega(\square)}^{a(\square)-\frac{1}{2}}
        }
        {
        \sqrt{\abs{\det\bigl[\mathbf{L_*}\diag\{\cot\arg\omega(\square\}_\square\mathbf{L_*^\transp}\bigr]}}
        }
    \exp\left[i\tfrac{\pi}{4}\bigl(N_+ - N_- +\Sigma(\omega)\bigr)\right]
    \sqrt{2\pi\kappa}\,
    e^{i\kappa\Vol(\omega)}.
\end{align}

We remark that when all tetrahedra are positively oriented, then
Lemma~5.3 of Futer-Gu\'eritaud \cite{futer-gueritaud} implies that $\omega$
must be a maximum point of the volume function, so that $\Sigma(\omega)=1-N$, $N_+=N$, and $N_-=0$.
Therefore, the integer $N_+ - N_- + \Sigma(\omega)$ is equal to $1$.
Likewise, when all tetrahedra are negatively oriented, we find $N_+ - N_- + \Sigma(\omega) = -1$.
\begin{question}\label{mysterious-integer}
Assuming that $\omega$ is a non-degenerate smooth critical point of the volume function,
does there exist an expression for the integer $N_+ - N_- + \Sigma(\omega)$ in terms of the
representation given by the shape parameters \eqref{from-omega-to-z}?
\end{question}

We are now going to shorten the expression \eqref{messy-tau} by introducing more compact notation.
Since the shape parameters $z(\square)$ defined by \eqref{from-omega-to-z}
satisfy complex gluing equations, we may write down \eqref{messy-tau}
starting with an algebraic solution~$z\in(\CC\setminus\RR)^{Q(\triang)}$
whose associated circle-valued angle
structure $\omega=\omega_z$, given by $\omega(\square)=z(\square)/|z(\square)|$,
lies in $\Omega$.
We therefore define
\begin{equation}\label{first-formula-for-tau}
    \tau = \tau(z,\alpha)
         := \frac{\prod\limits_{\square\in Q(\triang)}
        |\Imag \omega(\square)|^{a(\square)-\frac{1}{2}}}
        {\sqrt{\left|\det\bigl(\mathbf{L_*}
            \diag\{\frac{\Real\omega(\square)}{\Imag\omega(\square)}\}
        \mathbf{L_*^\transp} \bigr) \right|}},
    \quad a(\square) = \frac{\alpha(\square)}{\pi}.
\end{equation}
With the above notation, we have
\begin{equation}\label{contribution-of-smooth-crit}
    I_{\omega,\eps}(\kappa)\sim
    d_\triang\,\tau\,\exp\left[i\tfrac{\pi}{4}\bigl(N_+ - N_- +\Sigma(\omega)\bigr)\right]
    \sqrt{2\pi\kappa}\,e^{i\kappa\Vol(\omega)}.
\end{equation}
Although the approximation \eqref{heuristic-smooth-integral} is heuristic,
we do expect the expression in \eqref{contribution-of-smooth-crit} to describe the contribution
of a smooth non-degenerate critical point~$\omega$ of the function
\[
    \Vol\bigr|_{(\aholM,\aholL)}:\SAS_{(\exp(i\aholM),\exp(i\aholL))}(\triang)\to\RR
\]
to the $\kappa\to+\infty$ asymptotics of $\I{\triang}{\aholM}{\aholL}(-1/\kappa)$.
Numerical evidence for this claim is provided in Section~\ref{sec:numerical} below.
See also Section~\ref{sec:example-4_1} below for an explicit example.

%--------------------------------------------------------------------------------------------------
\subsection{\texorpdfstring{The contribution of $\ZZ_2$-taut angle structures}%
{The contribution of Z2-taut angle structures}}\label{sec:derivation-of-MB}
In this section, we assume that $\omega\in\Omega$
is a \Ztaut\ angle structure.
In other words, in every tetrahedron of $\triang$, $\omega$ equals $+1$ on
exactly two normal quadrilaterals and $-1$ on the remaining one.
Our goal is to study the contribution of a small neighbourhood of $\omega$ to the
$\kappa\to+\infty$ asymptotics of $\I{\triang}{\aholM}{\aholL}(-1/\kappa)$.
Observe that when $\omega=e^{i\alpha}$ for $\alpha\in\Areal(\triang)$,
then $\alpha(\square)\in\pi\ZZ$ for all $\square\in Q(\triang)$.
Therefore, the discussion in this section is relevant to $\aholM,\aholL\in\pi\ZZ$,
the most important case being that of the vanishing angle-holonomy, $\aholM=\aholL=0$.

For a fixed tetrahedron $\tet$ of $\triang$, we define the sets
\begin{alignat*}{2}
    Q_+(\triang) &= \{\square\in Q(\triang) : \omega(\square)=1\},\qquad&
    Q_-(\triang) &= \{\square\in Q(\triang) : \omega(\square)=-1\};\\
    Q_+(\tet) &= \{\square\in Q_+(\triang): \square\subset\tet\},\qquad&
    Q_-(\tet) &= \{\square\in Q_-(\triang): \square\subset\tet\}.
\end{alignat*}
Observe that $|Q_+(\tet)|=2$ for every tetrahedron, so condition~\eqref{cond:angle=1}
of Definition~\ref{definition-suppressing} is satisfied in every tetrahedron
and $\omega$ is never exponentially suppressed.

\begin{lemma}\label{beta-expansion-lemma}
Consider arbitrary complex numbers $A$, $A'$ and~$A''$, with positive real parts,
satisfying $A+A'+A''=1$.
When $\kappa$ is a positive real parameter and when $\hbar = -1/\kappa$ and $q=e^\hbar$,
the quantity
\begin{equation}\label{factor-at-Ztaut}
    f(\kappa)
    =
    c_q
    G_q\bigl(-e^{A\hbar}\bigr)
    G_q\bigl(e^{A'\hbar}\bigr)
    G_q\bigl(e^{A''\hbar}\bigr)
\end{equation}
has the asymptotic approximation
\begin{equation}\label{Beta-expansion}
    f(\kappa)
    =\kappa
    \Beta(A',A'')\bigl(1+O(\kappa^{-2})\bigr)\ \text{as}\ \kappa\to+\infty,
\end{equation}
where
\begin{equation}\label{definition-Euler-Beta}
    \Beta(z_1, z_2) = \frac{\Gamma(z_1)\Gamma(z_2)}{\Gamma(z_1+z_2)}
\end{equation}
denotes Euler's beta function.
Moreover, \eqref{Beta-expansion} holds uniformly
for $A$, $A'$ (and therefore $A''$) on compact
subsets of the right half-plane~$\{z\in\CC: \Real z>0\}$.
\end{lemma}

Suppose that $\tet$ is a tetrahedron of $\triang$ and that
$\omega$ is a fixed \Ztaut\ angle structure on $\triang$.
Setting $\{\square\}=Q_-(\tet)$ and putting
$A=a(\square)$, $A'=a(\square')$, $A''=a(\square'')$,
we see that term $\eqref{factor-at-Ztaut}$
is the value of the Boltzmann weight of $\Delta$ in
the state integral of the meromorphic \ind\ at the point $\omega$.

\begin{proof}
Using Theorem~\ref{th:asym-c_q-and-G_q} and
the expression $1-A=A'+A''$, we compute
\begin{align*}
    f(\kappa)
    &=
    \exp\Bigl[
        \log\kappa
        - \log\frac{\Gamma\bigl(1-A\bigr)}{\sqrt{2\pi}}
        + \log\frac{\Gamma\bigl(A'\bigr)}{\sqrt{2\pi}}
        + \log\frac{\Gamma\bigl(A''\bigr)}{\sqrt{2\pi}}\\
    &\phantom{=\exp\Bigl[}
        + \log\sqrt{4\pi} - \tfrac{1}{2}\log2 + O(\kappa^{-2})
    \Bigr]\\
    &=\kappa
    \frac{
        \Gamma\bigl(A'\bigr)\Gamma\bigl(A''\bigr)
    }
    {
         \Gamma\bigl(1-A\bigr)
    }
    \bigl(1+O(\kappa^{-2})\bigr)\\
    &= \kappa
    \Beta\bigl(A',A''\bigr)
    \bigl(1+O(\kappa^{-2})\bigr),\ \text{as}\ \kappa\to+\infty,
\end{align*}
where we have used the rule $\exp\bigl[O\bigl(\kappa^{-2}\bigr)\bigr]=1+O\bigl(\kappa^{-2}\bigr)$.
\end{proof}

In order to study the contribution of a small neighbourhood of the \Ztaut\ structure
$\omega$ to the state integral defining the meromorphic \ind, we consider again
a localised integral of the form
\[
    I_{\omega,\eps}(\kappa)
    = \frac{c_q^N}{(2\pi)^{N-1}} \int\limits_{\abs{t}\leq\eps}
    \prod_{\square\in Q(\triang)}
        G_q\left(
        e^{a(\square)\hbar} \omega(\square) e^{i(t\mathbf{L}_*)(\square)}
    \right) dt_1\dotsi dt_{N-1},
\]
where $t=(t_1,\dotsc,t_{N-1})$.
In analogy to \eqref{epsilon-integral}, we assume that $\eps>0$ is small enough so that
the set
\(
    \{\square\mapsto \omega(\square) e^{i (t\mathbf{L}_*)(\square)}: |t|\leq\eps\}\subset\Omega
\)
does not contain any \Ztaut\ structures other than $\omega$, nor any
smooth critical points of the volume function.
We perform the change of variables $t_m = r_m/\kappa$ for all $m\in\{1,\dotsc,N-1\}$,
obtaining
\begin{alignat*}{2}
    I_{\omega,\eps}(\kappa)
    &= \frac{c_q^N}{(2\pi\kappa)^{N-1}}
        &\int\limits_{|r|\leq\kappa\eps}
    &\prod_{\square\in Q(\triang)}
        G_q\Bigl(\omega(\square)
        \exp\bigl[ a(\square)\hbar +
            i (r\mathbf{L}_*)(\square)/\kappa\bigr]
        \Bigr)
    \,dr_1\dotsi dr_{N-1}\\
    &= \frac{c_q^N}{(2\pi\kappa)^{N-1}}
    &\int\limits_{|r|\leq\kappa\eps}
    \Biggl[
    &\prod_{\square\in Q_+(\triang)}
    G_q\bigl(e^{\hbar(a(\square) - i(r\mathbf{L}_*)(\square))}\bigr)\\
    &\quad&\times
    &\prod_{\square\in Q_-(\triang)}
    G_q\bigl(-e^{\hbar(a(\square) - i(r\mathbf{L}_*)(\square))}\bigr)\Biggr]
    \,dr_1\dotsi dr_{N-1}.
\end{alignat*}
Since the numbers~$A(\square) := a(\square) -i(r\mathbf{L}_*)(\square)$
satisfy $A(\square)+A(\square')+A(\square'')=1$ in every tetrahedron,
we may apply Lemma~\ref{beta-expansion-lemma} for each of the tetrahedral weights, obtaining
the approximation
\begin{align}
    \widetilde{I_{\omega,\eps}}(\kappa)
    &=
    \frac{1}{(2\pi\kappa)^{N-1}}
        \int\limits_{|r|\leq\kappa\eps}
        \kappa^N
        \frac{\prod_{\square\in Q^+(\triang)} \Gamma\bigl(A(\square)\bigr)}
        {\prod_{\square\in Q^-(\triang)} \Gamma\bigl(1-A(\square)\bigr)}
        \,dr\notag\\
    &=
    \frac{\kappa}{(2\pi)^{N-1}}
        \int\limits_{|r|\leq\kappa\eps}
        \prod_{\substack{\tet\ \text{in}\ \triang,\\ Q_+(\tet)=\{\square,\square'\}}}
        \mkern-10mu
        \Beta\bigl(a(\square)-i(r\mathbf{L}_*)(\square),\;a(\square')-i(r\mathbf{L}_*)(\square')\bigr)
        \,dr\label{beta-integral-on-compact}.
\end{align}
We stress that this approximation is heuristic, since we have replaced the integrand
with its asymptotic approximation without estimating the resulting error terms.
Given that we are interested in an asymptotic approximation as $\kappa\to\infty$, we would now like
to extend the limits of integration to infinity.
Observe that the product occurring under the integral~\eqref{beta-integral-on-compact} can
be rewritten as follows.
When $\tet$ is a tetrahedron of $\triang$ containing normal quads~$\square,\square'\subset\tet$
satisfying $\omega(\square)=\omega(\square')=+1$, we set
\begin{equation}\label{beta-weight-minus}
    \Betaweight(s_1,\dotsc,s_{N-1})
    = \Beta\left(\tfrac{\alpha(\square)}{\pi}  - \sum_{m=1}^{N-1} s_m l_m(\square),\;
                  \tfrac{\alpha(\square')}{\pi} - \sum_{m=1}^{N-1} s_m l_m(\square')\right).
\end{equation}
We define $\MB(\triang,\omega,\alpha)$ as Cauchy's principal
value of the following multivariate contour integral of
Mellin--Barnes type:
\begin{equation}\label{def:MB}
    \MB(\triang,\omega,\alpha) =
    \left(\frac{1}{2\pi i}\right)^{N-1}
    \pvint_{-i\infty}^{i\infty}
    \dotsi
    \int_{-i\infty}^{i\infty}
        \prod_{\tet\ \text{in}\ \triang}
        \Betaweight(s_1,\dotsc,s_{N-1})
        \,ds_1\dotsi ds_{N-1}.
\end{equation}
Observe that the contour of integration is symmetric with respect to complex
conjugation and that $\Beta(\bar{x},\bar{y})=\overline{\Beta(x,y)}$.
Hence, if the integral in \eqref{def:MB} converges at least in
the sense of Cauchy's principal value,
then $\MB(\triang,\omega,\alpha)$ is a real number
equal to the $\kappa\to\infty$ limit of the integral occurring in
\eqref{beta-integral-on-compact}.
In this case, we expect the localised integral $\widetilde{I_{\omega,\eps}}(\kappa)$
to have the linear asymptotic approximation
\(
    \widetilde{I_{\omega,\eps}}(\kappa)
    \sim \kappa\cdot\MB(\triang,\omega,\alpha)
\) as $\kappa\to\infty$.
In particular, this linear term does not depend on the chosen $\eps$.
As before, taking the multiplicity $d_\triang$ into account,
we expect that the contribution of an infinitesimal neighbourhood of $\omega$ to the
asymptotics of the meromorphic \ind\ is given by
\begin{equation}\label{MB-linear-approx}
    d_\triang\MB(\triang,\omega,\alpha)\cdot\kappa
    \quad\text{as}\ \kappa\to\infty.
\end{equation}
Note that Mellin--Barnes integrals similar to those of \eqref{def:MB} were
previously hinted at by Kashaev, Luo and Vartanov~\cite{kashaev-luo-vartanov}.
We shall study these integrals in more detail in
Section~\ref{sec:Mellin-Barnes}, where we discuss their gauge invariance properties.

%--------------------------------------------------------------------------------------------------
\subsection{Example: the figure-eight knot complement}\label{sec:example-4_1}
We are going to illustrate the results of this section
in the case of the  standard two-tetrahedron triangulation~$\triang$ of the figure-eight knot
complement~$M = S^3\setminus 4_1$.
The full gluing matrix is given by SnapPy~\cite{snappy} as
\[
    \begin{bmatrix}\mathbf{G}\\\mathbf{G_\del}\end{bmatrix} =
    \begin{bmatrix}G & G' & G''\\G_\del & G'_\del & G''_\del\end{bmatrix} =
    \left[\begin{array}{cc|cc|cc}
        2 & 2 & 1 & 1 & 0 & 0 \\
        0 & 0 & 1 & 1 &    2 & 2 \\
        -1& 0 & 0 & 0 & 0 & 1 \\
        -1& -1&-1 & 1 &-1 & 3
    \end{array}\right].
\]
The geometric solution of Thurston's gluing equations recovering the complete
hyperbolic structure of finite volume is given by $z_1=z_2=e^{i\pi/3}$.
The associated angle structure $\alpha^\geo\in\Astrict_0(\triang)$ has
$\alpha^\geo(\square)=\pi/3$ for all $\square\in Q(\triang)$.
Using
\[
    \mathbf{L_*} = [L_*\;|\;L'_*\;|\;L''_*] = [{-1}\ {-1}\;|\; 2\ 2 \;|\; {-1}\ {-1}],
\]
we obtain the parametrisation
\begin{align*}
    [-\pi, \pi] &\to \Omega\subset\SAS_0(\triang),\\
    t &\mapsto \omega_t = \Bigl(
    e^{i(\frac{\pi}{3}-t)},\;
    e^{i(\frac{\pi}{3}-t)}\;\Bigl|\;
    e^{i(\frac{\pi}{3}+2t)},\;
    e^{i(\frac{\pi}{3}+2t)}\;\Bigr|\;
    e^{i(\frac{\pi}{3}-t)},\;
    e^{i(\frac{\pi}{3}-t)}
    \Bigr)
\end{align*}
of degree $d_\triang = |H_1(\hat{M};\FF_2)| = 1$.
In terms of the parameter~$t$, the volume function has the expression
\begin{align*}\label{4_1-volume}
    V(t) = \Vol(\omega_t) &= 4\Lob\bigl(\tfrac{\pi}{3}-t\bigr) + 2\Lob\bigl(\tfrac{\pi}{3}+2t\bigr)\\
    &= 4\Lob\bigl(\tfrac{\pi}{3}-t\bigr) - 2\Lob\left(2\bigl(\tfrac{\pi}{3}-t\bigr)\right)\\
    &= 4\Lob\bigl(\tfrac{\pi}{6}+t\bigr),
\end{align*}
where the last equality follows from Lemma~7.1.4 of~\cite{thurston-notes}.
We see that $V(t)$ is differentiable whenever
$t\not\in\bigl\{-\tfrac{\pi}{6}, \tfrac{5\pi}{6}\bigr\}$ with derivative
\(
    V'(t) = -4\log\abs{2\sin\bigl(t+\tfrac{\pi}{6}\bigr)}
\).
Setting $V'(t)=0$, we find two non-degenerate critical points: $t=0$ and $t=\frac{2\pi}{3}$.

At $t=0$, we recover the (exponentiated) hyperbolic angle
structure~$\omega^\geo = \exp(i\alpha^\geo)$ and consequently
\[
    \tau = \frac{\prod\limits_{\square\in Q(\triang)}
            |\sin\tfrac{\pi}{3}|^{\frac{1}{3}-\frac{1}{2}}}
            {\sqrt{\left|\det\bigl(\mathbf{L_*}
                \diag\{\cot\tfrac{\pi}{3}\}
            \mathbf{L_*^\transp} \bigr) \right|}}
    = \frac{|\sin\tfrac{\pi}{3}|^{-1}}{
    \sqrt{\abs{\cot\tfrac{\pi}{3}\det\bigl(\mathbf{L_*}\mathbf{L_*}^\transp\bigr)}}}
    = 3^{-3/4}.
\]
Since this point is the global maximum of $V$, we have $\Sigma=-1$, $N_+=2$, $N_-=0$, so that
the approximation~\eqref{contribution-of-smooth-crit} becomes
\begin{equation}\label{4_1:positive_contrib}
    3^{-3/4}\,\exp\bigl(i\tfrac{\pi}{4}\bigr)
    \sqrt{2\pi\kappa}\,e^{i\kappa\Volhyp(M)},
\end{equation}
where $\Volhyp(M)=2.0298832128193\dots$ is the volume of $M$ in the complete hyperbolic metric.

At $t=\frac{2\pi}{3}$, we see the angle structure $\overline{\omega}^\geo = \exp(-i\alpha^\geo)$
which also corresponds to the complete hyperbolic structure, but with the negative orientation.
Since $\tau$ is the same in this case, but $\Sigma=1$, $N_+=0$, $N_-=2$ instead, we obtain the term
\begin{equation}\label{4_1:negative_contrib}
    3^{-3/4}\,\exp\bigl(-i\tfrac{\pi}{4}\bigr)
    \sqrt{2\pi\kappa}\,e^{-i\kappa\Volhyp(M)}.
\end{equation}
Therefore, the expected total contribution of smooth critical points to the $\kappa\to\infty$
asymptotics of $\I{M}{0}{0}(-1/\kappa)$ is
\begin{equation}\label{4_1:total_contrib}
    \frac{2\sqrt{2\pi}}{\sqrt[4]{27}}\;\sqrt{\kappa}\cos\bigl(\kappa\Volhyp(M)+\tfrac{\pi}{4}\bigr).
\end{equation}

We now search for angle structures with values in $\{+1,-1\}$.
Since $\exp\bigl(i(\frac{\pi}{3}-t)\bigr)=\pm1$ only when $t\in\bigl\{-\frac23\pi, \frac{\pi}{3}\bigr\}$,
we find two such structures:
\[
    \omega_{-2\pi/3} = \bigl(-1, -1\;\bigl|\; -1, -1\;\bigr|\; -1, -1\bigr),\quad
    \omega_{\pi/3} = \bigl(1, 1\;\bigl|\; -1, -1\;\bigl|\; 1, 1\bigr)
\]
of which the first, $\omega_{-2\pi/3}$, is exponentially suppressed by
virtue of Proposition~\ref{proposition-suppression}.
The second structure~$\omega=\omega_{\pi/3}$ is \Ztaut.
According to \eqref{beta-weight-minus}, we have
$\Betaweight(s)    = \Beta\bigl(\tfrac{1}{3} + s ,\,\tfrac{1}{3} + s\bigr)$
for both tetrahedra $\tet$ of $\triang$.
The associated Mellin--Barnes integral is therefore
\begin{align*}
    \MB(\triang,\omega,\alpha)
    &= \frac{1}{2\pi i}
        \int_{-i\infty}^{i\infty}
        \Beta^2\bigl(\tfrac{1}{3} + s ,\; \tfrac{1}{3} + s\bigr)
        \,ds
            = \frac{1}{2\pi i}
        \int_{\tfrac{1}{3}-i\infty}^{\tfrac{1}{3}+i\infty}
        \Beta^2\bigl(s ,\; s\bigr)
        \,ds.
\end{align*}
In fact,  since the integrand is holomorphic when $\Real s>0$,
we have, for any $c>0$, 
\[
    \MB(\triang,\omega,\alpha)
    =  \frac{1}{2\pi i}\int_{c-i\infty}^{c+i\infty}
    \Beta^2\bigl(s ,\;  s\bigr)
        \,ds.
\]
We now wish to push the contour of integration to the right, by letting $c\to+\infty$.
Stirling's formula for the gamma function leads to the asymptotic expansion
\[
    \Beta\bigl(s ,\;  s\bigr) = 2 \pi^{1/2} 4^{-s} s^{-1/2} (1+O(s^{-1})),
    \ \text{as}\ |s|\to\infty\ \text{with}\ |\Arg s| \le \pi - \epsilon,
\]
for any small $\epsilon >0$ (see equation (2.2.4) in Paris--Kaminski~\cite{paris-kaminski}).
Thus
\begin{equation}\label{explicit-MB-4_1}
    {2\pi i}\,  \MB(\triang,\omega,\alpha)
    =  4\pi   \int_{c-i\infty}^{c+i\infty}
        16^{-s} s^{-1} \left(1 + O\bigl(s^{-1}\bigr)\right)\,ds
\end{equation}
as $c=\Real s\to+\infty$.
For any real $x$ with $0 < x < 1$, equation~(3.3.20) in \cite{paris-kaminski} yields
\begin{equation*}
    \int_{c-i\infty}^{c+i\infty} x^{s} s^{-1}\,ds = 0,
    \quad\text{for}\ c>0;
\end{equation*}
hence
\[
    \int_{c-i\infty}^{c+i\infty}16^{-s}s^{-1}\,ds = 0.
\]
In order to estimate the error term in \eqref{explicit-MB-4_1}, we note that 
\[
    \abs{\int_{c-i\infty}^{c+i\infty} 16^{-s} s^{-2} \, ds}
    \le \int_{-\infty}^{\infty} \frac{16^{-c}}{c^2+y^2} dy
    = \frac{16^{-c} \pi}{c}\to0\ \text{as}\ c\to\infty.
\]
So we conclude that $\MB(\triang,\omega_{\pi/3},\alpha)=0$.

Thanks to the vanishing of the above Mellin--Barnes integral,
the expected asymptotic approximation of $\I{M}{0}{0}(-1/\kappa)$
as $\kappa\to+\infty$ reduces to the expression given in \eqref{4_1:total_contrib}.
We refer to Section~\ref{sec:numerix} for a comparison of this approximation
with values found numerically.

\subsection{Example: an infinite family of boundary-parabolic representations}
\label{example-infinite-family}
While the asymptotic analysis described above works for isolated, non-degenerate critical
points of the volume function, there may also exist infinite families of conjugacy classes 
of boundary-parabolic \nword{\PSLC}{representations} with the geometric obstruction class.
Such families may give rise to higher-dimensional critical loci of the volume function
on $\SAS_0(\triang)$.

To construct an explicit example of such infinite family, consider the free product
$\ZZ_3 * \ZZ\cong\langle x,y\;|\; x^3\rangle$ and the family of representations
$\psi_s:\ZZ_3*\ZZ\to\PSLC$ given in terms of the parameter $s\neq0$ by
\[
    \psi_s(x) = \pm\begin{bmatrix}
        -s^{-1} & \frac{1 + s -3s^2 - P(s)}{s^3 Q(s)} \\
        s\frac{s^2+s-3-P(s)}{Q(s)} & \frac{1-s}{s}
    \end{bmatrix},
    \quad
    \psi_s(y) = \pm\begin{bmatrix}
        0 & -s^{-1}\\
        s & \frac{-Q(s)}{2s} \\
    \end{bmatrix},
\]
where $P(s)^2 = s^4 + 2s^3 - 13s^2 + 2s + 1$ and $Q(s)^2 = 2\bigl(s^2 + s + 1 - P(s)\bigr)$.

Using SnapPy \cite{snappy}, we find the fundamental group of the census manifold
$M = \mathtt{t07934}$ to have the presentation
$\pi_1(M) \cong \langle a, b, c\;|\; acabc^{-1}bc^{-3},\;a^2c^3b^2\rangle$ and peripheral subgroup generated by $cb^{-1}a^{-1}c^2$ and $cb^{-1}cb^{-3}$.
Hence, we have an epimorphism $\varphi:\pi_1(M)\to\ZZ_3*\ZZ$ given by
$\varphi(a) = xy^{-1}$, $\varphi(b) = yx^{-1}$, $\varphi(c) = yx^{-1}y^{-1}$.
Then $\rho_s := \psi_s\circ\varphi$ is an infinite family of boundary-parabolic representations of $\pi_1(M)$. 
Since the trace of $\psi_s(y)$ is non-constant, these representations represent
infinitely many conjugacy classes. 
Now $\rho_s$ is an irreducible
representation if and only if trace $[\psi_s(x)$, $\psi_s(y)]\ne 2$ (see, for example, Heusener--Porti \cite[Lemma~3.4]{heusener-porti}); this happens for a Zariski-open set of values of the parameter $s$. 
Furthermore, the obstruction class to lifting $\rho_s$ to a boundary-unipotent
\nword{\SLC}{representation} coincides with the geometric obstruction class. 
We thank Matthias G\"orner for providing data on triangulations with positive dimensional Ptolemy varieties, which we used to help find this example.

%==============================================================================
% Section 6: The tau-invariant
%==============================================================================
\section{\texorpdfstring{The $\tau$-invariant}{The tau-invariant}}

This section is devoted to the study of the quantity
$\tau$ introduced in \eqref{first-formula-for-tau}.
Throughout the section, we assume that $M$ has a single toroidal boundary
component and that $\triang$ is an ideal triangulation of $M$.
Assume that $z\in(\CC\setminus\RR)^{Q(\triang)}$ is an algebraic solution
of the edge consistency equations~\eqref{edge-consistency-equations}.
The shapes $z$ determine an \Sval\ angle structure $\omega$
by the usual formula
\begin{equation}\label{omega-from-z}
    \omega(\square) = \frac{z(\square)}{\abs{z(\square)}}.
\end{equation}
This angle structure $\omega$ is in fact an element of the
open subspace
\[
    \SASnonflat(\triang) = \bigl\{\omega\in\SAS(\triang) : \omega(\square)\not\in\{-1,1\}
    \ \text{for all}\ \square\in Q(\triang)\bigr\}
\]
consisting of circle-valued angle structures free of any ``flat angles''.
When considering a neighbourhood of a point $\omega\in\SASnonflat(\triang)$,
it is convenient to use real-valued pseudo-angles given by
$\alpha(\square) = \Arg\omega(\square)\in(-\pi,0)\cup(0,\pi)$,
but it must be borne in mind that these pseudo-angles only satisfy edge equations
modulo $2\pi$.

In this notation, $\tau=\tau(z)$ is defined as
\begin{equation}\label{tau-recalled}
    \tau = \tau(z,\beta) = \frac{\prod\limits_{\square\in Q(\triang)}
            \abs{\sin\alpha(\square)}^{\beta(\square)/\pi-\frac{1}{2}}}
            {\sqrt{\left|\det\bigl(\mathbf{L_*}
                \diag\{\cot\alpha(\square)\}
            \mathbf{L_*^\transp} \bigr) \right|}},
\end{equation}
where $\beta\in\Areal(\triang)$ is a generalised angle structure for which
$\exp(i\beta)$ and $\omega$ lie in the same connected component of
$\SAS_{(\exp(i\aholM),\;\exp(i\aholL))}(\triang)$.
In the above equation, $\mathbf{L}_*$ is the matrix of \eqref{definition-L_star}.
In Section~\ref{tau-exponent-invariance}, we shall show that the value of $\tau(z,\beta)$ does
not depend on the choice of $\beta$.

The key observation which underlies most of the results of this section is the
interpretation of the matrix~$\mathbf{L_*}\diag\{\cot\alpha(\square)\}\mathbf{L_*^\transp}$
in terms of the Jacobian matrix of the \emph{shearing displacements} of an angle structure~$\alpha$
along the edges of the triangulation~$\triang$.
Therefore, we start by introducing the concept of shearing singularities.

\subsection{Shearing singularities of angle structures}
For every flat-free angle structure~$\omega\in\SASnonflat(\triang)$
it is possible to invert the relationship~\eqref{omega-from-z}
by defining the complex shapes via
\begin{equation}\label{shapes-from-angles}
    z(\square) = \abs{\frac{\Imag\omega(\square')}{\Imag\omega(\square'')}}\omega(\square)
    = \abs{\frac{\sin\alpha(\square')}{\sin\alpha(\square'')}}e^{i\alpha(\square)}.
\end{equation}
Note that the above shapes may not satisfy Thurston's gluing equations:
even though their product around every edge of $\triang$ is a real number, it may fail to equal $1$. 
Geometrically, the sum of the angles around each edge is a multiple of $2\pi$, but there may be {\em shearing type singularities} along edges where the holonomy around the edge is a non-trivial translation along the edge.
The failure of the shapes~\eqref{shapes-from-angles} to satisfy complex edge consistency equations
is measured by these translations along the edges.

\begin{definition}\label{def-shear}
The \emph{shearing displacement}~$\Shear_j$ of $\omega\in\SASnonflat(\triang)$
along the edge~$E_j\in\edges(\triang)$ is given by
\[
    \Shear_j(\omega)
    = \sum_{\square\in Q(\triang)} G(E_j,\square) \log\left|\frac{\sin\alpha(\square')}{\sin\alpha(\square'')}\right|
    = -\sum_{\square\in Q(\triang)}l_j(\square)\log|\sin\alpha(\square)|,
    \quad 1\leq j \leq N
\]
where $G(E_j, \square)$ is the number of times the quad $\square$ faces the edge~$E_j$.
\end{definition}
It is easy to see that the shape assignment~$z=z(\omega)$ of \eqref{shapes-from-angles} satisfies the edge consistency
equations~\eqref{edge-consistency-equations} if and only if $\Shear_j(\omega)=0$ for all $j\in\{1,\dotsc,N\}$.

Since $\triang$ has a single toroidal end, the shearing displacements~$\Shear_j$ defined
above always satisfy $\sum_{j=1}^N \Shear_j = 0$,
and this is the only nontrivial linear relation among them.
Omitting the last edge~$E_N$, we obtain the reduced \emph{shear map}
\begin{equation}\label{shear-map}
    \Shear\colon \SASnonflat(\triang) \to \RR^{N-1},\quad
    \Shear(\omega) = \bigl(\Shear_1(\omega),\dotsc,\Shear_{N-1}(\omega)\bigr).
\end{equation}

We shall now fix the peripheral angle-holonomy $(x,y)\in S^1\times S^1$ arbitrarily and consider
the restriction of the reduced shear map~$\Shear$ to $\SASnonflat_{(x,y)}(\triang)$.
Observe that the derivative of this restriction,
$D\Shear\colon TAS_0(\triang)\to\RR^{N-1}$ can be described as
a Jacobian matrix with respect to the basis of leading-trailing deformations
$l_1,\dotsc,l_{N-1}\in TAS_0(\triang)=T_\omega\bigl(\SASnonflat_{(x,y)}(\triang)\bigr)$.
Using the \nword{\alpha}{coordinates}, we calculate explicitly,
for any $i,j\in\{1,\dotsc,N-1\}$,
\begin{align}
    \nabla_{l_j}\Shear_i(\alpha)
    &= \sum_{\square\in Q(\triang)}l_j(\square)\frac{\partial\Shear_i(\alpha)}{\partial\alpha(\square)}\notag\\
    &= -\sum_{\square\in Q(\triang)}l_j(\square)l_i(\square)\frac{\partial}{\partial\alpha(\square)}
        \log\abs{\sin\alpha(\square)}\notag\\
    &= -\sum_{\square\in Q(\triang)}l_j(\square)l_i(\square)\cot\alpha(\square)\notag\\
    &= -[\mathbf{L}_*\diag\{\cot\alpha(\square)\}\mathbf{L}_*^\transp]_{i,j}.\label{explicit-Jacobian-of-shears}
\end{align}
In other words, $-\mathbf{L}_*\diag\{\cot\alpha(\square)\}\mathbf{L}_*^\transp$ is the Jacobian matrix of $\Shear$
in the coordinates given by the basis of leading-trailing deformations.
We define
\begin{align}
    \label{def-tau_1}
    \tau_1(z) &= \abs{\det\bigl(\mathbf{L}_*\diag\{\cot\alpha(\square)\} \mathbf{L}_*^\transp\bigr)}^{-1},\\
    \tau_2(z,\beta) &= \prod_{\square\in Q(\triang)} \abs{\sin\alpha(\square)}^{2\beta(\square)/\pi-1},
    \label{def-tau_2}
\end{align}
so that $\tau(z,\beta) = \sqrt{\tau_1(z)\tau_2(z,\beta)}$.
Then \eqref{explicit-Jacobian-of-shears} says that
\begin{equation}\label{tau_1-shear-Jacobian}
    \tau_1(z)=\abs{\det\bigl(D\Shear\bigr|_{\TAS_0(\triang)}\bigr)}^{-1}
\end{equation}
with respect to the basis of $\TAS_0(\triang)$ given by the rows of $\mathbf{L}_*$.

\subsection{\texorpdfstring{Invariance of $\tau$ under the change of exponents}%
{Invariance of tau under the change of exponents}}\label{tau-exponent-invariance}
We shall now show that $\tau(z,\beta)$ does not depend on the choice of the
generalised angle structure $\beta$.

\begin{proposition}\label{tau-independent-of-exponents}
Suppose that $\beta$ and $\tilde{\beta}$ are two real-valued angle structures with
identical peripheral angle-holonomies. Then
$\tau(z,\beta)=\tau(z,\tilde{\beta})$.
\end{proposition}
\begin{proof}
It suffices to study the factor~$\tau_2(z,\beta)$ given by eq.~\eqref{def-tau_2}.
Since $\beta$ and $\tilde{\beta}$ have identical peripheral angle-holonomies,
the difference $\beta-\tilde{\beta}$ is an element of $\TAS_0(\triang)$
and therefore $\beta = \tilde{\beta} + \sum_{j=1}^N c_j l_j$ for some real coefficients $c_j$.
Hence,
\[
    \frac{\tau_2(z,\beta)}{\tau_2(z,\tilde{\beta})}
    = \prod_{\square\in Q(\triang)}
             |\sin\alpha(\square)|^{2(\beta(\square) - \tilde{\beta}(\square))/\pi}
    = \prod_{j=1}^N \exp\left[ \frac{2c_j}{\pi}
        \sum_{\square\in Q(\triang)} l_j(\square) \log\abs{\sin\alpha(\square)}\right].
\]
Since the shapes~$z$ satisfy the gluing equations, for every $j\in\{1,\dotsc,N\}$ we have
\[
    0 = -\Shear_j(\alpha)
    = \sum_{\square\in Q(\triang)} l_j(\square) \log\abs{\sin\alpha(\square)}.
\]
This implies $\tau_2(z,\beta)=\tau_2(z,\tilde{\beta})$.
\end{proof}
%--------------------------------------------------------------------------------------------------
\subsection{\texorpdfstring{Invariance of $\tau$ under Pachner~$2$--$3$ moves}%
{Invariance of tau under Pachner 2-3 moves}}
\label{pachner-invariance}
Our goal is now to prove that $\tau(z)=\sqrt{\tau_1(z)\tau_2(z)}$
is invariant under shaped Pachner~$2$--$3$ moves on ideal triangulations.

\begin{figure}
    \centering
    \begin{tikzpicture}
        \node[anchor=south west,inner sep=0] (image) at (0,0)
            {\includegraphics[width=0.7\columnwidth]{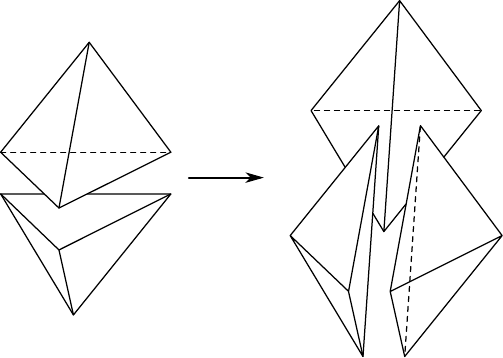}};
        \begin{scope}[x={(image.south east)},y={(image.north west)},
                every node/.style={inner sep=0, outer sep=0}]
            \node[anchor=south] at (0.18, 0.89)  {$0$};
            \node[anchor=east] at (0, 0.57)  {$1$};
            \node[anchor=west] at (0.345, 0.57)  {$3$};
            \node[anchor=east] at (0, 0.45)  {$1$};
            \node[anchor=west] at (0.345, 0.45)  {$3$};
            \node[anchor=north] at (0.12, 0.405)  {$2$};
            \node[anchor=south] at (0.12, 0.31)  {$2$};
            \node[anchor=north] at (0.15, 0.105)  {$4$};
            \node[anchor=south] at (0.795, 1.01)  {$0$};
            \node[anchor=east] at (0.61, 0.69)  {$1$};
            \node[anchor=west] at (0.965, 0.69)  {$3$};
            \node[anchor=north] at (0.77, 0.35)  {$4$};
            \node[anchor=west] at (1.005, 0.34)  {$3$};
            \node[anchor=north] at (0.75, 0.2)  {$2$};
            \node at (0.76, 0)  {$4$};
            \node[anchor=east] at (0.57, 0.34)  {$1$};

            \node[anchor=west] at (0.05, 0.9) {\normalsize$\mathcal{T}_2$};
            \node[anchor=east] at (0.65, 0.95) {\normalsize$\mathcal{T}_3$};
            % Edge labels left, tetrahedron 4
            \node[anchor=south west] at (0.26, 0.74) {$z_4$};
            \node[anchor=south east] at (0.094, 0.75) {$z'_4$};
            \node[anchor=north west] at (0.16, 0.72) {$z''_4$};
            \node[anchor=south west] at (0.07, 0.485) {$z_4$};
            \node[anchor=south east] at (0.2, 0.485) {$z'_4$};
            \node[anchor=south] at (0.08, 0.58) {$z''_4$};
            % Edge labels left, tetrahedron 0
            \node[anchor=north] at (0.23, 0.442) {$z_0$};
            \node[anchor=south east] at (0.196, 0.36) {$z'_0$};
            \node[anchor=south west] at (0.055, 0.38) {$z''_0$};
            \node[anchor=south west] at (0.137, 0.22) {$z_0$};
            \node[anchor=north east] at (0.07, 0.28) {$z'_0$};
            \node[anchor=north west] at (0.24, 0.28) {$z''_0$};
            % Edge labels right, tetrahedron 2
            \node[anchor=west] at (0.796, 0.8) {$z_2$};
            \node[anchor=south west] at (0.89, 0.83) {$z'_2$};
            \node[anchor=south east] at (0.69, 0.83) {$z''_2$};
            % Edge labels right, tetrahedron 3
            \node[anchor=south west] at (0.64, 0.255) {$z_3$};
            \node[anchor=south east] at (0.63, 0.44) {$z'_3$};
            \node[anchor=north west] at (0.613, 0.17) {$z''_3$};
            % Edge labels right, tetrahedron 1
            \node[anchor=south] at (0.9, 0.287) {$z_1$};
            \node[anchor=north west] at (0.91, 0.18) {$z'_1$};
            \node[anchor=south west] at (0.95, 0.44) {$z''_1$};
        \end{scope}
    \end{tikzpicture}
    \caption{%
    \label{fig:shaped-Pachner}
    The shaped Pachner \twothree\ move on a labelled bipyramid.
    The triangulations $\triang_2$ and $\triang_3$ are identical
    outside of the bipyramid.
    }
\end{figure}
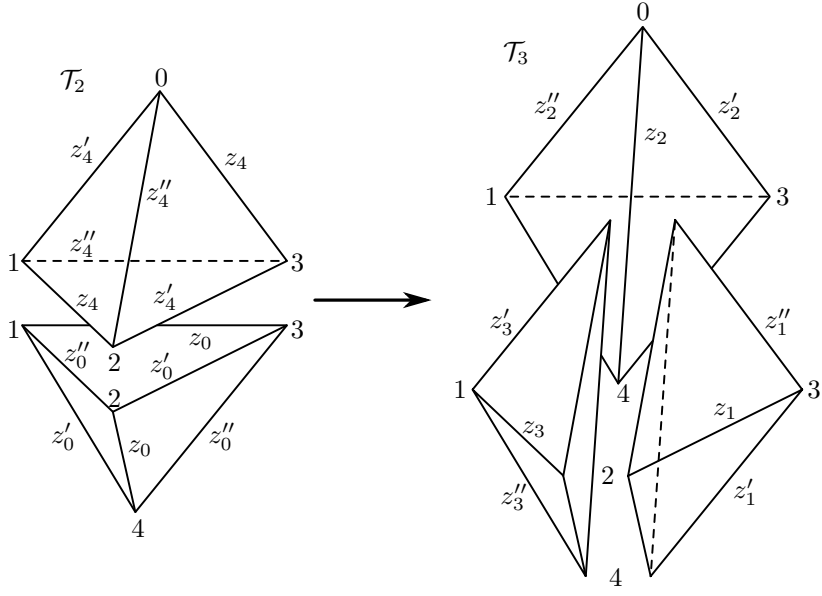

Assume that the triangulations $\triang_2$ and $\triang_3$ are related by a Pachner 2-3 move as
depicted in Figure~\ref{fig:shaped-Pachner}.
The shape parameters $z_1$, $z_2$, $z_3$ on the three-tetrahedron side
are uniquely determined by the two shapes $z_0$, $z_4$; explicitly,
\begin{equation}\label{shape-formula-23}
    z_1 = z'_0z'_4, \quad z_2 = z_0z''_4,\quad z_3 = z''_0z_4.
\end{equation}
Note that \eqref{shape-formula-23} implies the gluing equation~$z_1z_2z_3=1$ along the central
edge~$[04]$ of the bipyramid.
Denote by $\tau^{(2)}(z)$ and $\tau^{(3)}(z)$ be the \nword{\tau}{invariants} of
$\triang_2$ and $\triang_3$, respectively.

\begin{theorem}\label{invariance-of-tau}
Assume that $\triang_2$ and $\triang_3$ are two ideal triangulations equipped with shapes
in $\CC\setminus\RR$ related by a shaped Pachner $2$--$3$ move.
Then $\tau^{(2)}(z)=\tau^{(3)}(z)$.
\end{theorem}

In order to prove the above theorem, we first need to establish a number of lemmas; the proof
itself is contained in Section~\ref{tau-invariance-proofsection}.

\subsubsection{Angled Pachner moves and their derivatives}
Throughout the section, we assume that $\triang_2$ and $\triang_3$ are two triangulations
related by a shaped Pachner $2$--$3$ move as in Figure~\ref{fig:shaped-Pachner}.
If $\triang_2$ has $N$ tetrahedra and edges~$E_1,\dotsc,E_N$,
then $\triang_3$ contains $N+1$ tetrahedra and in addition to the edges already present in $\triang_2$,
it has an additional edge $E_*$ which arises as the central edge of the bipyramid on which
the Pachner move is performed (edge $[04]$ in Figure~\ref{fig:shaped-Pachner}).

In the case of a $3\mathord\rightarrow2$ move, an angle structure $\omega\in\SAS_{(x,y)}(\triang_3)$
uniquely determines an angle structure $\PachD(\omega)\in\SAS_{(x,y)}(\triang_2)$.
Namely, $\PachD(\omega)$ is identical to $\omega$ outside of the bipyramid of the Pachner move,
and the angles in the two tetrahedra of the bipyramid in $\triang_2$ are given uniquely
as products in $S^1$ of respective angles on $\triang_3$.
We can view this `angled Pachner move' as a map between the spaces of angle
structures of the two triangulations,
\begin{equation}\label{def-PachD}
    \PachD\colon\SAS_{(x,y)}(\triang_3)\to\SAS_{(x,y)}(\triang_2).
\end{equation}
We shall now construct a canonical, partially defined section
$\PachA\colon\SASnonflat_{(x,y)}(\triang_2)\to\SAS_{(x,y)}(\triang_3)$ of the map $\PachD$.
This map between angle spaces is determined uniquely by the additional requirement that it preserves
the shearing displacements along the edges~$E_1,\dotsc,E_N$ common to the two triangulations.

\begin{definition}
Given an angle structure~$\omega\in\SASnonflat(\triang_2)$, consider
the shape parameters on $\triang_2$ given by \eqref{shapes-from-angles}.
In particular, the shapes~$z_0$, $z_4$ on $\triang_2$ uniquely determine shapes $z_1$, $z_2$, $z_3$
on $\triang_3$ via~\eqref{shape-formula-23}.
Let $z$ denote the resulting shape assignment on $\triang_3$.
We define
\[
    \PachA\colon\SASnonflat(\triang_2)\to\SAS(\triang_3), \quad
    \left[\PachA(\omega)\right](\square) = \tfrac{z(\square)}{|z(\square)|}
    \ \text{for all}\ \square\in Q(\triang_3).
\]
\end{definition}
It is easy to see that $\PachD\circ\PachA$ is the identity on $\SASnonflat(\triang_2)$.
We now characterise the effect of $\PachA$ on the shearing displacements.

\begin{lemma}\label{shear-preserving}
    Let $\Shear^2$ and $\Shear^3$ be the shear maps of the triangulations $\triang_2$ and $\triang_3$,
    respectively.
    If $\omega\in\SASnonflat(\triang_2)$ is such that $\PachA(\omega)\in\SASnonflat(\triang_3)$,
    then $\Shear^3_j\bigl(\PachA(\omega)\bigr) = \Shear^2_j(\omega)$ for all $1\leq j \leq N$.
    Moreover, the shearing displacement of $\PachA(\omega)$ along $E_*$ is zero.
\end{lemma}
\begin{proof}
If an edge $E_j$ is disjoint from the bipyramid involved in the Pachner move, there is nothing to prove,
since $\PachA$ does not affect the angles of any of the tetrahedra meeting $E_j$.
On the other hand, whenever $E_j$ is incident to one of the external edges of the bipyramid,
the product of the shapes labelling tetrahedral edges
in $\triang_2$ incident to $E_j$ is equal to the product of shapes in $\triang_3$ incident to $E_j$.
In particular, the moduli of these products are equal, which shows that $\omega$ and $\PachA(\omega)$ have
the same contribution to the shearing displacement at $E_j$.
Finally, since the equation $z_1z_2z_3=1$ follows directly from \eqref{shape-formula-23}, the angle
structure $\PachA(\omega)$ has no shear along the edge $E_*$.
\end{proof}

\begin{remark}\label{rem:angled-Pachner-moves}
The angled Pachner move map $\PachD$ of \eqref{def-PachD} is a multiplicative
analogue of the well-known angled Pachner move with additive angles.
More specifically, every generalised angle structure $\alpha_3\in\Areal(\triang_3)$
determines a unique generalised angle structure $\alpha_2\in\Areal(\triang_2)$ in which
the external angles of the bipyramid of the Pachner move are the same as in $\alpha_3$
and which is identical to $\alpha_3$ outside the bipyramid.
By exponentiating generalised angle structures, we obtain \Sval\ angle structures
in the geometric components of the respective angle spaces:
$\omega_3=\exp(i\alpha_3)\in\SASdistinguished(\triang_3)$
and $\omega_2=\exp(i\alpha_2)\in\SASdistinguished(\triang_2)$.
By construction, we have $\PachD(\omega_3)=\omega_2$.
This shows that the map $\PachD$ respects the geometric components.
Since the angled Pachner moves preserve peripheral angle-holonomies,
for any $(x,y)\in S^1\times S^1$ we have
$\PachD\bigl(\SASdistinguished_{(x,y)}(\triang_3)\bigr)=\SASdistinguished_{(x,y)}(\triang_2)$.
\end{remark}

For arbitrarily fixed $(x,y)\in S^1\times S^1$, we shall now
compute the derivatives of the restrictions of $\PachD$ and $\PachA$ to
the spaces of circle-valued angles with the peripheral angle-holonomy~$(x,y)$.
We shall treat these derivatives as maps between the tangent spaces,
\[
    D\PachD : \TAS_0(\triang_3)\to \TAS_0(\triang_2),
    \qquad
    D\PachA : \TAS_0(\triang_2)\to \TAS_0(\triang_3),
\]
and study their Jacobian matrices with respect to the bases
of leading-trailing deformations
\begin{equation}\label{basis-choice}
    \TAS_0(\triang_2) = \linspan\{ l^{(2)}_1,\dotsc,l^{(2)}_{N-1}\},\qquad
    \TAS_0(\triang_3) = \linspan\{ l^{(3)}_1,\dotsc,l^{(3)}_{N-1};\; l_*\},
\end{equation}
where we assume that the ignored edge $E_N$ is not the central edge~$E_*$
of the bipyramid and that the numbering of the edges is consistent between the two triangulations.

\begin{lemma}
With respect to the bases \eqref{basis-choice}, the Jacobian matrices have the form
\begin{equation}\label{Jacobian-matrices}
    D\PachD=
        \left[
            \begin{array}{c|c}
                {} & 0\\
                \hspace{1.5em}I \hspace{1.5em} & \vdots\\
                {} & 0
            \end{array}
        \right],\qquad
    D\PachA=
        \left[
            \begin{array}{c}
                {}\\[0.2ex] I \\[1.5ex]
                \tikz[baseline=({w})]{%
                    \node (w) at (0,0) {w};
                    \draw (w) -- ++(1,0);
                    \draw (w) -- ++(-1,0);
                }
            \end{array}
        \right],\quad
\end{equation}
where $I$ is the $(N-1)\times(N-1)$ identity matrix and $w=w(\omega)$ is a certain row vector.
\end{lemma}
\begin{proof}
The second equality follows from the first.
Indeed, assume that $D\PachD$ has the matrix given by the first equality and that
$D\PachA=\begin{bmatrix}W\\w\end{bmatrix}$ for an indeterminate square matrix~$W$
and some row vector~$w$.
Since $\PachD\circ\PachA\equiv\Id$, we have
\[
    I = \bigl[D\PachD\bigr]\cdot\bigl[D\PachA\bigr] = IW + 0 = W.
\]

We shall now prove the first equality of \eqref{Jacobian-matrices}.
To keep things simple, we will work in the angle coordinates
$\alpha(\square)=\Arg\omega(\square)\in(-\pi,0)\cup(0,\pi)$
which are unambiguously defined in a neighbourhood of any circle-valued
angle structure~$\omega\in\SASnonflat_{(x,y)}(\triang_3)$.

For an arbitrarily fixed $j\in\{1,\dotsc,N-1\}$,
the basis vector $l^{(3)}_j\in TAS_0(\triang_3)$ is tangent to the
one-parameter family of angle structures
\begin{equation}\label{one-parameter-family-of-angle-structures}
 t\mapsto\left(\alpha(\square)+tl^{(3)}_j(\square)\right)_{\square\in Q(\triang_3)},\quad t\in(-\eps,\eps),
\end{equation}
where $\eps>0$ is small enough so that $\omega = e^{i\alpha}\in\SASnonflat_{(x,y)}(\triang_3)$.
The image of this family under $\PachD$ is a family of angle structures on $\triang_2$
whose infinitesimal generator can be expressed in the basis $\{l^{(2)}_1,\dotsc,l^{(2)}_{N-1}\}$.
By definition, the \nword{j}{th} column of the Jacobian $D\PachD$ contains the coefficients occurring in this expression.
We split our reasoning into a few cases depending on the position of the edge~$E_j$ relative to the bipyramid
of the Pachner move.

\Case{1:~$E_j$ disjoint from the bipyramid of the Pachner move.}
Since none of the angles in any of the tetrahedra incident to $E_j$ change under $\PachD$,
we have
$\PachD\bigl(\alpha+tl^{(3)}_j\bigr) = \PachD(\alpha) + tl^{(2)}_j$,
so that $D\PachD\bigl(l^{(3)}_j\bigr) = l^{(2)}_j$.

\Case{2:~$E_j$ is an external edge of the bipyramid emanating from either the
top or the bottom apex.}
Since $\PachD$ leaves all angles outside of the bipyramid unchanged,
it suffices to consider the angles on the two triangulations
of the bipyramid related by the Pachner move.
When $E_j$ is incident to the boldened edge in Figure~\ref{semi-meridinal-case}, the
one-parameter family~\eqref{one-parameter-family-of-angle-structures}
transforms under $\PachD$ as depicted in the figure.
\begin{figure}[tbh]
    \centering
    \begin{tikzpicture}
        \node[anchor=south west,inner sep=0] (image) at (0,0)
            {\includegraphics[width=0.9\textwidth]{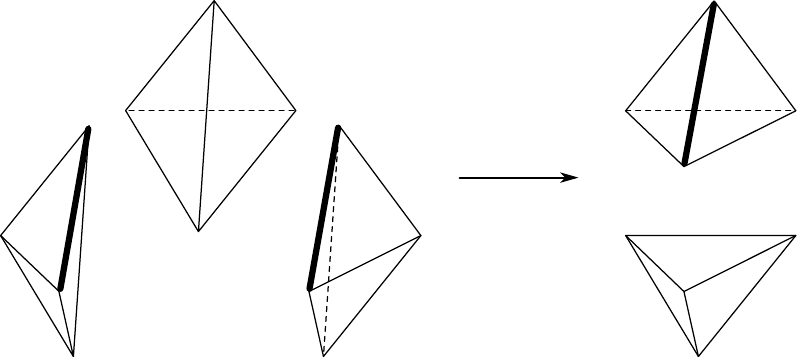}};
        \begin{scope}[x={(image.south east)},y={(image.north west)}]
            % Left panel
            \node[anchor=south west, inner sep=0] at (0.0, 0.53) {$\alpha'_2+t$};
            \node[anchor=south east, inner sep=0] at (0.09, 0.39) {$\alpha''_2$};
            \node[anchor=north west, inner sep=0] at (0.105, 0.3) {$\alpha_2-t$};
            \node[anchor=south east, inner sep=0] at (0.4, 0.39) {$\alpha'_3$};
            \node[anchor=south east, inner sep=0] at (0.49, 0.29) {$\alpha_3+t$};
            \node[anchor=east, inner sep=0] at (0.39, 0.1) {$\alpha''_3-t$};
            \node[anchor=west, inner sep=0] at (0.26, 0.6) {$\alpha_1$};
            \node[anchor=south east, inner sep=0] at (0.21, 0.85) {$\alpha''_1$};
            \node[anchor=south west, inner sep=0] at (0.32, 0.85) {$\alpha'_1$};
            % Right panel
            \node[anchor=south west, inner sep=0] at (0.92, 0.91) {$\alpha'_1+\alpha''_3-t$};
            \node[anchor=south east, inner sep=0] at (0.87, 0.91) {$\alpha''_1+\alpha'_2+t$};
            \node[anchor=south west, inner sep=0] at (0.88, 0.7) {$\alpha''_2+\alpha'_3$};
            \node[anchor=south, inner sep=0] at (0.88, 0.355) {$\alpha'_2+\alpha''_3$};
            \node[anchor=north west, inner sep=0] at (0.92, 0.1) {$\alpha''_1+\alpha'_3$};
            \node[anchor=north east, inner sep=0] at (0.84, 0.1) {$\alpha''_2+\alpha'_1$};
            % Arrow label
            \node[anchor=south] at (0.65, 0.5) {$\PachD$};
        \end{scope}
    \end{tikzpicture}
    \caption{An angled Pachner \threetwo\ move on a one-parameter family of angle structures
    tangent to the leading-trailing deformation of an external edge of the bipyramid
    emanating from an apex.}
    \label{semi-meridinal-case}
\end{figure}
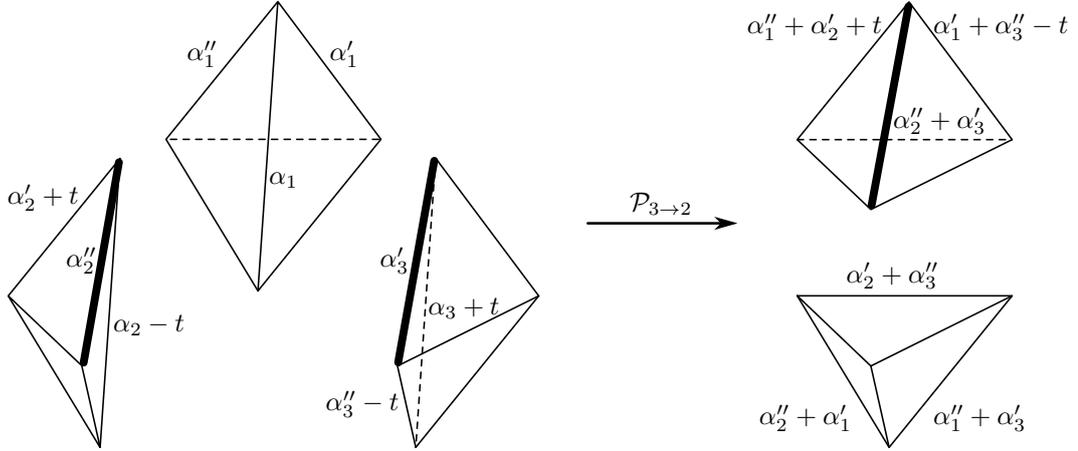
It is therefore evident that the one-parameter family of angles depicted on the right panel
of the figure has $l^{(2)}_j$ as its infinitesimal generator.
In other words, $D\PachD\bigl(l^{(3)}_j\bigr) = l^{(2)}_j$.

\Case{3:~$E_j$ incident to an equatorial edge of the bipyramid.}
This situation is depicted in Figure~\ref{equatorial-case}.
We consider the bold edge in the figure and compute the transformation of
the corresponding one-parameter family of angles.
\begin{figure}[tbh]
    \centering
    \begin{tikzpicture}
        \node[anchor=south west,inner sep=0] (image) at (0,0)
            {\includegraphics[width=0.9\textwidth]{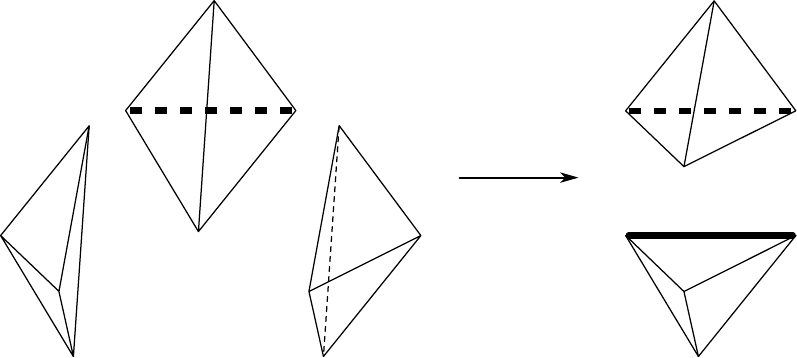}};
        \begin{scope}[x={(image.south east)},y={(image.north west)}]
            % Left panel
            \node[anchor=south east, inner sep=0] at (0.06, 0.52) {$\alpha'_2$};
            \node[anchor=south east, inner sep=0] at (0.09, 0.39) {$\alpha''_2$};
            \node[anchor=north west, inner sep=0] at (0.105, 0.3) {$\alpha_2$};
            \node[anchor=south east, inner sep=0] at (0.4, 0.39) {$\alpha'_3$};
            \node[anchor=south east, inner sep=0] at (0.48, 0.29) {$\alpha_3$};
            \node[anchor=east, inner sep=0] at (0.39, 0.1) {$\alpha''_3$};
            \node[anchor=west, inner sep=0] at (0.26, 0.6) {$\alpha_1$};
            \node[anchor=south east, inner sep=0] at (0.21, 0.85) {$\alpha''_1+t$};
            \node[anchor=south west, inner sep=0] at (0.32, 0.85) {$\alpha'_1-t$};
            % Right panel
            \node[anchor=south west, inner sep=0] at (0.92, 0.91) {$\alpha'_1+\alpha''_3-t$};
            \node[anchor=south east, inner sep=0] at (0.87, 0.91) {$\alpha''_1+\alpha'_2+t$};
            \node[anchor=south west, inner sep=0] at (0.88, 0.7) {$\alpha''_2+\alpha'_3$};
            \node[anchor=south, inner sep=0] at (0.88, 0.355) {$\alpha'_2+\alpha''_3$};
            \node[anchor=north west, inner sep=0] at (0.92, 0.1) {$\alpha'_3+\alpha''_1+t$};
            \node[anchor=north east, inner sep=0] at (0.84, 0.1) {$\alpha''_2+\alpha'_1-t$};
            % Arrow label
            \node[anchor=south] at (0.65, 0.5) {$\PachD$};
        \end{scope}
    \end{tikzpicture}
    \caption{An angled Pachner \threetwo\ move on a one-parameter family of angle structures
    tangent to the leading-trailing deformation of an equatorial edge of the bipyramid.}
    \label{equatorial-case}
\end{figure}
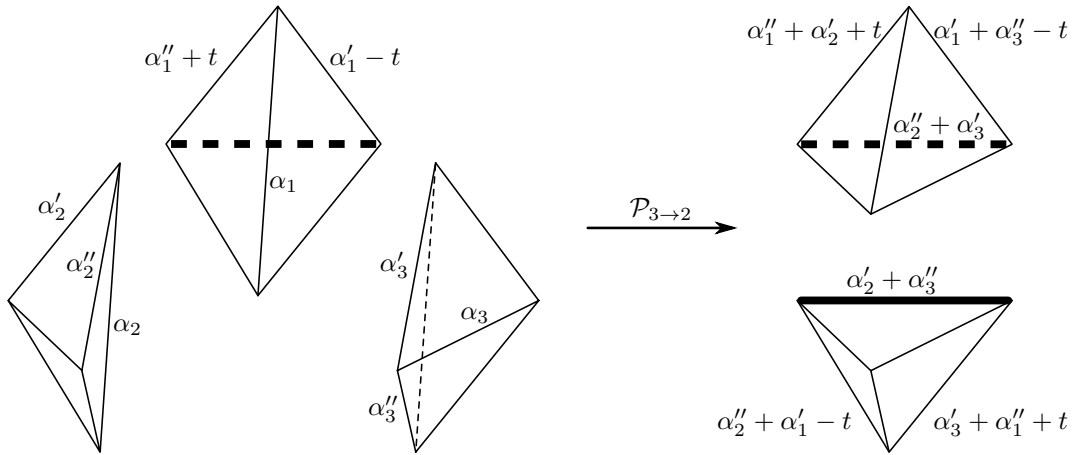
As before, the one-parameter family on the right-hand side of the figure
evidently has $l^{(2)}_j$ as its infinitesimal generator,
so that $D\PachD\bigl(l^{(3)}_j\bigr) = l^{(2)}_j$ once again.

\Case{4:~General case.}
In general, the external edges of the bipyramid may be glued to one another.
However, these gluings will be identical in $\triang_2$ and $\triang_3$.
Treating the edge class $E_j$ as an equivalence class of tetrahedral edges,
the calculation reduces to the cases 1--3 considered above.

Finally, we consider the central edge~$E_*$ of the bipyramid.
Since the leading--trailing deformation $l_*$ keeps the total external angles of the bipyramid constant
and vanishes around the bipyramid's equator,
we have $\PachD(\alpha+tl_*)\equiv\PachD(\alpha)$, and consequently $D\PachD(l_*) = 0$.
We have shown that the matrix of $D\PachD$ is exactly as given in~\eqref{Jacobian-matrices} and the proof
is finished.
\end{proof}

Let $\Shear^2$ and $\Shear^3=(\Shear_1,\dotsc,\Shear_{N-1}; \Shear_*)$ be as in Lemma~\ref{shear-preserving}.
Fix a point $\omega\in\SASnonflat(\triang_2)$ such that $\PachA(\omega)\in\SASnonflat(\triang_3)$.
Since these points are smooth points of $\Shear^2$ and $\Shear^3$, respectively,
we obtain the following diagram in the category of finite-dimensional real vector spaces:
\begin{equation}\label{big-diagram}
    \begin{tikzcd}[column sep=7em]
        \TAS_0(\triang_2)\arrow{r}{D\PachA}\arrow{d}{D\Shear^2} &  \TAS_0(\triang_3)\arrow{d}{D\Shear^3}\\
        \RR^{N-1}\arrow{r}{\iota_1} &  \RR^{N-1}\oplus\RR,
    \end{tikzcd}
\end{equation}
where $\iota_1$ is the inclusion onto the first component of $\RR^{N-1}\oplus\RR$.
This component corresponds to the edges of $\triang_3$ already present in $\triang_2$, with the same numbering.
Proposition~\ref{shear-preserving} says that $\PachA$ preserves the shearing
displacements along edges $E_1,\dotsc,E_{N-1}$,
implying that the diagram~\eqref{big-diagram} is commutative.
Moreover, the vertical arrows are isomorphisms.

The newly introduced edge $E_*$ furnishes a new tangent vector $l_*\in TAS_0(\triang_3)$
which we deal with by considering the rather trivial commutative square
\begin{equation}\label{small-diagram}
    \begin{tikzcd}[column sep=7em]
        \RR \arrow{r}{1\mapsto l_*}\arrow[equal]{d}
              & \TAS_0(\triang_3)\arrow{d}{D\Shear^3}\\
        \RR\arrow{r}{1\mapsto D\Shear^3(l_*) } &\RR^{N-1}\oplus\RR.
    \end{tikzcd}
\end{equation}
The direct sum of the commutative squares \eqref{big-diagram} and \eqref{small-diagram} now induces the diagram
\begin{equation}\label{combined-diagram}
    \begin{tikzcd}[column sep=3em]
        \TAS_0(\triang_2)\oplus\RR  \arrow{r}{R}\arrow{d}{D\Shear^2\oplus\Id}
             &  \TAS_0(\triang_3)\arrow{d}{D\Shear^3}\\
        \RR^{N-1}\oplus\RR\arrow{r}{S} & \RR^{N-1}\oplus\RR,
    \end{tikzcd}
\end{equation}
where
\begin{equation*}
    R(a\oplus b) = D\PachA(a) + bl_*,\qquad
    S(a\oplus b) = a\oplus0 + bD\Shear^3(l_*).
\end{equation*}
Although the commutativity of \eqref{combined-diagram} follows from general properties of direct sums,
we can verify it explicitly by calculating
\begin{align*}
    D\Shear^3(R(a\oplus b))
    &= D\Shear^3(D\PachA(a)) + D\Shear^3(bl_*) \\
    &= D\Shear^2(a)\oplus0 + bD\Shear^3(l_*)\\
    &= S\left(D\Shear^2(a)\oplus b\right)\\
    &= \bigl(S\circ(D\Shear^2\oplus\Id)\bigr)(a\oplus b),
\end{align*}
where we used Lemma~\ref{shear-preserving} upon passing from the first to the second line.

We equip the spaces in the diagram~\eqref{combined-diagram} with the bases given in
\eqref{basis-choice} and, where needed, the standard basis of $\RR^n$.
Moreover, direct sums are equipped with unions of images of the respective bases under the canonical inclusions.

\begin{lemma}\label{det-RS}
With respect to the bases described above,
$\det R = 1$ and $\displaystyle\det S = \frac{\partial\Shear^3_*}{\partial l_*}$.
\end{lemma}
\begin{proof}
    Using the matrix expression for $D\PachA$ given in~\eqref{Jacobian-matrices},
    we see that the matrix of $R$ has the block form
    \[
        R =
        \left[
            \begin{array}{c|c}
                {} & 0\\
                {} & 0\\
                \hspace{1.1em}D\PachA\hspace{1em} & \vdots\\
                {} & 0\\
                {} & 1
            \end{array}
        \right]
        =
            \begin{bmatrix}
                I & 0\\
                w & 1
            \end{bmatrix}.
    \]
    Regardless of the value of the row vector $w$, we therefore have $\det R=1$.

    At the same time, we have
    \[
        \det S =
        \det
        \begin{bmatrix}
            1 & 0 & \dots & 0 & \frac{\partial\Shear^3_1}{\partial l_*} \\
            0 & 1 &  & \vdots & \frac{\partial\Shear^3_2}{\partial l_*} \\
            \vdots &  & \ddots & 0 & \vdots \\
            0 & \dots & 0 & 1 & \frac{\partial\Shear^3_{N-1}}{\partial l_*} \\
            0 & 0    & \dots & 0 & \frac{\partial\Shear^3_*}{\partial l_*}
        \end{bmatrix}
        =\frac{\partial\Shear^3_*}{\partial l_*}.\qedhere
    \]
\end{proof}

\subsubsection{Proof of the invariance theorem}\label{tau-invariance-proofsection}
We are now going to prove the invariance of $\tau(z)$ under shaped Pachner moves.
Throughout the section, we assume that $\triang_2$ are $\triang_3$ are two ideal triangulations
related by a shaped Pachner move as depicted in Figure~\ref{fig:shaped-Pachner}.

\begin{proof}[Proof of Theorem~\ref{invariance-of-tau}]
We shall study the factors $\tau_1$ and $\tau_2$ given in \eqref{def-tau_1}--\eqref{def-tau_2}
separately.
We begin by comparing the values $\tau_1^{(2)}$, $\tau_1^{(3)}$ of the factor~$\tau_1$
on the triangulations $\triang_2$ and $\triang_3$, respectively.
Recall from \eqref{tau_1-shear-Jacobian} that
\[
    \left|\det \bigl(D\Shear^2\oplus\Id\bigr)\right|
    = |\det D\Shear^2|
    = \frac{1}{\tau_1^{(2)}},\quad\text{and}\quad
    |\det D\Shear^3| = \frac{1}{\tau_1^{(3)}}.
\]
Hence, using the commutativity of the square~\eqref{combined-diagram} and Lemma~\ref{det-RS}, we get
\begin{equation}
    \frac{\tau_1^{(3)}}{\tau_1^{(2)}}
    = \left|\frac{\det (D\Shear^2\oplus\Id)}{\det D\Shear^3}\right|
    = \left|\frac{\det R}{\det S}\right|
    = \left|\frac{\partial\Shear^3_*}{\partial l_*}\right|^{-1}.
\end{equation}
Finally, referring to Figure~\ref{fig:shaped-Pachner} and formula~\eqref{shear-map},
we compute
\[
    \frac{\partial\Shear^3_*}{\partial l_*} = -\sum_{j=1}^3 \bigl(\cot\alpha'_j + \cot\alpha''_j\bigr),
\]
so that
\begin{equation}\label{ratio-of-tau1s}
    \frac{\tau_1^{(3)}}{\tau_1^{(2)}}
    =\Bigl|\sum_{j=1}^3 \bigl(\cot\alpha'_j + \cot\alpha''_j\bigr)\Bigr|^{-1}.
\end{equation}

It remains to establish a transformation law for $\tau_2$ which relates $\tau_2^{(2)}$ to $\tau_2^{(3)}$.
To start with, observe that \eqref{shape-formula-23} implies
\begin{align}
    |z_0| &= |z'_3 z''_1| = |(1-z_3)^{-1}(1-z_1)z_1^{-1}|\notag\\
    |1-z_0| &= |z'_0|^{-1} = |z'_2 z''_3|^{-1} = |(1-z_2)(1-z_3)^{-1}z_3|.\notag\\
    |z_4| &= |z'_2 z''_1| = |(1-z_2)^{-1}(1-z_1)z_1^{-1}|\label{s=s(z)}\\
    |1-z_4| &= |z'_4|^{-1} = |z''_2 z'_3|^{-1} = |(1-z_2)^{-1}z_2(1-z_3)|.\notag
\end{align}
Moreover, a generalised angle structure~$\beta\in\Areal(\triang_3)$
determines a unique angle structure on $\triang_2$
whose values on the bipyramid are given below (in the notations of Figure~\ref{fig:shaped-Pachner}):
\begin{alignat}{3}\label{Pachner-angle-formulas}
    \beta_0   &= \beta''_1 + \beta'_3,\quad&
    \beta'_0  &= \beta'_2 + \beta''_3,\quad&
    \beta''_0 &= \beta'_1 + \beta''_2,\\
    \beta_4   &= \beta''_1 + \beta'_2,\quad&
    \beta'_4  &= \beta''_2 + \beta'_3,\quad&
    \beta''_4 &= \beta'_1 + \beta''_3 \notag
\end{alignat}
and which is identical to $\beta$ outside of the bipyramid.
Let $b(\square) = \beta(\square)/\pi$.
According to \eqref{def-tau_2}, we can write
\begin{equation}\label{initial-ratio-of-tau_2s}
    \frac{\tau_2^{(3)}}{\tau_2^{(2)}}
    = \frac{
        \prod_{j\in\{1,2,3\}}
        \Bigl(
        \abs{\sin\alpha_j}  ^{2b_j  - 1}
        \abs{\sin\alpha'_j} ^{2b'_j - 1}
        \abs{\sin\alpha''_j}^{2b''_j- 1}
        \Bigr)
    }{
        \prod_{j\in\{0,4\}}
        \Bigl(
        \abs{\sin\alpha_j}  ^{2b_j  - 1}
        \abs{\sin\alpha'_j} ^{2b'_j - 1}
        \abs{\sin\alpha''_j}^{2b''_j- 1}
        \Bigr)
    }.
\end{equation}

In order to simplify this expression, suppose that $z\in\CC$ has a positive imaginary part.
Taking
\begin{equation}\label{definition-of-zetas}
    \zeta = \frac{1}{z},\quad
    \zeta' = \frac{1}{1-z},\quad
    \zeta'' = \frac{1}{z(z-1)},
\end{equation}
we may construct a triangle in $\CC$ whose oriented sides are $\zeta, \zeta', \zeta''$.
This triangle is similar to the triangle with vertices $\{0,1,z\}$ and is depicted in
Figure~\ref{zetas}.
Setting $\alpha=\arg z$, $\alpha' = \arg z'$ and $\alpha''=\arg z''$,
the law of sines yields
\begin{equation}\label{sines-from-zetas}
    \sin\alpha   = |\zeta|   \Imag z,\qquad
    \sin\alpha'  = |\zeta'|  \Imag z,\qquad
    \sin\alpha'' = |\zeta''| \Imag z.
\end{equation}
When $z$ has a negative imaginary part, the above equations can be deduced by
replacing $z$ with $\bar{z}$ and $\alpha$ with $-\alpha$.
\begin{figure}
    \centering\hfil
    \begin{tikzpicture}[baseline=(BASE)]
        \node[anchor=south west,inner sep=0] (image) at (0,0)
            {\includegraphics[height=20mm]{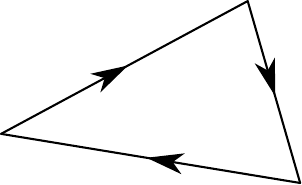}};
        \begin{scope}[x={(image.south east)},y={(image.north west)},every node/.style={inner sep=1pt}]
            \node[anchor=west] at (0.12, 0.31) {$z$};
            \node[anchor=south east] at (0.95, 0.08) {$z'$};
            \node[anchor=north] at (0.76, 0.9) {$z''$};
            \node[anchor=west] at (0.95, 0.5) {$\zeta$};
            \node[anchor=south east] at (0.4, 0.65) {$\zeta'$};
            \node[anchor=north] at (0.5, 0.1) {$\zeta''$};
            \node (BASE) at (1,1) {};
        \end{scope}
    \end{tikzpicture}
    \hfil
    \parbox[t][][t]{0.5\textwidth}{%
    \caption[A similarity class of a Euclidian triangle]{
        A Euclidean triangle with similarity class $z$ and sides $\zeta$, $\zeta'$ and $\zeta''$
        given by equation~\eqref{definition-of-zetas}.
    }\label{zetas}%
    }\hfil
\end{figure}
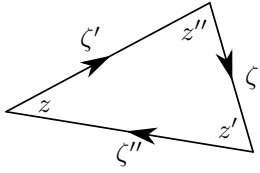%
Substituting \eqref{sines-from-zetas} into \eqref{initial-ratio-of-tau_2s}
and using equation~\eqref{Pachner-angle-formulas}, we obtain
\begin{align*}
\frac{\tau_2^{(3)}}{\tau_2^{(2)}}
&= \left(\frac{
        |\zeta_1|^{b_1} |\zeta'_1|^{b'_1} |\zeta''_1|^{b''_1}
        |\zeta_2|^{b_2} |\zeta'_2|^{b'_2} |\zeta''_2|^{b''_2}
        |\zeta_3|^{b_3} |\zeta'_3|^{b'_3} |\zeta''_3|^{b''_3}
    }
    {
        |\zeta_0|^{b'_3+b''_1} |\zeta'_0|^{b'_2+b''_3} |\zeta''_0|^{b''_2+b'_1}
        |\zeta_4|^{b'_2+b''_1} |\zeta'_4|^{b''_2+b'_3} |\zeta''_4|^{b''_3 + b'_1}
    }\right)^2\\
    &\phantom{{}={}} \times
    \left|
    \frac{ [z_1 z_2 z_3 (1-z_1)(1-z_2)(1-z_3)]^2 \Imag z_0 \Imag z_4 }
        { [z_0 z_4 (1-z_0)(1-z_4)]^2 \Imag z_1 \Imag z_2 \Imag z_3}
    \right|,
\end{align*}
since $b_j + b'_j + b''_j=1$ for all $0\leq j \leq4$.
Using the equations~\eqref{definition-of-zetas}, we
arrive at the following expression in terms of shapes:
\begin{align*}
\frac{\tau_2^{(3)}}{\tau_2^{(2)}}
    &= \left(\frac{
        |z_0|^{b'_3+b''_1+b''_2+b'_1-1}|1-z_0|^{b'_2+b''_3+b''_2+b'_1-1}
        |z_4|^{b'_2+b''_1+b''_3+b'_1-1}|1-z_4|^{b''_2+b'_3+b''_3+b'_1-1}
    }{
        |z_1|^{b_1+b''_1}|1-z_1|^{b'_1+b''_1-1}
        |z_2|^{b_2+b''_2}|1-z_2|^{b'_2+b''_2-1}
        |z_3|^{b_3+b''_3}|1-z_3|^{b'_3+b''_3-1}
    }\right)^2\\
    &\phantom{{}={}} \times
    \left| \frac{\Imag z_0 \Imag z_4}{\Imag z_1 \Imag z_2 \Imag z_3} \right|.
\end{align*}
The relations~\eqref{s=s(z)} enable us to express $|z_0|$, $|1-z_0|$, $|z_4|$
and $|1-z_4|$ in terms of the remaining shape variables.
The resulting formula can be simplified by combining the powers
of $|z_j|$ and $|1-z_j|$ ($1\leq j \leq 3$) and simplifying the exponents
with the help of the angle equation~$b_1+b_2+b_3=2$ along
the central edge of the bipyramid, as well as the relations~$b_j+b'_j+b''_j=1$.
After these simplifications, we obtain
\begin{align}
\frac{\tau_2^{(3)}}{\tau_2^{(2)}}
&=
    \left|z_1z_2z_3\right|^{2(-b_2-b_3+b'_1)}
    \left|z_1(1-z_2)(1-z_3)\right|^2
    \frac{|\Imag(z'_2z''_1)\Imag(z'_3z''_1)|}{|\Imag z_1 \Imag z_2 \Imag z_3|}
    \notag\\
&=    \frac{|\Imag(z'_2z''_1)\Imag(z'_3z''_1)|}
    {
        \frac{|z_1-1|}{|z_1(1-z_2)|} \cdot \frac{|z_1-1|}{|z_1(1-z_3)|}
    }\;\cdot\;
    \frac{|z_1 z_2 z_3|}
    {|\Imag z_1 \Imag z_2 \Imag z_3|}\label{last-rational}
    \;\cdot\;|z_1-1|^2|z_2-1|\,|z_3-1|,
\end{align}
where we have used the gluing equation $z_1z_2z_3=1$ along the central edge.
Note that in general, for $\alpha = \arg z$, we have
\[
    \sin\alpha = \frac{\Imag z}{|z|} \quad
    \text{and}\quad
    |z-1| = \frac{1}{|z'|} = \frac{\sin\alpha}{\sin\alpha''}.
\]
Applying the first of these identities to the first two factors of \eqref{last-rational}
and the second identity to the last factor puts the expression in the trigonometric form
\begin{equation}\label{first-trigonometric}
    \frac{\tau_2^{(3)}}{\tau_2^{(2)}}
    = \left|
        \frac{
            \sin(\alpha'_2+\alpha''_1)\sin(\alpha'_3+\alpha''_1)
            \sin\alpha_2\sin\alpha_3\sin^2\alpha_1
        }{
            \sin\alpha_2 \sin\alpha_3 \sin\alpha_1 \sin\alpha''_2 \sin\alpha''_3\sin^2\alpha''_1
        }\right|.
\end{equation}
Note that the complex gluing equation $z_1z_2z_3=1$ implies that
\begin{equation}\label{single-double-prime-sines}
    \sin\alpha'_1\sin\alpha'_2\sin\alpha'_3 = \pm\sin\alpha''_1\sin\alpha''_2\sin\alpha''_3.
\end{equation}
Using the above relation, we can further simplify \eqref{first-trigonometric} to
\begin{equation}\label{first-sinal}
    \frac{\tau_2^{(3)}}{\tau_2^{(2)}}
    = \left|
        \frac{
            \sin(\alpha'_2+\alpha''_1)\sin(\alpha'_3+\alpha''_1)\sin\alpha_1
        }{
            \sin\alpha'_2 \sin\alpha'_3\sin\alpha'_1\sin\alpha''_1
        }\right|.
\end{equation}
Recall the trigonometric identities
\begin{equation}\label{trig-identities}
    \sin(\alpha+\beta) = \cos\alpha\sin\beta + \cos\beta\sin\alpha,
    \quad\text{whence}\quad
    \frac{\sin(\alpha + \beta)}{\sin\alpha\; \sin\beta} = \cot\alpha + \cot\beta.
\end{equation}
We can now simplify the trigonometric expression occurring in \eqref{first-sinal} by repeatedly
applying \eqref{trig-identities}:
\begin{align*}
    &\phantom{{}={}}\frac{
        \sin(\alpha'_2+\alpha''_1)\sin(\alpha'_3+\alpha''_1)\sin\alpha_1
    }{
        \sin\alpha'_2 \sin\alpha'_3\sin\alpha'_1\sin\alpha''_1
    }\\
    &=
    \frac{
        (\sin\alpha'_2\cos\alpha''_1+\sin\alpha''_1\cos\alpha'_2)
        (\sin\alpha'_3\cos\alpha''_1+\sin\alpha''_1\cos\alpha'_3)\sin\alpha_1
    }{
        \sin\alpha'_2\sin\alpha'_3\sin\alpha'_1\sin\alpha''_1
    }
\\
    &=
    \frac{
        \bigl(\sin\alpha'_2 \sin\alpha'_3 (1-\sin^2\alpha''_1)
        + \sin\alpha'_2\cos\alpha'_3\cos\alpha''_1\sin\alpha''_1
        \bigr)\sin\alpha_1
    }
    {
        \sin\alpha'_2\sin\alpha'_3\sin\alpha'_1\sin\alpha''_1
    }
\\
    &\phantom{{}={}}+ \frac{
        \bigl(
        \sin\alpha'_3\cos\alpha'_2\sin\alpha''_1\cos\alpha''_1
        + \sin^2\alpha''_1\cos\alpha'_2\cos\alpha'_3\bigr)\sin\alpha_1
    }
    {
        \sin\alpha'_2\sin\alpha'_3\sin\alpha'_1\sin\alpha''_1
    }
\\
&=
    \frac{ \sin\alpha'_2 \sin\alpha'_3 \sin\alpha_1}{\sin\alpha'_2\sin\alpha'_3\sin\alpha'_1\sin\alpha''_1}
    +
    \frac{\bigl(
        -\sin\alpha'_2\sin\alpha'_3\sin^2\alpha''_1 + \sin^2\alpha''_1\cos\alpha'_2\cos\alpha'_3\bigr)\sin\alpha_1
    }{\sin\alpha'_2\sin\alpha'_3\sin\alpha'_1\sin\alpha''_1}\\
    &\phantom{{}={}}+ \frac{\sin\alpha_1\cos\alpha''_1 (\sin\alpha'_2\cos\alpha'_3+\sin\alpha'_3\cos\alpha'_2)}
            {\sin\alpha'_2\sin\alpha'_3\sin\alpha'_1}
\\
    &=
    \frac{\sin\alpha_1}{\sin\alpha'_1\sin\alpha''_1}
    +
    \frac{\sin\alpha_1\sin\alpha''_1(\cos\alpha'_2\cos\alpha'_3 - \sin\alpha'_2\sin\alpha'_3)}
    {\sin\alpha'_2\sin\alpha'_3\sin\alpha'_1}\\
    &\phantom{{}={}}+ \frac{
    \bigl(\sin(\alpha_1+\alpha''_1)-\cos\alpha_1\sin\alpha''_1\bigr)\sin(\alpha'_2+\alpha'_3)}
    {\sin\alpha'_2\sin\alpha'_3\sin\alpha'_1}
\\
    &=
    \frac{\sin(\alpha'_1+\alpha''_1)}{\sin\alpha'_1\sin\alpha''_1}
    +\frac{\sin\alpha_1\sin\alpha''_1\cos(\alpha'_2+\alpha'_3)}{\sin\alpha''_2\sin\alpha''_3\sin\alpha''_1}\\
    &\phantom{{}={}}+\frac{\sin(\pi-\alpha'_1)\sin(\alpha'_2+\alpha'_3)}{\sin\alpha'_2\sin\alpha'_3\sin\alpha'_1}
    -\frac{\cos\alpha_1\sin\alpha''_1\sin(\alpha'_2+\alpha'_3)}{\sin\alpha''_2\sin\alpha''_3\sin\alpha''_1}
\\
    &=
    \frac{\sin(\alpha'_1+\alpha''_1)}{\sin\alpha'_1\sin\alpha''_1}
    +\frac{\sin(\alpha'_2+\alpha'_3)}{\sin\alpha'_2\sin\alpha'_3}
    +\frac{\sin\alpha_1\cos(\alpha'_2+\alpha'_3)-\cos\alpha_1\sin(\alpha'_2+\alpha'_3)}{\sin\alpha''_2\sin\alpha''_3}.
\end{align*}
The angle equation $\alpha_1+\alpha_2+\alpha_3=2\pi$ implies that the numerator of the last term equals
\begin{align*}
    \sin\alpha_1\cos(\alpha'_2+\alpha'_3)-\cos\alpha_1\sin(\alpha'_2+\alpha'_3)
    &= \sin(\alpha_1-\alpha'_2-\alpha'_3) \\
    &= \sin(2\pi - \alpha_2 - \alpha_3 - \alpha'_2 - \alpha'_3)\\
    &= \sin(\alpha''_2 + \alpha''_3).
\end{align*}
Applying \eqref{trig-identities} one last time, we obtain the following result:
\begin{equation}\label{symmetric-sine-anomaly}
    \frac{\tau_2^{(3)}}{\tau_2^{(2)}}
    = \left|
        \frac{\sin(\alpha'_1+\alpha''_1)}{\sin\alpha'_1\sin\alpha''_1}
        +\frac{\sin(\alpha'_2+\alpha'_3)}{\sin\alpha'_2\sin\alpha'_3}
        +\frac{\sin(\alpha''_2+\alpha''_3)}{\sin\alpha''_2\sin\alpha''_3}
    \right|
    = \left|\sum_{j=1}^3(\cot\alpha'_j + \cot\alpha''_j)\right|.
\end{equation}
Combining \eqref{symmetric-sine-anomaly} with \eqref{ratio-of-tau1s}
shows that $\tau=\sqrt{\tau_1\tau_2}$ is invariant under shaped Pachner moves.
\end{proof}

\begin{remark}
In the special case when $z=z^\geo$ is a geometric solution
recovering the complete hyperbolic structure on $M$,
the associated strict angle structure~$\alpha=\Arg z$
is a smooth non-degenerate maximum point of the volume function~$\Vol:\Astrict_0(\triang)\to\RR$,
so that the Hessian of $\Vol$ is negative-definite at $\alpha$.
In this situation, $\tau(z)$ can be interpreted as a
state integral of Turaev--Viro type with a Gaussian Boltzmann weight.
To make this precise, we consider a labelled tetrahedron $\tet$ with angles~$\alpha=\alpha_{01}$,
$\alpha'=\alpha_{02}$, $\alpha''=\alpha_{03}$, as depicted in Figure~\ref{fig:tetrahedron}.
We associate to $\tet$ the charged Boltzmann weight
\begin{equation}\label{gaussian-boltzmann-weight}
    B_\alpha^\tet\!
    \left(\begin{smallmatrix} x_{01} & x_{02} & x_{03}\\
                          x_{23} & x_{13} & x_{12}\end{smallmatrix}\right)
    = (\sin\alpha)^{\frac{\alpha}{\pi}  -\frac{1}{2}}
     (\sin\alpha')^{\frac{\alpha'}{\pi} -\frac{1}{2}}
    (\sin\alpha'')^{\frac{\alpha''}{\pi}-\frac{1}{2}}
    e^{-\pi b(x_{01}+x_{23},\,x_{02}+x_{13},\,x_{03}+x_{12})},
\end{equation}
where the variables $x_{mn}$ are real and the quadratic form $b$ is given explicitly by
\begin{align*}
    b(t,t',t'') &= (\cot\alpha' + \cot\alpha'')t^2 +
                  (\cot\alpha + \cot\alpha'')(t')^2 +
                  (\cot\alpha + \cot\alpha')(t'')^2\\
    &\phantom{{}={}} - 2(\cot\alpha\,t't'' + \cot\alpha'\,tt'' + \cot\alpha''\,tt').
\end{align*}
In this way,
\[
    \tau(z,\alpha) = \int_{\RR^{N-1}}
    \prod_{\tet\ \text{in}\ \triang}
    B_\alpha^\tet\!
    \Bigl|_{
        x_{mn}=y_j\ \text{iff}\ \Delta([mn])\in E_j\in\edges(\triang)
        }
    dy_1\dotsi dy_{N-1},
\]
where $y_N=-\sum_{j=1}^{N-1}y_j$.
In the context of the Teichm\"uller TQFT \cite{kashaev-andersen-TQFT},
a Gaussian tetrahedral weight similar to \eqref{gaussian-boltzmann-weight} was also
found by Kashaev~\cite{kashaev-private}.
\end{remark}

%==============================================================================
% Section 7: Relationship with the 1-loop invariant
%==============================================================================

\section{Relationship with the 1-loop invariant of Dimofte--Garoufalidis}
\label{sec:oneloop}
In \cite{tudor-stavros}, Dimofte and Garoufalidis described a non-rigorous
construction of a state integral invariant associated to an ideal triangulation
with hyperbolic shapes.
Using Feynman diagrams, they found explicit expressions for the expected
formal asymptotic coefficients of the proposed invariant, called
the \emph{\nword{n}{loop} terms}, for $n\geq 0$.
In this section, we focus on the `\nword{1}{loop} invariant'
which is determined by the \nword{1}{loop} term.
We start by reviewing the definition of the \nword{1}{loop} invariant, based
in part on the work of the third author~\cite{siejakowski-infinitesimal}.

To make the formulae more explicit, we number both the edges and the tetrahedra in $\triang$
with integers $\{1,\dotsc,N\}$ and refer to the notations of Section~\ref{sec:gluing-equations}.
We suppose that $\triang$ has a single ideal vertex, with a toroidal link,
and we choose a peripheral curve $\gamma$ in normal position
with respect to $\triang$.
We shall denote by $G(\gamma,\square)$ the incidence numbers of \eqref{loghol}.
Since one of the edge consistency equations~\eqref{edge-consistency-equations} is always redundant,
we may eliminate the last equation ($j=N$) and replace it with the
completeness equation~\eqref{completeness-equations} along $\gamma$.
Explicitly, the matrices of exponents in this modified system of equations are given by
\begin{equation}\label{definition-G-hat}
    \mathbf{\hat{G}} = [\,\hat{G}\;|\;\hat{G}'\;|\;\hat{G}''\,],\qquad
    [\mathbf{\hat{G}}]_{j,\square} =
        \begin{cases}
            \mathbf{G}_{j,\square}=G(E_j,\square), &\text{when}\ j < N,\\
            G(\gamma,\square), &\text{when}\ j = N.
        \end{cases}
\end{equation}
From now on, we assume that $z\in \gluvar$ is a complex solution of Thurston's edge consistency
equations on $\triang$ which satisfies either of the two conditions:
\begin{enumerate}[(i)]
    \item
    the solution $z$ recovers the holonomy representation $\rho^\geo$
    of the complete hyperbolic structure of finite volume, which means
    that $z$ also satisfies the completeness equation $u=0$; or
    \item
    the solution $z$ does not satisfy the completeness equation
    but remains in the range of the parameter $u$,
    so that there is a curve in $\gluvar$ connecting $z$ to a solution for $\rho^\geo$.
\end{enumerate}
In terms of a solution $z$ as above,
the \nword{1}{loop} term of Dimofte--Garoufalidis~\cite{tudor-stavros} is then
defined as follows.

\begin{definition}
Let $\gamma$ be an arbitrarily chosen peripheral curve in normal position with respect to $\triang$.
For all $j\in\{1,\dotsc,N\}$, we define $\zeta_j$, $\zeta'_j$ and $\zeta''_j$
in terms of the shape parameter $z_j$ as in equation \eqref{definition-of-zetas}.
Moreover, let $f$ be any peripherally trivial integer-valued angle structure (see Definition~\ref{def:Zvals}).
Then the \emph{\nword{1}{loop} invariant} associated to the solution $z$ and
the curve $\gamma$ is given by the formula
\[
    \oneloop(z, \gamma) = \pm\frac{
        \det\bigl(
              \hat{G}\diag\{\zeta_j\} +
              \hat{G}'\diag\{\zeta'_j\} +
              \hat{G}''\diag\{\zeta''_j\} \bigr)
        }{
        2\prod_{j=1}^{N} {\zeta_j}^{f_j}{\zeta'_j}^{f'_j}{\zeta''_j}^{f''_j}
    },
\]
where
 $(f_j, f'_j, f''_j)=\bigl(f(\square_j), f(\square'_j),f(\square''_j)\bigr)$.
\end{definition}

Viewing the \nword{1}{loop} invariant as a complex number defined up to sign only,
Dimofte and Garoufalidis prove in \cite[Section~3]{tudor-stavros} that $\oneloop(z, \gamma)$
does not depend on the choice of the integer-valued angle structure~$f$ or on any of the
labelling choices.
For geometric triangulations, the value of $\oneloop(z, \gamma)$ depends only on the homology
class of $\gamma$ in $H_1(\partial M;\ZZ)$ \cite{siejakowski-infinitesimal}.

The `\nword{1}{loop} conjecture' \cite[Conjecture~1.8]{tudor-stavros} postulates that
\begin{equation}\label{oneloop-conjecture-eq}
    \oneloop(z, \gamma) = \pm\torsion(M,\rho_z, \gamma),
\end{equation}
where $\torsion(M,\rho_z,\gamma)\in\CC^*/\{\pm1\}$ is the adjoint Reidemeister
torsion introduced by Porti~\cite{porti1997}, corresponding to the conjugacy class of
a representation~$\rho_z$ determined by the solution~$z$.
Porti's construction of $\torsion$ requires the choice of non-trivial element of
$H_1(\partial M;\CC)$, which is accomplished
here by taking the (integral) homology class of the oriented curve~$\gamma$.
Thus, the \nword{1}loop conjecture predicts an explicit formula for the torsion
in terms of hyperbolic shapes.
For geometric ideal triangulations,
both sides of the conjectural equality~\eqref{oneloop-conjecture-eq}
transform in the same way under the change of the curve~$\gamma$;
cf. \cite[Corollary~5.2]{siejakowski-infinitesimal}.
In other words, given two peripheral curves $\gamma$, $\tilde{\gamma}$,
we have
\begin{equation}\label{torsion-oneloop-transform-alike}
    \frac{\oneloop(z, \tilde{\gamma})}{\oneloop(z, \gamma)}
    = \pm\frac{\torsion(M,\rho_z,\tilde{\gamma})}{\torsion(M,\rho_z,\gamma)}.
\end{equation}
Using a result of Porti \cite[Corollaire~4.3]{porti1997} we see that
when $\rho_z=\rho^\geo$ is boundary-parabolic, the ratio of torsions on the
right-hand side is given by
\begin{equation}\label{Porti-changement-de-courbes}
    \frac{\torsion(M,\rho,\tilde{\gamma})}{\torsion(M,\rho,\gamma)}
    = \pm\sigma(\tilde{\gamma}, \gamma),
\end{equation}
where $\sigma(\tilde{\gamma}, \gamma)$ is the relative modulus of the curves
$\tilde{\gamma}$ and $\gamma$, introduced in equation~\eqref{def-relative-modulus}.
%------------------------------------------------- -------------------------------------------------
\subsection{\texorpdfstring{The 1-loop conjecture for $\tau$}{The 1-loop conjecture for tau}}
\label{tau-vs-oneloop}
The following theorem provides an identification of $\tau$ with a
certain normalisation of the \nword{1}{loop} invariant
for the representations $\rho^\geo$ and $\overline{\rho}^\geo$.

\begin{theorem}\label{th:tau=oneloop}
    Suppose that $z$ is a geometric solution of the
    edge consistency and completeness equations on $\triang$.
    Consider the conformal class of non-degenerate
    Euclidean structures on the torus~$\partial M$ determined by the complete hyperbolic structure on $M$.
    Choose a scale so that the Euclidean structure on $\partial M$ has area~$A>0$
    and denote by $L(\gamma)>0$ the length of a geodesic representative
    of the free homotopy class of $\gamma$ in this Euclidean metric.
    Then
    \begin{equation}\label{eq:tau=oneloop}
        \tau(z) = \frac{L(\gamma)}{\sqrt{2A}}\left|\oneloop(z,\gamma)\right|^{-1}.
    \end{equation}
\end{theorem}
\begin{remark}
Porti's change-of-curves formula~\eqref{Porti-changement-de-courbes} immediately implies that
the right-hand side of \eqref{eq:tau=oneloop} is independent of the choice of $\gamma$.
Indeed, suppose that $\tilde{\gamma}$ is another curve in $\partial M$
with a Euclidean geodesic representative of length $L(\tilde{\gamma})$.
By combining equations \eqref{torsion-oneloop-transform-alike}
and \eqref{Porti-changement-de-courbes}, we obtain
\[
    \frac{
        L(\tilde{\gamma})\left|\oneloop(z,\tilde{\gamma})\right|^{-1}
    }{
        L(\gamma)\left|\oneloop(z,\gamma)\right|^{-1}
    }
    =
    \frac{L(\tilde{\gamma})}{L(\gamma)\abs{\sigma(\tilde{\gamma},\gamma)}}
    = 1.
\]
\end{remark}

\begin{remark}
Observe that the definition of the \nword{1}{loop} invariant~$\oneloop(z,\gamma)$
only requires non-degenerate shape parameters $z(\square)\in\CC\setminus\{0,1\}$, for all $\square\in Q(\triang)$,
whereas the definition of $\tau(z)$ requires the stronger condition $z(\square)\in\CC\setminus\RR$.
Hence, the right-hand side of \eqref{eq:tau=oneloop} provides an extension of $\tau$ to shaped triangulations
containing non-degenerate zero-volume tetrahedra.
\end{remark}

\begin{remark}
    We remark that the equality \eqref{torsion-oneloop-transform-alike}
    suggests a new way of making the adjoint Reidemeister torsion of \cite{porti1997} independent
    of the choice of the peripheral curve~$\gamma$.
    More precisely, if $\torsion(M,\rho^\geo,\gamma)$ is the
    adjoint Reidemeister torsion associated to the holonomy representation of the complete
    hyperbolic structure,
    then the quantity
    \begin{equation}\label{torsion-normalised}
        \torsionNormalised(M,\rho^\geo)
        \eqdef
        \frac{\abs{\torsion(M,\rho^\geo,\gamma)}\sqrt{A}}%
        {L(\gamma)}
        \in \RR_+
    \end{equation}
    does not depend on the choice of the curve~$\gamma$.
\end{remark}

We will next prove Theorem~\ref{th:tau=oneloop}.
This calculation requires the following technical result which will let us make a peripherally-trivial adjustment of a tangent vector to the space of angle structures to a tangent vector which can be lifted to a path in the gluing variety.
\begin{proposition}
\label{construction-of-lhat}
Let $\triang$ be an ideal triangulation with a single ideal vertex and
assume that $z_0$ is a geometric solution of edge consistency and completeness equations
so that $\alpha_0 = \Arg z_0$ is a peripherally trivial strict angle structure.
For an arbitrarily fixed peripheral curve $\mu$, we consider $l_\mu\in\TAS(\triang)$
given in \eqref{peripheral-LTD}.
Then there exists a unique vector $\hat{l}_\mu\in\TAS(\triang)$ such that
\begin{enumerate}[(i)]
    \item \label{difference-in-TAS0}
        $\hat{l}_\mu - l_\mu\in\TAS_0(\triang)$, and
    \item \label{tangent-to-gluvar}
        There exists an $\varepsilon>0$ and a one-parameter
        family~$t\mapsto z_t\in\gluvar$ ($-\varepsilon<t<\varepsilon$)
        for which $\hat{l}_\mu$ is the initial tangent vector to the family of angles, i.e.,
        $\hat{l}_\mu = \frac{d}{dt}\bigr|_{t=0}(\Arg z_t)$.
\end{enumerate}
\end{proposition}
We will defer the proof of this proposition to the end of this section.

%--------------------------------------------------------------------------------------------------
\begin{proof}[Proof of Theorem~\ref{th:tau=oneloop}]
Let $f$ be a peripherally trivial integer-valued angle structure in the sense of
Definition~\ref{def:Zvals}.
In view of Proposition~\ref{tau-independent-of-exponents}, we may
always set $\beta=\pi f$ in \eqref{tau-recalled}.
We split the proof into several steps.
\Step{1: Eliminating the trigonometric functions.}
We start by inserting the expressions \eqref{sines-from-zetas} into
the numerator of \eqref{tau-recalled}, obtaining
\begin{align}
    \prod_{\square\in Q(\triang)}
            |\sin\alpha(\square)|^{f(\square)-\frac{1}{2}}
    &=
    \prod_{j=1}^N \frac{|\Imag z_j|^{(f_j + f'_j + f''_j)-3/2}
            \left|{\zeta_j}^{f_j} {\zeta'_j}^{f'_j} {\zeta''_j}^{f''_j}\right|}
            {\sqrt{|\zeta_j \zeta'_j \zeta''_j|}}\notag\\
        &= \prod_{j=1}^N \left( \left|{\zeta_j}^{f_j} {\zeta'_j}^{f'_j} {\zeta''_j}^{f''_j}\right|
        \frac{|z_j(1-z_j)|}{\sqrt{|\Imag z_j|}}\right),\label{tau-numerator-without-sines}
\end{align}
where we have used the fact that $f_j+f'_j+f''_j=1$ for every $j\in\{1,\dotsc,N\}$.

For each shape parameter $z(\square)$, choose a logarithm $Z(\square)$
so that $Z(\square)+Z(\square')+Z(\square'')=i\pi$ in every tetrahedron.
We customarily denote these log-parameters by $Z_n$, $Z'_n$, $Z''_n$ subject
to our numbering of the tetrahedra ($1\leq n \leq N$) and the choice of
normal quadrilateral types.
For $1\leq j \leq N-1$, consider the function $X_j=X_j(Z)$
given by the log-holonomy around the edge $E_j$ and write $\mathbf{X}=(X_1,\dotsc,X_{N-1})$.
Additionally, we denote by $u=u(Z)$ the log-holonomy parameter along
the curve $\gamma$ and set $\mathbf{\hat{X}}=(\mathbf{X},u)$.
In terms of the matrices $\hat{G}$, $\hat{G}'$ and $\hat{G}''$, we have
\begin{equation}\label{enlarged-X}
    \mathbf{\hat{X}}
    = \coldotsvec{X_1}{X_{N-1}\\ u}
    =   \hat{G}  \coldotsvec{Z_1}{Z_N}
      + \hat{G}' \coldotsvec{Z'_1}{Z'_N}
      + \hat{G}''\coldotsvec{Z''_1}{Z''_N}.
\end{equation}
Moreover, we have
\begin{equation}\label{log-param-in-angles}
    Z(\square) \equiv \log\left|\frac{\sin\alpha(\square')}{\sin\alpha(\square'')}\right| + i\alpha(\square)
    \mod 2\pi i.
\end{equation}

Denote by $l_1,\dotsc,l_{N-1}$ the leading-trailing deformation vectors along the respective
edges of $\triang$ and let $t=(t_1,\dotsc,t_{N-1})\in\RR^{N-1}$.
In a small neighbourhood of $t=0$, we may treat $X_j$ as the composition
$t\mapsto X_j(Z_t)$, where the log-parameters $Z_t$ result from substituting a deformed angle
structure~$\alpha_t$ into \eqref{log-param-in-angles}, given by
$\alpha_t(\square) = \alpha(\square) + \sum_{m=1}^{N-1}t_ml_m(\square)$.
Using \eqref{explicit-Jacobian-of-shears}, for any $j,m\in\{1,\dotsc,N-1\}$ we compute
\begin{align}
    \frac{\partial\Real X_j}{\partial t_m}\biggr|_{t_m=0}
    &= -\frac{\partial\Shear_j}{\partial t_m}\biggr|_{t_m=0}\notag\\
    &= [\mathbf{L_*}\diag\{\cot\alpha(\square)\}\mathbf{L_*}^\transp]_{j,m}.\label{LdiagcotL-from-X}
\end{align}
Moreover, since the imaginary part of $X_j$ is the
angle sum around the edge $E_j$ and is therefore preserved by all leading-trailing
deformations, we have $\dfrac{\partial\Imag X_j}{\partial t_m}=0$.
Hence, the matrix \eqref{LdiagcotL-from-X}
is in fact the Jacobian matrix of the complex-valued function~$\mathbf{X}(t)$
at $t=0$.
Using \eqref{tau-numerator-without-sines}, we obtain the following formula for $\tau$:
\begin{equation}\label{tau-with-Jacobian-of-X}
    \tau(z)^2
    =
    \prod_{j=1}^N \left[ \left|{\zeta_j}^{f_j} {\zeta'_j}^{f'_j} {\zeta''_j}^{f''_j}\right|^2
        \frac{|z_j(1-z_j)|^2}{|\Imag z_j|}\right]\times
        \left|\det\left[ \frac{\partial X_j}{\partial t_m}\right]_{j,m}\right|^{-1}.
\end{equation}

\Step{2: Decomposition of the Jacobian.}
Since we are treating $\mathbf{X}(t)$ as the composition
$t\mapsto \mathbf{X}(Z_t)$
of differentiable maps, the Jacobian~$D\mathbf{X}(0)$ is the composition of two
derivatives which we now compute.
In order to simplify notation, we denote the three normal
quads in the \nword{n}{th} tetrahedron by $\square_n$, $\square'_n$
and $\square''_n$ and we write $\alpha_n=\alpha(\square_n)$,
$\alpha'_n=\alpha(\square'_n)$, $\alpha''_n=\alpha(\square''_n)$.
Using \eqref{log-param-in-angles} with $\alpha_t$ we find, for any $1\leq n \leq N$,
\begin{align*}
    \frac{\partial Z_n}{\partial t_m}
    &= \left[  l_m(\square_n)  \frac{\partial}{\partial\alpha_n} +
              l_m(\square'_n) \frac{\partial}{\partial\alpha'_n} +
              l_m(\square''_n)\frac{\partial}{\partial\alpha''_n}
      \right]
      \left(\log\left|\frac{\sin\alpha'_n}{\sin\alpha''_n}\right| + i\alpha_n\right)\\
    &= l_m(\square_n)i
    + l_m(\square'_n)\cot\alpha'_n
    - l_m(\square''_n)\cot\alpha''_n.
\end{align*}
Secondly, since $z_n=\exp(Z_n)$ does not depend on $Z_m$ for $m\neq n$, we may write
\begin{align*}
    \frac{\partial X_j}{\partial Z_n}
    = \frac{\partial X_j}{\partial z_n}\frac{\partial z_n}{\partial Z_n}
    &= \left( G_{jn}\frac{d Z_n}{d z_n}
           + G'_{jn}\frac{d Z'_n}{d z_n}
           + G''_{jn}\frac{d Z''_n}{d z_n}\right) \frac{d(e^{Z_n})}{\partial Z_n}\\
    &= \left( G_{jn}\zeta_n + G'_{jn}\zeta'_n + G''_{jn}\zeta''_n\right)z_n.
\end{align*}
Therefore, the Jacobian matrix occurring in \eqref{tau-with-Jacobian-of-X} can be rewritten as the product
\begin{align}
    \left[\frac{\partial X_j}{\partial t_m}\right]_{j,m}
    &=
    \left[
        \frac{\partial X_j}{\partial Z_n}
    \right]_{j,n}\times
    \left[
        \frac{\partial Z_n}{\partial t_m}
    \right]_{n,m} \label{explicit-jacobian-X}\\
    &=      \bigl(
            G_*\diag\{\zeta_n\}_n + G'_*\diag\{\zeta'_n\}_n + G''_*\diag\{\zeta''_n\}_n
        \bigr)\diag\{z_n\}_n
    \\
    &\phantom{{}={}}\times
    \left(
        iL_*^\transp + \diag\{\cot\alpha'_n\}_n {L'_*}^\transp - \diag\{\cot\alpha''_n\}_n {L''_*}^\transp
    \right),\notag
\end{align}
where the matrices $G_*$, $G'_*$, $G''_*$ consist of rows $1,\dotsc,N-1$ of the
matrices $G$, $G'$ and $G''$, respectively.

\Step{3: Enlarging the matrices.}
We now wish to replace the rectangular matrices $G_*$, $G'_*$, $G''_*$ with the square matrices
$\hat{G}$, $\hat{G}'$, $\hat{G}''$ containing additionally
the coefficients of the completeness equation along $\gamma$.
Correspondingly, we use the vector $\hat{l}_\gamma\in\TAS(\triang)$ of
Proposition~\ref{construction-of-lhat} in the last row of
the enlarged $N\times N$ matrices $\hat{L}$, $\hat{L}'$, $\hat{L}''$, setting
\begin{equation}\label{definition-L-hat}
    \hat{L}_{m,n} = \begin{cases}
                          [L_*]_{m,n}, &\text{when}\ m < N,\\
                          \hat{l}_\gamma(\square_n), &\text{when}\ m = N,
                   \end{cases}
\end{equation}
and analogously for $\hat{L}'$ and $\hat{L}''$.
By construction, the vector $\hat{l}_\gamma$ is tangent to
the gluing variety of $\triang$.
Denote by $\nabla_\gamma = \nabla_{\hat{l}_\gamma}$ the directional derivative
in the direction of $\hat{l}_\gamma$.
In particular, we have $\nabla_\gamma X_j = 0$ for all $1\leq j \leq N-1$.

We may now view the quantity $\mathbf{\hat{X}}=(X_1,\dotsc,X_{N-1}; u)$ of \eqref{enlarged-X}
as a function of the variables $\hat{t}=(t_1,\dotsc,t_{N-1}; t_\gamma)$,
where the additional coordinate $t_\gamma$ is tangent to $\hat{l}_\gamma$ at zero.
In this way, the enlarged Jacobian~$D\mathbf{\hat{X}}(\hat{t})$ at $\hat{t}=0$
has the block form
\[
    D\mathbf{\hat{X}}(0) = \left[
    \begin{tikzpicture}[baseline={([yshift=-\the\dimexpr\fontdimen22\textfont2\relax]
                    current bounding box.center)}]
        \node[inner xsep=1em, inner ysep=8ex] (w) at (0,0)
        {\(\displaystyle
            \left[\frac{\partial X_j}{\partial t_m}\right]_{j,m\in\{1,\dotsc,N-1\}}
        \)};
        \draw ({w.south west}) rectangle ({w.north east});
        \node[yshift=-0.8ex, anchor=north west] (b1) at ({w.south west})
            {$\frac{\partial u}{\partial t_1}$};
        \node[yshift=-2ex, anchor=north] (b2) at ({w.south}) {$\dots$};
        \node[yshift=-0.8ex, anchor=north east] (b3) at ({w.south east})
            {$\frac{\partial u}{\partial t_{N-1}}$};
        \node[xshift=0.6em, anchor=north west] (r1) at ({w.north east})
            {$\frac{\partial X_1}{\partial t_\gamma}$};
        \node[xshift=1.2em, anchor=west] (r2) at ({w.east}) {$\vdots$};
        \node[xshift=0.6em, anchor=south west] (r1) at ({w.south east})
            {$\frac{\partial X_{N-1}}{\partial t_\gamma}$};
        \node[yshift=-0.8ex, xshift=1.2em, anchor=north west] (c) at ({w.south east})
            {$\frac{\partial u}{\partial t_\gamma}$};
    \end{tikzpicture}
    \right]_{\mathbf{\hat{t}}=0}
    =
    \begin{bmatrix}
    D\mathbf{X}(0) & 0\\
    * & \nabla_\gamma u
    \end{bmatrix}.
\]
Since $\det D\mathbf{\hat{X}}(0)=\det D\mathbf{X}(0)\,\cdot\,\nabla_\gamma u$, we may
rewrite \eqref{tau-with-Jacobian-of-X} as
\begin{equation}\label{tau-with-enlarged-jacobian}
    \tau(z)^2
            =
        \frac{\abs{\nabla_\gamma u}}
        {\left|\det D\mathbf{\hat{X}}(0)\right|}
        \prod_{j=1}^N \left[ \left|{\zeta_j}^{f_j} {\zeta'_j}^{f'_j} {\zeta''_j}^{f''_j}\right|^2
        \frac{|z_j(1-z_j)|^2}{|\Imag z_j|}\right].
\end{equation}

\Step{4: Computing $\det D\mathbf{\hat{X}}(0)$.}
We now compute the Jacobian determinant $\det D\mathbf{\hat{X}}(0)$ explicitly
in terms of the shape parameters and the enlarged matrices of \eqref{definition-G-hat}
and \eqref{definition-L-hat}.
To this end, we substitute $\cot\alpha'_j = \frac{1-\Real z_j}{\Imag z_j}$ and
$\cot\alpha''_j = \frac{|z_j|^2-\Real z_j}{\Imag z_j}$ into \eqref{explicit-jacobian-X},
obtaining
\begin{align}
\det D\mathbf{\hat{X}}(0) &=
    \det\left(
            \hat{G}\diag\{\zeta_j\} + \hat{G}'\diag\{\zeta'_j\} + \hat{G}''\diag\{\zeta''_j\}
        \right)\det\diag\{z_j\}_j\notag\\
    &\phantom{{}={}}\times\det
    \left(
        i{\hat{L}}^\transp + \diag\left\{\frac{1-\Real z_j}{\Imag z_j}\right\}_j\left.\hat{L}'\right.^\transp -
        \diag\left\{\frac{|z_j|^2-\Real z_j}{\Imag z_j}\right\}_j\left.\hat{L}''\right.^\transp
    \right).\label{second-line-junk}
\end{align}
As per \eqref{definition-L-hat}, the last rows of the matrices $\hat{L}$, $\hat{L}'$, $\hat{L}''$
are given by $\hat{l}_\gamma$.
By Part~\eqref{difference-in-TAS0} of Proposition~\ref{construction-of-lhat},
the difference $\hat{l}_\gamma - l_\gamma$ is a linear combination of rows of $\mathbf{L_*}$.
Hence, replacing $\hat{l}_\gamma$ by $l_\gamma$ in the last rows
will not affect the determinant in the second line~\eqref{second-line-junk}.
Therefore,
\begin{align}
    &\det\left(
        i{\hat{L}}^\transp + \diag\left\{\frac{1-\Real z_j}{\Imag z_j}\right\}_j\left.\hat{L}'\right.^\transp -
        \diag\left\{\frac{|z_j|^2-\Real z_j}{\Imag z_j}\right\}_j\left.\hat{L}''\right.^\transp
    \right)\notag\\
    = &\det\left(\notag
    (\hat{G}''-\hat{G}')i + (\hat{G}-\hat{G}'')\diag\left\{\tfrac{1-\Real z_j}{\Imag z_j}\right\}_j
    + (\hat{G}'-\hat{G})\diag\left\{\tfrac{\Real z_j - |z_j|^2}{\Imag z_j}\right\}_j\right)\\
    = &\det\Bigl[\bigl(
             \hat{G}\diag\{\bar{\zeta}_j\}
           + \hat{G}'\diag\{\bar{\zeta}'_j\}
           + \hat{G}''\diag\{\bar{\zeta}''_j\}
    \bigr)
    \diag\left\{\tfrac{\bar{z}_j}{\Imag z_j} (1-z_j)(1-\bar{z}_j)\right\}_j\Bigr],
    \label{last-messy}
\end{align}
where the second equality is trivially verified with the help of~\eqref{definition-of-zetas}.
In Corollary~\ref{derivative-of-u-wrt-lhat} 
below 
we obtain the identity
$\abs{\nabla_\gamma u} = 2L(\gamma)^2/A$.
Inserting the expression~\eqref{last-messy} back into the line~\eqref{second-line-junk} and then
into \eqref{tau-with-enlarged-jacobian}, we obtain
\begin{align*}
    \tau(z)^2
            &= \abs{\nabla_\gamma u}
        \prod_{j=1}^N \left[ \left|{\zeta_j}^{f_j} {\zeta'_j}^{f'_j} {\zeta''_j}^{f''_j}\right|^2
        \frac{|z_j(1-z_j)|^2}{|\Imag z_j|}\right]
            \times
        \left|\det D\mathbf{\hat{X}}(0)\right|^{-1}\\
    &= \frac{2L(\gamma)^2}{A}\frac{
        \prod_{j=1}^{N} \left|{\zeta_j}^{f_j}{\zeta'_j}^{f'_j}{\zeta''_j}^{f''_j}\right|^2
    }
    {
        \left|\det\bigl(
              \hat{G}\diag\{\zeta_j\} +
              \hat{G}'\diag\{\zeta'_j\} +
              \hat{G}''\diag\{\zeta''_j\} \bigr)\right|^2
        }\\
    &= \frac{L(\gamma)^2}{2A}\left|\oneloop(z,\gamma)\right|^{-2}.
\end{align*}
Taking square roots now produces the desired equation~\eqref{eq:tau=oneloop}.
\end{proof}

\begin{proof}[Proof of Proposition \ref{construction-of-lhat}]
There always exists a dual curve~$\lambda\subset\partial\overline{M}$
satisfying $\intersection(\mu,\lambda)=1$.
Let $u=\log\hol(\mu)$ and $v=\log\hol(\lambda)$ be the corresponding local parameters
on $\gluvar$ in a neighbourhood of $z_0$, so that $u=v=0$ at $z_0$.
Neumann and Zagier proved \cite[Lemma~4.1]{neumann-zagier} that $v=u\tau(u)$ for a
certain even holomorphic function~$\tau$ defined in a neighbourhood of the origin.
In particular, $\tau(0)$ is the \emph{cusp shape} associated to the complete
hyperbolic structure at $u=0$ and for this reason $\Imag\tau(0)\neq0$.
Therefore, there exists an $\varepsilon>0$ such that for all
$t\in(-\varepsilon,\varepsilon)\subset\RR$ the parameter
\begin{equation}\label{u=u(t)}
    u(t) = \frac{2t}{\Imag \tau(0)}
\end{equation}
remains in the neighbourhood on which $u$ is a coordinate and where $\tau$ is defined.
Let $z_t$ be the shape coordinates of a point in $\gluvar$ corresponding to $u(t)$.
In particular, $z_t$ is a nontrivial deformation of the initial point $z_0$.
For $\alpha_t = \Arg z_t$, we may compute the angle-holonomies along $\mu$ and $\lambda$ as
\begin{align*}
    \ahol{\alpha_t}(\mu)&=\Imag u(t) \equiv 0,\\
    \ahol{\alpha_t}(\lambda)&=\Imag v(t)
    = \Imag\bigl(u(t)\tau(u(t))\bigr)
    = u(t) \Imag\tau(u(t)),
\end{align*}
because $u(t)\in\RR$ by \eqref{u=u(t)}.
Since $\tau$ is an even function, we have in particular
\begin{equation}\label{two-as-desired}
    \tfrac{d}{dt}\bigr|_{t=0} \ahol{\alpha_t}(\lambda)
    = u'(0) \Imag\tau(u(0)) 
    = 2.
\end{equation}
We define the vector $\hat{l}_\mu\in\TAS(\triang)$ as the infinitesimal generator
of the family of angles $\alpha_t=\Arg z_t$.
In this way, $\hat{l}_\mu$ satisfies \eqref{tangent-to-gluvar} automatically.
In order to prove \eqref{difference-in-TAS0}, we
consider the vector $w=\hat{l}_\mu - l_\mu\in\TAS(\triang)$
and analyse its expression in the basis of leading-trailing deformations,
\[
    w = c_\mu l_\mu + c_\lambda l_\lambda + \sum_{m=1}^{N-1}c_m l_m,
    \quad (c_\mu, c_\lambda, c_1,\dotsc,c_{N-1} \in \RR).
\]
By \eqref{futer-gueritaud-peripheral} and \eqref{two-as-desired}, we have
$\nabla_w \ahol{}(\mu) = \nabla_w\ahol{}(\lambda) = 0$,
so that $c_\mu=c_\lambda=0$, i.e., $w\in\TAS_0(\triang)$.

It remains to prove the uniqueness of $\hat{l}_\mu$.
Suppose that $\hat{l}'_\mu$ is another element
satisfying \eqref{difference-in-TAS0} and \eqref{tangent-to-gluvar}
and set $r = \hat{l}_\mu-\hat{l}'_\mu$.
Then we can regard $r$ as a tangent vector to $\gluvar$ at $z_0$ lying in $\TAS_0(\triang)$, so
$\nabla_r\Imag u = \nabla_r\Imag v = 0$. Thus $\nabla_r u \in \RR$ and 
 $0=\nabla_r\Imag v = \Imag \tau(0) \, \nabla_r  u$. Since $\Imag \tau(0) \ne 0$ this implies that  $\nabla_r  u=0$.
Since $r\in T_{z_0}\gluvar$ and $u$ is a local coordinate on  $\gluvar$ near $z_0$, we must have $r=0$. \end{proof}

The quantity $\tau(0)$ occurring in \eqref{u=u(t)} can be interpreted
as the \emph{relative modulus} of the curves $\mu$ and $\lambda$ in
the boundary-parabolic representation~$\rho_z:\pi_1(M)\to\PSLC$ defined
by the shapes $z$ (up to conjugation).
More precisely, we may conjugate $\rho_z$ so that its restriction to
$\pi_1(\partial\overline{M})$ satisfies
\begin{equation}\label{def-relative-modulus}
    \rho_z(\mu)=\pm\begin{pmatrix}
    1 & 1\\
    0 & 1
    \end{pmatrix},\quad\text{and}\quad
    \rho_z(\lambda)=\pm\begin{pmatrix}
    1 & \sigma(\lambda, \mu)\\
    0 & 1
    \end{pmatrix},
\end{equation}
where $\sigma(\lambda, \mu)=\tau(0)$.
When $\rho_z=\rho^\geo$ is the holonomy representation of a complete hyperbolic
structure, the relative modulus~$\sigma(\lambda, \mu)$ also has a geometric interpretation
as the \emph{cusp shape}.
Specifically, $\sigma(\lambda, \mu)\in\CC_{\Imag>0}$ corresponds
to the conformal structure common to all horospherical cross-sections of the cusp under
the usual identification of the Teichm\"uller space of the torus with the upper half-plane.

An elementary way to see this is to interpret equation~\eqref{def-relative-modulus}
by scaling the Euclidean structure on the cusp
torus so that the length~$L(\mu)$ of a geodesic representative of $\mu$ equals $1$.
Then the parallelogram in $\CC$ with vertices~$\{0,1,\sigma(\lambda, \mu),\sigma(\lambda, \mu)+1\}$ is
a fundamental polygon for the Euclidean structure on $\partial\overline{M}$.
In this normalisation, $\Imag\sigma(\lambda,\mu)$ is equal, up to sign, to the Euclidean area~$A$
of the torus. In other words,
\begin{equation}\label{introducing-cusp-area}
    |\Imag\sigma(\lambda,\mu)| = \frac{A}{L(\mu)^2},
\end{equation}
where $L(\mu)$ is the length of any geodesic representative of $\mu$.
Note that the right-hand side does not depend on the chosen scale.

\begin{corollary}\label{derivative-of-u-wrt-lhat}
With the hypotheses and notations of Proposition~\ref{construction-of-lhat},
for any peripheral curve $\mu$, we have
\[
    \abs{\nabla_{\hat{l}_\mu} \log\hol(\mu)} = \frac{2L(\mu)^2}{A},
\]
where the area~$A$ and the geodesic length~$L(\mu)$ are computed in any
Euclidean metric representing the conformality class of $\partial\overline{M}$.
\end{corollary}
\begin{proof}
Let $u=\log\hol(\mu)$ and let $\lambda$ be a peripheral curve such
that $\intersection(\mu,\lambda)=1$.
The one-parameter family $u=u(t)$ defined in $\eqref{u=u(t)}$ has $\hat{l}_\mu$
as its infinitesimal generator, so that
\[
    \abs{\nabla_{\hat{l}_\mu} \log\hol(\mu)}
    = \abs{\frac{d}{dt} \left(\frac{2t}{\Imag \tau(0)}\right)}
    = \frac{2}{\abs{\Imag \sigma(\lambda,\mu)}}
    = \frac{2L(\mu)^2}{A},
\]
where the last equality uses \eqref{introducing-cusp-area}.
\end{proof}

In light of the above theorem and of the \nword{1}{loop} conjecture~\eqref{oneloop-conjecture-eq},
we expect that the quantity $\tau(z)$ is related to the adjoint Reidemeister torsion as follows.
\begin{conjecture}
Let $z^\geo$ be a geometric solution of edge consistency and completeness equations on $\triang$
endowing $M$ with the complete hyperbolic structure of finite volume and let $\rho^\geo$
be the holonomy representation of this hyperbolic structure.
Then
\[
        \torsionNormalised(M,\rho^\geo) = \bigl(\sqrt{2}\;\tau(z^\geo)\bigr)^{-1},
\]
where $\torsionNormalised(M,\rho^\geo)$ is the normalised
adjoint Reidemeister torsion given by equation~\eqref{torsion-normalised}.
\end{conjecture}

Moreover, we expect similar relationships to hold for other irreducible
boundary-parabolic representations with the geometric obstruction class.

%==============================================================================
% Section 7: Mellin-Barnes integrals
%==============================================================================

\section{Mellin--Barnes integrals associated to ideal triangulations}
\label{sec:Mellin-Barnes}

The goal of this section is a more systematic study of the
multivariate Mellin--Barnes integral~$\MB(\triang,\omega,\alpha)$
which first appeared in equation~\eqref{def:MB} and which we define more
carefully below.
Subsequently, we show that the sum of the integrals $\MB(\triang,\omega,\alpha)$
over a suitable finite set of \Ztaut\ angle structures $\omega$ constitutes a topological
invariant of the underlying 3-manifold.

We still assume that $\triang$ is a topological ideal \nword{N}{tetrahedron}
triangulation of a non-compact, connected, orientable 3-manifold $M$ 
whose ideal boundary is a torus.
If $M$ admits a complete hyperbolic structure with a cusp,
we shall denote by $\SASdistinguished_0(\triang)$
the ``geometric component'' of the space $\SAS_0(\triang)$ of \Sval\ angle
structures with trivial peripheral holonomy.
In terms of the map $\Phi_0$ of Theorem~\ref{obstruction-theorem}, we have
$\SASdistinguished_0(\triang)=\Phi_0^{-1}\bigl(\Obs_0(\rho^\geo)\bigr)$,
where $\Obs_0(\rho^\geo)\in H^2(M,\partial M;\FF_2)$ is the obstruction class
to lifting the holonomy representation of the complete hyperbolic structure to a boundary-unipotent
$\SLC$-representation.

Even without the assumption of hyperbolicity, Lemma~6.1 of \cite{neumann1990} shows 
that $\triang$ always admits a generalised (real-valued) angle structure with
vanishing peripheral angle-holonomy, and since the set of all such structures 
is connected, its image under the exponential map $\alpha\mapsto\exp(i\alpha)$ 
provides a general definition of the distinguished component $\SASdistinguished_0(\triang)$,
regardless of whether $\triang$ admits a strict angle structure.

For the time being, we shall assume that $\triang$ is equipped with a strict angle
structure~$\alpha$ of vanishing peripheral angle-holonomy,
but Section~\ref{sec:MB-extension} below explains how to do away with this assumption.

\begin{definition}\label{definition-MB-integral-at-omega}
Suppose that $\omega\in\SASdistinguished_0(\triang)$ is a \Ztaut\ angle structure.
When $\tet$ is a tetrahedron of $\triang$ with quads $\square,\square'\subset\tet$
satisfying $\omega(\square)=\omega(\square')=+1$, we define the
\emph{beta Boltzmann weight} of $\tet$ to be the analytic function
\begin{equation}\label{beta-weights}
    \Betaweight(s_1,\dotsc,s_{N-1})
    = \Beta\left(\frac{\alpha(\square)}{\pi}  - \sum_{m=1}^{N-1} s_m l_m(\square),\;
                 \frac{\alpha(\square')}{\pi} - \sum_{m=1}^{N-1} s_m l_m(\square')\right),
\end{equation}
where $\Beta(z_1,z_2)=\Gamma(z_1)\Gamma(z_2)/\Gamma(z_1+z_2)$ is Euler's beta function
and where the variables $s_m$ are contained in a vertical strip $0\leq\Real s_m\leq\eps_m$
which is narrow enough to ensure that both arguments of the beta function
have positive real parts.
Then, the \emph{Mellin--Barnes integral} associated to the \Ztaut\ structure~$\omega$ is defined
by the formula
\begin{align}\label{MB-rigorous}
    \MB(\triang,\omega,\alpha)
    &=
    \frac{1}{(2\pi i)^{N-1}}\;
     \pvint_{(i\RR)^{N-1}}
         \Bigl(\prod_{\tet\ \text{in}\ \triang}
         \Betaweight(s)\Bigr)\,ds\\
    &= \frac{1}{(2\pi)^{N-1}}\lim_{R\to+\infty}
     \int\limits_{\substack{x\in\RR^{N-1},\\|x|\leq R}}
         \Bigl(\prod_{\tet\ \text{in}\ \triang}\label{MB-explicit-vp}
         \Betaweight(ix)\Bigr)\,dx,
\end{align}
where the symbol ``$\pvint$'' stands for Cauchy's principal value of the improper integral,
in the sense made precise by the second line.
\end{definition}
\begin{remark}
    \begin{enumerate}[(i)]
    \item
    Note that the region 
    of integration in \eqref{MB-rigorous}
    is symmetric with respect to complex conjugation 
    so the integral is real when it exists. 
    Further, the integral does not depend on the chosen orientation of the triangulation $\triang$,
    since changing the orientation merely changes the signs of the leading-trailing deformations.
    \item
    The integral \eqref{MB-rigorous} does not depend on the choice of the edge
    $E_N\in\edges(\triang)$ whose integration variable is omitted.
    Indeed, as stated in Section~\ref{subsec:TAS}, we have $-l_N = l_1+\dotsb+l_{N-1}$,
    so changing the removed edge amounts to a linear substitution in the integral.
    The determinant of this substitution is always equal to $\pm1$.
    \item
    The integral \eqref{MB-explicit-vp} is over a norm ball $\{|x|\leq R\}$.
    Since all norms on $\RR^{N-1}$ are equivalent, the value of the limit does not
    depend on the choice of the norm.
    \end{enumerate}
\end{remark}

The following proposition was first stated as a part of Theorem~1 in \cite{kashaev-luo-vartanov}
by Kashaev, Luo and Vartanov, in a slightly different context.

\begin{proposition}\label{gauge-invariance-of-Mellin-Barnes}
Let $\triang$ be an ideal triangulation of a 3-manifold in which the link of the ideal
vertex is a torus and let $\omega\in\SASdistinguished_0(\triang)$ be an arbitrary \Ztaut\ angle structure.
Assume that $\alpha_0, \alpha_1$ are two peripherally trivial strict angle structures on the triangulation $\triang$.
Then the integral $\MB(\triang,\omega,\alpha_1)$ converges in the sense of Cauchy's principal value
whenever $\MB(\triang,\omega,\alpha_0)$ does and these two principal values are equal.
\end{proposition}
\begin{proof}
The vector $\alpha_1-\alpha_0\in\TAS_0(\triang)$ can be written as a linear combination
$\alpha_1-\alpha_0 = \sum_{m=1}^{N-1} c_m l_m$ with unique
coefficients $c_m\in\RR$, $1\leq m \leq N-1$.
Consider the one-parameter family~$\mathcal{C}_t$ of integration contours given by
\[
    \mathcal{C}_t:\quad s_m \in i\RR + \frac{c_m}{\pi}t,\quad
    t\in[0,1],\quad 1\leq m\leq N-1.
\]
By convexity, the angle
structures~$\alpha_t(\square) = \alpha_0(\square) + t\sum_{m=1}^{N-1} c_m l_m(\square)$
are all strict, so that
$\Real\bigl[\tfrac{\alpha_0(\square)}{\pi}  + \sum_{m=1}^{N-1} s_m l_m(\square)\bigr]>0$
all $t\in[0,1]$ when $s\in\mathcal{C}_t$.
Since the Euler beta function~$\Beta(z_1,z_2)$ is holomorphic for $\Real z_1, \Real z_2>0$,
the sliding contour~$\mathcal{C}_t$ never leaves the domain of holomorphicity of the integrand,
implying that 
\begin{align*}
    (2\pi i)^{N-1} \MB(\triang,\omega,\alpha_0)
    &=
    \pvint_{\mathcal{C}_0}
        \prod_{\tet}\customBetaweight{\alpha_0}(s)\,ds\\
    &=
    \pvint_{\mathcal{C}_1}
        \prod_{\tet}\customBetaweight{\alpha_0}(s)\,ds\\
    &=
    \pvint_{(i\RR)^{N-1}}
        \prod_{\tet}\customBetaweight{\alpha_1}(s)\,ds\\
    &=
    (2\pi i)^{N-1} \MB(\triang,\omega,\alpha_1).
\end{align*}
In other words, the Mellin--Barnes integral~$\MB(\triang,\omega,\alpha)$
does not depend on the choice of the strict angle structure~$\alpha$.
\end{proof}

%--------------------------------------------------------------------------------------------------
\subsection{Extension of the Mellin--Barnes integral}\label{sec:MB-extension}
The Mellin--Barnes integral~$\MB(\triang,\omega,\alpha)$ can be extended
to triangulations which do not admit strict angle structures.
Since this extension is in full analogy with the constructions
of Section~\ref{3D-index-general},
we limit ourselves to a brief outline.

For simplicity, we assume that the peripheral curves $\mu$ and $\lambda$ are arbitrarily fixed
and we use the notations of Section~\ref{3D-index-general}.
Suppose that $(\aholM,\aholL)\in\RR^2$ are such that
the set $\preAsub{(\aholM,\aholL)}(\triang)$ of positive pre-angle structures
with angle-holonomy~$(\aholM,\aholL)$ is nonempty.
In particular, the definitions in the preceding section use $(\aholM,\aholL)=(0,0)$.

Choose a basepoint $\alpha_*\in\preAsub{(\aholM,\aholL)}(\triang)$
and consider the corresponding coordinates $(r,\epsilon)$ on $\preAsub{(\aholM,\aholL)}(\triang)$, 
where $r_j$ are coordinates on $\TAS_0(\triang)$ and $\epsilon_i$ are angle excesses. 
Using any $\alpha=\alpha(r,\epsilon)\in\preAsub{(\aholM,\aholL)}(\triang)$ in the definition
of the beta Boltzmann weights~\eqref{beta-weights},
we obtain the contour integral $\preMB\bigl(\triang,\omega,\alpha(r,\epsilon)\bigr)$.
If the integral converges, then
an argument essentially identical to Proposition~\ref{gauge-invariance-of-Mellin-Barnes}
shows that $\preMB\bigl(\triang,\omega,\alpha(r,\epsilon)\bigr)$ does not depend on $r$ and thus
defines a germ of an analytic function in the variables $\epsilon=(\epsilon_1,\dotsc,\epsilon_{N-1})$ only.
The Mellin--Barnes integral $\MB(\triang,\omega)$ can then be defined as
the analytic continuation of this germ to $\epsilon=0$.
Note that this analytic continuation may have a pole at $\epsilon=0$, in which case
$\MB(\triang,\omega)=\infty$.

%--------------------------------------------------------------------------------------------------
\subsection{The beta invariant}
In order to simplify our discussion, we consider only the case of $(\aholM,\aholL)=(0,0)$
and focus on the distinguished component $\SASdistinguished_0(\triang)$.
\begin{definition}\label{def:total-MB}
Let $\Ztautset=\{\omega\in\SASdistinguished_0(\triang) : \omega\ \text{is \Ztaut}\}$.
We define the \emph{beta invariant of the triangulation~$\triang$} to be the
finite sum
\begin{equation}\label{total-MB}
    \hyperinv(\triang) = \sum_{\omega\in\Ztautset} \MB(\triang,\omega),
\end{equation}
provided that all summands are finite, so that $\hyperinv(\triang)\in\RR$.
In case any of the summands is either given by an integral
divergent in the sense of Cauchy's principal value or is obtained via analytic continuation
to a singular point, we shall say that $\hyperinv(\triang)$ is undefined.
\end{definition}
We remark that our beta invariant $\hyperinv$ can be seen as a precise formulation
of the candidate invariant whose existence was pondered by Kashaev in \cite[p.~2]{kashaev-beta}.

\begin{remark}
More generally, we may construct a beta invariant $\hyperinv_{(\aholM,\aholL)}(\triang)$
given as a function of the peripheral angle-holonomies $(\aholM,\aholL)$.
This function is defined initially for peripheral angle-holonomies realised on the set 
of positive pre-angle structures, but can be subsequently extended by analytic continuation.
In what follows, we focus on $\hyperinv(\triang)=\hyperinv_{(0,0)}(\triang)$,
although many of our results hold for general $(\aholM,\aholL)$.
\end{remark}

In the context of the asymptotics of the meromorphic \ind\ for a cusped hyperbolic manifold $M$,
we expect the boundary-parabolic $\PSLR$-representations of geometric obstruction class
to contribute linear terms with slopes given by Mellin--Barnes integrals.
For this reason, it is not possible to discern these individual slopes from the total
predicted asymptotics of $\I{M}{0}{0}(-1/\kappa)$.
In other words, the overall coefficient of the \nword{\kappa^1}{term} in our conjectural
approximation~\eqref{general-approximation} is given by the beta invariant~$\hyperinv(\triang)$.
Thus, the discussion of Section~\ref{sec:derivation-of-MB} suggests
that the beta invariant is determined by $\I{M}{0}{0}(-1/\kappa)$, by the means of the limit
\begin{equation}\label{MB-determined-by-3Dindex}
    \hyperinv(\triang) = \lim_{\kappa\to\infty}\frac{\I{M}{0}{0}(-1/\kappa)}{|H_1(\hat{M};\FF_2)|\;\kappa}.
\end{equation}
Although we do not have a rigorous asymptotic argument justifying the above equality,
the theorem stated below guarantees that the
beta invariant $\hyperinv(\triang)$ constitutes a topological invariant
of 3-manifolds with torus boundary.
The topological invariance is a consequence of the connectedness of ideal triangulations 
with at least two tetrahedra under 2-3 and 3-2 Pachner moves (see, e.g., \cite{RST}).

\begin{theorem}\label{beta-invariance-theorem}
Suppose that $M$ is an oriented, non-compact 3-manifold whose ideal boundary
is a torus and that $\triang_2$, $\triang_3$ are two ideal triangulations of
$M$ related by a Pachner 2-3 move. 
Then the beta invariant is defined on $\triang_2$ whenever it is defined on $\triang_3$ and we have
$\hyperinv(\triang_2) = \hyperinv(\triang_3)$.
\end{theorem}

\begin{conjecture}\label{MB-convergence-invariance-conjecture}
Whenever $M$ admits an ideal triangulation $\triang$ with a strict angle structure,  
all Mellin--Barnes integrals in the sum \eqref{total-MB} 
are finite and therefore $\hyperinv(M)=\hyperinv(\triang) \in \RR$ is defined.
\end{conjecture}

%--------------------------------------------------------------------------------------------------
\subsection{Behaviour of the beta invariant under Pachner 2-3 moves}
We shall now study the behaviour of the beta invariant
under a Pachner~$2$-$3$ move on an ideal triangulation, with the aim of establishing
Theorem~\ref{beta-invariance-theorem}.
Throughout this section, we assume that the triangulation $\triang_3$ is
the result of a Pachner \twothree\ move on a triangulation $\triang_2$.
Since the beta invariants of $\triang_2$ and $\triang_3$ are defined
as sums over certain finite sets of \Ztaut\ angle structures on these triangulations,
we can establish Theorem~\ref{beta-invariance-theorem} by studying,
for each \Ztaut\ angle structure $\omega_2\in\SASdistinguished_0(\triang_2)$,
the \Ztaut\ angle structures on $\triang_3$ which map to $\omega_2$
under the angled Pachner move map \eqref{def-PachD}.

Typically, the topological invariance of state integrals of Turaev--Viro type
on ideal triangulations is established through the use of ``pentagon identities''
\cite{kashaev-beta,kashaev-andersen-new}.
A crucial part of this reasoning is the Fubini Theorem which allows us to focus only on
the part of the state integral associated to the bipyramid of the Pachner move.
We remark that although the contour integrals considered here may only converge conditionally
and hence may not satisfy the assumptions of Fubini's theorem,
our use of Cauchy's principal value allows us to sidestep this difficulty.
More explicitly, one may restrict the domain of integration to
a large box $\{\|x\|_\infty \leq R\}$ (cf. eq.~\eqref{MB-explicit-vp}),
where Fubini's theorem holds, and then pass to the limit as $R\to\infty$.
Therefore, all integrals considered from now on will be understood in the
sense of Cauchy's principal value.

\begin{figure}%
    \definecolor{darkred}{RGB}{140,0,0}%
    \centering%
    \begin{tikzpicture}
        \node[anchor=south west,inner sep=0] (image) at (0,0)
            {\includegraphics[width=0.2\columnwidth]{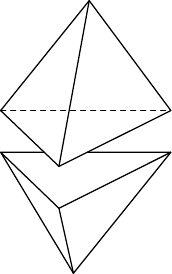}};
        \begin{scope}[x={(image.south east)},y={(image.north west)},
                every node/.style={inner sep=0, outer sep=0}]
            \node[anchor=south east] at (0.7, 0.5)  {\textcolor{darkred}{$-1$}};
            \node[anchor=south west] at (0.17, 0.5)  {$+1$};
            \node[anchor=south] at (0.6, 0.61) {$+1$};
            \node[anchor=south east] at (0.7, 0.35)  {$+1$};
            \node[anchor=south west] at (0.178, 0.32)  {\textcolor{darkred}{$-1$}};
            \node[anchor=south] at (0.85, 0.45) {$+1$};
        \end{scope}
    \end{tikzpicture}%
    \parbox[b][][t]{\dimexpr0.6\columnwidth-2mm\relax}{\caption{%
        \label{fig:twoZ2tautPachners}
        Two possible behaviours of a \Ztaut\ angle structure on a pair of tetrahedra
        sharing a common face.\\
        \textsc{Left:} Non-adjacent case. 
        The tetrahedral edges with \Ztaut\ angles of $-1$ do not
        touch one another along the equator of the bipyramid.\\
        \textsc{Right:} Adjacent case. 
	Along the equator, there are two
        adjacent tetrahedral edges with angles of $-1$.
    }}%
    \begin{tikzpicture}
        \node[anchor=south west,inner sep=0] (image) at (0,0)
            {\includegraphics[width=0.2\columnwidth]{bipyr2tet.pdf}};
        \begin{scope}[x={(image.south east)},y={(image.north west)},
                every node/.style={inner sep=0, outer sep=0}]
            \node[anchor=south east] at (0.7, 0.5)  {\textcolor{darkred}{$-1$}};
            \node[anchor=south west] at (0.17, 0.5)  {$+1$};
            \node[anchor=south] at (0.6, 0.61) {$+1$};
            \node[anchor=south east] at (0.7, 0.35)  {\textcolor{darkred}{$-1$}};
            \node[anchor=south west] at (0.178, 0.32)  {$+1$};
            \node[anchor=south] at (0.85, 0.45) {$+1$};
        \end{scope}
    \end{tikzpicture}%
\end{figure}
Suppose that $\omega_2\in\SASdistinguished_0(\triang_2)$ is a \Ztaut\ angle structure
and let $\PachD$ be the angled Pachner move map of \eqref{def-PachD}.
By Remark~\ref{rem:angled-Pachner-moves}, the pre-image of $\omega_2$ under $\PachD$
is contained in $\SASdistinguished_0(\triang_3)$.
For this reason, in order to show that $\hyperinv(\triang_2)$ equals $\hyperinv(\triang_3)$, 
it suffices to prove that $\MB(\triang_2,\omega_2)=\sum_{\omega\in\Omega}\MB(\triang_3,\omega)$,
where $\Omega$ is the (finite) set of \Ztaut\ angle structures contained
in the pre-image $\PachD^{-1}(\{\omega_2\})$.
It turns out that the cardinality of $\Omega$ is always either 1 or 2.
Thus, we distinguish two cases, depicted in Figure~\ref{fig:twoZ2tautPachners}, depending on
the values of $\omega_2$ around the equator of the bipyramid of the Pachner move.

In what follows, we shall say that the Pachner move is of the \emph{non-adjacent} type
when the bipyramid in the \Ztaut-angled triangulation
$(\triang_2,\omega_2)$ looks like the left panel of Figure~\ref{fig:twoZ2tautPachners}
or its mirror image.
Otherwise, the Pachner move is of \emph{adjacent} type, as shown in the right panel of the figure.

For concreteness, we label the vertices of the bipyramid with the integers $\{0,1,2,3,4\}$,
as shown in Figure~\ref{fig:shaped-Pachner}, and denote by $\sigma_n$ the ideal simplex with
vertices $\{0,\dotsc,4\}\setminus\{n\}$.
When $m,n$ are two distinct ideal vertices of a simplex $\sigma$, we shall denote
by $\sigma([mn])$ the tetrahedral edge of $\sigma$ connecting $m$ and $n$.

\subsubsection{The non-adjacent case}
\begin{figure}%
    \centering%
    \def\minusone{$-1$}%
    \parbox[b][][t]{0.5\columnwidth}%
    {%
        \begin{tikzpicture}
            \node[anchor=south west,inner sep=0] (image) at (0,0)
                {\includegraphics[width=0.5\columnwidth]{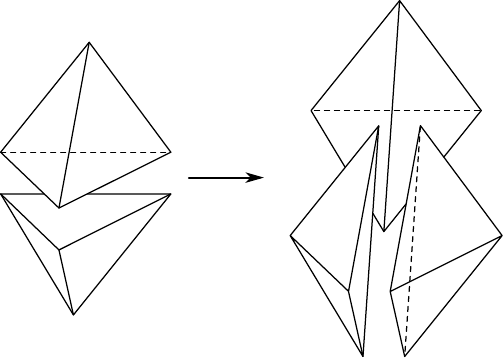}};
            \begin{scope}[x={(image.south east)},y={(image.north west)},
                    every node/.style={inner sep=1pt, outer sep=0}]
                %\draw[help lines,step=0.1] (0,0) grid (1,1);
                \node[anchor=south east] at (0.25, 0.5) {\minusone};
                \node[anchor=south west] at (0.06, 0.35) {\minusone};
                \node[anchor=south east] at (0.1, 0.8) {$\omega_2$};
                \node[anchor=south west] at (0.9, 0.9) {$\omega_3$};
                \node[anchor=south east] at (0.91, 0.28) {\minusone};
                \node[anchor=south] at (0.65, 0.26) {\minusone};
                \node[anchor=south east] at (0.7, 0.83) {\minusone};
            \end{scope}
        \end{tikzpicture}%
    }\hspace{-8mm}
    \parbox[b][][t]{\dimexpr0.5\columnwidth + 5mm\relax}%
    {%
        \caption{
        \label{fig:nonadjacent-Pachner}
        The Pachner \protect\twothree\ move on a \protect\Ztaut\ angle
        structure~$\omega_2$ without adjacent angles of $-1$ induces
        a unique \protect\Ztaut\ angle structure~$\omega_3$ on $\protect\triang_3$.
        For clarity, only one angle of $-1$ has been marked in each tetrahedron.
        }%
    }
\end{figure}
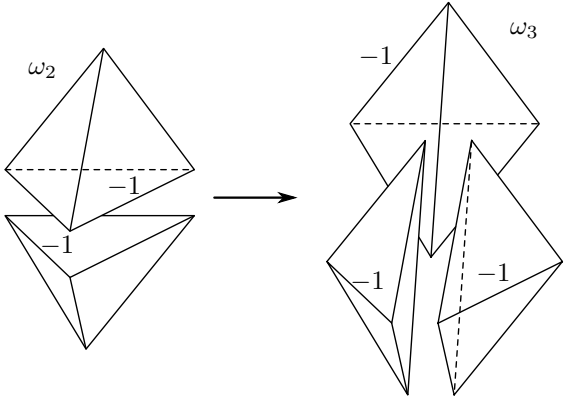

In the non-adjacent case,
the \Ztaut\ angle structure~$\omega_2$ on $\triang_2$ induces a unique
\Ztaut\ angle structure~$\omega_3$ on $\triang_3$.
Explicitly, the rules for finding $\omega_3$ are a multiplicative analogue of
the additive angle relations \eqref{Pachner-angle-formulas}.
For instance, if $\omega_2$ satisfies
$\omega_2\bigl(\sigma_4([23]))\bigr)=\omega_2\bigl(\sigma_0([12]))\bigr)=-1$,
then $\omega_3$ must take the value $-1$ on the tetrahedral edges $\sigma_1([23])$,
$\sigma_2([01])$ as well as $\sigma_3([12])$ and these values 
determine $\omega_3$ completely, as shown in Figure~\ref{fig:nonadjacent-Pachner}.  

As a consequence, it suffices to focus solely on the contributions of the \Ztaut\ angle structures
$\omega_2$ and $\omega_3$ to the
beta invariants of $\triang_2$ and $\triang_3$, respectively.

\begin{lemma}
Let $\omega_3$ be the unique \Ztaut\ angle structure on $\triang_3$ corresponding under
a non-adjacent Pachner \twothree\ move to a \Ztaut\ angle structure $\omega_2$.
\label{noninvmellin}
Suppose that $\alpha_3\in\Astrict_0(\triang_3)$ is a fixed, peripherally trivial strict angle structure
satisfying the inequalities
\begin{alignat}{2}
    \alpha_3\bigl(\sigma_1([03])\bigr) + \alpha_3\bigl(\sigma_2([03])\bigr) &< \pi, &\qquad
    \alpha_3\bigl(\sigma_1([34])\bigr) + \alpha_3\bigl(\sigma_2([34])\bigr) &< \pi,\notag\\
    \alpha_3\bigl(\sigma_2([01])\bigr) + \alpha_3\bigl(\sigma_3([01])\bigr) &< \pi, &\qquad
    \alpha_3\bigl(\sigma_2([14])\bigr) + \alpha_3\bigl(\sigma_3([14])\bigr) &< \pi,\label{angle-strict-ineqs}\\
    \alpha_3\bigl(\sigma_1([02])\bigr) + \alpha_3\bigl(\sigma_3([02])\bigr) &< \pi, &\qquad
    \alpha_3\bigl(\sigma_1([24])\bigr) + \alpha_3\bigl(\sigma_3([24])\bigr) &< \pi,\notag
\end{alignat}
so that 
the above angle sums determine a unique \emph{strict} angle structure~$\alpha_2$ on $\triang_2$.
Then
\[
    \MB(\triang_2, \omega_2, \alpha_2) = \MB(\triang_3, \omega_3, \alpha_3).
\]
%\label{noninvmellin}
\end{lemma}

\begin{proof}
It suffices to consider the \Ztaut\ structure~$\omega_2$
with $\omega_2\bigl(\sigma_4([23]))\bigr)=\omega_2\bigl(\sigma_0([12]))\bigr)=-1$.
We compute the integral $\MB(\triang_3, \omega_3, \alpha_3)$ by integrating over 
variables associated to $N-1$ of the $N$ edges in $\triang_2$ 
together with an extra variable for the new central edge $[04]$ in the bipyramid.
In terms of the edge variables $z_{mn}\in i\RR$ ($0\leq m<n\leq 4$),
the tetrahedral weights associated to the simplices $\sigma_1$, $\sigma_2$, $\sigma_3$
are given by the following formulae:
\begin{align*}
    \Betasigma{\sigma_1}_{\omega_3,\alpha_3} &= \Beta\biggl(
        \frac{\alpha_3\bigl(\sigma_1([02])\bigr)}{\pi}
        + z_{23} + z_{04} - z_{24} - z_{03},\;
        \frac{\alpha_3\bigl(\sigma_1([03])\bigr)}{\pi}
        + z_{02} + z_{34} - z_{04} - z_{23}
    \biggr),\\
    \Betasigma{\sigma_2}_{\omega_3,\alpha_3} &= \Beta\biggl(
        \frac{\alpha_3\bigl(\sigma_2([13])\bigr)}{\pi}
        + z_{01} + z_{34} - z_{03} - z_{14},\;
        \frac{\alpha_3\bigl(\sigma_2([03])\bigr)}{\pi}
        + z_{04} + z_{13} - z_{01} - z_{34}
    \biggr),\\
    \Betasigma{\sigma_3}_{\omega_3,\alpha_3} &= \Beta\biggl(
        \frac{\alpha_3\bigl(\sigma_3([01])\bigr)}{\pi}
        + z_{04} + z_{12} - z_{02} - z_{14},\;
        \frac{\alpha_3\bigl(\sigma_3([02])\bigr)}{\pi}
        + z_{01} + z_{24} - z_{04} - z_{12}
    \biggr).
\end{align*}
We are going to use the substitutions:
\begin{alignat}{2}
    A_1 &= \frac{\alpha_3\bigl(\sigma_1([02])\bigr)}{\pi} + z_{23} - z_{24} - z_{03}, &\qquad
    B_1 &= \frac{\alpha_3\bigl(\sigma_1([03])\bigr)}{\pi} + z_{02} + z_{34} - z_{23}, \notag\\
    A_2 &= \frac{\alpha_3\bigl(\sigma_3([01])\bigr)}{\pi} + z_{12} - z_{02} - z_{14}, &\qquad
    B_2 &= \frac{\alpha_3\bigl(\sigma_3([02])\bigr)}{\pi} + z_{01} + z_{24} - z_{12},
        \label{Beta-pentagon-substitutions} \\
    A_3 &= \frac{\alpha_3\bigl(\sigma_2([03])\bigr)}{\pi} + z_{13} - z_{01} - z_{34}, &\qquad
    s &= z_{04}.\notag
\end{alignat}
Since $\alpha_3$ satisfies the angle equation around the central edge, we have
\[
    A_1 + A_2 + B_1 + B_2
    = \frac{\alpha_3\bigl(\sigma_2([13])\bigr)}{\pi} + z_{01} + z_{34} - z_{03} - z_{14}.
\]
Inserting the new variables, we obtain
\begin{multline}
    \frac{1}{2\pi i} \int_{-i\infty}^{i\infty}
    \Betasigma{\sigma_1}_{\omega_3,\alpha_3}
    \Betasigma{\sigma_2}_{\omega_3,\alpha_3}
    \Betasigma{\sigma_3}_{\omega_3,\alpha_3}dz_{04}\\
    =
    \frac{1}{2\pi i}\int_{-i\infty}^{i\infty}
    \Beta(A_1+s, B_1-s)
    \Beta(A_2+s, B_2-s)
    \Beta(A_3+s, A_1+A_2+B_1+B_2)\,ds.
\end{multline}
We now apply the integral identity
\begin{multline}\label{Barnes-KLV-identity}
    \frac{1}{2\pi i}\int_{-i\infty}^{i\infty}
    \Beta(A_1+s,\;B_1-s)
    \Beta(A_2+s,\;B_2-s)
    \Beta(A_3+s,\;A_1+A_2+B_1+B_2)\,ds \\
    = \Beta(A_2+B_1,\;A_3+B_2)\,\Beta(A_1+B_2,\;A_3+B_1)
\end{multline}
which is stated as equation~(34) by Kashaev, Luo and Vartanov~\cite{kashaev-luo-vartanov};
see also Remark~\ref{remark-Barnes} below for a direct proof.
Undoing the substitutions~\eqref{Beta-pentagon-substitutions} and
expressing the angles of the angle structure~$\alpha_2$ in terms of those of $\alpha_3$,
we transform the right-hand side into
\begin{align*}
    &\phantom{{}={}}\Beta(A_2+B_1,\;A_3+B_2)\,\Beta(A_1+B_2,\;A_3+B_1)\\
    &=\Beta\biggl(
        \frac{\alpha_2\bigl(\sigma_0([24])\bigr)}{\pi}
        + z_{12} + z_{34} - z_{23} - z_{14},\;
        \frac{\alpha_2\bigl(\sigma_0([14])\bigr)}{\pi}
        + z_{24} + z_{13} - z_{12} - z_{34}
    \biggr)\\
    &\phantom{{}={}}\times\Beta\biggl(
        \frac{\alpha_2\bigl(\sigma_4([02])\bigr)}{\pi}
        + z_{01} + z_{23} - z_{03} - z_{12},\;
        \frac{\alpha_2\bigl(\sigma_4([03])\bigr)}{\pi}
        + z_{02} + z_{13} - z_{01} - z_{23}
    \biggr)\\
    &=
    \Betasigma{\sigma_0}_{\omega_2,\alpha_2}
    \Betasigma{\sigma_4}_{\omega_2,\alpha_2}.
\end{align*}
This establishes the required ``charged pentagon identity'' 
\begin{equation}\label{charged-beta-pentagon-nonadjacent}
    \frac{1}{2\pi i} \int_{-i\infty}^{i\infty}
    \Betasigma{\sigma_1}_{\omega_3,\alpha_3}
    \Betasigma{\sigma_2}_{\omega_3,\alpha_3}
    \Betasigma{\sigma_3}_{\omega_3,\alpha_3}dz_{04}
    =
    \Betasigma{\sigma_0}_{\omega_2,\alpha_2}
    \Betasigma{\sigma_4}_{\omega_2,\alpha_2}.
    \qedhere
\end{equation}
\end{proof}

\begin{remark}\label{remark-Barnes}
The integral identity \eqref{Barnes-KLV-identity} is essentially equivalent to
a well-known theorem of Barnes~\cite{barnes1910}.
A modern statement of this celebrated result, known as Barnes' Second Lemma, can be found
in Paris and Kaminski~\cite[Lemma~3.6]{paris-kaminski}.
It implies the equality
\begin{multline}
\frac{1}{2\pi i}\int_{-i\infty}^{i\infty}
\frac{\Gamma(a-\eps+s)\Gamma(b-\eps+s)\Gamma(c-\eps+s)\Gamma(d+\eps-s)\Gamma(\eps-s)}
    {\Gamma(a+b+c+d-\eps+s)}\,ds\\
    =
    \frac{
        \Gamma(a)\Gamma(b)\Gamma(c)\Gamma(a+d)\Gamma(b+d)\Gamma(c+d)
    }{
        \Gamma(b+c+d)\Gamma(a+c+d)\Gamma(a+b+d)
    },
\end{multline}
in which $a,b,c,d$ are all real and positive and $0<\eps<\min(a,b,c)$.
Introducing new parameters $A_1$, $A_2$, $A_3$, $B_1$ and $B_2$
given by
\begin{alignat*}{3}
    a &= A_1 + B_2, &\quad c &= A_3+B_2, & \quad &\\
    b &= A_2 + B_2, &\quad d &= B_1-B_2, &\quad \text{and}\ \eps &= B_2
\end{alignat*}
and observing that
\(
    a+b+c+d-\eps = A_1+A_2+A_3+B_1+B_2
\),
we arrive at the identity
\begin{multline*}
    \frac{1}{2\pi i}\int_{-i\infty}^{i\infty}
    \frac{\Gamma(A_1+s)\Gamma(A_2+s)\Gamma(A_3+s)\Gamma(B_1-s)\Gamma(B_2-s)}
    {\Gamma(A_1+A_2+A_3+B_1+B_2+s)}\,ds\\
    =
    \frac{\Gamma(A_1+B_2) \Gamma(A_2+B_2) \Gamma(A_3+B_2) \Gamma(A_1 + B_1) \Gamma(A_2+B_1) \Gamma(A_3+B_1)}
    {\Gamma(A_2+A_3+B_1+B_2)\Gamma(A_1+A_3+B_1+B_2)\Gamma(A_1+A_2+B_1+B_2)}.
\end{multline*}
Since both sides depend analytically on the parameters, the above identity can be analytically continued to
all values of $A_j$, $B_j$ in the right half-plane~$\{z\in\CC:\Real z>0\}$.
Somewhat miraculously, one may rearrange the terms and organise them
into values of the Euler beta function, as was first observed
in \cite{kashaev-luo-vartanov}.
In this way, we obtain exactly the ``pentagon identity''~\eqref{Barnes-KLV-identity}.
\end{remark}

\subsubsection{The adjacent case}
We now turn to the case of a \Ztaut\ angle structure~$\omega_2$ on the triangulation~$\triang_2$
which assigns the angles of $-1$ to two adjacent tetrahedral edges along the equator of the bipyramid
of the Pachner move.
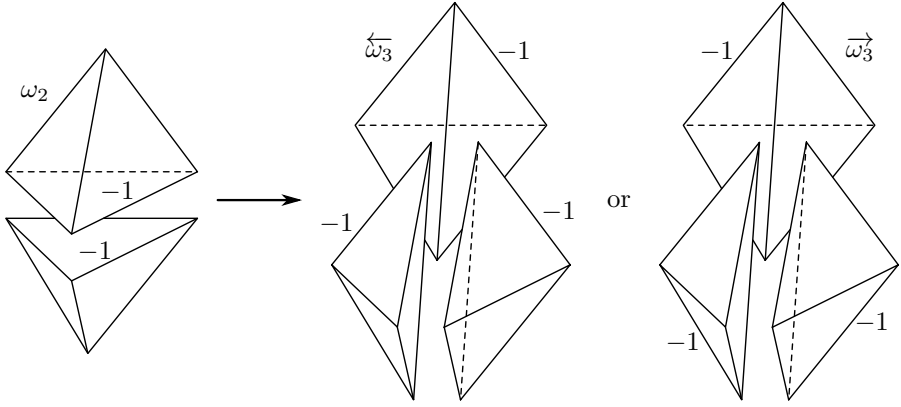
\begin{figure}
    \centering
    \def\minusone{$-1$}
    \begin{tikzpicture}
        \node[anchor=south west,inner sep=0] (image) at (0,0)
            {\includegraphics[width=0.8\columnwidth]{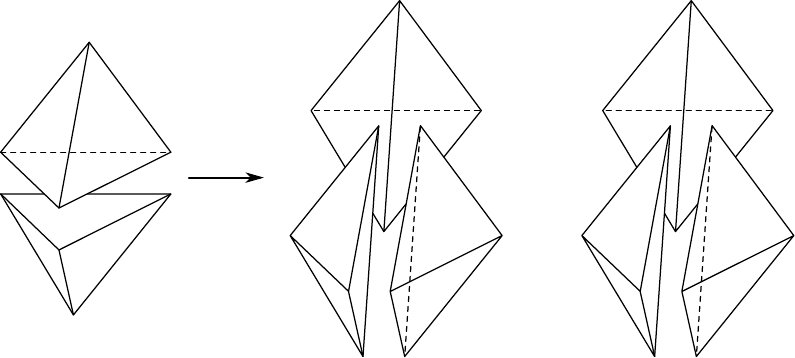}};
        \begin{scope}[x={(image.south east)},y={(image.north west)},
                every node/.style={inner sep=0, outer sep=0}]
            %\draw[help lines,step=0.1] (0,0) grid (1,1);
            \node[anchor=south east] at (0.145, 0.5) {\minusone};
            \node[anchor=south east] at (0.12, 0.35) {\minusone};
            \node[anchor=south east] at (0.05, 0.75) {$\omega_2$};
            \node[anchor=south west] at (0.55, 0.85) {\minusone};
            \node[anchor=south west] at (0.6, 0.45) {\minusone};
            \node[anchor=south east] at (0.39, 0.42) {\minusone};
            \node[anchor=south east] at (0.435, 0.85) {$\overleftarrow{\omega_3}$};
            \node[anchor=south west] at (0.94, 0.85) {$\overrightarrow{\omega_3}$};
            \node[anchor=east] at (0.7, 0.5) {or};
            \node[anchor=south east] at (0.815, 0.85) {\minusone};
            \node[anchor=north west] at (0.95, 0.22) {\minusone};
            \node[anchor=north east] at (0.777, 0.17) {\minusone};
        \end{scope}
    \end{tikzpicture}
    \caption{
    The Pachner \protect\twothree\ move on a \protect\Ztaut\ angle
    structure~$\omega_2$ with two adjacent angles of $-1$ induces
    two distinct \protect\Ztaut\ angle structures~$\protect\overleftarrow{\omega_3}$,
    $\protect\overrightarrow{\omega_3}$ on $\protect\triang_3$.
    For clarity, only one angle of $-1$ has been marked in each tetrahedron.
    }\label{fig:adjacent-Pachner}
\end{figure}
In this case, illustrated in Figure~\ref{fig:adjacent-Pachner}, the 
preimage of the angle structure $\omega_2$ under the map $\PachD$ of \eqref{def-PachD}
contains two different \Ztaut\ angle structures on $\triang_3$,
which we denote by $\overleftarrow{\omega_3}$ and $\overrightarrow{\omega_3}$, since their values
on the bipyramid are mirror images of one another.
In order to ensure that $\hyperinv(\triang_2)=\hyperinv(\triang_3)$ in the adjacent case,
we need to establish the following ``bifurcated pentagon identity'':
\begin{equation}\label{general-double-pentagon}
    \frac{1}{2\pi i} \int_{-i\infty}^{i\infty}
    \Betasigma{\sigma_1}_{\overleftarrow{\omega_3},\alpha_3}
    \Betasigma{\sigma_2}_{\overleftarrow{\omega_3},\alpha_3}
    \Betasigma{\sigma_3}_{\overleftarrow{\omega_3},\alpha_3}
    \,dz_{04}
    +
    \frac{1}{2\pi i} \int_{-i\infty}^{i\infty}
    \Betasigma{\sigma_1}_{\overrightarrow{\omega_3},\alpha_3}
    \Betasigma{\sigma_2}_{\overrightarrow{\omega_3},\alpha_3}
    \Betasigma{\sigma_3}_{\overrightarrow{\omega_3},\alpha_3}
    \,dz_{04}
    =
    \Betasigma{\sigma_0}_{\omega_2,\alpha_2}
    \Betasigma{\sigma_4}_{\omega_2,\alpha_2}
\end{equation}
for all positive angle assignments $\alpha_3$ satisfying
the inequalities \eqref{angle-strict-ineqs}.
In order to simplify our notation, we make the substitutions \eqref{Beta-pentagon-substitutions}
once again, obtaining
\begin{align}\label{adjacent-beta-weights-begin}
    \Betasigma{\sigma_1}_{\overleftarrow{\omega_3},\alpha_3}
    \Betasigma{\sigma_2}_{\overleftarrow{\omega_3},\alpha_3}
    \Betasigma{\sigma_3}_{\overleftarrow{\omega_3},\alpha_3}
    &=
    \Beta(A_1 + s,\, 1 - A_1 - B_1)\,
    \Beta(B_2 - s,\, 1 - A_2 - B_2)\\
    &\pheq\times
    \Beta(B_3 - s,\, 1 - A_3 - B_3),\notag
\\
    \Betasigma{\sigma_1}_{\overrightarrow{\omega_3},\alpha_3}
    \Betasigma{\sigma_2}_{\overrightarrow{\omega_3},\alpha_3}
    \Betasigma{\sigma_3}_{\overrightarrow{\omega_3},\alpha_3}
    &=
    \Beta(B_1 - s,\,1 - A_1 - B_1)\,
    \Beta(A_2 + s,\,1 - A_2 - B_2)\notag\\
    &\pheq\times
    \Beta(A_3 + s,\,1 - A_3 - B_3),\label{adjacent-beta-weights-end}\\
    \text{and}\ 
    \Betasigma{\sigma_0}_{\omega_2,\alpha_2}
    \Betasigma{\sigma_4}_{\omega_2,\alpha_2} &= 
    \Beta(A_1+B_3,A_2+B_1)\, \Beta(A_1+B_2,A_3+B_1),\notag
\end{align}
where $B_3 := 1-A_1-A_2-A_3-B_1-B_2$.

Observe that only one of the arguments of each of the beta functions in
\eqref{adjacent-beta-weights-begin}--\eqref{adjacent-beta-weights-end} depends on the variable~$s$,
so these beta functions reduce to ratios of gamma functions with constant prefactors 
$\Gamma(1-A_1-B_1)$, $\Gamma(1-A_2-B_2)$, $\Gamma(1-A_3-B_3)$.
Therefore, setting
\begin{equation}\label{def-Is}
    \begin{aligned}
        I_1 &= \frac{1}{2\pi i} \int_{-i\infty}^{i\infty}
            \frac{
                \Gamma(A_1 + s)\Gamma(B_2-s)\Gamma(B_3-s)  
            }{
                \Gamma(1-B_1+s)\Gamma(1-A_2-s)\Gamma(1-A_3-s)
            }\,ds,\\
        I_2 &= \frac{1}{2\pi i} \int_{-i\infty}^{i\infty}
            \frac{
                \Gamma(B_1-s)\Gamma(A_2+s)\Gamma(A_3+s)
            }{
                \Gamma(1-A_1-s)\Gamma(1-B_2+s)\Gamma(1-B_3+s)
            }\,ds,
    \end{aligned}
\end{equation}
allows us to rewrite the identity \eqref{general-double-pentagon} in the following equivalent form:
\begin{equation}\label{concrete-bifurcated}
    \Gamma(1-A_1-B_1)\,(I_1 + I_2) =
    \frac{
        \Gamma(A_1+B_2)\Gamma(A_1+B_3)\Gamma(A_2+B_1)\Gamma(A_3+B_1)
    }{
        \Gamma(1-A_2-B_2)\Gamma(1-A_2-B_3)\Gamma(1-A_3-B_2)\Gamma(1-A_3-B_3)
    }.
\end{equation}
The remainder of this section is devoted to the proof of the above equation.
\begin{lemma}\label{lemma-adjacent}
Suppose that the complex parameters $A_j$, $B_j$, $1\leq j\leq3$, have positive real parts and satisfy
\begin{equation}\label{parameter-relation}
    A_1 + A_2 + A_3 + B_1 + B_2 + B_3 = 1.
\end{equation}
Then equation \eqref{concrete-bifurcated} holds.
Consequently, we have
\[
    \MB(\triang_2, \omega_2, \alpha_2)
    = \MB(\triang_3, \overleftarrow{\omega_3}, \alpha_3)
    + \MB(\triang_3, \overrightarrow{\omega_3}, \alpha_3).
\]
\end{lemma}

In order to prove the above lemma, we are going to express the integrals $I_1$ and $I_2$ of \eqref{def-Is}
in terms of hypergeometric series. 
Recall that the hypergeometric function $\hypersymbol{p}{q}$ is defined by
\[
    \pFq{p}{q}{a_1\ a_2\ \dotsc\ a_p}{b_1\ b_2\ \dotsc\ b_q}{z}
    = \frac{
        \prod_{r=1}^q \Gamma(b_r)
    }{
        \prod_{r=1}^p \Gamma(a_r)
    }
    \sum_{n=0}^\infty
    \frac{
        \prod_{r=1}^p \Gamma(a_r+n)
    }{
        \prod_{r=1}^q \Gamma(b_r+n)
    }
    \frac{z^n}{n!},
\]
cf. \cite[Ch.~IV]{erdelyi-bateman}, \cite[Ch.~2]{slater}.
According to the results discussed by Slater in \cite[\S4.6]{slater},
and in particular equation (4.6.2.6) therein, we may express the Mellin--Barnes
integrals $I_1$ and $I_2$ as follows:
\begin{equation}\label{expressions-of-integrals-with-3F2}
    \begin{aligned}
        I_1 &= \frac{
            \Gamma(A_1 + B_2)\Gamma(A_1 + B_3)
            \iiiFii{A_1+B_1}{A_1+B_2}{A_1+B_3}{1+A_1-A_2}{1+A_1-A_3}
        }{
            \Gamma(1-A_1-B_1)\Gamma(1+A_1-A_2)\Gamma(1+A_1-A_3)
        },\\
        I_2 &= \frac{
            \Gamma(A_2 + B_1)\Gamma(A_3 + B_1)  
            \iiiFii{A_1+B_1}{A_2+B_1}{A_3+B_1}{1+B_1-B_2}{1+B_1-B_3}
        }{
            \Gamma(1-A_1-B_1)\Gamma(1+B_1-B_2)\Gamma(1+B_1-B_3)
        }.\\        
    \end{aligned}
\end{equation}
Thanks to the above expressions, we will be able to derive a proof
of Lemma~\ref{lemma-adjacent} from the theory of transformation
identities for the the hypergeometric series $\hypersymbol32(1)$.
We refer to Slater \cite{slater} for more details and an excellent 
historical account of this theory.

\begin{lemma}\label{lemma-simplified-identity}
If $a, b, c, d$ are complex parameters whose real parts lie 
in the interval $(0,1)$ and $a+b+c+d-1$ has positive real part, then 
we have
\begin{equation}\label{simplified-identity}
    \begin{split}
          ab \iiiFii{1-a}{1-b}{1}{1+c}{1+d}
        + cd \iiiFii{1-c}{1-d}{1}{1+a}{1+b}
        \\=
        \frac{
            \Gamma(a+b+c+d-1)\Gamma(a+1)\Gamma(b+1)\Gamma(c+1)\Gamma(d+1)
        }{
            \Gamma(a+c)\Gamma(a+d)\Gamma(b+c)\Gamma(b+d)
        }.
    \end{split}
\end{equation}
\end{lemma}
\begin{proof}
The equation (4.3.4) of Slater \cite{slater} says that
\begin{align*}
    \pheq&
    \iiiFii{a_1}{a_2}{a_3}{b_1}{b_2} 
    \\
    &=
    \frac{
        \Gamma(1-a_1)\Gamma(b_1)\Gamma(b_2)\Gamma(a_3-a_2)
    }{
        \Gamma(b_1-a_2)\Gamma(b_2-a_2)\Gamma(1+a_2-a_1)\Gamma(a_3)
    }
    \iiiFii{a_2}{1+a_2-b_1}{1+a_2-b_2}{1+a_2-a_3}{1+a_2-a_1}\\
    &\pheq+
    \frac{
        \Gamma(1-a_1)\Gamma(b_1)\Gamma(b_2)\Gamma(a_2-a_3)
    }{
        \Gamma(b_2-a_3)\Gamma(b_1-a_3)\Gamma(1+a_3-a_1)\Gamma(a_2)
    }
    \iiiFii{a_3}{1+a_3-b_2}{1+a_3-b_1}{1+a_3-a_2}{1+a_3-a_1}.
\end{align*}
Setting $a_1 = 1-a$, $a_2 = 1-b$, $a_3 = 1$, $b_1 = 1+c$ and
$b_2 = 1+d$, we obtain
\begin{align}
    &\pheq\notag
    \iiiFii{1-a}{1-b}{1}{1+c}{1+d}\\
    &=
    \frac{
        \Gamma(a)\Gamma(b)\Gamma(c+1)\Gamma(d+1)
    }{
        \Gamma(b+d)\Gamma(c+b)\Gamma(1+a-b)
    }\notag
    \iiiFii{1-b}{1-b-c}{1-b-d}{1-b}{1+a-b}\\
    &\pheq+
    \frac{
        \Gamma(a)\Gamma(-b)\Gamma(c+1)\Gamma(d+1)
    }{
        \Gamma(a+1)\Gamma(1-b)\Gamma(c)\Gamma(d)
    }\notag
    \iiiFii{1-c}{1-d}{1}{1+a}{1+b}\\
    &=
    \frac{
        \Gamma(a)\Gamma(b)\Gamma(c+1)\Gamma(d+1)
    }{
        \Gamma(b+c)\Gamma(b+d)\Gamma(1+a-b)
    }\label{apply-Gauss-here}
    \pFq{2}{1}{1-b-c\mkern20mu 1-b-d}{1+a-b}{1}\\
    &\pheq+
    \frac{cd}{a(-b)}\notag
    \iiiFii{1-c}{1-d}{1}{1+a}{1+b}.
\end{align}
Applying the summation theorem of Gauss,
\[
    \pFq{2}{1}{\alpha \mkern20mu \beta}{\gamma}{1} = \frac{\Gamma(\gamma)\Gamma(\gamma-\alpha-\beta)}{\Gamma(\gamma-\alpha)\Gamma(\gamma-\beta)}
\]
(cf. \cite[eq.~2.1.3.(14)]{erdelyi-bateman}; \cite[eq.~(1.1.5)]{slater})
to the line \eqref{apply-Gauss-here}, multiplying by $ab$ and rearranging terms gives equation \eqref{simplified-identity}.
\end{proof}

\begin{remark}
The equality \eqref{simplified-identity} can be extended by analytic continuation
in the parameters $a,b,c,d$, so that it holds whenever both sides are defined.
\end{remark}
\begin{proof}[Proof of Lemma~\ref{lemma-adjacent}]
Suppose that the parameters $A_j$, $B_j$, $1\leq j \leq3$ have positive
real parts and satisfy \eqref{parameter-relation}.
We define
\[
    \left\{
        \begin{aligned}
            a &= A_1 + B_2\\
            b &= A_1 + B_3\\
            c &= A_2 + B_1\\
            d &= A_3 + B_1
        \end{aligned}
    \right.
    \quad
    \text{which implies}
    \quad
    \left\{
        \begin{aligned}
            a + c &= 1 - A_3 - B_3\\
            a + d &= 1 - A_2 - B_3\\
            b + c &= 1 - A_3 - B_2\\
            b + d &= 1 - A_2 - B_2\\
        \end{aligned}
    \right.,
\]
as well as $a+b+c+d-1 = A_1 + B_1$,
and we substitute these values in \eqref{simplified-identity}.
In this way, we arrive at the identity
\begin{align}
    &\pheq
    (A_1 + B_2)(A_1 + B_3)\label{first-hypergeom}
    \iiiFii{1-A_1-B_2}{1-A_1-B_3}{1}{1+A_2+B_1}{1+A_3+B_1}\\
    &\pheq+
    (A_2 + B_1)(A_3 + B_1)\label{second-hypergeom}
    \iiiFii{1-A_2-B_1}{1-A_3-B_1}{1}{1+A_1+B_2}{1+A_1+B_3}\\
    &=
    \frac{
        \Gamma(A_1+B_1)\Gamma(1+A_1+B_2)\Gamma(1+A_1+B_3)\Gamma(1+A_2+B_1)\Gamma(1+A_3+B_1)
    }{
        \Gamma(1-A_2-B_2)\Gamma(1-A_2-B_3)\Gamma(1-A_3-B_2)\Gamma(1-A_3-B_3)
    }.\label{halfbaked-equality}
\end{align}
Observe that each of the above $\hypersymbol32(1)$ series can be rewritten
with the help of the formula
\begin{equation}\label{Saalschutz}
    \iiiFii{b_1-a_1}{b_2-a_1}{1}{a_2+1}{a_3+1}
    =
    \frac{
        \Gamma(a_1)\Gamma(a_2+1)\Gamma(a_3+1)
    }{
        \Gamma(b_1)\Gamma(b_2)
    }
    \iiiFii{a_1}{a_2}{a_3}{b_1}{b_2}, 
\end{equation}
where the parameters satisfy $b_1 + b_2 - a_1 - a_2 - a_3 = 1$ 
(Saalsch\"utzian case). 
The above formula is the result of setting $s=1$ in equation (4.3.1) 
of Slater \cite{slater} and rearranging terms; cf. also \cite[eq.~4.4.(4)]{erdelyi-bateman}.
Explicitly, we set
\begin{align*}
    (a_1, a_2, a_3, b_1, b_2)
    &=
    (A_1+B_1, A_2+B_1, A_3+B_1, 1+B_1-B_2, 1+B_1-B_3),\\
    (a_1, a_2, a_3, b_1, b_2)
    &=
    (A_1+B_1, A_1+B_2, A_1+B_3, 1+A_1-A_2, 1+A_1-A_3)
\end{align*}
when inserting \eqref{Saalschutz} into the lines \eqref{first-hypergeom} 
and \eqref{second-hypergeom}, respectively.
These operations transform the equality \eqref{halfbaked-equality} into a rather
formidable equation, which nonetheless may be simplified with the help
of the functional equation $\Gamma(z+1)=z\Gamma(z)$ and through the subsequent division of both
sides by the common factor 
$(A_1+B_2)(A_1+B_3)(A_2+B_1)(A_3+B_1)\Gamma(A_1+B_1)$.
After these simplifications, we arrive at the identity
\begin{align*}
    &\pheq
    \frac{
        \Gamma(A_2 + B_1)\Gamma(A_3 + B_1)
    }{
        \Gamma(1 + B_1 - B_2)\Gamma(1 + B_1 - B_3)
    }
    \iiiFii{A_1+B_1}{A_2+B_1}{A_3+B_1}{1+B_1-B_2}{1+B_1-B_3}\\
    &\pheq+
    \frac{
        \Gamma(A_1 + B_2)\Gamma(A_1 + B_3)
    }{
        \Gamma(1 + A_1 - A_2)\Gamma(1 + A_1 - A_3)
    }
    \iiiFii{A_1+B_1}{A_1+B_2}{A_1+B_3}{1+A_1-A_2}{1+A_1-A_3}\\
    &=
    \frac{
        \Gamma(A_1 + B_2)\Gamma(A_1 + B_3)\Gamma(A_2 + B_1)\Gamma(A_3 + B_1)
    }{
        \Gamma(1 - A_2 - B_2)\Gamma(1 - A_2 - B_3)\Gamma(1 - A_3 - B_2)\Gamma(1 - A_3 - B_3)
    }.
\end{align*}
With the help of \eqref{expressions-of-integrals-with-3F2}, we may rewrite the above as
\[
    \Gamma(1-A_1-B_1)\,(I_1 + I_2) =
    \frac{
        \Gamma(A_1+B_2)\Gamma(A_1+B_3)\Gamma(A_2+B_1)\Gamma(A_3+B_1)
    }{
        \Gamma(1-A_2-B_2)\Gamma(1-A_2-B_3)\Gamma(1-A_3-B_2)\Gamma(1-A_3-B_3)
    },
\]
thus establishing equation \eqref{concrete-bifurcated}.
Consequently, the bifurcated pentagon identity \eqref{general-double-pentagon}
holds and the proof is finished.
\end{proof}

\begin{proof}[Proof of Theorem~\ref{beta-invariance-theorem}]
Suppose, as before, that $\triang_2$ and $\triang_3$ are triangulations related by a Pachner 2-3 move.
Although we assumed in Lemma~\ref{noninvmellin} that $\alpha_3$ is an angle structure
on $\triang_3$, we only needed to consider the angle equation about the 
central edge~$E_*\in\edges(\triang_3)$ of the bipyramid of the Pachner move.
Hence, Lemma~\ref{noninvmellin} holds whenever $\alpha_3$ is a positive \emph{pre-angle
structure} with angle sum $2\pi$ along $E_*$.
Likewise, Lemma~\ref{lemma-adjacent} holds for all positive pre-angle structures.

The analytic continuation construction described in Section~\ref{sec:MB-extension} defines
the Mellin--Barnes integral associated to a \Ztaut\ angle structure as the analytic
continuation to the point $\epsilon=0$ of the meromorphic germ determined by positive pre-angle structures,
where $\epsilon$ is the vector of their \emph{angle excesses}.
Therefore, the beta invariant $\hyperinv(\triang_3)$ can be viewed as the analytic
continuation to $\epsilon=0$ of a quatity $\hyperinv(\triang_3)(\epsilon)$ defined
as a finite sum of Mellin--Barnes integrals \eqref{total-MB}, where we may
assume that all of the integrals in the sum are expressed in terms of the same pre-angle 
structure $\alpha_3=\alpha_3(r,\epsilon)$ on $\triang_3$.
We may order the edges in such a way that 
$\epsilon=(\epsilon_1,\dotsc,\epsilon_{N-1}, \epsilon_*)$, where $N$ is the number of tetrahedra in $\triang_2$
and where the variable $\epsilon_*$ represents the angle excess along $E_*$.

We can now analytically continue the multivariate germ $\hyperinv(\triang_3)(\epsilon)$
with respect to the angle excess variable $\epsilon_*$ to the value $\epsilon_*=0$.
In this way, we obtain a germ of a meromorphic function 
$\hyperinv(\triang_3)(\epsilon_1,\dotsc,\epsilon_{N-1}, 0)$ in the remaining variables 
$\epsilon_1,\dotsc,\epsilon_{N-1}$.
In order to ensure that the pre-angle structure remains positive, it may be necessary
to adjust some of the remaining angle excesses, which may modify
the domain of definition of the new \nword{(N-1)}{variable} germ.
Crucially, the pre-angle structures obtained in this way satisfy the angle equation along $E_*$, 
so Lemmas \ref{noninvmellin} and  \ref{lemma-adjacent} can be applied directly to 
the integrals whose sum defines $\hyperinv(\triang_3)(\epsilon_1,\dotsc,\epsilon_{N-1}, 0)$.
In this way, we obtain the equality 
$\hyperinv(\triang_3)(\epsilon_1,\dotsc,\epsilon_{N-1}, 0)=\hyperinv(\triang_2)(\epsilon_1,\dotsc,\epsilon_{N-1})$
of germs of meromorphic functions, where the angle excesses on the right-hand side are computed
along the edges common to $\triang_2$ and $\triang_3$, with a consistent numbering.
In particular, the analytic continuations of both germs to the fully balanced case 
of $\epsilon=0$ must coincide.
\end{proof}
%--------------------------------------------------------------------------------------------------
\subsection{Questions and expectations}
We finish this section by stating further conjectures regarding the expected behaviour
of Mellin--Barnes integrals associated to \Ztaut\ angle structures.
To this end, we assume that $M$ is a connected, non-compact, orientable 3-manifold with a
single toroidal end admitting a complete hyperbolic structure with a cusp.
Furthermore, we suppose that $\rho:\pi_1(M)\to\PSLR$ is an irreducible, boundary-parabolic
representation with the geometric obstruction class, i.e., $\Obs_0(\rho)=\Obs_0(\rho^\geo)$.

\begin{conjecture}
There exists a well-defined quantity $\MB(\rho)\in\RR$ which depends only on the
topology of $M$ and the conjugacy class of $\rho$ and moreover satisfies the following property.
Suppose that $\triang$ is an ideal triangulation of $M$ admitting real solutions
$z:Q(\triang)\to\RR\setminus\{0,1\}$ of Thurston's gluing equations for the representation
$\rho$ and consider the set $\Omega_\rho(\triang)\neq\varnothing$ consisting of all \Ztaut\ angle structures
on $\triang$ associated to real solutions for $\rho$.
Then $\MB(\rho) = \sum_{\omega\in\Omega_\rho(\triang)}\MB(\triang,\;\omega)$.
\end{conjecture}

Note that the above conjecture is natural in the context of Theorem~\ref{beta-invariance-theorem},
whose proof relied on the study of individual contributions of \Ztaut\ angle structures
to the beta invariant $\hyperinv(M)$.
When a given representation $\rho:\pi_1(M)\to\PSLR$ has just
one \Ztaut\ angle structure associated to it, then the above conjecture
postulates that the corresponding Mellin--Barnes integral is a well-defined
quantity depending on the conjugacy class of $\rho$.
In particular, it would be interesting to know whether $\MB(\rho)$
is related to other known quantities associated to such representations.

\begin{conjecture}\label{MB-vanishing-conjecture}
Let $\triang$ be an ideal triangulation of $M$ and suppose that
$\omega\in\SASdistinguished_0(\triang)$ is a \Ztaut\ angle structure for
which $\MB(\triang,\omega)\neq0$.
Then there exists a real solution $z:Q(\triang)\to\RR\setminus\{0,1\}$ of Thurston's
edge consistency and completeness equations satisfying $\omega=\sgn z$.
\end{conjecture}

Note that the solution $z$ whose existence is postulated in Conjecture~\ref{MB-vanishing-conjecture}
must necessarily determine a conjugacy class of a boundary-parabolic $\PSLR$-representation
with the geometric obstruction class; cf. Theorem~\ref{obstruction-theorem}.
In other words, the above conjecture predicts that Mellin--Barnes integrals
vanish at all \Ztaut\ angle structures which \emph{do not} correspond to $\PSLR$-representations of
the fundamental group of $M$.
If true, this would point to an interpretation of the beta invariant $\hyperinv(M)$
in terms of conjugacy classes of boundary-parabolic $\PSLR$-representations of $\pi_1(M)$.

We remark that real-valued solutions of Thurston's gluing equations were studied in detail
by Luo in \cite[\S5]{luo-solving}.
In particular, he provided a variational framework for detecting which \Ztaut\ angle structures 
arise as signs of real-valued solutions.
For this reason, there is hope that Luo's methods might help
in the study of Conjecture~\ref{MB-vanishing-conjecture}.

%==============================================================================
% Section 9: The Asymptotic Conjecture for the meromorphic 3D-index
%==============================================================================

\section{Asymptotic Conjecture for the meromorphic 3D-index}
\label{sec:numerical}

In this section, we summarise our asymptotic analysis by formulating a precise
conjecture describing an $\hbar\to0^-$ asymptotic approximation of the
meromorphic \ind~$\I{M}{\aholM}{\aholL}(\hbar)$.
We will mostly focus on the peripherally trivial case $(\aholM,\aholL)=(0,0)$.
In this case, our proposed approximation can be written as a sum over
conjugacy classes of irreducible, boundary-parabolic $\PSLC$ representations
of $\pi_1(M)$ with the geometric obstruction class.
Subsequently, we present some numerical evidence for this conjecture.

As before, we assume that $M$ is a connected, orientable non-compact manifold
admitting a complete hyperbolic structure of finite volume with a single cusp.
We denote by $\rho^\geo:\pi_1(M)\to\PSLC$ a holonomy representation of this structure,
which is defined up to conjugation only.

Let $\mathcal{X}^\geo$ be the set of conjugacy classes of
irreducible, boundary-parabolic representations
$\rho:\pi_1(M)\to\PSLC$ satisfying $\Obs_0(\rho)=\Obs_0(\rho^\geo)$.
In general, it is possible for the set $\mathcal{X}^\geo$ to be infinite, 
as in the example of Section~\ref{example-infinite-family}.
However, we shall assume, for the rest of this section,
that \emph{$\mathcal{X}^\geo$ is finite for the given $M$}.

Observe that complex conjugation leaves $\mathcal{X}^\geo$ invariant
and fixes the set~$\mathcal{X}^\geo_\RR\subset\mathcal{X}^\geo$
of representations whose images can be conjugated to lie in $\PSLR$.
Recall that equation \eqref{contribution-of-smooth-crit} describes the asymptotic
contribution of a conjugacy class $[\rho]\in\mathcal{X}^\geo\setminus\mathcal{X}^\geo_\RR$
in terms of an associated algebraic solution~$z\in(\CC\setminus\RR)^{Q(\triang)}$
of edge consistency and completeness equations on the triangulation~$\triang$.
In this case, $\overline{z}$ is another algebraic solution and it corresponds to
the ``mirror image'' representation $\overline{\rho}$.
Since $\tau(z)=\tau({\overline{z}})$,
we may combine the contributions of the pair $\{[\rho], [\overline{\rho}]\}$ of complex conjugate
representations into a single asymptotic term:
\begin{equation}\label{combined-contribution}
    \contribnonreal([\rho],\kappa)+\contribnonreal([\overline{\rho}],\kappa)
    =
    2\,|H_1(\hat{M};\FF_2)|\,\tau(z)\,\sqrt{2\pi\kappa}
    \cos\left(\kappa\Vol(\rho) + \tfrac{\pi}{4}(N_+ - N_- + \Sigma(\omega))\right)\!,
\end{equation}
where the volume $\Vol(\rho)$ can be computed as the critical value $\Vol(\omega)$
of the volume function on $\SAS_0(\triang)$ and
where $\Sigma(\omega)$ is the signature of its Hessian at the smooth critical point $\omega=z/|z|$.
Moreover $N_+$ and $N_-$ are the numbers of positively and negatively oriented tetrahedra,
respectively, under the shape assignment $z$.

As stated in Question~\ref{mysterious-integer}, we expect the integer
$n_\rho := N_+ - N_- + \Sigma(\omega)$ to be fully determined by the conjugacy class of $\rho=\rho_z$.
In general, only the congruence class of $n_\rho$ modulo 8 affects the expression \eqref{combined-contribution}.
Note that $n_{\overline{\rho}} = -n_\rho$ and $\Vol(\overline{\rho})=-\Vol(\rho)$, so that
the right-hand side of \eqref{combined-contribution} does not depend on which
representation we choose from the complex conjugate pair $\{\rho,\overline{\rho}\}$.
Furthermore, we always have $n_{\rho^\geo} = 1$, and in fact all representations
we studied had $n_\rho=\sgn\Vol(\rho)$, although we do not know if this is true for all 1-cusped manifolds.

\begin{conjecture}[Main Conjecture]\label{VolConj-parabolic}
When $-1/\hbar = \kappa\to +\infty$, we have
\begin{align*}
    \I{M}{0}{0}(\hbar)
    &=
    |H_1(\hat{M}; \FF_2)| \Biggl(
    \hyperinv(M)\,\kappa
    \\&\phantom{{}={}} + \sqrt{8\pi}\mkern-15mu
    \sum_{\{[\rho],[\overline{\rho}]\}\subset\mathcal{X}^\geo\setminus\mathcal{X}^\geo_\RR}
        \mkern-15mu
        \tau(z)\,\sqrt{\kappa}
    \cos\Bigl(\kappa\Vol(\rho) + \tfrac{\pi}{4}\bigl(N_+ - N_- + \Sigma(\omega)\bigr)\Bigr)\Biggr)
    + o(1),
\end{align*}
where:
\begin{itemize}
\item $\hyperinv(M)$ is the beta invariant of $M$;
\item $\hat{M}$ is the end compactification of $M$;
\item $\tau(z)$ is defined by \eqref{tau-recalled};
\item $\Vol(\rho)=-\Vol(\overline\rho)$ is the algebraic volume of the
    representation $\rho$, cf. \cite{francaviglia,zickert-volrep}.
\end{itemize}
\end{conjecture}

\begin{remark}
\begin{enumerate}
    \item
    Since we expect the beta invariant $\hyperinv(M)$ to be a sum of contributions
    arising from the real locus $\mathcal{X}^\geo_\RR$, the overall asymptotic
    approximation postulated by the Conjecture
    can be viewed as a sum of contributions from \emph{all} elements
    of $\mathcal{X}^\geo$, as explained in the Introduction.
    More generally, for a manifold with $k$ cusps we expect the term coming from
    real representations to have degree $k$ in the variable $\kappa$.
    \item
    For the holonomy representation $\rho^\geo$ of the complete hyperbolic structure
    of finite volume, we expect
    \[
        \tau(z^\geo) = \frac{1}{\sqrt{2}\;\torsionNormalised(M,\rho^\geo)},
    \]
    where $\torsionNormalised(M,\rho^\geo)$ is the normalised adjoint Reidemeister torsion
    given in \eqref{torsion-normalised}.
    This expectation is equivalent to the 1-loop conjecture \cite{tudor-stavros,siejakowski-infinitesimal}.
    \item
    The integer prefactor $|H_1(\hat{M};\FF_2)|$, which is always a power of two,
    can also be expressed as the cardinality of the kernel of the natural
    map $H^2(M;\FF_2)\to H^2(\partial M;\FF_2)$ induced by the
    inclusion of the boundary in $M$.
    Hence, this integer also coincides with the number of connected components
    of $\SAS(\triang)$.
\end{enumerate}
\end{remark}

The reason for the restriction
to representations sharing their obstruction class with $\rho^\geo$ lies in
the interpretation in Section~\ref{state-integral-to-circle-valued} of the meromorphic \ind\ as a state-integral over a single connected component of the space $\SAS_0(\triang)$, and the identification via Theorem~\ref{obstruction-theorem} of that connected component as $\Phi_0^{-1}\bigl(\Obs_0(\rho^\geo)\bigr)\subset\SAS_0(\triang)$.

This motivates the following question.
\begin{question}
Can the state-integral of the meromorphic \ind\ be extended to
an integral over \emph{all} connected components of the space of \Sval\ angle
structures on suitable ideal triangulations in such a way that its
analytic continuation provides a topological invariant of 3-manifolds with
toroidal boundary?
\end{question}

We finish by briefly explaining how to extend Conjecture~\ref{VolConj-parabolic}
beyond the boundary-parabolic case of $(\aholM, \aholL)=(0,0)$.
Suppose that $(\aholM, \aholL)\in\RR^2$ satisfy $|\aholM|, |\aholL| < \frac{\pi}{2}$.
We assume that the complete hyperbolic structure, with holonomy $\rho^\geo$, can be deformed
into an incomplete structure whose holonomy representation $\rho^\geo(\aholM,\aholL)$
has \Sval\ peripheral angle-holonomies~$\bigl(\exp(i\aholM), \exp(i\aholL)\bigr)$.
Then we may study the obstruction class $\nu\in H^2(M,\del M;\FF_2)$ to lifting $\rho^\geo(\aholM,\aholL)$
to an $\SLC$ representation in which the images of the peripheral elements
have eigenvalues with positive real parts.
We expect that the asymptotic contributions to $\I{M}{\aholM}{\aholL}(\hbar)$
come from the set of conjugacy classes of irreducible $\PSLC$ representations of $\pi_1(M)$
with \Sval\ peripheral angle-holonomies~$\bigl(\exp(i\aholM), \exp(i\aholL)\bigr)$
which have $\nu$ as their obstruction class to lifting to $\SLC$ representations with positive
real parts of the peripheral eigenvalues.
Moreover, we conjecture that
each such contribution has the form given by \eqref{contribution-of-smooth-crit}
as $-1/\hbar = \kappa\to+\infty$.

\subsection{Numerical results}\label{sec:numerix}
We shall now compare the asymptotic approximations
of Conjecture~\ref{VolConj-parabolic} with values of the meromorphic
\ind\ computed by numerical quadrature of the state integral~\eqref{3Dindex-int-2}.
As the explicit state integral definition of $\I{M}{0}{0}(\hbar)$
requires using a triangulation~$\triang$ of $M$ admitting a strict angle structure,
we shall restrict our attention to the minimal geometric triangulations
from the SnapPea~\cite{snappy} orientable cusped census.

\subsubsection{The figure-eight knot complement}
\begin{figure}
    \centering
    \begin{tikzpicture}
        \node[anchor=south west,inner sep=0] (image) at (0,0)
            {\includegraphics[width=.8\columnwidth]{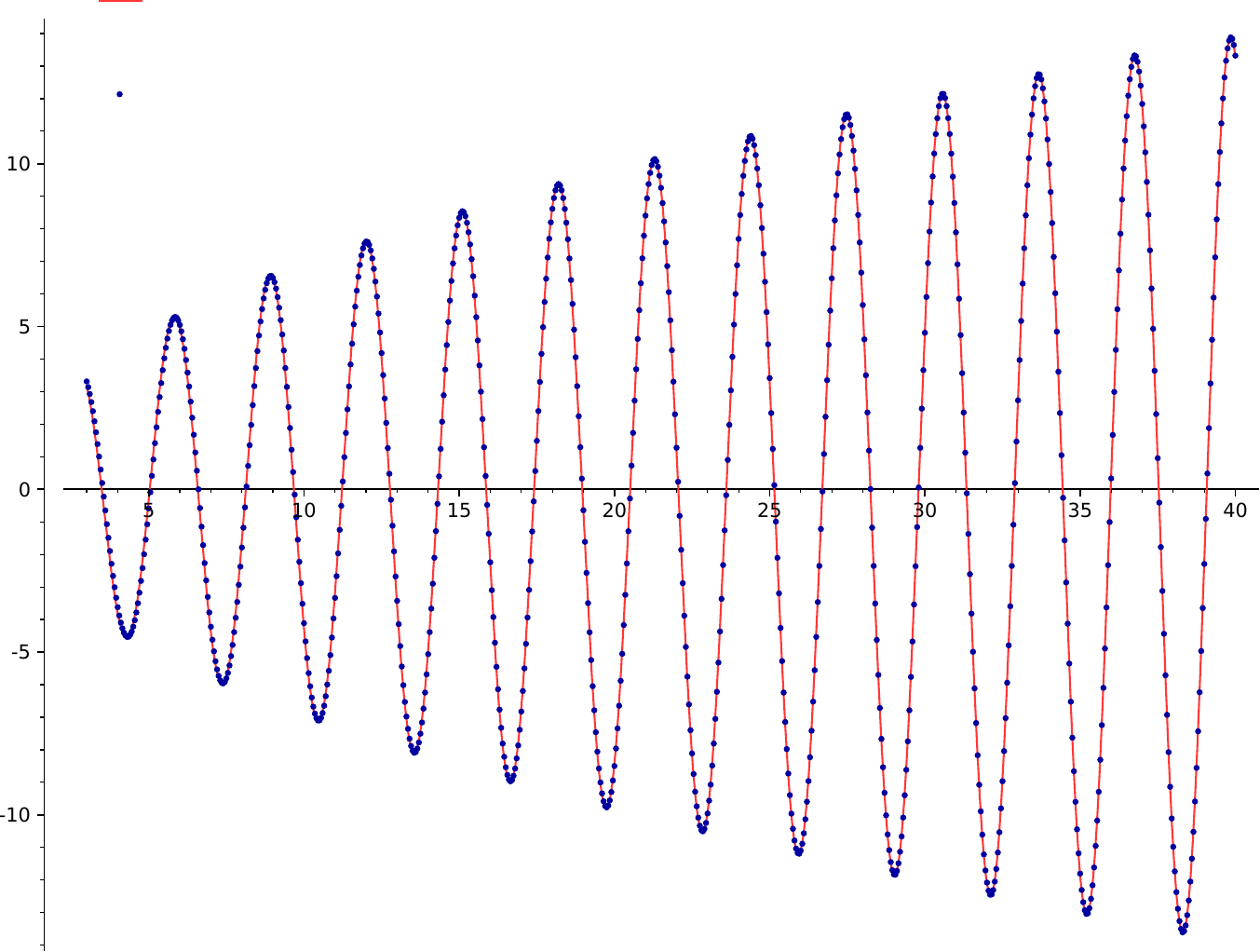}};
        \begin{scope}[x={(image.south east)},y={(image.north west)},
                every node/.style={inner sep=0, outer sep=0}]
            \node[anchor=west] at (1.01, 0.48)  {$\kappa$};
            \node[anchor=west] at (0.13, 1) {$\frac{2\sqrt{2\pi}}{\sqrt[4]{27}}\;
                \sqrt{\kappa}\cos\bigl(\kappa\Volhyp(M)+\tfrac{\pi}{4}\bigr)$};
            \node[anchor=west] at (0.13, 0.9) {$\I{M}{0}{0}(-1/\kappa)$ found numerically};
        \end{scope}
    \end{tikzpicture}
    \caption{%
    \label{plot:asym-m004-00}
    A comparison of the predicted asymptotic expansion \protect\eqref{postulated-approx-4_1}
    with the results of a numerical quadrature of the state integral
    for the figure-eight knot complement~$M$, at $(\aholM,\aholL)=(0,0)$.
    }
\end{figure}

Let $M=S^3\setminus4_1$ be the complement of the figure-eight knot in the \nword{3}{sphere}.
As discussed in detail in Section~\ref{sec:example-4_1},
the asymptotic approximation predicted by Conjecture~\ref{VolConj-parabolic} is
\begin{equation}\label{postulated-approx-4_1}
\I{M}{0}{0}(-1/\kappa)=
    \frac{2\sqrt{2\pi}}{\sqrt[4]{27}}\;\sqrt{\kappa}\cos\bigl(\kappa\Volhyp(M)+\tfrac{\pi}{4}\bigr)
    + o(1)\quad
    \text{as}\ \kappa\to+\infty.
\end{equation}
We computed the values of $\I{M}{0}{0}(-1/\kappa)$ numerically for
$\kappa\in\{3,3.1,3.2,\dotsc,39.9,40\}$ and
observed a very good agreement with the approximation~\eqref{postulated-approx-4_1},
as illustrated in Figure~\ref{plot:asym-m004-00}.

\subsubsection{\texorpdfstring{The complement of the knot $5_2$}{The complement
of the knot 5\char`_2}}
Let $M=S^3\setminus5_2$ be the complement of the knot~$5_2$
and let $\triang$ be the triangulation of $M$ with Regina~\cite{regina}
isomorphism signature~\texttt{dLQbcccdero}.
In one possible ordering of tetrahedra and edges, the gluing matrix of $\triang$ is given by
\begin{equation}\label{gluing-matrix-5_2}
    \begin{bmatrix}\mathbf{G}\\\mathbf{G_\del}\end{bmatrix} =
    \begin{bmatrix}G & G' & G''\\G_\del & G'_\del & G''_\del\end{bmatrix} =
        \left[\begin{array}{rrr|rrr|rrr}
        0 &  0  &  0 & 1 &  1 &  0 & 1 &  1  & 1 \\
        1 &  1  &  0 & 0 &  0 &  2 & 1 &  1  & 0 \\
        1 &  1  &  2 & 1 &  1 &  0 & 0 &  0  & 1 \\\hline
        1 &  -3 &  0 & 0 &  0 &  0 & 0 &  0  & 1 \\
        1 &  0  &  1 & 0 &  0 &  0 & 0 &  -1 & 0
\end{array}\right],
\end{equation}
whence we may take
\[
    \mathbf{L_*} = \left[\begin{array}{rrr|rrr|rrr}
        0 & 0 & 1  & -1 & -1 & -1 &  1 & 1  & 0\\
        1 & 1 & -2 & 0  &  0 & 0  & -1 & -1 & 2
        \end{array}\right].
\]
The complete hyperbolic structure on $M$ has volume~$\Volhyp(M)\approx2.8281220883$
and is recovered at the angle structure $\alpha^\geo$ with angles
\[
    (\alpha, \alpha', \alpha'') \approx (0.7038577213, 1.0300194897, 1.4077154426)
\]
in each of the three tetrahedra.
Using \eqref{tau-recalled}, we compute the \nword{\tau}{invariant} at $\alpha^\geo$,
which equals approximately $0.3527963427$.

We find three \Ztaut\ angle structures on $\triang$, which we write down below, in
the same ordering of the normal quads as that used in \eqref{gluing-matrix-5_2}:
\begin{align*}
    \omega_1 &= ( 1, 1, 1\ |\ -1,-1,-1\ |\  1, 1, 1),\\
    \omega_2 &= ( 1, 1,-1\ |\  1, 1, 1\ |\ -1,-1, 1),\\
    \omega_3 &= (-1,-1,-1\ |\  1, 1, 1\ |\  1, 1, 1).
\end{align*}
Using the strict angle structure
\(
    \alpha_*=\bigl(\tfrac{\pi}{6},\tfrac{\pi}{6},\tfrac{\pi}{6};
    \tfrac{\pi}{2},\tfrac{\pi}{2},\tfrac{\pi}{2};
    \tfrac{\pi}{3},\tfrac{\pi}{3},\tfrac{\pi}{3}\bigr)
\),
we write down the corresponding Mellin--Barnes integrals,
\begin{align*}
    \MB(\triang,\omega_1,\alpha_*) &= \frac{1}{(2\pi i)^2}
    \int_{(i\RR)^2}\Beta^2\bigl(\tfrac16 + s_2,\,\tfrac13+s_1-s_2\bigr)
        \Beta\bigl(\tfrac16+s_1-2s_2,\,\tfrac13+2s_2\bigr)\,ds_1 ds_2,\\
    \MB(\triang,\omega_2,\alpha_*) &= \frac{1}{(2\pi i)^2}
    \int_{(i\RR)^2}\Beta^2\bigl(\tfrac16 + s_2,\,\tfrac12-s_1\bigr)
        \Beta\bigl(\tfrac12-s_1,\,\tfrac13+2s_2\bigr)\,ds_1 ds_2,\\
    \MB(\triang,\omega_3,\alpha_*) &= \frac{1}{(2\pi i)^2}
    \int_{(i\RR)^2}\Beta^2\bigl(\tfrac12 - s_1,\,\tfrac13+s_1-s_2\bigr)
        \Beta\bigl(\tfrac12-s_1,\,\tfrac13+2s_2\bigr)\,ds_1 ds_2.
\end{align*}
Numerical evaluation shows that
$\MB(\triang,\omega_1)\approx\MB(\triang,\omega_2)\approx0$,
whereas $\MB(\triang,\omega_3)\approx0.534$.
This is in accordance with Conjecture~\ref{MB-vanishing-conjecture},
since $\omega_3$ is in fact the only \Ztaut\ angle structure which corresponds to
a \realrep\ representation.
This representation is listed in the Ptolemy database \cite{ptolemy} and
has a Chern-Simons invariant of approximately $-0.53148$.

Since $H_1(\hat{M};\FF_2)=0$, the expected asymptotic approximation of $\I{M}{0}{0}(-1/\kappa)$ is
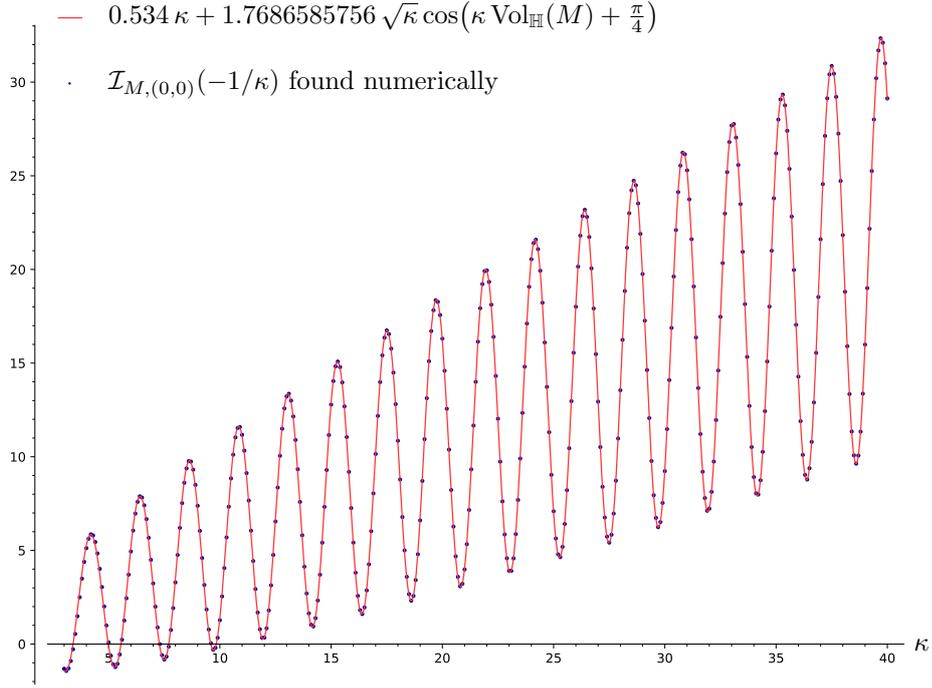
\begin{figure}
    \centering
    \begin{tikzpicture}
        \node[anchor=south west,inner sep=0] (image) at (0,0)
            {\includegraphics[width=.8\columnwidth]{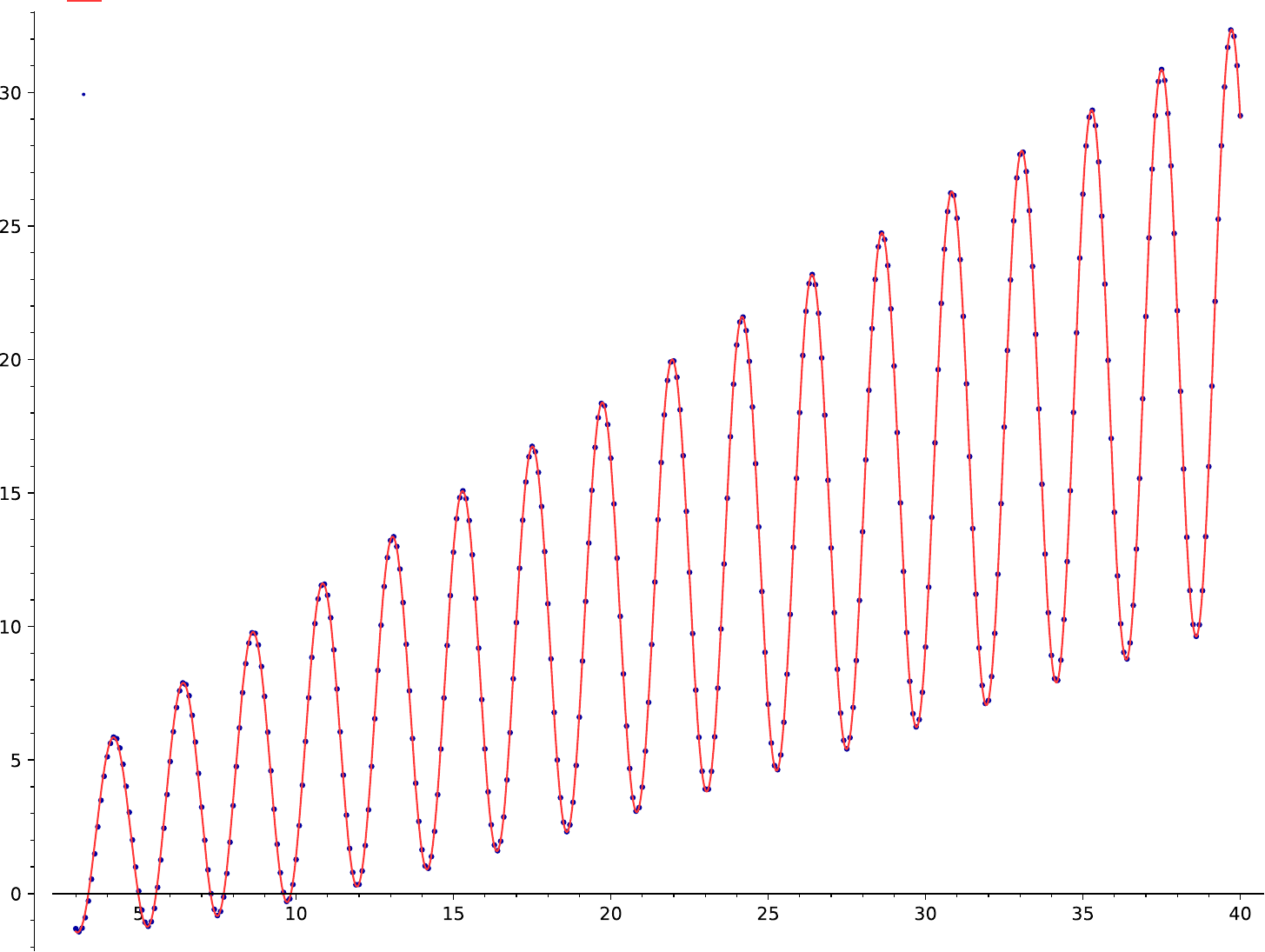}};
        \begin{scope}[x={(image.south east)},y={(image.north west)},
                every node/.style={inner sep=0, outer sep=0}]
            \node[anchor=west] at (1.01, 0.06)  {$\kappa$};
            \node[anchor=west] at (0.11, 1) {$0.534\,\kappa +
                1.7686585756\,\sqrt{\kappa}\cos\bigl(\kappa\Volhyp(M) + \tfrac{\pi}{4}\bigr)$};
            \node[anchor=west] at (0.11, 0.9) {$\I{M}{0}{0}(-1/\kappa)$ found numerically};
        \end{scope}
    \end{tikzpicture}
    \caption{%
    \label{plot:asym-m015-00}
    A comparison of the predicted asymptotic expansion \eqref{predicted:5_2}
    with the results of a numerical quadrature of the state integral
    for the $5_2$ knot complement~$M$, at $(\aholM,\aholL)=(0,0)$.
    }
\end{figure}
\begin{align}
    \I{M}{0}{0}(-1/\kappa) &=
    \hyperinv(M)\kappa + 2 \tau(z^\geo) \sqrt{2\pi\kappa} \cos\bigl(\kappa\Volhyp(M) + \tfrac{\pi}{4}\bigr)
    + o(1)\label{predicted:5_2}\\
    &\approx 0.534\,\kappa + 1.7686585756 \sqrt{\kappa}\cos\bigl(\kappa\Volhyp(M) + \tfrac{\pi}{4}\bigr)
    \ \text{as}\ \kappa\to+\infty.\notag
\end{align}
The above approximation and the values of the meromorphic \ind\ found numerically are plotted
together in Figure~\ref{plot:asym-m015-00}.

\subsubsection{The complement of the Stevedore knot}
Let $M=S^3\setminus6_1$ be the complement of the Stevedore knot~$6_1$.
The manifold $M$ has a geometric triangulation $\triang$
with the isomorphism signature \texttt{eLPkbcddddcwjb}, which
appears as \texttt{m032} in the oriented cusped census \cite{snappy}.
According to the entry in the Ptolemy database \cite{ptolemy}, there are two complex conjugate pairs
of irreducible boundary-parabolic \nword{\PSLC}{representations} with the geometric
obstruction class, up to conjugacy:
\begin{itemize}
    \item
    The pair $\{\rho^\geo, \overline{\rho^\geo}\}$ formed by holonomy representations
    of the complete hyperbolic structure,
    with volume $\Vol(\rho^\geo)=\Volhyp(M)\approx 3.16396$.
    \item
    A pair $\{\rho_1, \overline{\rho_1}\}$ of boundary-parabolic representations
    with algebraic volumes approximately equal to $\pm1.4151$.
\end{itemize}
The geometric representation $\rho^\geo$ can be constructed from the positively oriented
solution $z^\geo$ for which we compute the \nword{\tau}{invariant} numerically as
$\tau(z^\geo)\approx0.2445989349486$.
Since $H_1(\hat{M};\FF_2)=0$, the expected asymptotic contribution from the holonomy representation is
\begin{equation}\label{6_1:contrib-geo}
    \begin{aligned}
        \contribnonreal([\rho^\geo],\kappa)
        + \contribnonreal([\overline{\rho^\geo}],\kappa)
        &= \tau(z^\geo)\cdot 2\cos\left(\kappa\Volhyp(M) + \tfrac{\pi}{4}\right) \sqrt{2\pi\kappa}\\
        &\approx 1.22623721257\sqrt{\kappa}\cos\left(\kappa\Volhyp(M) + \tfrac{\pi}{4}\right)\!.
    \end{aligned}
\end{equation}

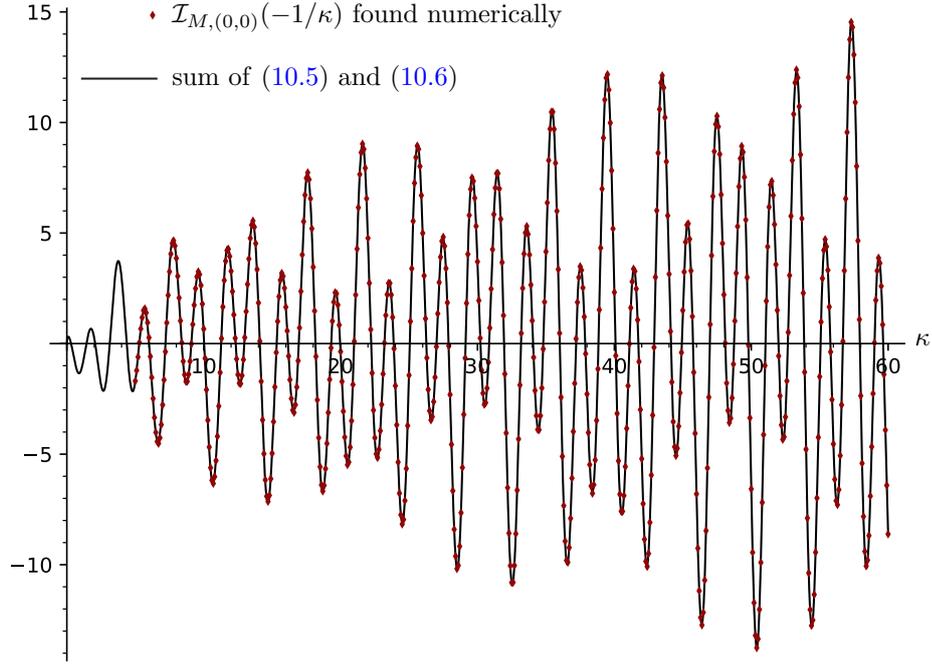
\begin{figure}
    \centering
    \begin{tikzpicture}
        \node[anchor=south west,inner sep=0] (image) at (0,0)
            {\includegraphics[width=.8\columnwidth]{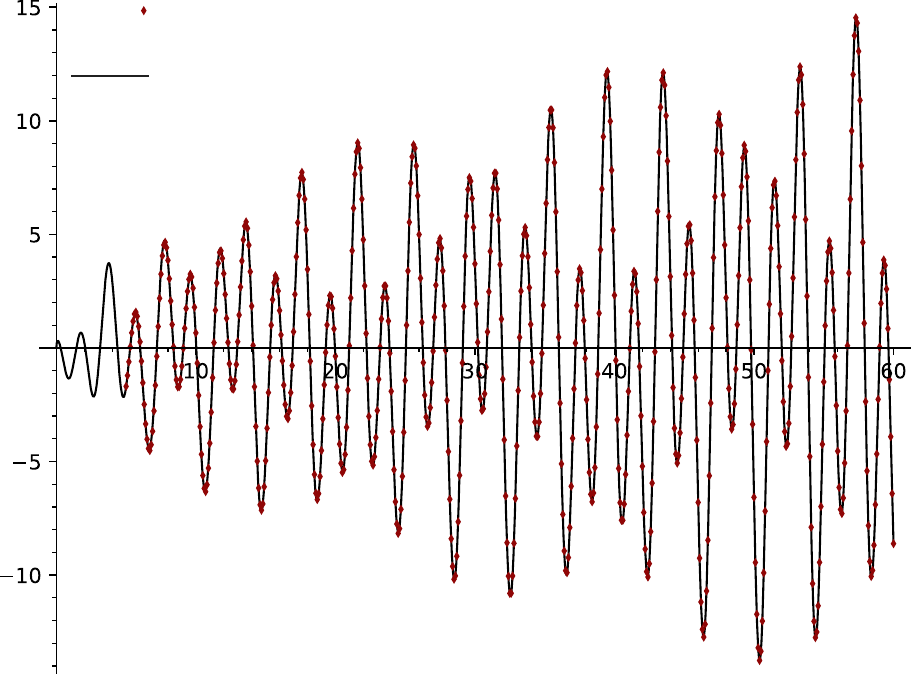}};
        \begin{scope}[x={(image.south east)},y={(image.north west)},
                every node/.style={inner sep=0, outer sep=0}]
            \node[anchor=west] at (1.01, 0.49)  {$\kappa$};
            \node[anchor=west] at (0.18, 0.887) {sum of \eqref{6_1:contrib-geo} and \eqref{6_1:contrib-nongeo}};
            \node[anchor=west] at (0.18, 0.98) {$\I{M}{0}{0}(-1/\kappa)$ found numerically};
        \end{scope}
    \end{tikzpicture}
    \caption{%
    \label{asym:6_1}
    A comparison of the meromorphic \ind\ of the Stevedore knot complement $M=S^3\setminus6_1$
    with the asymptotic approximation of Conjecture~\ref{VolConj-parabolic}.
    }
\end{figure}

The representation $\rho_1$ of volume $\Vol(\rho_1)\approx1.4151$ is described by an algebraic
solution $z_1$ of the gluing equations, which can be found with the help of the \texttt{ptolemy} module
available inside SnapPy \cite{snappy}.
When equipped with the shapes $z_1$, three tetrahedra of $\triang$ are positively oriented and
one tetrahedron is negatively oriented.
The Hessian \eqref{computation-of-Hessian} of the volume function has signature $\Sigma(\omega_1)=-1$ at the
angle structure $\omega_1(\square)=z_1(\square)/|z_1(\square)|$, so that
$N_+ - N_- + \Sigma(\omega_1) = 1$.
Moreover, we have $\tau(z_1) \approx 0.143746443341$,
so we obtain the contributions
\begin{equation}\label{6_1:contrib-nongeo}
    \begin{aligned}
        \contribnonreal([\rho_1],\kappa)
        + \contribnonreal([\overline{\rho_1}],\kappa)
        &= \tau(z_1)\cdot 2\cos\left(\kappa\Vol(\rho_1) + \tfrac{\pi}{4}\right) \sqrt{2\pi\kappa}\\
        &\approx 0.7206377985127918\sqrt{\kappa}\cos\left(\kappa\Vol(\rho_1) + \tfrac{\pi}{4}\right).
    \end{aligned}
\end{equation}
Therefore, the asymptotic approximation predicted by Conjecture~\ref{VolConj-parabolic}
is the sum of the contributions \eqref{6_1:contrib-geo} and \eqref{6_1:contrib-nongeo}.
In Figure~\ref{asym:6_1}, we compare this approximation to the values of the
meromorphic \ind\ found numerically for $\kappa\in\{5, 5.1, 5.2, \dotsc, 60\}$.

\subsubsection{\texorpdfstring{The complement of the knot $7_2$}{The complement
of the knot 7\char`_2}}
As our last example, we consider the manifold $M=S^3\setminus7_2$.
A geometric four-tetrahedron triangulation $\triang$ of $M$ is found in the SnapPea
oriented cusped census under the designation \texttt{m053} and has the
isomorphism signature \texttt{eLAkbccddmejln}.
According to the Ptolemy database \cite{ptolemy}, there are five conjugacy classes of boundary-parabolic
$\PSLC$-representations of $\pi_1(M)$ with the geometric obstruction class:
\begin{itemize}
    \item
    The pair $\{\rho^\geo, \overline{\rho^\geo}\}$ of holonomy representations of the complete structure
    with volume $\Vol(\rho^\geo)=\Volhyp(M)\approx 3.33174423$.
    \item
    Another complex conjugate pair $\{\rho_1, \overline{\rho_1}\}$, with algebraic volumes
    equal approximately $\pm 2.21397$.
    \item
    A \realrep\ representation $\rho_2$ with a Chern-Simons invariant of approximately
    $-0.762751$.
\end{itemize}

As in the previous examples, we can find shape parameter solutions $z^\geo, z_1, z_2$
for all of the above representations using the \texttt{ptolemy} functionality inside SnapPy \cite{snappy}.
From \eqref{tau-recalled}, we then calculate the associated \nword{\tau}{invariants}, finding
\[
    \tau(z^\geo) \approx 0.16264170908787,\quad
    \tau(z_1)    \approx 0.176853150149.
\]
With the algebraic solution $z_1$, the triangulation $\triang$ has three positively oriented tetrahedra
and one negatively oriented tetrahedron.
At the \Sval\ angle structure $\omega_1(\square) = z_1(\square)/|z_1(\square)|$, the Hessian of the volume
has signature $\Sigma(\omega_1) = -1$, so that $N_+ - N_- + \Sigma(\omega_1)=1$.
Since $H_1(\hat{M};\FF_2)=0$, the conjectured contribution of the four \nonrealrep\ representations is
\begin{equation}\label{7_2:nonrealcontrib}
    \begin{aligned}
        &\phantom{{}={}}\contribnonreal([\rho^\geo],\kappa) + \contribnonreal([\overline{\rho^\geo}],\kappa)
        + \contribnonreal([\rho_1],\kappa)  + \contribnonreal([\overline{\rho_1}],\kappa)\\
        &= \tau(z^\geo)\cdot 2\cos\left(\kappa\Volhyp(M)   + \tfrac{\pi}{4}\right) \sqrt{2\pi\kappa}
          +\tau(z_1)   \cdot 2\cos\left(\kappa\Vol(\rho_1) + \tfrac{\pi}{4}\right) \sqrt{2\pi\kappa}\\
        &\approx \sqrt{\kappa}\Bigl(
            0.8153646132679\,\cos\left(\kappa\Volhyp(M)   + \tfrac{\pi}{4}\right)\\
        &\phantom{{}\approx \sqrt{\kappa}\Bigl(}
        + 0.8866102132422\,\cos\left(\kappa\Vol(\rho_1) + \tfrac{\pi}{4}\right)\Bigr).
    \end{aligned}
\end{equation}

\begin{figure}
    \centering
    \begin{tikzpicture}
        \node[anchor=south west,inner sep=0] (image) at (0,0)
            {\includegraphics[width=.8\columnwidth]{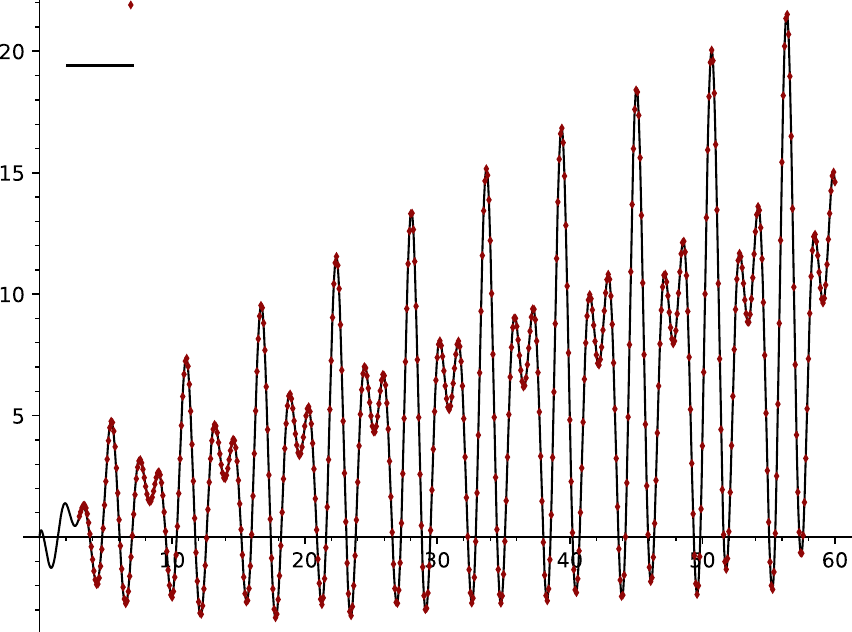}};
        \begin{scope}[x={(image.south east)},y={(image.north west)},
                every node/.style={inner sep=0, outer sep=0}]
            \node[anchor=west] at (0.99, 0.18)  {$\kappa$};
            \node[anchor=west] at (0.17, 0.9) {sum of $\contribreal([\rho_2], \kappa)$ and
            the expression \eqref{7_2:nonrealcontrib}};
            \node[anchor=west] at (0.17, 0.99) {$\I{M}{0}{0}(-1/\kappa)$ found numerically};
        \end{scope}
    \end{tikzpicture}
    \caption{%
    \label{asym:7_2}
    A comparison of the meromorphic \ind\ of $M=S^3\setminus7_2$
    with the asymptotic approximation of Conjecture~\ref{VolConj-parabolic}.
    }
\end{figure}
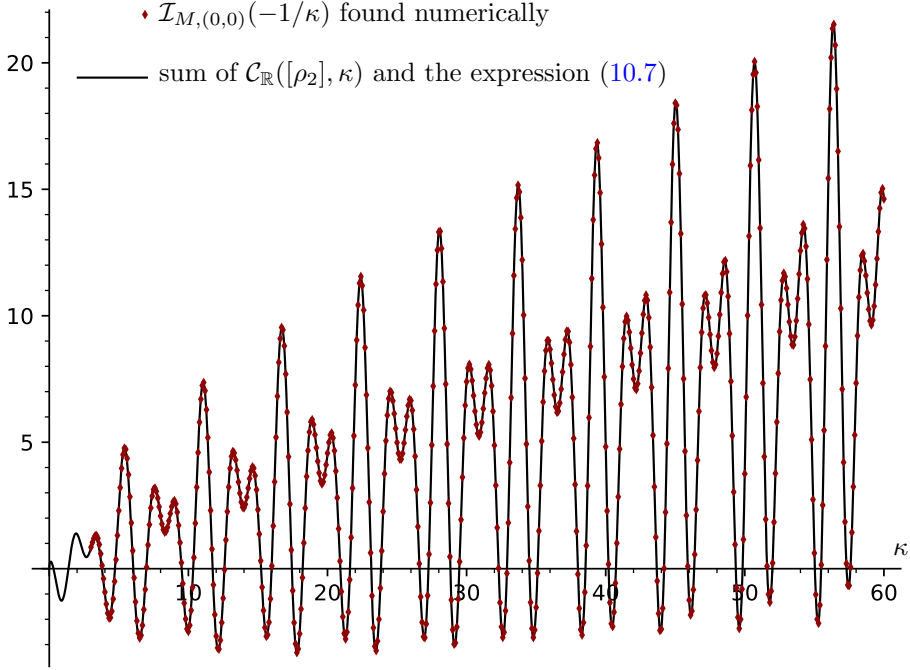

It remains to write down the contribution of the \realrep\ representation $\rho_2$ for which
we find an algebraic solution $z_2\in\left(\RR\setminus\{0,1\}\right)^{Q(\triang)}$
and the associated \Ztaut\ angle structure $\omega_2$.
In this way, we obtain the following triple Mellin--Barnes integral:
\begin{multline}
    \MB(\triang,\omega_2)=\frac{1}{(2\pi i)^3}
    \int_{-i\infty}^{i\infty}
    \int_{-i\infty}^{i\infty}
    \int_{-i\infty}^{i\infty}
    \Beta\bigl(\tfrac13+s_2-s_3,\; \tfrac23-2s_1+2s_3\bigr)\times\\
    \Beta^2\bigl(\tfrac12 + s_1,\; \tfrac13+s_2-s_3\bigr)
    \Beta\bigl(\tfrac13+2s_2,\;\tfrac12-2s_2+s_3\bigr)\,ds_1\,ds_2\,ds_3,
\end{multline}
which we believe to be the only non-vanishing contribution to the beta invariant.
Hence, $\contribreal([\rho_2], \kappa) = \MB(\triang,\omega_2)\,\kappa\approx 0.15698\,\kappa$.
The asymptotic approximation predicted by Conjecture~\ref{VolConj-parabolic} is the sum of these
five contributions and can be seen in Figure~\ref{asym:7_2}, plotted together with the values of
$\I{M}{0}{0}(-1/\kappa)$ found numerically for $\kappa\in\{3,\,3.1,\,3.2,\dotsc,60\}$.

\subsection{Methodology of numerical computations}
The values of the meromorphic \ind\ plotted in this section were found
with the help of a simple command line program called ``m3di'',
created by the last named author.
The source code of the program \cite{m3di} is publicly available under the terms
of the GNU General Public License.
We checked many of the numerical results produced by \texttt{m3di} against an independent implementation
developed in Matlab by the second named author, thus verifying the correctness of
the numerical results presented here.
In this section, we provide some details on how the program \texttt{m3di} works.

\subsubsection{Principle of operation}
The computation of the meromorphic \ind\ by
the program \texttt{m3di} is based on equation \eqref{3Dindex-int-2}, which can be rewritten
in the following form:
\begin{equation}\label{practical-formula}
    \I{\triang}{\aholM}{\aholL}(\hbar)
    = \frac{c_q^N}{(2\pi)^{N-1}}
    \int\limits_{[0,2\pi]^{N-1}}
    \prod_{\square\in Q(\triang)}
    G_q \Biggl(
    e^{(\hbar + i\pi)a(\square)}
    \exp\Bigl[i\sum_{j=1}^{N-1}t_j l_j(\square)\Bigr]
    \Biggr)\,dt.
\end{equation}
The above integral is approximated by a Riemann sum with a user-specified number
$S\in\NN$ of sample points per dimension, for a total of $S^{N-1}$
samples in the domain of integration.
We use evenly spaced sample points of the form
\[
    t_j = \frac{2\pi n}{S}, \quad n\in\{1,2,\dotsc,S\}, \quad 1\leq j \leq N-1.
\]
Since the leading-trailing deformations $l_j$ are integer vectors,
the \nword{G_q}{factor} corresponding to a quad~$\square$ will only be evaluated at points
of the form $\exp\bigl((\hbar + i\pi)a(\square)\bigr)\cdot\omega$,
where $\omega$ is an $S$-th root of unity.
For these reasons, the operation of \texttt{m3di} is split into two principal stages:
\begin{enumerate}
    \item%\label{algorithm:tabulation} 
    Tabulation: for every quad $\square\in Q(\triang)$,
    we compute the values
    \[
        \left\{G_q\left(\exp\bigl[(\hbar + i\pi)a(\square)\bigr]\cdot\omega\right) : \omega^S=1\right\}.
    \]
    This set contains the values of the respective factor of the integrand at all possible sample points.
    \item%\label{algorithm:integration} %commented label out to stop strange error message
    Integration: we compute the Riemann sum of the integral \eqref{practical-formula}.
    In order to evaluate the integrand at a sample point, we locate and multiply the values of the individual
    quad factors tabulated in the first stage.
\end{enumerate}

\subsubsection{Algorithmic complexity}
We shall now analyse the algorithmic complexity of the computations outlined above
in terms of the parameters $N$ and $S$.

We start with the tabulation stage.
Since there are $3N$ quads, each of which corresponds to a factor of the integrand that needs to be
evaluated at $S$ sample points, we see that the time complexity of the tabulation stage is $O(NS)$.
Moreover, the memory complexity is also $O(NS)$, because the tabulated values must be stored
in memory until the integration stage completes.

Meanwhile, the integration stage requires visiting all $S^{N-1}$ sample points in the integration domain,
and at each such point we must perform $O(N)$ multiplications to obtain the value of the integrand.
Hence, the integration stage has time complexity $O\bigl(NS^{N-1}\bigr)$.
The integration stage does not need much memory compared to the order $O(NS)$
amount already allocated for the tabulation.

As a consequence, the entire algorithm has exponential time complexity in the number of tetrahedra
of the triangulation.
The actual runtime of \texttt{m3di} depends strongly on the value of the parameter $\hbar$: the closer
$\Real\hbar$ is to zero, the slower the convergence of the infinite products \eqref{q-Pochhammer},
which lengthens the tabulation stage significantly.
With $N\leq 3$, \texttt{m3di} can be successfully run on an ordinary PC.
We used a supercomputing cluster at Nanyang Technological University
to study triangulations with $4$ tetrahedra and we believe that $N=5$
is also within reach for $S$ below $10^4$.

\subsubsection{Example computation}
Since \texttt{m3di} uses the formula \eqref{practical-formula}, we should supply the
following data as input:
\begin{itemize}
    \item The value of the parameter $\hbar\in\CC$ with $\Real\hbar < 0$;
    \item The desired number $S$ of samples per dimension;
    \item The number $N$ of tetrahedra in $\triang$;
    \item The matrix $\begin{bmatrix}\mathbf{L}\\\mathbf{L_\partial}\end{bmatrix}$
    of leading-trailing deformations of $\triang$ -- see Section~\ref{subsec:TAS};
    \item The vector $a\in\RR^{Q(\triang)}$ given by $a(\square)=\alpha(\square)/\pi$,
    where $\alpha$ is a strict angle structure with the desired peripheral angle-holonomy $(\aholM,\aholL)$.
    The entries of the vector $a$ should be ordered consistently with the columns of the matrix of
    leading-trailing deformations.
\end{itemize}

\begin{figure}%
    \centering
    \includegraphics[width=\columnwidth]{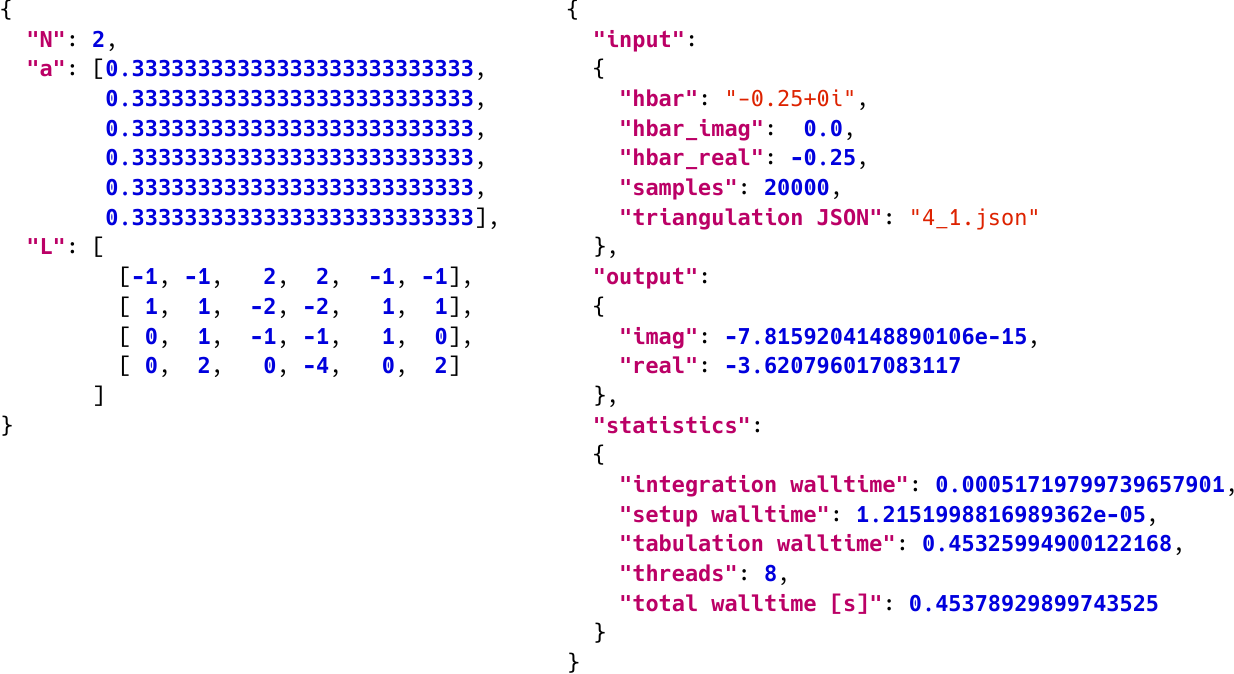}
    \caption{\label{fig:JSON}
       \textsc{Left:} A hypothetical JSON data file \texttt{4\char`_1.json} describing
       the two-tetrahedron triangulation of the figure-eight knot complement.
       \textsc{Right:} An example of possible output of the command \eqref{shell-command}.
    }
\end{figure}

We will now explain how to pass this data to \texttt{m3di} on the Linux command line,
although it ought to be possible to compile the program for other operating systems
using the source code \cite{m3di}.
As a concrete example, we will show how to calculate one of the values of the
meromorphic \ind\ of the figure-eight knot complement $M=S^3\setminus 4_1$ shown in
the plot of Figure~\ref{plot:asym-m004-00}.
We will use the standard two-tetrahedron triangulation $\triang$ and the geometric
strict angle structure $\alpha=\alpha^\geo$, which gives $a(\square)=\frac{1}{3}$ for
all $\square\in Q(\triang)$.
This structure is peripherally trivial, so that $(\aholM, \aholL)=(0,0)$.
In order to store the necessary input, we create a JSON data file
whose contents are shown on the left side of Figure~\ref{fig:JSON} and save it as
\texttt{4\char`_1.json}.
The remaining parameters $\hbar$ and $S$ can then be passed on the command line.
For example, to compute $\I{M}{0}{0}(-1/\kappa)$ for $\kappa=4$ (i.e., $\hbar = -0.25 + 0i$)
with $S=20\,000$ sample points, we could issue the following shell command:
\begin{equation}\label{shell-command}
    \text{\texttt{%
        m3di integrate 4\char`_1.json -0.25 0 20000
    }}
\end{equation}
whose output is another data structure in JSON format, shown on the right side
of Figure~\ref{fig:JSON}.

The output data contains three objects: \texttt{input}, \texttt{output} and \texttt{statistics}.
The object \texttt{input} records the input parameters for bookkeeping purposes, whereas
the object \texttt{statistics} contains technical information including the total run time of the program,
which was $0.45\,\mathrm{s}$ in our example.

The actual results are stored in the object \texttt{output}.
Since $\I{M}{0}{0}(-1/4)\in\RR$, the interesting information
is the number \texttt{output.real}, equal to $-3.620796017083117$ in our example.
Note that \texttt{output.imag} is of the order of $10^{-15}$, giving a sense of the magnitude of
numerical errors in this computation, albeit not a rigorous error estimate.

The plots in this paper were produced by running such computations for
different values of $\kappa$.
We used SnapPy \cite{snappy} to automatically create input data
files from census triangulations.
More information about \texttt{m3di} can be found in the `Read Me' file attached to the source
code of the program \cite{m3di}.

%==============================================================================
% Start the appendix
%==============================================================================
\appendix

%==============================================================================
% Appendix A: Relationship to the q-series 3D index
%==============================================================================

\section{\texorpdfstring{The $q$-series 3D index from the meromorphic 3D index}{The q-series 3D index from the meromorphic 3D index}}
\label{recover-qseries-appendix}

In this appendix we will state and prove a formula expressing the meromorphic \ind\ of Garoufalidis-Kashaev \cite{garoufalidis-kashaev} in terms of the $q$-series \ind\
of Dimofte-Gaiotto-Gukov \cite{dimofte-gaiotto-gukov 3D index}. The case of a knot complement was discovered and proven by Garoufalidis and Kashaev \cite{garoufalidis-kashaev}. In this appendix we will prove a version of this theorem for a general oriented ideal triangulation of a cusped 3-manifold carrying a strict angle structure. Our approach partially follows the proof of \cite{garoufalidis-kashaev}, but works with the definition of the $q$-series \ind\ which is symmetric with respect the quad-types, and we feel this allows us to clarify and generalize the deduction. Our proof relies on an interesting identity, Lemma~\ref{tripleJlemma}, which appears here for the first time. 

We begin by recalling the definition of the $q$-series \ind\ based on $Q$-normal surface theory. This approach was introduced in \cite{garoufalidis-hodgson-hoffman-rubinstein}. 

Let $M$ be an orientable cusped 3-manifold with $k\ge 1$ torus boundary components and let $\triang$ be an oriented ideal triangulation of $M$. 
As in Section~\ref{sec:gluing-equations}, once we choose orderings of the $N$ edges and the $N$ tetrahedra of the triangulation, and a distinguished quad-type in every tetrahedron, then we have gluing matrices. In order to conform to the notation of \cite{garoufalidis-kashaev} as much as possible in this appendix, we will denote the matrices $G,G',G''$ discussed in Section~\ref{sec:gluing-equations} here by $\overline{A}$, $\overline{B}$ and $\overline{C}$. A solution of the $Q$-normal equations is an assignment of an integer to every quad-type of the triangulation; here this is represented by a vector $S=(a_1,\ldots,a_N,b_1,\ldots,b_N,c_1,\ldots,c_N)\in \ZZ^{3N}$, such that
\[
[\, \overline{C}-\overline{B} \;|\; \overline{A}-\overline{C} \;|\; \overline{B}-\overline{A} \,]S^\transp=0.
\]
The group of integer solutions of the $Q$-normal equations is denoted $Q(\triang;\Z)\subset \ZZ^{3N}$. 

Solutions of the $Q$-normal equations have an interpretation related to spun normal surfaces in the  triangulation. Corresponding to this interpretation, there are two $\ZZ$-linear maps defined on the space of solutions, the \emph{homology map} $h:Q(\triang;\ZZ)\rightarrow H_2(M,\partial M;\ZZ_2)\oplus H_1(\partial M; \ZZ)$, and the \emph{Euler characteristic} $\chi:Q(\triang;\ZZ)\rightarrow \ZZ$. 

A \emph{tetrahedral solution} of the $Q$-normal equations is a solution which assigns $1$ to the three quad types of some selected tetrahedron, and zero to every other quad type. The subgroup of 
$Q(\triang;\ZZ)$ generated by tetrahedral solutions is denoted $\Tet$. Then $\Tet$ lies in the kernel of the homology map, so we get a well-defined map $h:Q(\triang;\ZZ)/\Tet\rightarrow H_2(M,\partial M;\ZZ_2)\oplus H_1(\partial M; \ZZ)$. Sometimes we will denote the projection to the second summand  by $\partial: Q(\triang;\ZZ)/\Tet\rightarrow H_1(\partial M; \ZZ)$; we call this the \emph{boundary map}.

The building block of the $q$-series \ind\ is the \emph{tetrahedral index} $I_\Delta:\ZZ^2\rightarrow \ZZ[[q^{\frac{1}{2}}]]$ defined by the formula
\[
I_\Delta(m,e)(q)=\sum_{n=\mathrm{max}\{0,-e\}}^\infty (-1)^n q^{\frac{1}{2}n(n+1)-(n+\frac{1}{2}e)m}\frac{1}{(q)_n(q)_{n+e}},
\]
where $(q)_n=\prod_{i=1}^n(1-q^i)$ denotes the $q$-Pochhammer symbol, with $(q)_0=1$.

We will use the more symmetric version of the tetrahedral index introduced by \cite{garoufalidis-hodgson-rubinstein-segerman}, which can be defined for $a,b,c\in\ZZ$ by
\[
J_\Delta(a,b,c)=(-q^{\frac{1}{2}})^{-b}I_\Delta(b-c,a-b).
\]
This function is invariant under all permutations of the variables \cite{garoufalidis-hodgson-rubinstein-segerman}. (An elegant way to see this is via our Lemma~\ref{tripleJlemma} below.) If we think of the variables as corresponding to the 3 quad-types of a tetrahedron, then this function has a particularly nice property with respect to the addition of a tetrahedral solution, namely
\begin{equation}
    J_\Delta(a+1,b+1,c+1)=(-q^{-\frac{1}{2}})J_\Delta(a,b,c).
\end{equation}
We extend this function to solutions of the $Q$-normal equations on the triangulation by defining
\[
J(S)=\prod_{i=1}^N J_\Delta(a_i,b_i,c_i)
\]
where $S=(a_1,\ldots,a_N,b_1,\ldots,b_N,c_1,\ldots,c_N)\in Q(\triang;\ZZ)$.

Now we review the definitions of several versions of the \ind\ following \cite{garoufalidis-hodgson-hoffman-rubinstein}. Their work generalizes and refines the 
original definitions using solutions of the $Q$-normal equations and their corresponding homology classes.

Let $(a,b)\in H_2(M,\partial M;\ZZ_2)\oplus H_1(\partial M; \ZZ)$ such that $\partial(a)=b$ modulo 2. Then
\[
I_\triang^{a}(b)(q)=\sum_{[S]\in Q(\triang;\ZZ)/\Tet, \, h([S])=(a,b)}(-q^{\frac{1}{2}})^{-\chi(S)}J(S).
\]
In other words, this is a sum over the elements of 
$Q(\triang;\ZZ)/\Tet$, 
where the term corresponding to some coset $x$ is calculated from a coset representative $S$, $[S]=x$.

Given a boundary homology class $\delta\in H_1(\partial M; \ZZ)$, the corresponding \emph{total} \ind\ is a sum over all $\Tet$-cosets of solutions to the $Q$-normal equations whose corresponding boundary homology class is $\delta$. That is, cosets $[S]$ such that $h([S])=(\omega,\delta)$ from some $\omega\in H_2(M,\partial M; \ZZ_2)$. This can be written:
\[
I_\triang^{\mathrm{tot}}(\delta)=\sum_{\omega\in H_2(M,\partial M; \ZZ_2),\partial(\omega)=\delta\,\mathrm{mod}\, 2} I^\omega_\triang(\delta)=
\sum_{[S]\in Q(\triang;\ZZ)/\Tet, \, \partial([S])=\delta} (-q^{\frac{1}{2}})^{-\chi(S)}J(S).
\]

Note that for a given boundary class $\delta\in H_1(\partial M;\ZZ)$, there will exist $\omega\in H_2(M,\partial M; \ZZ_2)$ with $\partial(\omega)=\delta$ mod 2 if and only if $\delta\in \ker(H_1(\bd M;\ZZ) \to H_1(M;\ZZ/2\ZZ))$. This follows from the relative homology exact sequence. In the special case that $M$ is a knot complement this subgroup equals
\[
\ker(H_1(\bd M;\ZZ) \to H_1(M;\ZZ/2\ZZ)) = 
\left\{ 2x\mu+y\lambda : x,y\in \ZZ\right\}
\]
where $\mu$ and $\lambda$ are elements representing a choice of meridian and longitude.

Now we turn to the main theorem of this appendix. Let
\[
\theta = (\alpha_1,\ldots,\alpha_N,\beta_1,\ldots,\beta_N,\gamma_1,\ldots,\gamma_N)\in \RR^{3N}
\]
denote a strict angle structure on the triangulation $\triang$. Corresponding to this angle structure there is the peripheral angular holonomy, $\rho_\theta: H_1(\partial M, \ZZ)\rightarrow \RR$, also denoted $\ahol{\theta}$ in \eqref{def-angle-holonomy}.
The state integral of the meromorphic \ind\ evaluated at this angle structure
(see equation~\ref{state-integral-formula}) only depends on $\rho_\theta$. This is due to \cite{garoufalidis-kashaev} and also follows elegantly from the formula stated in our main theorem below. Here we will write the state integral $I_{\triang,\rho_\theta}(q)$. In the notations of Section~\ref{sec:3D-index-with-strict-angles}, when we have fixed a meridian $\mu$ and longitude $\lambda$ this equals the meromorphic \ind\ $I_{\triang, (\rho_\theta(\mu),\rho_\theta(\lambda))}(q)$.

\medskip\noindent
\begin{theorem}
%{\bf Claim:} 
With the notations as above, if $\theta$ is a strict angle structure on the ideal triangulation $\triang$, then 
the meromorphic 3D-index is given by the formula 
\begin{equation}\label{meromseries}
I_{\triang,\, \rho_\theta}(q)
= \sum_{\delta \in K} \, (-q)^{ \frac{1}{2\pi} \, \rho_\theta(\delta)} \, I_\triang^{\mathrm{tot}}(\delta)(q^2),
\end{equation}
where $K=\ker(H_1(\bd M;\ZZ) \to H_1(M;\ZZ/2\ZZ))$.
\end{theorem}

\begin{proof} 
In the notation of \cite[(53)]{garoufalidis-kashaev}, the state integral can be expressed
$$I_{\triang,\, \rho_\theta}(q)=
\int_{\T^N} \prod_{i=1}^N  B(T_i,x,\theta) \, d\mu(x),$$
where 
$d\mu(x) = \frac{1}{(2\pi i)^N} \prod_{i=1}^N \frac{dx_i}{x_i}$ is the normalised Haar measure on $\T^N = (S^1)^N$ (denoted $\rm{dHaar}$ in (3.5)),
and
\cite[(52b)]{garoufalidis-kashaev} 
$$B(T_i,x,\theta)=c_q \, G_q\left( (-q)^{\frac{\alpha_i}{\pi}} x^{(\oC-\oB)_i} \right) 
G_q\left( (-q)^{\frac{\beta_i}{\pi}} x^{(\oA-\oC)_i} \right) G_q\left( (-q)^{\frac{\gamma_i}{\pi}} x^{(\oB-\oA)_i} \right)
$$ 
is the  Boltzmann weight for the $i$th tetrahedron $T_i$. The function $G_q(z)$ is given here by equation~\ref{def:G_q}. In this expression, the notation $A_i$ denotes the $i$-th column of the matrix $A$, and if $x =(x_1,\ldots,x_n) \in \CC_*^N$ and $v=(v_1,\ldots,v_N)\in\ZZ^N$ then the notation
$x^v$ denotes $\prod_{i=1}^N x_i^{v_i}$.

\medskip
Now from \cite[Lemma 2.2(c)]{garoufalidis-kashaev}, $G_q(z)$ has a convergent Laurent series expansion for $0<|z|<1$:
$$G_q(z)=\sum_{n\in\Z} J(n)(q) z^n,$$
where \cite[(10)]{garoufalidis-kashaev}
$$J(n)(q):= \sum_{k=(-n)_+}^\infty  \frac{q^{k(k+1)/2}}{(q)_k (q)_{n+k}},  \text{ with } (n)_+=\max(n,0).$$
Alternatively,
$$J(n)(q):= \sum_{k\in\Z}  \frac{q^{k(k+1)/2}}{(q)_k (q)_{n+k}},$$
if we take $\frac{1}{(q)_m}=0$ for $m<0$.

Given a strict 
angle structure $\theta$, we can use 
the Laurent series expansion for all $x \in \T^n$ to obtain
$$G_q\left( (-q)^{\frac{\alpha_i}{\pi}} x^{(\oC-\oB)_i} \right)  
=\sum_{a_i\in\Z} J(a_i)(q) (-q)^{\frac{a_i \alpha_i}{\pi}} x^{(\oC-\oB)_i a_i} ,$$
$$G_q\left( (-q)^{\frac{\beta_i}{\pi}} x^{(\oA-\oC)_i} \right)  
=\sum_{b_i\in\Z} J(b_i)(q) (-q)^{\frac{b_i \beta_i}{\pi}} x^{(\oA-\oC)_i b_i} , $$
$$G_q\left( (-q)^{\frac{\gamma_i}{\pi}} x^{(\oB-\oA)_i} \right)  
=\sum_{c_i\in\Z} J(c_i)(q) (-q)^{\frac{c_i \gamma_i}{\pi}} x^{(\oB-\oA)_i c_i} .$$

Hence
$$I_{\triang,\, \rho_\theta}(q)=
\int_{x \in \T^N} \sum_{a,b,c \in \Z^N} \prod_{i=1}^N c_q J(a_i)J(b_i) J(c_i)  
\, (-q)^{\frac{1}{\pi}\sum_{i=1}^N a_i\alpha_i + b_i \beta_i+c_i\gamma_i}
\, x^{(\oC-\oB)a^T + (\oA-\oC)b^T+(\oB-\oA)c^T} \, d\mu(x) 
$$
where $a=(a_1, \ldots, a_N)$, $b=(b_1, \ldots, b_N)$, $c=(c_1, \ldots, c_N)$.

Now orthogonailty of complex exponentials on $\T^N$ gives
$$
\int_{\T^N} x^v d\mu(x) =\int_{\T^N} x_1^{v_1} \ldots x_n^{v_N} \, d\mu(x) 
=\begin{cases} 1 \text{ if } v = 0,\\  0 \text{ otherwise }\end{cases}
$$
and 
$(\oC-\oB)a^T + (\oA-\oC)b^T+(\oB-\oA)c^T=0$ 
if and only if $S=(a,b,c) \in Q(\triang;\Z)$ is an integer solution to the Q-matching equations for normal surfaces in $\triang$. Note that $(a,b,c)$ here denotes $(a_1,\ldots,a_N,b_1,\ldots,b_N,c_{1},\ldots,c_N)\in \ZZ^{3N}$.

Hence
\begin{equation}\label{preindex1}
I_{\triang,\, \rho_\theta}(q)=
\sum_{S=(a,b,c) \in Q(\triang;\Z)} (-q)^{\frac{1}{\pi} \langle \theta ,S \rangle}\, \prod_{i=1}^N c_q J(a_i)J(b_i) J(c_i)  ,
\end{equation}
where
$$\langle \theta ,S \rangle := \sum_{i=1}^N a_i\alpha_i + b_i \beta_i+c_i\gamma_i$$
denotes the natural pairing between generalised angle structures $\theta \in \RR^{3N}$ and integral $Q$-normal classes $S \in \Z^{3N}$.

Now,  because $\theta$ is a generalised angle structure, a combinatorial Gauss-Bonnet formula from 
\cite[(39)]{garoufalidis-hodgson-hoffman-rubinstein} gives
\begin{equation}\label{GB}
\frac{1}{\pi} \langle \theta ,S \rangle = -\chi(S)+\frac{1}{2\pi} \rho_\theta(\bd S)
\end{equation}
where  $\chi(S)$ is the formal Euler characteristic of $S$ and 
$\rho_\theta(\bd S)$ denotes the sum of rotational holonomies of the boundary components of $S$, oriented as in $[\bd S]$.

Thus:
\begin{equation}\label{eq1}
I_{\triang,\, \rho_\theta}(q)=
\sum_{S=(a,b,c) \in Q(\triang;\Z)} (-q)^{ \frac{1}{2\pi} \, \rho_\theta(\bd S)}  \,(-q)^{-\chi(S)} \,  \prod_{i=1}^N c_q J(a_i)J(b_i) J(c_i).
\end{equation}

To finish we collect the terms of this sum into $\Tet$-cosets, then apply Lemma~\ref{tripleJlemma} to evaluate the corresponding sums. If some coset is represented by $S=(a,b,c)\in Q(\triang;\Z)$ then the sum over all the terms corresponding to that coset equals
\[
(-q)^{ \frac{1}{2\pi} \, \rho_\theta(\bd S)}(-q)^{-\chi(S)}
\sum_{(t_1,\ldots,t_N)\in \ZZ^N} 
\prod_{i=1}^N c_q (-q)^{t_i} J(a_i+t_i)J(b_i+t_i) J(c_i+t_i).
\]
This expression uses the fact that one of the generating tetrahedral solutions $X$, consisting of $1$'s in the quad-types of some tetrahedron and zeroes in every other tetrahedron, satisfies $h(X)=0$ and $\chi(X)=-1$.

Using Lemma~\ref{tripleJlemma} we obtain
\[
I_{\triang,\, \rho_\theta}(q)
=\sum_{[S]\in Q(\triang;\ZZ)/\Tet}
(-q)^{ \frac{1}{2\pi} \, \rho_\theta(\bd S)}(-q)^{-\chi(S)}J(S)(q^2).
\]
If we organize this sum into homology classes on the boundary, we finally obtain
\[
I_{\triang,\, \rho_\theta}(q)
=\sum_{\delta\in \ker(H_1(\bd M;\ZZ) \to H_1(M;\ZZ/2\ZZ)) }
(-q)^{ \frac{1}{2\pi} \, \rho_\theta(\delta)}   
I_\triang^{\mathrm{tot}}(\delta)(q^2),
\]
as required.
\end{proof}

\begin{lemma}\label{tripleJlemma} If $a,b,c \in \ZZ$ then
$$J_\Delta(a,b,c)(q^2) = c_q \sum_{t\in\ZZ}  (-q)^t J(a+t)(q) J(b+t)(q) J(c+t)(q)$$
\end{lemma}

\begin{proof}
Following [GK,(21)] and the definition of $J_\Delta$  [GHHR, (5)] we have
$$I^\Delta(m,e)(q)=(-q)^e I_\Delta(m,e)(q^2)=(-q)^e J_\Delta(-m,e,0)(q^2).$$
Equating (18) and (24) in [GK] gives
$$I^\Delta(m,e)(q)=c_q \sum_{\substack{k_1,k_2,k_3\in\ZZ\\k_3-k_2=m\\k_1-k_3=e}} (-q)^{k_1}J(k_1)(q)J(k_2)(q)J(k_2)(q).$$
Combining the last two equations gives 
$$
J_\Delta(-m,e,0)(q^2) = (-q)^{-e} c_q \sum_{\substack{k_1,k_2,k_3\in\ZZ\\k_3-k_2=m\\k_1-k_3=e}} (-q)^{k_1}J(k_1)(q)J(k_2)(q)J(k_2)(q).
$$
Now for $t \in \ZZ$, [GHHR,(6)] gives
$$J_\Delta(-m+t,e+t,t)(q^2)=(-q)^{-t} J_\Delta(-m,e,0)(q^2),$$
and we can rewrite this with $a=-m+t, b=e+t, c=t$ by taking $t=c,m=c-a,e=b-c$.
This gives 
$$J_\Delta(a,b,c)(q^2)=
(-q)^{-b} c_q  \sum_{\substack{k_1,k_2,k_3\in\ZZ\\k_3-k_2=c-a\\k_1-k_3=b-c}} (-q)^{k_1}J(k_1)(q)J(k_2)(q)J(k_2)(q).
$$
The linear equations for $k_1,k_2,k_3$ have general solution $k_1=b+t,k_2=a+t,k_3=c+t$ where $t\in\ZZ$, and 
$-b+k_1=t$, 
 so this can be rewritten in the simple symmetric form
$$J_\Delta(a,b,c)(q^2)=
c_q  \sum_{t\in\ZZ} (-q)^{t}J(a+t)(q)J(b+t)(q)J(c+t)(q).
$$
This proves the Lemma.
\end{proof}

% %%%%%%%%
\begin{example}
\label{appendix-fig-eight-example}
For the exterior $M$ of the figure eight knot, we use 
Thurston's ideal triangulation $\triang$  with 2 tetrahedra, and the standard meridian and longitude $\mu,\lambda \in H_1(\bd M;\Z)$.
Then $H_2(M,\bd M;\Z/2\Z)\cong H^1(M;\Z/2\Z)=0$, 
$K=\ker(H_1(\bd M;\Z) \to H_1(M;\Z/2\Z))$ is spanned by $\{2\mu,\lambda\}$, and
for $x,y\in\Z$, we have
 
\begin{align*} 
I_\triang^{\mathrm{tot}}(2x \mu + y\lambda)(q)=I_\triang^0(2x \mu + y\lambda)(q)
& = \sum_{k \in \Z} \, \ID(k-x,k) \ID(k+y,k-x+y)
\end{align*}
using \cite[Example 4.1, Remark 8.4] {garoufalidis-hodgson-hoffman-rubinstein}.

Hence, if  $\theta$ is a strict angle structure with trivial peripheral rotational holonomy $\rho_\theta$ then the meromorphic 3D-index is given by
\begin{align*} I_{\triang,\, \rho_\theta}(q) &= \sum_{x,y \in \Z} \, I_\triang^0(2x \mu + y\lambda)(q^2) \\
&=1-4 q-8 q^2-8 q^3-q^4+4 q^5+48 q^6+24 q^7+110 q^8+36 q^9+144 q^{10} + \ldots
\end{align*}

\noindent{\bf Remark:} In \cite{garoufalidis-wheeler}, the authors consider the asymptotics of a related, but different, $q$-series given by
\begin{align*} &\sum_{x,y \in \Z} \, I_\triang^0(2x \mu + 2y\lambda)(q) \\
&=1-4 q-q^2+36 q^3+70 q^4+100 q^5+34 q^6-116 q^7-410 q^8-808 q^9-1140 q^{10} + \ldots .
\end{align*}

\end{example}

%==============================================================================
% Appendix B: The asymptotic expansion
%==============================================================================

\section{\texorpdfstring{The asymptotic expansion of the $q$-dilogarithm}{The asymptotic expansion of the q-dilogarithm}}
\label{asymptotic-appendix}

In this appendix we will explain how to obtain the asymptotic expansion given in
equation~\eqref{Andrews-expansion} using the method of Mellin transforms. See \cite{mellin} for a clear exposition of this method.

We fix arbitrarily the parameter $\omega\in\CC$ satisfying $|\omega|=1$ and $\omega\neq 1$,
as well as the parameter $a\in\CC$ with $\Real a>0$.
In terms of the negative real parameter $\hbar < 0$, the expansion is
\[
    \Li(\omega e^{a\hbar}; e^\hbar) =
    -\Li(\omega)\hbar^{-1}
    + (a-{\textstyle\frac12})\log(1-\omega)
    - \sum_{k=1}^{\infty}\frac{B_{k+1}(a)}{(k+1)!}\Lis{1-k}(\omega)\,\hbar^k
\]
as $\hbar\to0^-$.
The analogous formula for the case $\omega=1$, equation~\eqref{Zhang},
was obtained by Zhang using a short and conceptually attractive Mellin transform
argument\footnote{Note that the formula itself appears earlier in McIntosh~\cite{mcintosh-qfactorials}.}.
The formula is stated as Theorem~2 in \cite{zhang-qexp} and is clearly explained there.
This approach will be convenient for generalisation to the case $\omega\neq 1$.

To match the discussion in \cite{zhang-qexp} we set $\hbar=-x$, for $x>0$,
and consider $\Li(\omega e^{-ax}; e^{-x})$ as $x\rightarrow 0^+$.
The heart of the computation of the case $\omega=1$ is equation~(3.1) of \cite{zhang-qexp} which
presents the computation of the Mellin transform of the function $\Li(e^{-ax}; e^{-x})$,
\[
    \int_0^\infty \Li(e^{-ax}; e^{-x})x^{s-1}\;dx=\Gamma(s)\zeta(s+1)\zeta(s,a).
\]
Here
\begin{itemize}
    \item $\Gamma(s)$ is the Euler gamma function,
    \item $\zeta(s)$ is the Riemann zeta function $\zeta(s)=\sum_{n=1}^{\infty}\frac{1}{n^s}$,
    \item and $\zeta(s,a)=\sum_{n=0}^\infty\frac{1}{(a+n)^s}$ is the Hurwitz zeta function.
\end{itemize}

The calculation leading to this formula is Section~4.2 of \cite{zhang-qexp}.
When the phase term $\omega$ is introduced into this calculation, it is straightforward to compute
that the only difference in the result is that the Riemann zeta function factor $\zeta(s+1)$ is
generalised to the \emph{Lerch zeta function}:
\[
    \zeta_\omega(s+1)=\sum_{n=1}^{\infty}\frac{\omega^n}{n^{s+1}}.
\]
In other words:
\begin{equation}\label{eq:Mellin-transform}
    \int_0^\infty \Li(e^{-ax}\omega; e^{-x})x^{s-1}\;dx=\Gamma(s)\zeta_\omega(s+1)\zeta(s,a).
\end{equation}

To be able to apply the inverse Mellin transform along some vertical contour in the $s$-plane,
$\gamma:\RR\rightarrow \CC$, $\gamma(t)=c+it$, for some fixed $c>1$, it is sufficient to know
that the Mellin transform \eqref{eq:Mellin-transform} is integrable along that contour
(see Theorem~2 of \cite{mellin}).
Zhang checks this for the $\omega=1$ case $\Gamma(s)\zeta(s+1)\zeta(s,a)$ using the
asymptotics of $\Gamma(c+it)$ and $\zeta(c+it,a)$ as $t\rightarrow\pm\infty$ (see Lemma~5 of \cite{zhang-qexp}).
When $\zeta(s)$ is replaced by $\zeta_\omega(s)$ one can see that there is no change needed to this argument because
\[
    |\zeta_\omega(c+it)| \leq \zeta(c)
\]
and hence the Lerch zeta factor is bounded on the contour.
Thus the asymtptotics given in \cite{zhang-qexp} for $\Gamma(c+it)$ and $\zeta(c+it,a)$ still imply
integrability in the case $\omega\neq1$.
As a consequence,
\begin{equation}\label{eq:inverse-Mellin-transform}
    \Li(e^{-ax}\omega; e^{-x})
    =
    \frac{1}{2\pi i}
    \int_{c-i\infty}^{c+i\infty} \Gamma(s)\zeta_\omega(s+1)\zeta(s,a) \frac{1}{x^s}\;ds.
\end{equation}
Following the method of Mellin transforms \cite{mellin},
we will now apply residue calculus to convert the integral \eqref{eq:inverse-Mellin-transform}
into an asymptotic expansion.
Note that the integral is of the order $O(x^{-c})$ and that the integrand
can be continued meromorphically onto the whole \nword{s}{plane}.

The method is to shift the horizontal position of the vertical contour to the left.
When the contour crosses a pole, the expression picks up a term from the residue at that pole.
So if $d<c$, and $P$ denotes the set of poles of
$\Gamma(s)\zeta_\omega(s+1)\zeta(s,a) \frac{1}{x^s}$ lying between the vertical contours at $d$ and $c$, then
\begin{equation}\label{residue-expansion}
    \Li(e^{-ax}\omega; e^{-x})
    =\sum_{p\in P} \Res_{s=p}
    \left(
    \Gamma(s)\zeta_\omega(s+1)\zeta(s,a) \frac{1}{x^s}
    \right)
    + O(x^{-d}).
\end{equation}

It remains for us to understand the poles of $\Gamma(s)\zeta_\omega(s+1)\zeta(s,a) \frac{1}{x^s}$
and the residues at those poles under the assumption that $\omega\neq 1$:
\begin{itemize}
    \item When $\omega\neq 1$, the Lerch zeta function $\zeta_\omega(s+1)$ is entire and has no poles \cite{apostol}.
    \item The Hurwitz zeta function $\zeta(s,a)$ has a simple pole at $s=1$ with the
    corresponding residue equal to $1$ \cite{zhang-qexp}.
    \item The gamma function $\Gamma(s)$ has simple poles at all non-positive integers
    $s=-k$ with corresponding residues $\frac{(-1)^k}{k!}$. (See \cite{erdelyi-bateman,zhang-qexp}.)
\end{itemize}
In summary, the Mellin transform \eqref{eq:Mellin-transform} has simple poles at
$s\in\{1,0\}\cup \{-k,\;k\in \NN \}$.

To determine the corresponding residues, we will need the following formula for the values of the
Lerch zeta function at the non-positive integers $m\leq 0$:
\begin{equation}\label{lerche-equals-dilog}
    \zeta_\omega(m)=\Lis{m}(\omega).
\end{equation}
This formula can be obtained from the analysis in Apostol \cite{apostol}.
To relate our notations to the formulae in \cite{apostol}, note that for a real $t$
such that $\omega = \exp(2\pi it)$ and a complex $s$, we have
\[
    \zeta_\omega(s) = \zeta_{\exp(2\pi it)}(s) = e^{2\pi it}\phi(t,1,s),
\]
where $\phi$ is the function defined in equation~(1.1) of \cite{apostol}, also known
as the \emph{Lerch transcendent.}
Then the case $m=0$ of equation~\eqref{lerche-equals-dilog} is given at the start of Section~3
in \cite{apostol}:
\[
    \zeta_{\exp(2\pi it)}(0)=e^{2\pi i t}\phi(t,1,s)=e^{2\pi i t}\frac{1}{1-e^{2\pi i t}}=\Lis{0}(\omega),
\]
using equation~\eqref{definition-polylogarithm}.
Furthermore, the differential equation
\[
    \frac{\partial}{\partial t}\phi(t,a,s)+2\pi i a \phi(t,a,s)=2\pi i \phi(t,a,s-1),
\]
given in \cite[eq.~(2.6)]{apostol},
can be used to derive a recursion relation as follows:
\begin{align*}
    \zeta_{\omega}(s-1)
    &=  e^{2\pi i t}\phi(t,1,s-1) \\
    &= \frac{1}{2\pi i}\left(e^{2\pi i t} \frac{\partial}{\partial t}\phi(t,1,s)+2\pi i e^{2\pi i t}\phi(t,1,s)\right) \\
    &= \frac{1}{2\pi i}\frac{\partial}{\partial t}\left(e^{2\pi i t}\phi(t,1,s) \right) \\
    &= \frac{1}{2\pi i}\frac{\partial}{\partial t} \zeta_{\exp(2\pi it)}(s) \\
    &= \frac{1}{2\pi i}\left(\frac{\partial}{\partial t}e^{2\pi i t}\right)
    \frac{\partial}{\partial \omega}\zeta_\omega(s)
    =\omega \frac{\partial}{\partial \omega}\zeta_\omega(s).
\end{align*}
Thus the two sides of \eqref{lerche-equals-dilog} coincide at $m=0$ and satisfy the same recursion relation, so the equation is established for all $m\leq0$.

Now we can turn to computing the residues of $\Gamma(s)\zeta_\omega(s+1)\zeta(s,a) \frac{1}{x^s}$:
\begin{itemize}
    \item The residue at the pole at $s=1$ equals:
    \[
    \Gamma(1)\zeta_\omega(1+1)
    \Res_{s=1}\bigl(\zeta(s,a)\bigr)\frac{1}{x^1}=1\cdot \left(\sum_{n=1}^\infty\frac{\omega^n}{n^2}\right)\cdot 1\cdot \frac{1}{x}=\Li(\omega)\frac{1}{x}.
    \]
    \item The residue at the pole at $s=0$ equals:
    \[
    \Res_{s=0}\bigl(\Gamma(s)\bigr)\zeta_\omega(1)\zeta(0,a)\frac{1}{x^0} = 1\cdot \left(\sum_{n=1}^\infty \frac{\omega^n}{n}\right)\cdot \bigl(\tfrac{1}{2}-a\bigr) = \left(a-\tfrac{1}{2}\right)\log(1-\omega)
    \]
    using the evaluation $\zeta(0,a)=\frac{1}{2} - a$; cf. \cite[eq.~1.10(8)]{erdelyi-bateman}.
    \item The residue at the pole at $s=-k$, where $k\in \NN$, equals:
    \[
        \Res_{s=-k}\bigl(\Gamma(s)\bigr)\zeta_\omega(1-k)\zeta(-k,a)x^k
        =\frac{(-1)^k}{k!}\cdot
        \zeta_\omega(1-k)
        \cdot\left(-\frac{B_{k+1}(a)}{k+1}\right)x^k.
    \]
By virtue of equation~\eqref{lerche-equals-dilog}, this equals
\[
    -\frac{1}{(k+1)!}B_{k+1}(a)\Lis{1-k}(\omega)(-x)^k.
\]
\end{itemize}
Substituting these calculations in equation~\eqref{residue-expansion}
and putting $\hbar = -x$ yields the asymptotic expansion given in equation~\eqref{Andrews-expansion}.

%==============================================================================
% The bibliography
%==============================================================================

%==============================================================================
% The End
%==============================================================================
\end{document}